# The Lie algebra $f_4(O_s)$ (split $f_4$) with *Mathematica*

**Pablo Alberca Bjerregaard**

*Department of Applied Mathematics*
*Escuela Técnica Superior de Ingeniería Industrial*
*Campus de El Ejido*
*Universidad de Málaga*
*Málaga –Spain*
*E–Mail:pgalberca@uma.es*

**Cándido Martín González**

*Department of Algebra, Geometry and Topology*
*Facultad de Ciencias*
*Universidad de Málaga*
*Campus de Teatinos*
*Málaga –Spain*
*E–Mail:candido@apncs.cie.uma.es*

2001.

**Abstract.** *We present in this paper all the details for a complete description of the Lie algebra $f_4$ in the split case at any characteristic. We finish with the determination of the expression of a generic element of this algebra. First of all is necessary to implement its quadratic Jordan structure (see the O. Loss' server http://mathematik.uibk.ac.at/jordan for our preprint). We write in this file only the computational details.*

The definition of *derivation* in the context of quadratic Jordan algebras $(J, U, T)$ is a linear map $D$ such that:

$$\text{and } D(y)) \qquad \forall\, x,y$$

We will use also the *Leibniz* identity $D(\{x, y, z\}) = \{D(x), y, z\} + \{x, D(y), z\} + \{x, y, D(z)\}$

## ■ The images of $E_i[1]$, $i = 1, 2, 3$.

```
Do[Print[U_{E_i[1]}[ONE] == E_i[1]], {i, 3}]
```

True

True

True

```
Leib[x_, y_, z_] := T[Δ[x], y, z] + T[x, Δ[y], z] + T[x, y, Δ[z]]

Δ[E_1[1]] = generic;
```

```
Δ[E₁[1]] - T[Δ[E₁[1]], ONE, E₁[1]] // Expand
```

$$\begin{pmatrix}
\begin{pmatrix} -\lambda_1 & \{0,0,0\} \\ \{0,0,0\} & -\lambda_1 \end{pmatrix} & \begin{pmatrix} 0 & \{0,0,0\} \\ \{0,0,0\} & 0 \end{pmatrix} & \begin{pmatrix} 0 & \{0,0,0\} \\ \{0,0,0\} & 0 \end{pmatrix} \\
\begin{pmatrix} 0 & \{0,0,0\} \\ \{0,0,0\} & 0 \end{pmatrix} & \begin{pmatrix} \lambda_2 & \{0,0,0\} \\ \{0,0,0\} & \lambda_2 \end{pmatrix} & \begin{pmatrix} \lambda_{12} & \{\lambda_{14},\lambda_{15},\lambda_{16}\} \\ \{\lambda_{17},\lambda_{18},\lambda_{19}\} & \lambda_{13} \end{pmatrix} \\
\begin{pmatrix} 0 & \{0,0,0\} \\ \{0,0,0\} & 0 \end{pmatrix} & \begin{pmatrix} \lambda_{13} & \{-\lambda_{14},-\lambda_{15},-\lambda_{16}\} \\ \{-\lambda_{17},-\lambda_{18},-\lambda_{19}\} & \lambda_{12} \end{pmatrix} & \begin{pmatrix} \lambda_3 & \{0,0,0\} \\ \{0,0,0\} & \lambda_3 \end{pmatrix}
\end{pmatrix}$$

```
Δ[E₁[1]] = Δ[E₁[1]] //. {λ₁ → 0, λ₂ → 0, λ₃ → 0, λ₁₂ → 0,
   λ₁₃ → 0, λ₁₄ → 0, λ₁₅ → 0, λ₁₆ → 0, λ₁₇ → 0, λ₁₈ → 0, λ₁₉ → 0}
```

$$\begin{pmatrix}
\begin{pmatrix} 0 & \{0,0,0\} \\ \{0,0,0\} & 0 \end{pmatrix} & \begin{pmatrix} \lambda_4 & \{\lambda_6,\lambda_7,\lambda_8\} \\ \{\lambda_9,\lambda_{10},\lambda_{11}\} & \lambda_5 \end{pmatrix} & \begin{pmatrix} \lambda_{20} & \{\lambda_{22},\lambda_{23},\lambda_{24}\} \\ \{\lambda_{25},\lambda_{26},\lambda_{27}\} & \lambda_{21} \end{pmatrix} \\
\begin{pmatrix} \lambda_5 & \{-\lambda_6,-\lambda_7,-\lambda_8\} \\ \{-\lambda_9,-\lambda_{10},-\lambda_{11}\} & \lambda_4 \end{pmatrix} & \begin{pmatrix} 0 & \{0,0,0\} \\ \{0,0,0\} & 0 \end{pmatrix} & \begin{pmatrix} 0 & \{0,0,0\} \\ \{0,0,0\} & 0 \end{pmatrix} \\
\begin{pmatrix} \lambda_{21} & \{-\lambda_{22},-\lambda_{23},-\lambda_{24}\} \\ \{-\lambda_{25},-\lambda_{26},-\lambda_{27}\} & \lambda_{20} \end{pmatrix} & \begin{pmatrix} 0 & \{0,0,0\} \\ \{0,0,0\} & 0 \end{pmatrix} & \begin{pmatrix} 0 & \{0,0,0\} \\ \{0,0,0\} & 0 \end{pmatrix}
\end{pmatrix}$$

```
Δ[E₁[1]] =
 Δ[E₁[1]] //. {λ₄ -> α₁, λ₅ -> α₂, λ₆ -> α₃, λ₇ -> α₄, λ₈ -> α₅,
   λ₉ -> α₆, λ₁₀ -> α₇, λ₁₁ -> α₈, λ₂₀ -> β₁, λ₂₁ -> β₂, λ₂₂ -> β₃,
   λ₂₃ -> β₄, λ₂₄ -> β₅, λ₂₅ -> β₆, λ₂₆ -> β₇, λ₂₇ -> β₈}
```

$$\begin{pmatrix}
\begin{pmatrix} 0 & \{0,0,0\} \\ \{0,0,0\} & 0 \end{pmatrix} & \begin{pmatrix} \alpha_1 & \{\alpha_3,\alpha_4,\alpha_5\} \\ \{\alpha_6,\alpha_7,\alpha_8\} & \alpha_2 \end{pmatrix} & \begin{pmatrix} \beta_1 & \{\beta_3,\beta_4,\beta_5\} \\ \{\beta_6,\beta_7,\beta_8\} & \beta_2 \end{pmatrix} \\
\begin{pmatrix} \alpha_2 & \{-\alpha_3,-\alpha_4,-\alpha_5\} \\ \{-\alpha_6,-\alpha_7,-\alpha_8\} & \alpha_1 \end{pmatrix} & \begin{pmatrix} 0 & \{0,0,0\} \\ \{0,0,0\} & 0 \end{pmatrix} & \begin{pmatrix} 0 & \{0,0,0\} \\ \{0,0,0\} & 0 \end{pmatrix} \\
\begin{pmatrix} \beta_2 & \{-\beta_3,-\beta_4,-\beta_5\} \\ \{-\beta_6,-\beta_7,-\beta_8\} & \beta_1 \end{pmatrix} & \begin{pmatrix} 0 & \{0,0,0\} \\ \{0,0,0\} & 0 \end{pmatrix} & \begin{pmatrix} 0 & \{0,0,0\} \\ \{0,0,0\} & 0 \end{pmatrix}
\end{pmatrix}$$

```
T[E₁[1], E₁[1], E₁[1]] == 2 E₁[1]
```

True

```
Leib[E₁[1], E₁[1], E₁[1]] == 2 Δ[E₁[1]]
```

True

```
Δ[E₂[1]] = generic;
```

```
Δ[E₂[1]] - T[Δ[E₂[1]], ONE, E₂[1]] // Expand
```

$$\begin{pmatrix}
\begin{pmatrix} \lambda_1 & \{0,0,0\} \\ \{0,0,0\} & \lambda_1 \end{pmatrix} & \begin{pmatrix} 0 & \{0,0,0\} \\ \{0,0,0\} & 0 \end{pmatrix} & \begin{pmatrix} \lambda_{20} & \{\lambda_{22},\lambda_{23},\lambda_{24}\} \\ \{\lambda_{25},\lambda_{26},\lambda_{27}\} & \lambda_{21} \end{pmatrix} \\
\begin{pmatrix} 0 & \{0,0,0\} \\ \{0,0,0\} & 0 \end{pmatrix} & \begin{pmatrix} -\lambda_2 & \{0,0,0\} \\ \{0,0,0\} & -\lambda_2 \end{pmatrix} & \begin{pmatrix} 0 & \{0,0,0\} \\ \{0,0,0\} & 0 \end{pmatrix} \\
\begin{pmatrix} \lambda_{21} & \{-\lambda_{22},-\lambda_{23},-\lambda_{24}\} \\ \{-\lambda_{25},-\lambda_{26},-\lambda_{27}\} & \lambda_{20} \end{pmatrix} & \begin{pmatrix} 0 & \{0,0,0\} \\ \{0,0,0\} & 0 \end{pmatrix} & \begin{pmatrix} \lambda_3 & \{0,0,0\} \\ \{0,0,0\} & \lambda_3 \end{pmatrix}
\end{pmatrix}$$

```
Δ[E₂[1]] = Δ[E₂[1]] //. {λ₁ → 0, λ₂ → 0, λ₃ → 0, λ₂₀ → 0,
    λ₂₁ → 0, λ₂₂ → 0, λ₂₃ → 0, λ₂₄ → 0, λ₂₅ → 0, λ₂₆ → 0, λ₂₇ → 0}
```

$$\begin{pmatrix} \begin{pmatrix} 0 & \{0,0,0\} \\ \{0,0,0\} & 0 \end{pmatrix} & \begin{pmatrix} \lambda_4 & \{\lambda_6,\lambda_7,\lambda_8\} \\ \{\lambda_9,\lambda_{10},\lambda_{11}\} & \lambda_5 \end{pmatrix} & \begin{pmatrix} 0 & \{0,0 \\ \{0,0,0\} & 0 \end{pmatrix} \\ \begin{pmatrix} \lambda_5 & \{-\lambda_6,-\lambda_7,-\lambda_8\} \\ \{-\lambda_9,-\lambda_{10},-\lambda_{11}\} & \lambda_4 \end{pmatrix} & \begin{pmatrix} 0 & \{0,0,0\} \\ \{0,0,0\} & 0 \end{pmatrix} & \begin{pmatrix} \lambda_{12} & \{\lambda_{14}, \\ \{\lambda_{17},\lambda_{18},\lambda_{19}\} & \end{pmatrix} \\ \begin{pmatrix} 0 & \{0,0,0\} \\ \{0,0,0\} & 0 \end{pmatrix} & \begin{pmatrix} \lambda_{13} & \{-\lambda_{14},-\lambda_{15},-\lambda_{16}\} \\ \{-\lambda_{17},-\lambda_{18},-\lambda_{19}\} & \lambda_{12} \end{pmatrix} & \begin{pmatrix} 0 & \{0,0 \\ \{0,0,0\} & 0 \end{pmatrix} \end{pmatrix}$$

**U_{E₁[1]}[E₂[1]] == ZERO**

True

**T[Δ[E₁[1]], E₂[1], E₁[1]] + U_{E₁[1]}[Δ[E₂[1]]] == ZERO**

True

**?Leib**

Global`Leib

Leib[x_, y_, z_] := T[Δ[x], y, z] + T[x, Δ[y], z] + T[x, y, Δ[z]]

**T[E₁[1], E₁[1], E₂[1]] == ZERO**

True

**Leib[E₁[1], E₁[1], E₂[1]]**

$$\begin{pmatrix} \begin{pmatrix} 0 & \{0,0,0\} \\ \{0,0,0\} & 0 \end{pmatrix} & \begin{pmatrix} \alpha_1+\lambda_4 & \{\alpha_3+\lambda_6,\alpha_4+\lambda_7,\alpha_5+\lambda_8\} \\ \{\alpha_6+\lambda_9,\alpha_7+\lambda_{10},\alpha_8+\lambda_{11}\} & \alpha_2+\lambda_5 \end{pmatrix} & \begin{pmatrix} \\ \{ \end{pmatrix} \\ \begin{pmatrix} \alpha_2+\lambda_5 & \{-\alpha_3-\lambda_6,-\alpha_4-\lambda_7,-\alpha_5-\lambda_8\} \\ \{-\alpha_6-\lambda_9,-\alpha_7-\lambda_{10},-\alpha_8-\lambda_{11}\} & \alpha_1+\lambda_4 \end{pmatrix} & \begin{pmatrix} 0 & \{0,0,0\} \\ \{0,0,0\} & 0 \end{pmatrix} & \begin{pmatrix} \\ \{ \end{pmatrix} \\ \begin{pmatrix} 0 & \{0,0,0\} \\ \{0,0,0\} & 0 \end{pmatrix} & \begin{pmatrix} 0 & \{0,0,0\} \\ \{0,0,0\} & 0 \end{pmatrix} & \begin{pmatrix} \\ \{ \end{pmatrix} \end{pmatrix}$$

```
Δ[E₂[1]] = Δ[E₂[1]] //. {λ₄ → -α₁, λ₅ → -α₂, λ₆ → -α₃,
    λ₇ → -α₄, λ₈ → -α₅, λ₉ → -α₆, λ₁₀ → -α₇, λ₁₁ → -α₈}
```

$$\begin{pmatrix} \begin{pmatrix} 0 & \{0,0,0\} \\ \{0,0,0\} & 0 \end{pmatrix} & \begin{pmatrix} -\alpha_1 & \{-\alpha_3,-\alpha_4,-\alpha_5\} \\ \{-\alpha_6,-\alpha_7,-\alpha_8\} & -\alpha_2 \end{pmatrix} & \begin{pmatrix} 0 & \{0,0,0\} \\ \{0,0,0\} & 0 \end{pmatrix} \\ \begin{pmatrix} -\alpha_2 & \{\alpha_3,\alpha_4,\alpha_5\} \\ \{\alpha_6,\alpha_7,\alpha_8\} & -\alpha_1 \end{pmatrix} & \begin{pmatrix} 0 & \{0,0,0\} \\ \{0,0,0\} & 0 \end{pmatrix} & \begin{pmatrix} \lambda_{12} & \{\lambda_{14},\lambda_{15},\lambda_{16}\} \\ \{\lambda_{17},\lambda_{18},\lambda_{19}\} & \lambda_{13} \end{pmatrix} \\ \begin{pmatrix} 0 & \{0,0,0\} \\ \{0,0,0\} & 0 \end{pmatrix} & \begin{pmatrix} \lambda_{13} & \{-\lambda_{14},-\lambda_{15},-\lambda_{16}\} \\ \{-\lambda_{17},-\lambda_{18},-\lambda_{19}\} & \lambda_{12} \end{pmatrix} & \begin{pmatrix} 0 & \{0,0,0\} \\ \{0,0,0\} & 0 \end{pmatrix} \end{pmatrix}$$

**T[E₁[1], E₂[1], E₁[1]] == ZERO**

True

**Leib[E₁[1], E₂[1], E₁[1]] == ZERO**

True

**T[E₁[1], E₂[1], E₂[1]] == ZERO**

True

**Leib[E₁[1], E₂[1], E₂[1]] == ZERO**

True

**T[E₂[1], E₁[1], E₂[1]]**

$$\begin{pmatrix} \begin{pmatrix} 0 & \{0,0,0\} \\ \{0,0,0\} & 0 \end{pmatrix} & \begin{pmatrix} 0 & \{0,0,0\} \\ \{0,0,0\} & 0 \end{pmatrix} & \begin{pmatrix} 0 & \{0,0,0\} \\ \{0,0,0\} & 0 \end{pmatrix} \\ \begin{pmatrix} 0 & \{0,0,0\} \\ \{0,0,0\} & 0 \end{pmatrix} & \begin{pmatrix} 0 & \{0,0,0\} \\ \{0,0,0\} & 0 \end{pmatrix} & \begin{pmatrix} 0 & \{0,0,0\} \\ \{0,0,0\} & 0 \end{pmatrix} \\ \begin{pmatrix} 0 & \{0,0,0\} \\ \{0,0,0\} & 0 \end{pmatrix} & \begin{pmatrix} 0 & \{0,0,0\} \\ \{0,0,0\} & 0 \end{pmatrix} & \begin{pmatrix} 0 & \{0,0,0\} \\ \{0,0,0\} & 0 \end{pmatrix} \end{pmatrix}$$

**Leib[E₂[1], E₁[1], E₂[1]] == ZERO**

True

**T[E₂[1], E₂[1], E₂[1]] == 2 E₂[1]**

True

**Leib[E₂[1], E₂[1], E₂[1]] == 2 Δ[E₂[1]]**

True

**Δ[E₂[1]] = Δ[E₂[1]] //. {λ₁₂ → γ₁, λ₁₃ → γ₂, λ₁₄ → γ₃, λ₁₅ → γ₄, λ₁₆ → γ₅, λ₁₇ → γ₆, λ₁₈ → γ₇, λ₁₉ → γ₈}**

$$\begin{pmatrix} \begin{pmatrix} 0 & \{0,0,0\} \\ \{0,0,0\} & 0 \end{pmatrix} & \begin{pmatrix} -\alpha_1 & \{-\alpha_3, -\alpha_4, -\alpha_5\} \\ \{-\alpha_6, -\alpha_7, -\alpha_8\} & -\alpha_2 \end{pmatrix} & \begin{pmatrix} 0 & \{0,0,0\} \\ \{0,0,0\} & 0 \end{pmatrix} \\ \begin{pmatrix} -\alpha_2 & \{\alpha_3, \alpha_4, \alpha_5\} \\ \{\alpha_6, \alpha_7, \alpha_8\} & -\alpha_1 \end{pmatrix} & \begin{pmatrix} 0 & \{0,0,0\} \\ \{0,0,0\} & 0 \end{pmatrix} & \begin{pmatrix} \gamma_1 & \{\gamma_3, \gamma_4, \gamma_5\} \\ \{\gamma_6, \gamma_7, \gamma_8\} & \gamma_2 \end{pmatrix} \\ \begin{pmatrix} 0 & \{0,0,0\} \\ \{0,0,0\} & 0 \end{pmatrix} & \begin{pmatrix} \gamma_2 & \{-\gamma_3, -\gamma_4, -\gamma_5\} \\ \{-\gamma_6, -\gamma_7, -\gamma_8\} & \gamma_1 \end{pmatrix} & \begin{pmatrix} 0 & \{0,0,0\} \\ \{0,0,0\} & 0 \end{pmatrix} \end{pmatrix}$$

As $E_1[1] + E_2[1] + E_3[1] == 1$ and $D(1) = 0$,

```
Δ[E₃[1]] = -Δ[E₂[1]] - Δ[E₁[1]]
```

$$\left( \begin{array}{ccc} \begin{pmatrix} 0 & \{0,0,0\} \\ \{0,0,0\} & 0 \end{pmatrix} & \begin{pmatrix} 0 & \{0,0,0\} \\ \{0,0,0\} & 0 \end{pmatrix} & \begin{pmatrix} -\beta_1 & \{-\beta_3,-\beta_4,-\beta_5\} \\ \{-\beta_6,-\beta_7,-\beta_8\} & -\beta_2 \end{pmatrix} \\ \begin{pmatrix} 0 & \{0,0,0\} \\ \{0,0,0\} & 0 \end{pmatrix} & \begin{pmatrix} 0 & \{0,0,0\} \\ \{0,0,0\} & 0 \end{pmatrix} & \begin{pmatrix} -\gamma_1 & \{-\gamma_3,-\gamma_4,-\gamma_5\} \\ \{-\gamma_6,-\gamma_7,-\gamma_8\} & -\gamma_2 \end{pmatrix} \\ \begin{pmatrix} -\beta_2 & \{\beta_3,\beta_4,\beta_5\} \\ \{\beta_6,\beta_7,\beta_8\} & -\beta_1 \end{pmatrix} & \begin{pmatrix} -\gamma_2 & \{\gamma_3,\gamma_4,\gamma_5\} \\ \{\gamma_6,\gamma_7,\gamma_8\} & -\gamma_1 \end{pmatrix} & \begin{pmatrix} 0 & \{0,0,0\} \\ \{0,0,0\} & 0 \end{pmatrix} \end{array} \right)$$

```
U_{E₃[1]}[E₁[1]] == ZERO
```

True

```
T[Δ[E₃[1]], E₁[1], E₃[1]] + U_{E₃[1]}[Δ[E₁[1]]] == ZERO
```

True

```
U_{E₃[1]}[E₂[1]] == ZERO
```

True

```
T[Δ[E₃[1]], E₂[1], E₃[1]] + U_{E₃[1]}[Δ[E₂[1]]] == ZERO
```

True

```
(Δ[E₃[1]] - T[Δ[E₃[1]], ONE, E₃[1]] // Expand) == ZERO
```

True

## ■ The images of $X_i[e_j]$.

```
Leib[x_, y_, z_] := T[Δ[x], y, z] + T[x, Δ[y], z] + T[x, y, Δ[z]]
```

$$O_1 = \begin{pmatrix} \alpha_1 & \{\alpha_3, \alpha_4, \alpha_5\} \\ \{\alpha_6, \alpha_7, \alpha_8\} & \alpha_2 \end{pmatrix};$$

$$O_2 = \begin{pmatrix} \beta_1 & \{\beta_3, \beta_4, \beta_5\} \\ \{\beta_6, \beta_7, \beta_8\} & \beta_2 \end{pmatrix};$$

$$O_3 = \begin{pmatrix} \gamma_1 & \{\gamma_3, \gamma_4, \gamma_5\} \\ \{\gamma_6, \gamma_7, \gamma_8\} & \gamma_2 \end{pmatrix};$$

$$\Delta[E_1[1]] = \begin{pmatrix} \text{zero} & O_1 & O_2 \\ \sigma[O_1] & \text{zero} & \text{zero} \\ \sigma[O_2] & \text{zero} & \text{zero} \end{pmatrix};$$

$$\Delta[E_2[1]] = \begin{pmatrix} \text{zero} & -O_1 & \text{zero} \\ -\sigma[O_1] & \text{zero} & O_3 \\ \text{zero} & \sigma[O_3] & \text{zero} \end{pmatrix};$$

```
Δ[E₃[1]] = -Δ[E₁[1]] - Δ[E₂[1]];

X₁[ξ_] := X_{1,2}[ξ]; X₂[ξ_] := X_{2,3}[ξ]; X₃[ξ_] := X_{1,3}[ξ];
```

- $X_1[e_1]$

    ```
    Δ[X₁[e₁]] = generic;
    ```

    ```
    U_{X₁[e₁]}[E₁[1]] == ZERO
    ```

    True

    ```
    defin[x_, y_] := T[Δ[x], y, x] + U_x[Δ[y]]
    ```

    ```
    defin[X₁[e₁], E₁[1]] // Expand
    ```

    $$\left( \begin{array}{ccc} \begin{pmatrix} 0 & \{0,0,0\} \\ \{0,0,0\} & 0 \end{pmatrix} & \begin{pmatrix} \alpha_2+\lambda_1 & \{0,0,0\} \\ \{0,0,0\} & 0 \end{pmatrix} & \begin{pmatrix} 0 & \{0,0,0\} \\ \{0,0,0\} & 0 \end{pmatrix} \\ \begin{pmatrix} 0 & \{0,0,0\} \\ \{0,0,0\} & \alpha_2+\lambda_1 \end{pmatrix} & \begin{pmatrix} \lambda_5 & \{0,0,0\} \\ \{0,0,0\} & \lambda_5 \end{pmatrix} & \begin{pmatrix} 0 & \{0,0,0\} \\ \{\lambda_{25},\lambda_{26},\lambda_{27}\} & \lambda_{21} \end{pmatrix} \\ \begin{pmatrix} 0 & \{0,0,0\} \\ \{0,0,0\} & 0 \end{pmatrix} & \begin{pmatrix} \lambda_{21} & \{0,0,0\} \\ \{-\lambda_{25},-\lambda_{26},-\lambda_{27}\} & 0 \end{pmatrix} & \begin{pmatrix} 0 & \{0,0,0\} \\ \{0,0,0\} & 0 \end{pmatrix} \end{array} \right)$$

    ```
    Δ[X₁[e₁]] = Δ[X₁[e₁]] //.
      {λ₁ -> -α₂, λ₅ -> 0, λ₂₁ -> 0, λ₂₅ -> 0, λ₂₆ -> 0, λ₂₇ -> 0}
    ```

    $$\left( \begin{array}{ccc} \begin{pmatrix} -\alpha_2 & \{0,0,0\} \\ \{0,0,0\} & -\alpha_2 \end{pmatrix} & \begin{pmatrix} \lambda_4 & \{\lambda_6,\lambda_7,\lambda_8\} \\ \{\lambda_9,\lambda_{10},\lambda_{11}\} & 0 \end{pmatrix} & \begin{pmatrix} \lambda_{20} & \{\lambda_{22},\lambda_{23},\lambda_{24}\} \\ \{0,0,0\} & 0 \end{pmatrix} \\ \begin{pmatrix} 0 & \{-\lambda_6,-\lambda_7,-\lambda_8\} \\ \{-\lambda_9,-\lambda_{10},-\lambda_{11}\} & \lambda_4 \end{pmatrix} & \begin{pmatrix} \lambda_2 & \{0,0,0\} \\ \{0,0,0\} & \lambda_2 \end{pmatrix} & \begin{pmatrix} \lambda_{12} & \{\lambda_{14},\lambda_{15},\lambda_{16}\} \\ \{\lambda_{17},\lambda_{18},\lambda_{19}\} & \lambda_{13} \end{pmatrix} \\ \begin{pmatrix} 0 & \{-\lambda_{22},-\lambda_{23},-\lambda_{24}\} \\ \{0,0,0\} & \lambda_{20} \end{pmatrix} & \begin{pmatrix} \lambda_{13} & \{-\lambda_{14},-\lambda_{15},-\lambda_{16}\} \\ \{-\lambda_{17},-\lambda_{18},-\lambda_{19}\} & \lambda_{12} \end{pmatrix} & \begin{pmatrix} \lambda_3 & \{0,0,0\} \\ \{0,0,0\} & \lambda_3 \end{pmatrix} \end{array} \right)$$

    ```
    U_{X₁[e₁]}[E₂[1]] == ZERO
    ```

    True

    ```
    defin[X₁[e₁], E₂[1]] // Expand
    ```

    $$\left( \begin{array}{ccc} \begin{pmatrix} 0 & \{0,0,0\} \\ \{0,0,0\} & 0 \end{pmatrix} & \begin{pmatrix} \lambda_2-\alpha_2 & \{0,0,0\} \\ \{0,0,0\} & 0 \end{pmatrix} & \begin{pmatrix} \lambda_{12} & \{\lambda_{14},\lambda_{15},\lambda_{16}\} \\ \{0,0,0\} & 0 \end{pmatrix} \\ \begin{pmatrix} 0 & \{0,0,0\} \\ \{0,0,0\} & \lambda_2-\alpha_2 \end{pmatrix} & \begin{pmatrix} 0 & \{0,0,0\} \\ \{0,0,0\} & 0 \end{pmatrix} & \begin{pmatrix} 0 & \{0,0,0\} \\ \{0,0,0\} & 0 \end{pmatrix} \\ \begin{pmatrix} 0 & \{-\lambda_{14},-\lambda_{15},-\lambda_{16}\} \\ \{0,0,0\} & \lambda_{12} \end{pmatrix} & \begin{pmatrix} 0 & \{0,0,0\} \\ \{0,0,0\} & 0 \end{pmatrix} & \begin{pmatrix} 0 & \{0,0,0\} \\ \{0,0,0\} & 0 \end{pmatrix} \end{array} \right)$$

    ```
    Δ[X₁[e₁]] =
     Δ[X₁[e₁]] //. {λ₂ -> α₂, λ₁₂ -> 0, λ₁₄ -> 0, λ₁₅ -> 0, λ₁₆ -> 0}
    ```

    $$\left( \begin{array}{ccc} \begin{pmatrix} -\alpha_2 & \{0,0,0\} \\ \{0,0,0\} & -\alpha_2 \end{pmatrix} & \begin{pmatrix} \lambda_4 & \{\lambda_6,\lambda_7,\lambda_8\} \\ \{\lambda_9,\lambda_{10},\lambda_{11}\} & 0 \end{pmatrix} & \begin{pmatrix} \lambda_{20} & \{\lambda_{22},\lambda_{23},\lambda_{24}\} \\ \{0,0,0\} & 0 \end{pmatrix} \\ \begin{pmatrix} 0 & \{-\lambda_6,-\lambda_7,-\lambda_8\} \\ \{-\lambda_9,-\lambda_{10},-\lambda_{11}\} & \lambda_4 \end{pmatrix} & \begin{pmatrix} \alpha_2 & \{0,0,0\} \\ \{0,0,0\} & \alpha_2 \end{pmatrix} & \begin{pmatrix} 0 & \{0,0,0\} \\ \{\lambda_{17},\lambda_{18},\lambda_{19}\} & \lambda_{13} \end{pmatrix} \\ \begin{pmatrix} 0 & \{-\lambda_{22},-\lambda_{23},-\lambda_{24}\} \\ \{0,0,0\} & \lambda_{20} \end{pmatrix} & \begin{pmatrix} \lambda_{13} & \{0,0,0\} \\ \{-\lambda_{17},-\lambda_{18},-\lambda_{19}\} & 0 \end{pmatrix} & \begin{pmatrix} \lambda_3 & \{0,0,0\} \\ \{0,0,0\} & \lambda_3 \end{pmatrix} \end{array} \right)$$

$U_{X_1[e_1]}[E_3[1]] == \text{ZERO}$

True

$(\text{defin}[X_1[e_1], E_3[1]] \;//\; \text{Expand}) == \text{ZERO}$

True

$U_{X_1[e_1]}[X_1[e_1]] == \text{ZERO}$

True

$(\text{defin}[X_1[e_1], X_1[e_1]] \;//\; \text{Expand}) == \text{ZERO}$

True

$T[E_1[1], E_1[1], X_1[e_1]] == X_1[e_1]$

True

$\Delta[X_1[e_1]] - \text{Leib}[E_1[1], E_1[1], X_1[e_1]] \;//\; \text{Expand}$

$$\begin{pmatrix}
\begin{pmatrix} 0 & \{0,0,0\} \\ \{0,0,0\} & 0 \end{pmatrix} & \begin{pmatrix} 0 & \{0,0,0\} \\ \{0,0,0\} & 0 \end{pmatrix} & \begin{pmatrix} 0 & \{0,0,0\} \\ \{0,0,0\} & 0 \end{pmatrix} \\
\begin{pmatrix} 0 & \{0,0,0\} \\ \{0,0,0\} & 0 \end{pmatrix} & \begin{pmatrix} 0 & \{0,0,0\} \\ \{0,0,0\} & 0 \end{pmatrix} & \begin{pmatrix} 0 & \{0,0,0\} \\ \{\lambda_{17}-\beta_6, \lambda_{18}-\beta_7, \lambda_{19}-\beta_8\} & \lambda_{13}-\beta_2 \end{pmatrix} \\
\begin{pmatrix} 0 & \{0,0,0\} \\ \{0,0,0\} & 0 \end{pmatrix} & \begin{pmatrix} \lambda_{13}-\beta_2 & \{0,0,0\} \\ \{\beta_6-\lambda_{17}, \beta_7-\lambda_{18}, \beta_8-\lambda_{19}\} & 0 \end{pmatrix} & \begin{pmatrix} \lambda_3 & \{0,0,0\} \\ \{0,0,0\} & \lambda_3 \end{pmatrix}
\end{pmatrix}$$

$\Delta[X_1[e_1]] =$
$\quad \Delta[X_1[e_1]] \;//.\; \{\lambda_3 \to 0, \lambda_{13} \to \beta_2, \lambda_{17} \to \beta_6, \lambda_{18} \to \beta_7, \lambda_{19} \to \beta_8\}$

$$\begin{pmatrix}
\begin{pmatrix} -\alpha_2 & \{0,0,0\} \\ \{0,0,0\} & -\alpha_2 \end{pmatrix} & \begin{pmatrix} \lambda_4 & \{\lambda_6, \lambda_7, \lambda_8\} \\ \{\lambda_9, \lambda_{10}, \lambda_{11}\} & 0 \end{pmatrix} & \begin{pmatrix} \lambda_{20} & \{\lambda_{22}, \lambda_{23}, \lambda_{24}\} \\ \{0,0,0\} & 0 \end{pmatrix} \\
\begin{pmatrix} 0 & \{-\lambda_6, -\lambda_7, -\lambda_8\} \\ \{-\lambda_9, -\lambda_{10}, -\lambda_{11}\} & \lambda_4 \end{pmatrix} & \begin{pmatrix} \alpha_2 & \{0,0,0\} \\ \{0,0,0\} & \alpha_2 \end{pmatrix} & \begin{pmatrix} 0 & \{0,0,0\} \\ \{\beta_6, \beta_7, \beta_8\} & \beta_2 \end{pmatrix} \\
\begin{pmatrix} 0 & \{-\lambda_{22}, -\lambda_{23}, -\lambda_{24}\} \\ \{0,0,0\} & \lambda_{20} \end{pmatrix} & \begin{pmatrix} \beta_2 & \{0,0,0\} \\ \{-\beta_6, -\beta_7, -\beta_8\} & 0 \end{pmatrix} & \begin{pmatrix} 0 & \{0,0,0\} \\ \{0,0,0\} & 0 \end{pmatrix}
\end{pmatrix}$$

$T[E_1[1], X_1[e_1], E_1[1]] == \text{ZERO}$

True

$\text{Leib}[E_1[1], X_1[e_1], E_1[1]] == \text{ZERO}$

True

$T[E_2[1], E_2[1], X_1[e_1]] == X_1[e_1]$

True

**Δ[X₁[e₁]] - Leib[E₂[1], E₂[1], X₁[e₁]] // Expand**

$$\begin{pmatrix} \begin{pmatrix} 0 & \{0,0,0\} \\ \{0,0,0\} & 0 \end{pmatrix} & \begin{pmatrix} 0 & \{0,0,0\} \\ \{0,0,0\} & 0 \end{pmatrix} & \begin{pmatrix} \lambda_{20}-\gamma_1 & \{\lambda_{22}-\gamma_3, \lambda_{23}-\gamma_4, \lambda_{24}-\gamma_5\} \\ \{0,0,0\} & 0 \end{pmatrix} \\ \begin{pmatrix} 0 & \{0,0,0\} \\ \{0,0,0\} & 0 \end{pmatrix} & \begin{pmatrix} 0 & \{0,0,0\} \\ \{0,0,0\} & 0 \end{pmatrix} & \begin{pmatrix} 0 & \{0,0,0\} \\ \{0,0,0\} & 0 \end{pmatrix} \\ \begin{pmatrix} 0 & \{\gamma_3-\lambda_{22}, \gamma_4-\lambda_{23}, \gamma_5-\lambda_{24}\} \\ \{0,0,0\} & \lambda_{20}-\gamma_1 \end{pmatrix} & \begin{pmatrix} 0 & \{0,0,0\} \\ \{0,0,0\} & 0 \end{pmatrix} & \begin{pmatrix} 0 & \{0,0,0\} \\ \{0,0,0\} & 0 \end{pmatrix} \end{pmatrix}$$

**Δ[X₁[e₁]] = Δ[X₁[e₁]] //. {λ₂₀ → γ₁, λ₂₂ → γ₃, λ₂₃ → γ₄, λ₂₄ → γ₅}**

$$\begin{pmatrix} \begin{pmatrix} -\alpha_2 & \{0,0,0\} \\ \{0,0,0\} & -\alpha_2 \end{pmatrix} & \begin{pmatrix} \lambda_4 & \{\lambda_6, \lambda_7, \lambda_8\} \\ \{\lambda_9, \lambda_{10}, \lambda_{11}\} & 0 \end{pmatrix} & \begin{pmatrix} \gamma_1 & \{\gamma_3, \gamma_4, \gamma_5\} \\ \{0,0,0\} & 0 \end{pmatrix} \\ \begin{pmatrix} 0 & \{-\lambda_6, -\lambda_7, -\lambda_8\} \\ \{-\lambda_9, -\lambda_{10}, -\lambda_{11}\} & \lambda_4 \end{pmatrix} & \begin{pmatrix} \alpha_2 & \{0,0,0\} \\ \{0,0,0\} & \alpha_2 \end{pmatrix} & \begin{pmatrix} 0 & \{0,0,0\} \\ \{\beta_6, \beta_7, \beta_8\} & \beta_2 \end{pmatrix} \\ \begin{pmatrix} 0 & \{-\gamma_3, -\gamma_4, -\gamma_5\} \\ \{0,0,0\} & \gamma_1 \end{pmatrix} & \begin{pmatrix} \beta_2 & \{0,0,0\} \\ \{-\beta_6, -\beta_7, -\beta_8\} & 0 \end{pmatrix} & \begin{pmatrix} 0 & \{0,0,0\} \\ \{0,0,0\} & 0 \end{pmatrix} \end{pmatrix}$$

**T[E₂[1], X₁[e₁], E₂[1]] == ZERO**

True

**Leib[E₂[1], X₁[e₁], E₂[1]] == ZERO**

True

**T[E₁[1], E₂[1], X₁[e₁]] == ZERO**

True

**(Leib[E₁[1], E₂[1], X₁[e₁]] // Expand) == ZERO**

True

**T[E₂[1], E₁[1], X₁[e₁]] == ZERO**

True

**(Leib[E₂[1], E₁[1], X₁[e₁]] // Expand) == ZERO**

True

**T[E₁[1], X₁[e₁], E₂[1]] == X₁[e₁]**

True

**Δ[X₁[e₁]] - Leib[E₁[1], X₁[e₁], E₂[1]] == ZERO**

True

```
T[X₁[e₁], X₁[e₁], X₁[e₁]] == ZERO
```

True

```
Leib[X₁[e₁], X₁[e₁], X₁[e₁]] == ZERO
```

True

```
T[X₁[e₁], E₁[1], X₁[e₁]] == ZERO
```

True

```
Leib[X₁[e₁], E₁[1], X₁[e₁]] == ZERO
```

True

```
T[X₁[e₁], E₂[1], X₁[e₁]] == ZERO
```

True

```
Leib[X₁[e₁], E₂[1], X₁[e₁]] == ZERO
```

True

```
T[E₁[1], X₁[e₁], X₁[e₁]] == ZERO
```

True

```
Leib[E₁[1], X₁[e₁], X₁[e₁]] == ZERO
```

True

```
T[E₂[1], X₁[e₁], X₁[e₁]] == ZERO
```

True

```
Leib[E₂[1], X₁[e₁], X₁[e₁]] == ZERO
```

True

```
Δ[X₁[e₁]] = Δ[X₁[e₁]] //.
  {λ₄ → δ₁, λ₆ → δ₂, λ₇ → δ₃, λ₈ → δ₄, λ₉ → δ₅, λ₁₀ → δ₆, λ₁₁ → δ₇}
```

$$\begin{pmatrix} \begin{pmatrix} -\alpha_2 & \{0,0,0\} \\ \{0,0,0\} & -\alpha_2 \end{pmatrix} & \begin{pmatrix} \delta_1 & \{\delta_2,\delta_3,\delta_4\} \\ \{\delta_5,\delta_6,\delta_7\} & 0 \end{pmatrix} & \begin{pmatrix} \gamma_1 & \{\gamma_3,\gamma_4,\gamma_5\} \\ \{0,0,0\} & 0 \end{pmatrix} \\ \begin{pmatrix} 0 & \{-\delta_2,-\delta_3,-\delta_4\} \\ \{-\delta_5,-\delta_6,-\delta_7\} & \delta_1 \end{pmatrix} & \begin{pmatrix} \alpha_2 & \{0,0,0\} \\ \{0,0,0\} & \alpha_2 \end{pmatrix} & \begin{pmatrix} 0 & \{0,0,0\} \\ \{\beta_6,\beta_7,\beta_8\} & \beta_2 \end{pmatrix} \\ \begin{pmatrix} 0 & \{-\gamma_3,-\gamma_4,-\gamma_5\} \\ \{0,0,0\} & \gamma_1 \end{pmatrix} & \begin{pmatrix} \beta_2 & \{0,0,0\} \\ \{-\beta_6,-\beta_7,-\beta_8\} & 0 \end{pmatrix} & \begin{pmatrix} 0 & \{0,0,0\} \\ \{0,0,0\} & 0 \end{pmatrix} \end{pmatrix}$$

- **X₁[e₂]**

    **Δ[X₁[e₂]] = generic;**

    **U_{X₁[e₂]}[E₁[1]] == ZERO**

    True

    **defin[X₁[e₂], E₁[1]] // Expand**

    $$\left(\begin{array}{ccc}
    \begin{pmatrix} 0 & \{0,0,0\} \\ \{0,0,0\} & 0 \end{pmatrix} &
    \begin{pmatrix} 0 & \{0,0,0\} \\ \{0,0,0\} & \alpha_1+\lambda_1 \end{pmatrix} &
    \begin{pmatrix} 0 & \{0,0,0\} \\ \{0,0,0\} & 0 \end{pmatrix} \\
    \begin{pmatrix} \alpha_1+\lambda_1 & \{0,0,0\} \\ \{0,0,0\} & 0 \end{pmatrix} &
    \begin{pmatrix} \lambda_4 & \{0,0,0\} \\ \{0,0,0\} & \lambda_4 \end{pmatrix} &
    \begin{pmatrix} \lambda_{20} & \{\lambda_{22},\lambda_{23},\lambda_{24}\} \\ \{0,0,0\} & 0 \end{pmatrix} \\
    \begin{pmatrix} 0 & \{0,0,0\} \\ \{0,0,0\} & 0 \end{pmatrix} &
    \begin{pmatrix} 0 & \{-\lambda_{22},-\lambda_{23},-\lambda_{24}\} \\ \{0,0,0\} & \lambda_{20} \end{pmatrix} &
    \begin{pmatrix} 0 & \{0,0,0\} \\ \{0,0,0\} & 0 \end{pmatrix}
    \end{array}\right)$$

    **Δ[X₁[e₂]] = Δ[X₁[e₂]] //. {λ₁ -> -α₁, λ₄ -> 0, λ₂₀ -> 0, λ₂₂ -> 0, λ₂₃ -> 0, λ₂₄ -> 0}**

    $$\left(\begin{array}{ccc}
    \begin{pmatrix} -\alpha_1 & \{0,0,0\} \\ \{0,0,0\} & -\alpha_1 \end{pmatrix} &
    \begin{pmatrix} 0 & \{\lambda_6,\lambda_7,\lambda_8\} \\ \{\lambda_9,\lambda_{10},\lambda_{11}\} & \lambda_5 \end{pmatrix} &
    \begin{pmatrix} 0 & \{\cdots \\ \{\lambda_{25},\lambda_{26},\lambda_{27}\} & \cdots \end{pmatrix} \\
    \begin{pmatrix} \lambda_5 & \{-\lambda_6,-\lambda_7,-\lambda_8\} \\ \{-\lambda_9,-\lambda_{10},-\lambda_{11}\} & 0 \end{pmatrix} &
    \begin{pmatrix} \lambda_2 & \{0,0,0\} \\ \{0,0,0\} & \lambda_2 \end{pmatrix} &
    \begin{pmatrix} \lambda_{12} & \{\lambda_{14},\cdots \\ \{\lambda_{17},\lambda_{18},\lambda_{19}\} & \cdots \end{pmatrix} \\
    \begin{pmatrix} \lambda_{21} & \{0,0,0\} \\ \{-\lambda_{25},-\lambda_{26},-\lambda_{27}\} & 0 \end{pmatrix} &
    \begin{pmatrix} \lambda_{13} & \{-\lambda_{14},-\lambda_{15},-\lambda_{16}\} \\ \{-\lambda_{17},-\lambda_{18},-\lambda_{19}\} & \lambda_{12} \end{pmatrix} &
    \begin{pmatrix} \lambda_3 & \{0,\cdots \\ \{0,0,0\} & \lambda_3 \end{pmatrix}
    \end{array}\right)$$

    **U_{X₁[e₂]}[E₂[1]] == ZERO**

    True

    **defin[X₁[e₂], E₂[1]] // Expand**

    $$\left(\begin{array}{ccc}
    \begin{pmatrix} 0 & \{0,0,0\} \\ \{0,0,0\} & 0 \end{pmatrix} &
    \begin{pmatrix} 0 & \{0,0,0\} \\ \{0,0,0\} & \lambda_2-\alpha_1 \end{pmatrix} &
    \begin{pmatrix} 0 & \{0,0,0\} \\ \{\lambda_{17},\lambda_{18},\lambda_{19}\} & \lambda_{13} \end{pmatrix} \\
    \begin{pmatrix} \lambda_2-\alpha_1 & \{0,0,0\} \\ \{0,0,0\} & 0 \end{pmatrix} &
    \begin{pmatrix} 0 & \{0,0,0\} \\ \{0,0,0\} & 0 \end{pmatrix} &
    \begin{pmatrix} 0 & \{0,0,0\} \\ \{0,0,0\} & 0 \end{pmatrix} \\
    \begin{pmatrix} \lambda_{13} & \{0,0,0\} \\ \{-\lambda_{17},-\lambda_{18},-\lambda_{19}\} & 0 \end{pmatrix} &
    \begin{pmatrix} 0 & \{0,0,0\} \\ \{0,0,0\} & 0 \end{pmatrix} &
    \begin{pmatrix} 0 & \{0,0,0\} \\ \{0,0,0\} & 0 \end{pmatrix}
    \end{array}\right)$$

    **Δ[X₁[e₂]] = Δ[X₁[e₂]] //. {λ₂ -> α₁, λ₁₃ -> 0, λ₁₇ -> 0, λ₁₈ -> 0, λ₁₉ -> 0}**

    $$\left(\begin{array}{ccc}
    \begin{pmatrix} -\alpha_1 & \{0,0,0\} \\ \{0,0,0\} & -\alpha_1 \end{pmatrix} &
    \begin{pmatrix} 0 & \{\lambda_6,\lambda_7,\lambda_8\} \\ \{\lambda_9,\lambda_{10},\lambda_{11}\} & \lambda_5 \end{pmatrix} &
    \begin{pmatrix} 0 & \{0,0,0\} \\ \{\lambda_{25},\lambda_{26},\lambda_{27}\} & \lambda_{21} \end{pmatrix} \\
    \begin{pmatrix} \lambda_5 & \{-\lambda_6,-\lambda_7,-\lambda_8\} \\ \{-\lambda_9,-\lambda_{10},-\lambda_{11}\} & 0 \end{pmatrix} &
    \begin{pmatrix} \alpha_1 & \{0,0,0\} \\ \{0,0,0\} & \alpha_1 \end{pmatrix} &
    \begin{pmatrix} \lambda_{12} & \{\lambda_{14},\lambda_{15},\lambda_{16}\} \\ \{0,0,0\} & 0 \end{pmatrix} \\
    \begin{pmatrix} \lambda_{21} & \{0,0,0\} \\ \{-\lambda_{25},-\lambda_{26},-\lambda_{27}\} & 0 \end{pmatrix} &
    \begin{pmatrix} 0 & \{-\lambda_{14},-\lambda_{15},-\lambda_{16}\} \\ \{0,0,0\} & \lambda_{12} \end{pmatrix} &
    \begin{pmatrix} \lambda_3 & \{0,0,0\} \\ \{0,0,0\} & \lambda_3 \end{pmatrix}
    \end{array}\right)$$

$U_{X_1[e_2]}[E_3[1]] ==$ ZERO

True

$(\text{defin}[X_1[e_2], E_3[1]] \ // \ \text{Expand}) ==$ ZERO

True

$U_{X_1[e_2]}[X_1[e_2]] ==$ ZERO

True

$(\text{defin}[X_1[e_2], X_1[e_2]] \ // \ \text{Expand}) ==$ ZERO

True

$U_{X_1[e_2]}[X_1[e_1]] == X_1[e_2]$

True

$(\Delta[X_1[e_2]] - \text{defin}[X_1[e_2], X_1[e_1]] \ // \ \text{Expand})$

$$\begin{pmatrix} \begin{pmatrix} 0 & \{0,0,0\} \\ \{0,0,0\} & 0 \end{pmatrix} & \begin{pmatrix} 0 & \{0,0,0\} \\ \{0,0,0\} & -\delta_1-\lambda_5 \end{pmatrix} & \begin{pmatrix} 0 & \{0,0,0\} \\ \{0,0,0\} & 0 \end{pmatrix} \\ \begin{pmatrix} -\delta_1-\lambda_5 & \{0,0,0\} \\ \{0,0,0\} & 0 \end{pmatrix} & \begin{pmatrix} 0 & \{0,0,0\} \\ \{0,0,0\} & 0 \end{pmatrix} & \begin{pmatrix} 0 & \{0,0,0\} \\ \{0,0,0\} & 0 \end{pmatrix} \\ \begin{pmatrix} 0 & \{0,0,0\} \\ \{0,0,0\} & 0 \end{pmatrix} & \begin{pmatrix} 0 & \{0,0,0\} \\ \{0,0,0\} & 0 \end{pmatrix} & \begin{pmatrix} \lambda_3 & \{0,0,0\} \\ \{0,0,0\} & \lambda_3 \end{pmatrix} \end{pmatrix}$$

$\Delta[X_1[e_2]] = \Delta[X_1[e_2]] \ //. \ \{\lambda_3 \to 0, \lambda_5 \to -\delta_1\}$

$$\begin{pmatrix} \begin{pmatrix} -\alpha_1 & \{0,0,0\} \\ \{0,0,0\} & -\alpha_1 \end{pmatrix} & \begin{pmatrix} 0 & \{\lambda_6,\lambda_7,\lambda_8\} \\ \{\lambda_9,\lambda_{10},\lambda_{11}\} & -\delta_1 \end{pmatrix} & \begin{pmatrix} 0 & \{0,0,0\} \\ \{\lambda_{25},\lambda_{26},\lambda_{27}\} & \lambda_{21} \end{pmatrix} \\ \begin{pmatrix} -\delta_1 & \{-\lambda_6,-\lambda_7,-\lambda_8\} \\ \{-\lambda_9,-\lambda_{10},-\lambda_{11}\} & 0 \end{pmatrix} & \begin{pmatrix} \alpha_1 & \{0,0,0\} \\ \{0,0,0\} & \alpha_1 \end{pmatrix} & \begin{pmatrix} \lambda_{12} & \{\lambda_{14},\lambda_{15},\lambda_{16}\} \\ \{0,0,0\} & 0 \end{pmatrix} \\ \begin{pmatrix} \lambda_{21} & \{0,0,0\} \\ \{-\lambda_{25},-\lambda_{26},-\lambda_{27}\} & 0 \end{pmatrix} & \begin{pmatrix} 0 & \{-\lambda_{14},-\lambda_{15},-\lambda_{16}\} \\ \{0,0,0\} & \lambda_{12} \end{pmatrix} & \begin{pmatrix} 0 & \{0,0,0\} \\ \{0,0,0\} & 0 \end{pmatrix} \end{pmatrix}$$

$T[E_1[1], E_1[1], X_1[e_2]] == X_1[e_2]$

True

$\Delta[X_1[e_2]] - \text{Leib}[E_1[1], E_1[1], X_1[e_2]] \ // \ \text{Expand}$

$$\begin{pmatrix} \begin{pmatrix} 0 & \{0,0,0\} \\ \{0,0,0\} & 0 \end{pmatrix} & \begin{pmatrix} 0 & \{0,0,0\} \\ \{0,0,0\} & 0 \end{pmatrix} & \begin{pmatrix} 0 & \{0,0,0\} \\ \{0,0,0\} & 0 \end{pmatrix} \\ \begin{pmatrix} 0 & \{0,0,0\} \\ \{0,0,0\} & 0 \end{pmatrix} & \begin{pmatrix} 0 & \{0,0,0\} \\ \{0,0,0\} & 0 \end{pmatrix} & \begin{pmatrix} \lambda_{12}-\beta_1 & \{\lambda_{14}-\beta_3,\lambda_{15}-\beta_4,\lambda_{16}-\beta_5\} \\ \{0,0,0\} & 0 \end{pmatrix} \\ \begin{pmatrix} 0 & \{0,0,0\} \\ \{0,0,0\} & 0 \end{pmatrix} & \begin{pmatrix} 0 & \{\beta_3-\lambda_{14},\beta_4-\lambda_{15},\beta_5-\lambda_{16}\} \\ \{0,0,0\} & \lambda_{12}-\beta_1 \end{pmatrix} & \begin{pmatrix} 0 & \{0,0,0\} \\ \{0,0,0\} & 0 \end{pmatrix} \end{pmatrix}$$

**Δ[X₁[e₂]] = Δ[X₁[e₂]] //. {λ₁₂ → β₁, λ₁₄ → β₃, λ₁₅ → β₄, λ₁₆ → β₅}**

$$\begin{pmatrix} \begin{pmatrix} -\alpha_1 & \{0,0,0\} \\ \{0,0,0\} & -\alpha_1 \end{pmatrix} & \begin{pmatrix} 0 & \{\lambda_6, \lambda_7, \lambda_8\} \\ \{\lambda_9, \lambda_{10}, \lambda_{11}\} & -\delta_1 \end{pmatrix} & \begin{pmatrix} 0 & \{0,0,0\} \\ \{\lambda_{25}, \lambda_{26}, \lambda_{27}\} & \lambda_{21} \end{pmatrix} \\ \begin{pmatrix} -\delta_1 & \{-\lambda_6, -\lambda_7, -\lambda_8\} \\ \{-\lambda_9, -\lambda_{10}, -\lambda_{11}\} & 0 \end{pmatrix} & \begin{pmatrix} \alpha_1 & \{0,0,0\} \\ \{0,0,0\} & \alpha_1 \end{pmatrix} & \begin{pmatrix} \beta_1 & \{\beta_3, \beta_4, \beta_5\} \\ \{0,0,0\} & 0 \end{pmatrix} \\ \begin{pmatrix} \lambda_{21} & \{0,0,0\} \\ \{-\lambda_{25}, -\lambda_{26}, -\lambda_{27}\} & 0 \end{pmatrix} & \begin{pmatrix} 0 & \{-\beta_3, -\beta_4, -\beta_5\} \\ \{0,0,0\} & \beta_1 \end{pmatrix} & \begin{pmatrix} 0 & \{0,0,0\} \\ \{0,0,0\} & 0 \end{pmatrix} \end{pmatrix}$$

**T[E₁[1], X₁[e₂], E₁[1]] == ZERO**

True

**Leib[E₁[1], X₁[e₂], E₁[1]] == ZERO**

True

**T[E₂[1], E₂[1], X₁[e₂]] == X₁[e₂]**

True

**Δ[X₁[e₂]] - Leib[E₂[1], E₂[1], X₁[e₂]] // Expand**

$$\begin{pmatrix} \begin{pmatrix} 0 & \{0,0,0\} \\ \{0,0,0\} & 0 \end{pmatrix} & \begin{pmatrix} 0 & \{0,0,0\} \\ \{0,0,0\} & 0 \end{pmatrix} & \begin{pmatrix} 0 & \{0,0,0\} \\ \{\lambda_{25}-\gamma_6, \lambda_{26}-\gamma_7, \lambda_{27}-\gamma_8\} & \lambda_{21}-\gamma_2 \end{pmatrix} \\ \begin{pmatrix} 0 & \{0,0,0\} \\ \{0,0,0\} & 0 \end{pmatrix} & \begin{pmatrix} 0 & \{0,0,0\} \\ \{0,0,0\} & 0 \end{pmatrix} & \begin{pmatrix} 0 & \{0,0,0\} \\ \{0,0,0\} & 0 \end{pmatrix} \\ \begin{pmatrix} \lambda_{21}-\gamma_2 & \{0,0,0\} \\ \{\gamma_6-\lambda_{25}, \gamma_7-\lambda_{26}, \gamma_8-\lambda_{27}\} & 0 \end{pmatrix} & \begin{pmatrix} 0 & \{0,0,0\} \\ \{0,0,0\} & 0 \end{pmatrix} & \begin{pmatrix} 0 & \{0,0,0\} \\ \{0,0,0\} & 0 \end{pmatrix} \end{pmatrix}$$

**Δ[X₁[e₂]] = Δ[X₁[e₂]] //. {λ₂₅ → γ₆, λ₂₆ → γ₇, λ₂₇ → γ₈, λ₂₁ → γ₂}**

$$\begin{pmatrix} \begin{pmatrix} -\alpha_1 & \{0,0,0\} \\ \{0,0,0\} & -\alpha_1 \end{pmatrix} & \begin{pmatrix} 0 & \{\lambda_6, \lambda_7, \lambda_8\} \\ \{\lambda_9, \lambda_{10}, \lambda_{11}\} & -\delta_1 \end{pmatrix} & \begin{pmatrix} 0 & \{0,0,0\} \\ \{\gamma_6, \gamma_7, \gamma_8\} & \gamma_2 \end{pmatrix} \\ \begin{pmatrix} -\delta_1 & \{-\lambda_6, -\lambda_7, -\lambda_8\} \\ \{-\lambda_9, -\lambda_{10}, -\lambda_{11}\} & 0 \end{pmatrix} & \begin{pmatrix} \alpha_1 & \{0,0,0\} \\ \{0,0,0\} & \alpha_1 \end{pmatrix} & \begin{pmatrix} \beta_1 & \{\beta_3, \beta_4, \beta_5\} \\ \{0,0,0\} & 0 \end{pmatrix} \\ \begin{pmatrix} \gamma_2 & \{0,0,0\} \\ \{-\gamma_6, -\gamma_7, -\gamma_8\} & 0 \end{pmatrix} & \begin{pmatrix} 0 & \{-\beta_3, -\beta_4, -\beta_5\} \\ \{0,0,0\} & \beta_1 \end{pmatrix} & \begin{pmatrix} 0 & \{0,0,0\} \\ \{0,0,0\} & 0 \end{pmatrix} \end{pmatrix}$$

**T[E₂[1], X₁[e₂], E₂[1]] == ZERO**

True

**Leib[E₂[1], X₁[e₂], E₂[1]] == ZERO**

True

**T[E₁[1], E₂[1], X₁[e₂]] == ZERO**

True

```
(Leib[E₁[1], E₂[1], X₁[e₂]] // Expand) == ZERO
```
True

```
T[E₂[1], E₁[1], X₁[e₂]] == ZERO
```
True

```
(Leib[E₂[1], E₁[1], X₁[e₂]] // Expand) == ZERO
```
True

```
T[E₁[1], X₁[e₂], E₂[1]] == X₁[e₂]
```
True

```
Δ[X₁[e₂]] - Leib[E₁[1], X₁[e₂], E₂[1]] == ZERO
```
True

```
U_{X₁[e₁]}[X₁[e₂]] == X₁[e₁]
```
True

```
(Δ[X₁[e₁]] - defin[X₁[e₁], X₁[e₂]] // Expand) == ZERO
```
True

```
T[E₁[1], X₁[e₁], X₁[e₂]] == E₁[1]
```
True

```
(Δ[E₁[1]] - Leib[E₁[1], X₁[e₁], X₁[e₂]] // Expand) == ZERO
```
True

```
T[X₁[e₁], E₁[1], X₁[e₂]] == E₂[1]
```
True

```
Δ[E₂[1]] - Leib[X₁[e₁], E₁[1], X₁[e₂]] == ZERO
```
True

```
T[E₂[1], X₁[e₁], X₁[e₂]] == E₂[1]
```
True

```
(Δ[E₂[1]] - Leib[E₂[1], X₁[e₁], X₁[e₂]] // Expand) == ZERO
```

True

```
T[X₁[e₂], X₁[e₁], X₁[e₁]] == ZERO
```

True

```
Leib[X₁[e₂], X₁[e₁], X₁[e₁]] == ZERO
```

True

```
T[X₁[e₁], X₁[e₂], X₁[e₁]] == 2 X₁[e₁]
```

True

```
(2 Δ[X₁[e₁]] - Leib[X₁[e₁], X₁[e₂], X₁[e₁]] // Expand) == ZERO
```

True

```
T[X₁[e₂], X₁[e₂], X₁[e₁]] == ZERO
```

True

```
Leib[X₁[e₂], X₁[e₂], X₁[e₁]] == ZERO
```

True

```
T[X₁[e₂], X₁[e₂], X₁[e₂]] == ZERO
```

True

```
Leib[X₁[e₂], X₁[e₂], X₁[e₂]] == ZERO
```

True

```
Δ[X₁[e₂]]
```

$$\begin{pmatrix} \begin{pmatrix} -\alpha_1 & \{0,0,0\} \\ \{0,0,0\} & -\alpha_1 \end{pmatrix} & \begin{pmatrix} 0 & \{\lambda_6, \lambda_7, \lambda_8\} \\ \{\lambda_9, \lambda_{10}, \lambda_{11}\} & -\delta_1 \end{pmatrix} & \begin{pmatrix} 0 & \{0,0,0\} \\ \{\gamma_6, \gamma_7, \gamma_8\} & \gamma_2 \end{pmatrix} \\ \begin{pmatrix} -\delta_1 & \{-\lambda_6, -\lambda_7, -\lambda_8\} \\ \{-\lambda_9, -\lambda_{10}, -\lambda_{11}\} & 0 \end{pmatrix} & \begin{pmatrix} \alpha_1 & \{0,0,0\} \\ \{0,0,0\} & \alpha_1 \end{pmatrix} & \begin{pmatrix} \beta_1 & \{\beta_3, \beta_4, \beta_5\} \\ \{0,0,0\} & 0 \end{pmatrix} \\ \begin{pmatrix} \gamma_2 & \{0,0,0\} \\ \{-\gamma_6, -\gamma_7, -\gamma_8\} & 0 \end{pmatrix} & \begin{pmatrix} 0 & \{-\beta_3, -\beta_4, -\beta_5\} \\ \{0,0,0\} & \beta_1 \end{pmatrix} & \begin{pmatrix} 0 & \{0,0,0\} \\ \{0,0,0\} & 0 \end{pmatrix} \end{pmatrix}$$

$\Delta[X_1[e_2]] =$
$\Delta[X_1[e_2]] //. \{\lambda_6 \to \rho_1, \lambda_7 \to \rho_2, \lambda_8 \to \rho_3, \lambda_9 \to \rho_4, \lambda_{10} \to \rho_5, \lambda_{11} \to \rho_6\}$

$$\begin{pmatrix} \begin{pmatrix} -\alpha_1 & \{0,0,0\} \\ \{0,0,0\} & -\alpha_1 \end{pmatrix} & \begin{pmatrix} 0 & \{\rho_1,\rho_2,\rho_3\} \\ \{\rho_4,\rho_5,\rho_6\} & -\delta_1 \end{pmatrix} & \begin{pmatrix} 0 & \{0,0,0\} \\ \{\gamma_6,\gamma_7,\gamma_8\} & \gamma_2 \end{pmatrix} \\ \begin{pmatrix} -\delta_1 & \{-\rho_1,-\rho_2,-\rho_3\} \\ \{-\rho_4,-\rho_5,-\rho_6\} & 0 \end{pmatrix} & \begin{pmatrix} \alpha_1 & \{0,0,0\} \\ \{0,0,0\} & \alpha_1 \end{pmatrix} & \begin{pmatrix} \beta_1 & \{\beta_3,\beta_4,\beta_5\} \\ \{0,0,0\} & 0 \end{pmatrix} \\ \begin{pmatrix} \gamma_2 & \{0,0,0\} \\ \{-\gamma_6,-\gamma_7,-\gamma_8\} & 0 \end{pmatrix} & \begin{pmatrix} 0 & \{-\beta_3,-\beta_4,-\beta_5\} \\ \{0,0,0\} & \beta_1 \end{pmatrix} & \begin{pmatrix} 0 & \{0,0,0\} \\ \{0,0,0\} & 0 \end{pmatrix} \end{pmatrix}$$

- $X_1[e_3]$

$\Delta[X_1[e_3]] = \text{generic};$

$U_{X_1[e_3]}[E_1[1]] == \text{ZERO}$

True

$\text{defin}[X_1[e_3], E_1[1]] // \text{Expand}$

$$\begin{pmatrix} \begin{pmatrix} 0 & \{0,0,0\} \\ \{0,0,0\} & 0 \end{pmatrix} & \begin{pmatrix} 0 & \{\lambda_1-\alpha_6,0,0\} \\ \{0,0,0\} & 0 \end{pmatrix} & \begin{pmatrix} 0 & \{0,0,0\} \\ \{0,0,0\} & 0 \end{pmatrix} \\ \begin{pmatrix} 0 & \{\alpha_6-\lambda_1,0,0\} \\ \{0,0,0\} & 0 \end{pmatrix} & \begin{pmatrix} -\lambda_9 & \{0,0,0\} \\ \{0,0,0\} & -\lambda_9 \end{pmatrix} & \begin{pmatrix} -\lambda_{25} & \{-\lambda_{21},0,0\} \\ \{0,\lambda_{24},-\lambda_{23}\} & 0 \end{pmatrix} \\ \begin{pmatrix} 0 & \{0,0,0\} \\ \{0,0,0\} & 0 \end{pmatrix} & \begin{pmatrix} 0 & \{\lambda_{21},0,0\} \\ \{0,-\lambda_{24},\lambda_{23}\} & -\lambda_{25} \end{pmatrix} & \begin{pmatrix} 0 & \{0,0,0\} \\ \{0,0,0\} & 0 \end{pmatrix} \end{pmatrix}$$

$\Delta[X_1[e_3]] = \Delta[X_1[e_3]] //.$
$\{\lambda_1 \to \alpha_6, \lambda_9 \to 0, \lambda_{21} \to 0, \lambda_{23} \to 0, \lambda_{24} \to 0, \lambda_{25} \to 0\}$

$$\begin{pmatrix} \begin{pmatrix} \alpha_6 & \{0,0,0\} \\ \{0,0,0\} & \alpha_6 \end{pmatrix} & \begin{pmatrix} \lambda_4 & \{\lambda_6,\lambda_7,\lambda_8\} \\ \{0,\lambda_{10},\lambda_{11}\} & \lambda_5 \end{pmatrix} & \begin{pmatrix} \lambda_{20} & \{\lambda_{22}, 0 \\ \{0,\lambda_{26},\lambda_{27}\} & 0 \end{pmatrix} \\ \begin{pmatrix} \lambda_5 & \{-\lambda_6,-\lambda_7,-\lambda_8\} \\ \{0,-\lambda_{10},-\lambda_{11}\} & \lambda_4 \end{pmatrix} & \begin{pmatrix} \lambda_2 & \{0,0,0\} \\ \{0,0,0\} & \lambda_2 \end{pmatrix} & \begin{pmatrix} \lambda_{12} & \{\lambda_{14},\lambda_1 \\ \{\lambda_{17},\lambda_{18},\lambda_{19}\} & \lambda_1 \end{pmatrix} \\ \begin{pmatrix} 0 & \{-\lambda_{22},0,0\} \\ \{0,-\lambda_{26},-\lambda_{27}\} & \lambda_{20} \end{pmatrix} & \begin{pmatrix} \lambda_{13} & \{-\lambda_{14},-\lambda_{15},-\lambda_{16}\} \\ \{-\lambda_{17},-\lambda_{18},-\lambda_{19}\} & \lambda_{12} \end{pmatrix} & \begin{pmatrix} \lambda_3 & \{0,0,0\} \\ \{0,0,0\} & \lambda_3 \end{pmatrix} \end{pmatrix}$$

$U_{X_1[e_3]}[E_2[1]] == \text{ZERO}$

True

$\text{defin}[X_1[e_3], E_2[1]] // \text{Expand}$

$$\begin{pmatrix} \begin{pmatrix} 0 & \{0,0,0\} \\ \{0,0,0\} & 0 \end{pmatrix} & \begin{pmatrix} 0 & \{\alpha_6+\lambda_2,0,0\} \\ \{0,0,0\} & 0 \end{pmatrix} & \begin{pmatrix} \lambda_{17} & \{\lambda_{13},0,0\} \\ \{0,-\lambda_{16},\lambda_{15}\} & 0 \end{pmatrix} \\ \begin{pmatrix} 0 & \{-\alpha_6-\lambda_2,0,0\} \\ \{0,0,0\} & 0 \end{pmatrix} & \begin{pmatrix} 0 & \{0,0,0\} \\ \{0,0,0\} & 0 \end{pmatrix} & \begin{pmatrix} 0 & \{0,0,0\} \\ \{0,0,0\} & 0 \end{pmatrix} \\ \begin{pmatrix} 0 & \{-\lambda_{13},0,0\} \\ \{0,\lambda_{16},-\lambda_{15}\} & \lambda_{17} \end{pmatrix} & \begin{pmatrix} 0 & \{0,0,0\} \\ \{0,0,0\} & 0 \end{pmatrix} & \begin{pmatrix} 0 & \{0,0,0\} \\ \{0,0,0\} & 0 \end{pmatrix} \end{pmatrix}$$

```
Δ[X₁[e₃]] =
 Δ[X₁[e₃]] //. {λ₂ -> -α₆, λ₁₃ -> 0, λ₁₅ -> 0, λ₁₆ -> 0, λ₁₇ -> 0}
```

$$\begin{pmatrix} \begin{pmatrix} \alpha_6 & \{0,0,0\} \\ \{0,0,0\} & \alpha_6 \end{pmatrix} & \begin{pmatrix} \lambda_4 & \{\lambda_6,\lambda_7,\lambda_8\} \\ \{0,\lambda_{10},\lambda_{11}\} & \lambda_5 \end{pmatrix} & \begin{pmatrix} \lambda_{20} & \{\lambda_{22},0,0\} \\ \{0,\lambda_{26},\lambda_{27}\} & 0 \end{pmatrix} \\ \begin{pmatrix} \lambda_5 & \{-\lambda_6,-\lambda_7,-\lambda_8\} \\ \{0,-\lambda_{10},-\lambda_{11}\} & \lambda_4 \end{pmatrix} & \begin{pmatrix} -\alpha_6 & \{0,0,0\} \\ \{0,0,0\} & -\alpha_6 \end{pmatrix} & \begin{pmatrix} \lambda_{12} & \{\lambda_{14},0,0\} \\ \{0,\lambda_{18},\lambda_{19}\} & 0 \end{pmatrix} \\ \begin{pmatrix} 0 & \{-\lambda_{22},0,0\} \\ \{0,-\lambda_{26},-\lambda_{27}\} & \lambda_{20} \end{pmatrix} & \begin{pmatrix} 0 & \{-\lambda_{14},0,0\} \\ \{0,-\lambda_{18},-\lambda_{19}\} & \lambda_{12} \end{pmatrix} & \begin{pmatrix} \lambda_3 & \{0,0,0\} \\ \{0,0,0\} & \lambda_3 \end{pmatrix} \end{pmatrix}$$

**U_{X₁[e₃]}[E₃[1]] == ZERO**

True

**(defin[X₁[e₃], E₃[1]] // Expand) == ZERO**

True

**U_{X₁[e₃]}[X₁[e₃]] == ZERO**

True

**(defin[X₁[e₃], X₁[e₃]] // Expand) == ZERO**

True

**U_{X₁[e₃]}[X₁[e₂]] == ZERO**

True

**defin[X₁[e₃], X₁[e₂]] // Expand**

$$\begin{pmatrix} \begin{pmatrix} 0 & \{0,0,0\} \\ \{0,0,0\} & 0 \end{pmatrix} & \begin{pmatrix} 0 & \{\lambda_4-\rho_4,0,0\} \\ \{0,0,0\} & 0 \end{pmatrix} & \begin{pmatrix} 0 & \{0,0,0\} \\ \{0,0,0\} & 0 \end{pmatrix} \\ \begin{pmatrix} 0 & \{\rho_4-\lambda_4,0,0\} \\ \{0,0,0\} & 0 \end{pmatrix} & \begin{pmatrix} 0 & \{0,0,0\} \\ \{0,0,0\} & 0 \end{pmatrix} & \begin{pmatrix} 0 & \{0,0,0\} \\ \{0,0,0\} & 0 \end{pmatrix} \\ \begin{pmatrix} 0 & \{0,0,0\} \\ \{0,0,0\} & 0 \end{pmatrix} & \begin{pmatrix} 0 & \{0,0,0\} \\ \{0,0,0\} & 0 \end{pmatrix} & \begin{pmatrix} 0 & \{0,0,0\} \\ \{0,0,0\} & 0 \end{pmatrix} \end{pmatrix}$$

**Δ[X₁[e₃]] = Δ[X₁[e₃]] //. {λ₄ → ρ₄}**

$$\begin{pmatrix} \begin{pmatrix} \alpha_6 & \{0,0,0\} \\ \{0,0,0\} & \alpha_6 \end{pmatrix} & \begin{pmatrix} \rho_4 & \{\lambda_6,\lambda_7,\lambda_8\} \\ \{0,\lambda_{10},\lambda_{11}\} & \lambda_5 \end{pmatrix} & \begin{pmatrix} \lambda_{20} & \{\lambda_{22},0,0\} \\ \{0,\lambda_{26},\lambda_{27}\} & 0 \end{pmatrix} \\ \begin{pmatrix} \lambda_5 & \{-\lambda_6,-\lambda_7,-\lambda_8\} \\ \{0,-\lambda_{10},-\lambda_{11}\} & \rho_4 \end{pmatrix} & \begin{pmatrix} -\alpha_6 & \{0,0,0\} \\ \{0,0,0\} & -\alpha_6 \end{pmatrix} & \begin{pmatrix} \lambda_{12} & \{\lambda_{14},0,0\} \\ \{0,\lambda_{18},\lambda_{19}\} & 0 \end{pmatrix} \\ \begin{pmatrix} 0 & \{-\lambda_{22},0,0\} \\ \{0,-\lambda_{26},-\lambda_{27}\} & \lambda_{20} \end{pmatrix} & \begin{pmatrix} 0 & \{-\lambda_{14},0,0\} \\ \{0,-\lambda_{18},-\lambda_{19}\} & \lambda_{12} \end{pmatrix} & \begin{pmatrix} \lambda_3 & \{0,0,0\} \\ \{0,0,0\} & \lambda_3 \end{pmatrix} \end{pmatrix}$$

**U$_{X_1[e_3]}$[X$_1$[e$_1$]] == ZERO**

True

**defin[X$_1$[e$_3$], X$_1$[e$_1$]] // Expand**

$$\left(\begin{array}{ccc} \begin{pmatrix} 0 & \{0,0,0\} \\ \{0,0,0\} & 0 \end{pmatrix} & \begin{pmatrix} 0 & \{\lambda_5-\delta_5,0,0\} \\ \{0,0,0\} & 0 \end{pmatrix} & \begin{pmatrix} 0 & \{0,0,0\} \\ \{0,0,0\} & 0 \end{pmatrix} \\ \begin{pmatrix} 0 & \{\delta_5-\lambda_5,0,0\} \\ \{0,0,0\} & 0 \end{pmatrix} & \begin{pmatrix} 0 & \{0,0,0\} \\ \{0,0,0\} & 0 \end{pmatrix} & \begin{pmatrix} 0 & \{0,0,0\} \\ \{0,0,0\} & 0 \end{pmatrix} \\ \begin{pmatrix} 0 & \{0,0,0\} \\ \{0,0,0\} & 0 \end{pmatrix} & \begin{pmatrix} 0 & \{0,0,0\} \\ \{0,0,0\} & 0 \end{pmatrix} & \begin{pmatrix} 0 & \{0,0,0\} \\ \{0,0,0\} & 0 \end{pmatrix} \end{array}\right)$$

**Δ[X$_1$[e$_3$]] = Δ[X$_1$[e$_3$]] //. {λ$_5$ → δ$_5$}**

$$\left(\begin{array}{ccc} \begin{pmatrix} \alpha_6 & \{0,0,0\} \\ \{0,0,0\} & \alpha_6 \end{pmatrix} & \begin{pmatrix} \rho_4 & \{\lambda_6,\lambda_7,\lambda_8\} \\ \{0,\lambda_{10},\lambda_{11}\} & \delta_5 \end{pmatrix} & \begin{pmatrix} \lambda_{20} & \{\lambda_{22},0,0\} \\ \{0,\lambda_{26},\lambda_{27}\} & 0 \end{pmatrix} \\ \begin{pmatrix} \delta_5 & \{-\lambda_6,-\lambda_7,-\lambda_8\} \\ \{0,-\lambda_{10},-\lambda_{11}\} & \rho_4 \end{pmatrix} & \begin{pmatrix} -\alpha_6 & \{0,0,0\} \\ \{0,0,0\} & -\alpha_6 \end{pmatrix} & \begin{pmatrix} \lambda_{12} & \{\lambda_{14},0,0\} \\ \{0,\lambda_{18},\lambda_{19}\} & 0 \end{pmatrix} \\ \begin{pmatrix} 0 & \{-\lambda_{22},0,0\} \\ \{0,-\lambda_{26},-\lambda_{27}\} & \lambda_{20} \end{pmatrix} & \begin{pmatrix} 0 & \{-\lambda_{14},0,0\} \\ \{0,-\lambda_{18},-\lambda_{19}\} & \lambda_{12} \end{pmatrix} & \begin{pmatrix} \lambda_3 & \{0,0,0\} \\ \{0,0,0\} & \lambda_3 \end{pmatrix} \end{array}\right)$$

**T[E$_1$[1], E$_1$[1], X$_1$[e$_3$]] == X$_1$[e$_3$]**

True

**Δ[X$_1$[e$_3$]] - Leib[E$_1$[1], E$_1$[1], X$_1$[e$_3$]] // Expand**

$$\left(\begin{array}{ccc} \begin{pmatrix} 0 & \{0,0,0\} \\ \{0,0,0\} & 0 \end{pmatrix} & \begin{pmatrix} 0 & \{0,0,0\} \\ \{0,0,0\} & 0 \end{pmatrix} & \begin{pmatrix} 0 & \{0,0,0\} \\ \{0,0,0\} & 0 \end{pmatrix} \\ \begin{pmatrix} 0 & \{0,0,0\} \\ \{0,0,0\} & 0 \end{pmatrix} & \begin{pmatrix} 0 & \{0,0,0\} \\ \{0,0,0\} & 0 \end{pmatrix} & \begin{pmatrix} \beta_6+\lambda_{12} & \{\beta_2+\lambda_{14},0, \\ \{0,\lambda_{18}-\beta_5,\beta_4+\lambda_{19}\} & 0 \end{pmatrix} \\ \begin{pmatrix} 0 & \{0,0,0\} \\ \{0,0,0\} & 0 \end{pmatrix} & \begin{pmatrix} 0 & \{-\beta_2-\lambda_{14},0,0\} \\ \{0,\beta_5-\lambda_{18},-\beta_4-\lambda_{19}\} & \beta_6+\lambda_{12} \end{pmatrix} & \begin{pmatrix} \lambda_3 & \{0,0,0\} \\ \{0,0,0\} & \lambda_3 \end{pmatrix} \end{array}\right)$$

**Δ[X$_1$[e$_3$]] = Δ[X$_1$[e$_3$]] //. {λ$_3$ → 0, λ$_{12}$ → -β$_6$, λ$_{14}$ → -β$_2$, λ$_{18}$ → β$_5$, λ$_{19}$ → -β$_4$}**

$$\left(\begin{array}{ccc} \begin{pmatrix} \alpha_6 & \{0,0,0\} \\ \{0,0,0\} & \alpha_6 \end{pmatrix} & \begin{pmatrix} \rho_4 & \{\lambda_6,\lambda_7,\lambda_8\} \\ \{0,\lambda_{10},\lambda_{11}\} & \delta_5 \end{pmatrix} & \begin{pmatrix} \lambda_{20} & \{\lambda_{22},0,0\} \\ \{0,\lambda_{26},\lambda_{27}\} & 0 \end{pmatrix} \\ \begin{pmatrix} \delta_5 & \{-\lambda_6,-\lambda_7,-\lambda_8\} \\ \{0,-\lambda_{10},-\lambda_{11}\} & \rho_4 \end{pmatrix} & \begin{pmatrix} -\alpha_6 & \{0,0,0\} \\ \{0,0,0\} & -\alpha_6 \end{pmatrix} & \begin{pmatrix} -\beta_6 & \{-\beta_2,0,0\} \\ \{0,\beta_5,-\beta_4\} & 0 \end{pmatrix} \\ \begin{pmatrix} 0 & \{-\lambda_{22},0,0\} \\ \{0,-\lambda_{26},-\lambda_{27}\} & \lambda_{20} \end{pmatrix} & \begin{pmatrix} 0 & \{\beta_2,0,0\} \\ \{0,-\beta_5,\beta_4\} & -\beta_6 \end{pmatrix} & \begin{pmatrix} 0 & \{0,0,0\} \\ \{0,0,0\} & 0 \end{pmatrix} \end{array}\right)$$

**T[E$_1$[1], X$_1$[e$_3$], E$_1$[1]] == ZERO**

True

```
Leib[E₁[1], X₁[e₃], E₁[1]] == ZERO
```

True

```
T[E₂[1], E₂[1], X₁[e₃]] == X₁[e₃]
```

True

```
Δ[X₁[e₃]] - Leib[E₂[1], E₂[1], X₁[e₃]] // Expand
```

$$\begin{pmatrix} \begin{pmatrix} 0 & \{0,0,0\} \\ \{0,0,0\} & 0 \end{pmatrix} & \begin{pmatrix} 0 & \{0,0,0\} \\ \{0,0,0\} & 0 \end{pmatrix} & \begin{pmatrix} \lambda_{20}-\gamma_6 & \{\lambda_{22}-\gamma_2,0,0\} \\ \{0,\gamma_5+\lambda_{26},\lambda_{27}-\gamma_4\} & 0 \end{pmatrix} \\ \begin{pmatrix} 0 & \{0,0,0\} \\ \{0,0,0\} & 0 \end{pmatrix} & \begin{pmatrix} 0 & \{0,0,0\} \\ \{0,0,0\} & 0 \end{pmatrix} & \begin{pmatrix} 0 & \{0,0,0\} \\ \{0,0,0\} & 0 \end{pmatrix} \\ \begin{pmatrix} 0 & \{\gamma_2-\lambda_{22},0,0\} \\ \{0,-\gamma_5-\lambda_{26},\gamma_4-\lambda_{27}\} & \lambda_{20}-\gamma_6 \end{pmatrix} & \begin{pmatrix} 0 & \{0,0,0\} \\ \{0,0,0\} & 0 \end{pmatrix} & \begin{pmatrix} 0 & \{0,0,0\} \\ \{0,0,0\} & 0 \end{pmatrix} \end{pmatrix}$$

```
Δ[X₁[e₃]] = Δ[X₁[e₃]] //. {λ₂₀ → γ₆, λ₂₂ → γ₂, λ₂₆ → -γ₅, λ₂₇ → γ₄}
```

$$\begin{pmatrix} \begin{pmatrix} \alpha_6 & \{0,0,0\} \\ \{0,0,0\} & \alpha_6 \end{pmatrix} & \begin{pmatrix} \rho_4 & \{\lambda_6,\lambda_7,\lambda_8\} \\ \{0,\lambda_{10},\lambda_{11}\} & \delta_5 \end{pmatrix} & \begin{pmatrix} \gamma_6 & \{\gamma_2,0,0\} \\ \{0,-\gamma_5,\gamma_4\} & 0 \end{pmatrix} \\ \begin{pmatrix} \delta_5 & \{-\lambda_6,-\lambda_7,-\lambda_8\} \\ \{0,-\lambda_{10},-\lambda_{11}\} & \rho_4 \end{pmatrix} & \begin{pmatrix} -\alpha_6 & \{0,0,0\} \\ \{0,0,0\} & -\alpha_6 \end{pmatrix} & \begin{pmatrix} -\beta_6 & \{-\beta_2,0,0\} \\ \{0,\beta_5,-\beta_4\} & 0 \end{pmatrix} \\ \begin{pmatrix} 0 & \{-\gamma_2,0,0\} \\ \{0,\gamma_5,-\gamma_4\} & \gamma_6 \end{pmatrix} & \begin{pmatrix} 0 & \{\beta_2,0,0\} \\ \{0,-\beta_5,\beta_4\} & -\beta_6 \end{pmatrix} & \begin{pmatrix} 0 & \{0,0,0\} \\ \{0,0,0\} & 0 \end{pmatrix} \end{pmatrix}$$

```
T[E₂[1], X₁[e₃], E₂[1]] == ZERO
```

True

```
Leib[E₂[1], X₁[e₃], E₂[1]] == ZERO
```

True

```
T[E₁[1], E₂[1], X₁[e₃]] == ZERO
```

True

```
(Leib[E₁[1], E₂[1], X₁[e₃]] // Expand) == ZERO
```

True

```
T[E₂[1], E₁[1], X₁[e₃]] == ZERO
```

True

```
(Leib[E₂[1], E₁[1], X₁[e₃]] // Expand) == ZERO
```

True

**T[E₁[1], X₁[e₃], E₂[1]] == X₁[e₃]**

True

**Δ[X₁[e₃]] - Leib[E₁[1], X₁[e₃], E₂[1]] == ZERO**

True

**U_{X₁[e₁]}[X₁[e₂]] == X₁[e₁]**

True

**(Δ[X₁[e₁]] - defin[X₁[e₁], X₁[e₂]] // Expand) == ZERO**

True

**T[E₁[1], X₁[e₁], X₁[e₃]] == ZERO**

True

**(Leib[E₁[1], X₁[e₁], X₁[e₃]] // Expand) == ZERO**

True

**T[X₁[e₁], E₁[1], X₁[e₃]] == ZERO**

True

**(Leib[X₁[e₁], E₁[1], X₁[e₃]] // Expand) == ZERO**

True

**T[E₂[1], X₁[e₁], X₁[e₃]] == ZERO**

True

**(Leib[E₂[1], X₁[e₁], X₁[e₃]] // Expand) == ZERO**

True

**T[X₁[e₃], X₁[e₁], X₁[e₁]] == ZERO**

True

**(Leib[X₁[e₃], X₁[e₁], X₁[e₁]] // Expand) == ZERO**

True

**T[X$_1$[e$_1$], X$_1$[e$_3$], X$_1$[e$_1$]] == ZERO**

True

**(Leib[X$_1$[e$_1$], X$_1$[e$_3$], X$_1$[e$_1$]] // Expand) == ZERO**

True

**T[X$_1$[e$_3$], X$_1$[e$_3$], X$_1$[e$_1$]] == ZERO**

True

**(Leib[X$_1$[e$_3$], X$_1$[e$_3$], X$_1$[e$_1$]] // Expand) == ZERO**

True

**T[X$_1$[e$_3$], X$_1$[e$_3$], X$_1$[e$_2$]] == ZERO**

True

**(Leib[X$_1$[e$_3$], X$_1$[e$_3$], X$_1$[e$_2$]] // Expand) == ZERO**

True

**T[X$_1$[e$_3$], X$_1$[e$_3$], X$_1$[e$_3$]] == ZERO**

True

**(Leib[X$_1$[e$_3$], X$_1$[e$_3$], X$_1$[e$_3$]] // Expand) == ZERO**

True

**T[X$_1$[e$_3$], X$_1$[e$_1$], X$_1$[e$_3$]] == ZERO**

True

**(Leib[X$_1$[e$_3$], X$_1$[e$_1$], X$_1$[e$_3$]] // Expand) == ZERO**

True

**Δ[X$_1$[e$_3$]] =
  Δ[X$_1$[e$_3$]] //. {λ$_6$ → η$_1$, λ$_7$ → η$_2$, λ$_8$ → η$_3$, λ$_{10}$ → η$_4$, λ$_{11}$ → η$_5$}**

$$\begin{pmatrix} \begin{pmatrix} \alpha_6 & \{0,0,0\} \\ \{0,0,0\} & \alpha_6 \end{pmatrix} & \begin{pmatrix} \rho_4 & \{\eta_1,\eta_2,\eta_3\} \\ \{0,\eta_4,\eta_5\} & \delta_5 \end{pmatrix} & \begin{pmatrix} \gamma_6 & \{\gamma_2,0,0\} \\ \{0,-\gamma_5,\gamma_4\} & 0 \end{pmatrix} \\ \begin{pmatrix} \delta_5 & \{-\eta_1,-\eta_2,-\eta_3\} \\ \{0,-\eta_4,-\eta_5\} & \rho_4 \end{pmatrix} & \begin{pmatrix} -\alpha_6 & \{0,0,0\} \\ \{0,0,0\} & -\alpha_6 \end{pmatrix} & \begin{pmatrix} -\beta_6 & \{-\beta_2,0,0\} \\ \{0,\beta_5,-\beta_4\} & 0 \end{pmatrix} \\ \begin{pmatrix} 0 & \{-\gamma_2,0,0\} \\ \{0,\gamma_5,-\gamma_4\} & \gamma_6 \end{pmatrix} & \begin{pmatrix} 0 & \{\beta_2,0,0\} \\ \{0,-\beta_5,\beta_4\} & -\beta_6 \end{pmatrix} & \begin{pmatrix} 0 & \{0,0,0\} \\ \{0,0,0\} & 0 \end{pmatrix} \end{pmatrix}$$

- $X_1[e_4]$

  $\Delta[X_1[e_4]]$ = generic;

  $U_{X_1[e_4]}[E_1[1]]$ == ZERO

  True

  defin[$X_1[e_4]$, $E_1[1]$] // Expand

  $$\begin{pmatrix} \begin{pmatrix} 0 & \{0,0,0\} \\ \{0,0,0\} & 0 \end{pmatrix} & \begin{pmatrix} 0 & \{0, \lambda_1 - \alpha_7, 0\} \\ \{0,0,0\} & 0 \end{pmatrix} & \begin{pmatrix} 0 & \{0,0,0\} \\ \{0,0,0\} & 0 \end{pmatrix} \\ \begin{pmatrix} 0 & \{0, \alpha_7 - \lambda_1, 0\} \\ \{0,0,0\} & 0 \end{pmatrix} & \begin{pmatrix} -\lambda_{10} & \{0,0,0\} \\ \{0,0,0\} & -\lambda_{10} \end{pmatrix} & \begin{pmatrix} -\lambda_{26} & \{0, -\lambda_{21}, 0\} \\ \{-\lambda_{24}, 0, \lambda_{22}\} & 0 \end{pmatrix} \\ \begin{pmatrix} 0 & \{0,0,0\} \\ \{0,0,0\} & 0 \end{pmatrix} & \begin{pmatrix} 0 & \{0, \lambda_{21}, 0\} \\ \{\lambda_{24}, 0, -\lambda_{22}\} & -\lambda_{26} \end{pmatrix} & \begin{pmatrix} 0 & \{0,0,0\} \\ \{0,0,0\} & 0 \end{pmatrix} \end{pmatrix}$$

  $\Delta[X_1[e_4]]$ =
   $\Delta[X_1[e_4]]$ //. {$\lambda_1 \to \alpha_7$, $\lambda_{10} \to 0$, $\lambda_{21} \to 0$, $\lambda_{22} \to 0$, $\lambda_{24} \to 0$, $\lambda_{26} \to 0$}

  $$\begin{pmatrix} \begin{pmatrix} \alpha_7 & \{0,0,0\} \\ \{0,0,0\} & \alpha_7 \end{pmatrix} & \begin{pmatrix} \lambda_4 & \{\lambda_6, \lambda_7, \lambda_8\} \\ \{\lambda_9, 0, \lambda_{11}\} & \lambda_5 \end{pmatrix} & \begin{pmatrix} \lambda_{20} & \{0, \lambda_{23}, \\ \{\lambda_{25}, 0, \lambda_{27}\} & 0 \end{pmatrix} \\ \begin{pmatrix} \lambda_5 & \{-\lambda_6, -\lambda_7, -\lambda_8\} \\ \{-\lambda_9, 0, -\lambda_{11}\} & \lambda_4 \end{pmatrix} & \begin{pmatrix} \lambda_2 & \{0,0,0\} \\ \{0,0,0\} & \lambda_2 \end{pmatrix} & \begin{pmatrix} \lambda_{12} & \{\lambda_{14}, \lambda_{15} \\ \{\lambda_{17}, \lambda_{18}, \lambda_{19}\} & \lambda_{13} \end{pmatrix} \\ \begin{pmatrix} 0 & \{0, -\lambda_{23}, 0\} \\ \{-\lambda_{25}, 0, -\lambda_{27}\} & \lambda_{20} \end{pmatrix} & \begin{pmatrix} \lambda_{13} & \{-\lambda_{14}, -\lambda_{15}, -\lambda_{16}\} \\ \{-\lambda_{17}, -\lambda_{18}, -\lambda_{19}\} & \lambda_{12} \end{pmatrix} & \begin{pmatrix} \lambda_3 & \{0,0,0\} \\ \{0,0,0\} & \lambda_3 \end{pmatrix} \end{pmatrix}$$

  $U_{X_1[e_4]}[E_2[1]]$ == ZERO

  True

  defin[$X_1[e_4]$, $E_2[1]$] // Expand

  $$\begin{pmatrix} \begin{pmatrix} 0 & \{0,0,0\} \\ \{0,0,0\} & 0 \end{pmatrix} & \begin{pmatrix} 0 & \{0, \alpha_7 + \lambda_2, 0\} \\ \{0,0,0\} & 0 \end{pmatrix} & \begin{pmatrix} \lambda_{18} & \{0, \lambda_{13}, 0\} \\ \{\lambda_{16}, 0, -\lambda_{14}\} & 0 \end{pmatrix} \\ \begin{pmatrix} 0 & \{0, -\alpha_7 - \lambda_2, 0\} \\ \{0,0,0\} & 0 \end{pmatrix} & \begin{pmatrix} 0 & \{0,0,0\} \\ \{0,0,0\} & 0 \end{pmatrix} & \begin{pmatrix} 0 & \{0,0,0\} \\ \{0,0,0\} & 0 \end{pmatrix} \\ \begin{pmatrix} 0 & \{0, -\lambda_{13}, 0\} \\ \{-\lambda_{16}, 0, \lambda_{14}\} & \lambda_{18} \end{pmatrix} & \begin{pmatrix} 0 & \{0,0,0\} \\ \{0,0,0\} & 0 \end{pmatrix} & \begin{pmatrix} 0 & \{0,0,0\} \\ \{0,0,0\} & 0 \end{pmatrix} \end{pmatrix}$$

  $\Delta[X_1[e_4]]$ = $\Delta[X_1[e_4]]$ //. {$\lambda_2 \to -\alpha_7$, $\lambda_{13} \to 0$, $\lambda_{14} \to 0$, $\lambda_{16} \to 0$, $\lambda_{18} \to 0$}

  $$\begin{pmatrix} \begin{pmatrix} \alpha_7 & \{0,0,0\} \\ \{0,0,0\} & \alpha_7 \end{pmatrix} & \begin{pmatrix} \lambda_4 & \{\lambda_6, \lambda_7, \lambda_8\} \\ \{\lambda_9, 0, \lambda_{11}\} & \lambda_5 \end{pmatrix} & \begin{pmatrix} \lambda_{20} & \{0, \lambda_{23}, 0\} \\ \{\lambda_{25}, 0, \lambda_{27}\} & 0 \end{pmatrix} \\ \begin{pmatrix} \lambda_5 & \{-\lambda_6, -\lambda_7, -\lambda_8\} \\ \{-\lambda_9, 0, -\lambda_{11}\} & \lambda_4 \end{pmatrix} & \begin{pmatrix} -\alpha_7 & \{0,0,0\} \\ \{0,0,0\} & -\alpha_7 \end{pmatrix} & \begin{pmatrix} \lambda_{12} & \{0, \lambda_{15}, 0\} \\ \{\lambda_{17}, 0, \lambda_{19}\} & 0 \end{pmatrix} \\ \begin{pmatrix} 0 & \{0, -\lambda_{23}, 0\} \\ \{-\lambda_{25}, 0, -\lambda_{27}\} & \lambda_{20} \end{pmatrix} & \begin{pmatrix} 0 & \{0, -\lambda_{15}, 0\} \\ \{-\lambda_{17}, 0, -\lambda_{19}\} & \lambda_{12} \end{pmatrix} & \begin{pmatrix} \lambda_3 & \{0,0,0\} \\ \{0,0,0\} & \lambda_3 \end{pmatrix} \end{pmatrix}$$

**U$_{X_1[e_4]}$ [E$_3$[1]] == ZERO**

True

**(defin[X$_1$[e$_4$], E$_3$[1]] // Expand) == ZERO**

True

**U$_{X_1[e_4]}$ [X$_1$[e$_4$]] == ZERO**

True

**(defin[X$_1$[e$_4$], X$_1$[e$_4$]] // Expand) == ZERO**

True

**U$_{X_1[e_4]}$ [X$_1$[e$_3$]] == ZERO**

True

**defin[X$_1$[e$_4$], X$_1$[e$_3$]] // Expand**

$$\begin{pmatrix} \begin{pmatrix} 0 & \{0,0,0\} \\ \{0,0,0\} & 0 \end{pmatrix} & \begin{pmatrix} 0 & \{0,-\eta_4-\lambda_9,0\} \\ \{0,0,0\} & 0 \end{pmatrix} & \begin{pmatrix} 0 & \{0,0,0\} \\ \{0,0,0\} & 0 \end{pmatrix} \\ \begin{pmatrix} 0 & \{0,\eta_4+\lambda_9,0\} \\ \{0,0,0\} & 0 \end{pmatrix} & \begin{pmatrix} 0 & \{0,0,0\} \\ \{0,0,0\} & 0 \end{pmatrix} & \begin{pmatrix} 0 & \{0,0,0\} \\ \{0,0,0\} & 0 \end{pmatrix} \\ \begin{pmatrix} 0 & \{0,0,0\} \\ \{0,0,0\} & 0 \end{pmatrix} & \begin{pmatrix} 0 & \{0,0,0\} \\ \{0,0,0\} & 0 \end{pmatrix} & \begin{pmatrix} 0 & \{0,0,0\} \\ \{0,0,0\} & 0 \end{pmatrix} \end{pmatrix}$$

**Δ[X$_1$[e$_4$]] = Δ[X$_1$[e$_4$]] //. {λ$_9$ → −η$_4$}**

$$\begin{pmatrix} \begin{pmatrix} \alpha_7 & \{0,0,0\} \\ \{0,0,0\} & \alpha_7 \end{pmatrix} & \begin{pmatrix} \lambda_4 & \{\lambda_6,\lambda_7,\lambda_8\} \\ \{-\eta_4,0,\lambda_{11}\} & \lambda_5 \end{pmatrix} & \begin{pmatrix} \lambda_{20} & \{0,\lambda_{23},0\} \\ \{\lambda_{25},0,\lambda_{27}\} & 0 \end{pmatrix} \\ \begin{pmatrix} \lambda_5 & \{-\lambda_6,-\lambda_7,-\lambda_8\} \\ \{\eta_4,0,-\lambda_{11}\} & \lambda_4 \end{pmatrix} & \begin{pmatrix} -\alpha_7 & \{0,0,0\} \\ \{0,0,0\} & -\alpha_7 \end{pmatrix} & \begin{pmatrix} \lambda_{12} & \{0,\lambda_{15},0\} \\ \{\lambda_{17},0,\lambda_{19}\} & 0 \end{pmatrix} \\ \begin{pmatrix} 0 & \{0,-\lambda_{23},0\} \\ \{-\lambda_{25},0,-\lambda_{27}\} & \lambda_{20} \end{pmatrix} & \begin{pmatrix} 0 & \{0,-\lambda_{15},0\} \\ \{-\lambda_{17},0,-\lambda_{19}\} & \lambda_{12} \end{pmatrix} & \begin{pmatrix} \lambda_3 & \{0,0,0\} \\ \{0,0,0\} & \lambda_3 \end{pmatrix} \end{pmatrix}$$

**U$_{X_1[e_4]}$ [X$_1$[e$_2$]] == ZERO**

True

**defin[X$_1$[e$_4$], X$_1$[e$_2$]] // Expand**

$$\begin{pmatrix} \begin{pmatrix} 0 & \{0,0,0\} \\ \{0,0,0\} & 0 \end{pmatrix} & \begin{pmatrix} 0 & \{0,\lambda_4-\rho_5,0\} \\ \{0,0,0\} & 0 \end{pmatrix} & \begin{pmatrix} 0 & \{0,0,0\} \\ \{0,0,0\} & 0 \end{pmatrix} \\ \begin{pmatrix} 0 & \{0,\rho_5-\lambda_4,0\} \\ \{0,0,0\} & 0 \end{pmatrix} & \begin{pmatrix} 0 & \{0,0,0\} \\ \{0,0,0\} & 0 \end{pmatrix} & \begin{pmatrix} 0 & \{0,0,0\} \\ \{0,0,0\} & 0 \end{pmatrix} \\ \begin{pmatrix} 0 & \{0,0,0\} \\ \{0,0,0\} & 0 \end{pmatrix} & \begin{pmatrix} 0 & \{0,0,0\} \\ \{0,0,0\} & 0 \end{pmatrix} & \begin{pmatrix} 0 & \{0,0,0\} \\ \{0,0,0\} & 0 \end{pmatrix} \end{pmatrix}$$

**Δ[X₁[e₄]] = Δ[X₁[e₄]] //. {λ₄ → ρ₅}**

$$\left( \begin{array}{ccc} \begin{pmatrix} \alpha_7 & \{0,0,0\} \\ \{0,0,0\} & \alpha_7 \end{pmatrix} & \begin{pmatrix} \rho_5 & \{\lambda_6, \lambda_7, \lambda_8\} \\ \{-\eta_4, 0, \lambda_{11}\} & \lambda_5 \end{pmatrix} & \begin{pmatrix} \lambda_{20} & \{0, \lambda_{23}, 0\} \\ \{\lambda_{25}, 0, \lambda_{27}\} & 0 \end{pmatrix} \\ \begin{pmatrix} \lambda_5 & \{-\lambda_6, -\lambda_7, -\lambda_8\} \\ \{\eta_4, 0, -\lambda_{11}\} & \rho_5 \end{pmatrix} & \begin{pmatrix} -\alpha_7 & \{0,0,0\} \\ \{0,0,0\} & -\alpha_7 \end{pmatrix} & \begin{pmatrix} \lambda_{12} & \{0, \lambda_{15}, 0\} \\ \{\lambda_{17}, 0, \lambda_{19}\} & 0 \end{pmatrix} \\ \begin{pmatrix} 0 & \{0, -\lambda_{23}, 0\} \\ \{-\lambda_{25}, 0, -\lambda_{27}\} & \lambda_{20} \end{pmatrix} & \begin{pmatrix} 0 & \{0, -\lambda_{15}, 0\} \\ \{-\lambda_{17}, 0, -\lambda_{19}\} & \lambda_{12} \end{pmatrix} & \begin{pmatrix} \lambda_3 & \{0,0,0\} \\ \{0,0,0\} & \lambda_3 \end{pmatrix} \end{array} \right)$$

**U_{X₁[e₄]}[X₁[e₁]] == ZERO**

True

**defin[X₁[e₄], X₁[e₁]] // Expand**

$$\left( \begin{array}{ccc} \begin{pmatrix} 0 & \{0,0,0\} \\ \{0,0,0\} & 0 \end{pmatrix} & \begin{pmatrix} 0 & \{0, \lambda_5 - \delta_6, 0\} \\ \{0,0,0\} & 0 \end{pmatrix} & \begin{pmatrix} 0 & \{0,0,0\} \\ \{0,0,0\} & 0 \end{pmatrix} \\ \begin{pmatrix} 0 & \{0, \delta_6 - \lambda_5, 0\} \\ \{0,0,0\} & 0 \end{pmatrix} & \begin{pmatrix} 0 & \{0,0,0\} \\ \{0,0,0\} & 0 \end{pmatrix} & \begin{pmatrix} 0 & \{0,0,0\} \\ \{0,0,0\} & 0 \end{pmatrix} \\ \begin{pmatrix} 0 & \{0,0,0\} \\ \{0,0,0\} & 0 \end{pmatrix} & \begin{pmatrix} 0 & \{0,0,0\} \\ \{0,0,0\} & 0 \end{pmatrix} & \begin{pmatrix} 0 & \{0,0,0\} \\ \{0,0,0\} & 0 \end{pmatrix} \end{array} \right)$$

**Δ[X₁[e₄]] = Δ[X₁[e₄]] //. {λ₅ → δ₆}**

$$\left( \begin{array}{ccc} \begin{pmatrix} \alpha_7 & \{0,0,0\} \\ \{0,0,0\} & \alpha_7 \end{pmatrix} & \begin{pmatrix} \rho_5 & \{\lambda_6, \lambda_7, \lambda_8\} \\ \{-\eta_4, 0, \lambda_{11}\} & \delta_6 \end{pmatrix} & \begin{pmatrix} \lambda_{20} & \{0, \lambda_{23}, 0\} \\ \{\lambda_{25}, 0, \lambda_{27}\} & 0 \end{pmatrix} \\ \begin{pmatrix} \delta_6 & \{-\lambda_6, -\lambda_7, -\lambda_8\} \\ \{\eta_4, 0, -\lambda_{11}\} & \rho_5 \end{pmatrix} & \begin{pmatrix} -\alpha_7 & \{0,0,0\} \\ \{0,0,0\} & -\alpha_7 \end{pmatrix} & \begin{pmatrix} \lambda_{12} & \{0, \lambda_{15}, 0\} \\ \{\lambda_{17}, 0, \lambda_{19}\} & 0 \end{pmatrix} \\ \begin{pmatrix} 0 & \{0, -\lambda_{23}, 0\} \\ \{-\lambda_{25}, 0, -\lambda_{27}\} & \lambda_{20} \end{pmatrix} & \begin{pmatrix} 0 & \{0, -\lambda_{15}, 0\} \\ \{-\lambda_{17}, 0, -\lambda_{19}\} & \lambda_{12} \end{pmatrix} & \begin{pmatrix} \lambda_3 & \{0,0,0\} \\ \{0,0,0\} & \lambda_3 \end{pmatrix} \end{array} \right)$$

**T[E₁[1], E₁[1], X₁[e₄]] == X₁[e₄]**

True

**Δ[X₁[e₄]] - Leib[E₁[1], E₁[1], X₁[e₄]] // Expand**

$$\left( \begin{array}{ccc} \begin{pmatrix} 0 & \{0,0,0\} \\ \{0,0,0\} & 0 \end{pmatrix} & \begin{pmatrix} 0 & \{0,0,0\} \\ \{0,0,0\} & 0 \end{pmatrix} & \begin{pmatrix} 0 & \{0,0,0\} \\ \{0,0,0\} & 0 \end{pmatrix} \\ \begin{pmatrix} 0 & \{0,0,0\} \\ \{0,0,0\} & 0 \end{pmatrix} & \begin{pmatrix} 0 & \{0,0,0\} \\ \{0,0,0\} & 0 \end{pmatrix} & \begin{pmatrix} \beta_7 + \lambda_{12} & \{0, \beta_2 + \lambda_{15}, \\ \{\beta_5 + \lambda_{17}, 0, \lambda_{19} - \beta_3\} & 0 \end{pmatrix} \\ \begin{pmatrix} 0 & \{0,0,0\} \\ \{0,0,0\} & 0 \end{pmatrix} & \begin{pmatrix} 0 & \{0, -\beta_2 - \lambda_{15}, 0\} \\ \{-\beta_5 - \lambda_{17}, 0, \beta_3 - \lambda_{19}\} & \beta_7 + \lambda_{12} \end{pmatrix} & \begin{pmatrix} \lambda_3 & \{0,0,0\} \\ \{0,0,0\} & \lambda_3 \end{pmatrix} \end{array} \right)$$

```
Δ[X₁[e₄]] =
 Δ[X₁[e₄]] //. {λ₃ → 0, λ₁₂ → -β₇, λ₁₅ → -β₂, λ₁₇ → -β₅, λ₁₉ → β₃}
```

$$\begin{pmatrix} \begin{pmatrix} \alpha_7 & \{0,0,0\} \\ \{0,0,0\} & \alpha_7 \end{pmatrix} & \begin{pmatrix} \rho_5 & \{\lambda_6,\lambda_7,\lambda_8\} \\ \{-\eta_4,0,\lambda_{11}\} & \delta_6 \end{pmatrix} & \begin{pmatrix} \lambda_{20} & \{0,\lambda_{23},0\} \\ \{\lambda_{25},0,\lambda_{27}\} & 0 \end{pmatrix} \\ \begin{pmatrix} \delta_6 & \{-\lambda_6,-\lambda_7,-\lambda_8\} \\ \{\eta_4,0,-\lambda_{11}\} & \rho_5 \end{pmatrix} & \begin{pmatrix} -\alpha_7 & \{0,0,0\} \\ \{0,0,0\} & -\alpha_7 \end{pmatrix} & \begin{pmatrix} -\beta_7 & \{0,-\beta_2,0\} \\ \{-\beta_5,0,\beta_3\} & 0 \end{pmatrix} \\ \begin{pmatrix} 0 & \{0,-\lambda_{23},0\} \\ \{-\lambda_{25},0,-\lambda_{27}\} & \lambda_{20} \end{pmatrix} & \begin{pmatrix} 0 & \{0,\beta_2,0\} \\ \{\beta_5,0,-\beta_3\} & -\beta_7 \end{pmatrix} & \begin{pmatrix} 0 & \{0,0,0\} \\ \{0,0,0\} & 0 \end{pmatrix} \end{pmatrix}$$

```
T[E₁[1], X₁[e₄], E₁[1]] == ZERO
```

True

```
Leib[E₁[1], X₁[e₄], E₁[1]] == ZERO
```

True

```
T[E₂[1], E₂[1], X₁[e₄]] == X₁[e₄]
```

True

```
Δ[X₁[e₄]] - Leib[E₂[1], E₂[1], X₁[e₄]] // Expand
```

$$\begin{pmatrix} \begin{pmatrix} 0 & \{0,0,0\} \\ \{0,0,0\} & 0 \end{pmatrix} & \begin{pmatrix} 0 & \{0,0,0\} \\ \{0,0,0\} & 0 \end{pmatrix} & \begin{pmatrix} \lambda_{20}-\gamma_7 & \{0,\lambda_{23}-\gamma_2,0\} \\ \{\lambda_{25}-\gamma_5,0,\gamma_3+\lambda_{27}\} & 0 \end{pmatrix} \\ \begin{pmatrix} 0 & \{0,0,0\} \\ \{0,0,0\} & 0 \end{pmatrix} & \begin{pmatrix} 0 & \{0,0,0\} \\ \{0,0,0\} & 0 \end{pmatrix} & \begin{pmatrix} 0 & \{0,0,0\} \\ \{0,0,0\} & 0 \end{pmatrix} \\ \begin{pmatrix} 0 & \{0,\gamma_2-\lambda_{23},0\} \\ \{\gamma_5-\lambda_{25},0,-\gamma_3-\lambda_{27}\} & \lambda_{20}-\gamma_7 \end{pmatrix} & \begin{pmatrix} 0 & \{0,0,0\} \\ \{0,0,0\} & 0 \end{pmatrix} & \begin{pmatrix} 0 & \{0,0,0\} \\ \{0,0,0\} & 0 \end{pmatrix} \end{pmatrix}$$

```
Δ[X₁[e₄]] = Δ[X₁[e₄]] //. {λ₂₀ → γ₇, λ₂₃ → γ₂, λ₂₅ → γ₅, λ₂₇ → -γ₃}
```

$$\begin{pmatrix} \begin{pmatrix} \alpha_7 & \{0,0,0\} \\ \{0,0,0\} & \alpha_7 \end{pmatrix} & \begin{pmatrix} \rho_5 & \{\lambda_6,\lambda_7,\lambda_8\} \\ \{-\eta_4,0,\lambda_{11}\} & \delta_6 \end{pmatrix} & \begin{pmatrix} \gamma_7 & \{0,\gamma_2,0\} \\ \{\gamma_5,0,-\gamma_3\} & 0 \end{pmatrix} \\ \begin{pmatrix} \delta_6 & \{-\lambda_6,-\lambda_7,-\lambda_8\} \\ \{\eta_4,0,-\lambda_{11}\} & \rho_5 \end{pmatrix} & \begin{pmatrix} -\alpha_7 & \{0,0,0\} \\ \{0,0,0\} & -\alpha_7 \end{pmatrix} & \begin{pmatrix} -\beta_7 & \{0,-\beta_2,0\} \\ \{-\beta_5,0,\beta_3\} & 0 \end{pmatrix} \\ \begin{pmatrix} 0 & \{0,-\gamma_2,0\} \\ \{-\gamma_5,0,\gamma_3\} & \gamma_7 \end{pmatrix} & \begin{pmatrix} 0 & \{0,\beta_2,0\} \\ \{\beta_5,0,-\beta_3\} & -\beta_7 \end{pmatrix} & \begin{pmatrix} 0 & \{0,0,0\} \\ \{0,0,0\} & 0 \end{pmatrix} \end{pmatrix}$$

```
T[E₂[1], X₁[e₄], E₂[1]] == ZERO
```

True

```
Leib[E₂[1], X₁[e₄], E₂[1]] == ZERO
```

True

```
T[E₁[1], E₂[1], X₁[e₄]] == ZERO
```

True

```
(Leib[E₁[1], E₂[1], X₁[e₄]] // Expand) == ZERO
```

True

```
T[E₂[1], E₁[1], X₁[e₄]] == ZERO
```

True

```
(Leib[E₂[1], E₁[1], X₁[e₄]] // Expand) == ZERO
```

True

```
T[E₁[1], X₁[e₄], E₂[1]] == X₁[e₄]
```

True

```
Δ[X₁[e₄]] - Leib[E₁[1], X₁[e₄], E₂[1]] == ZERO
```

True

```
T[E₁[1], X₁[e₁], X₁[e₄]] == ZERO
```

True

```
(Leib[E₁[1], X₁[e₁], X₁[e₄]] // Expand) == ZERO
```

True

```
T[X₁[e₁], E₁[1], X₁[e₄]] == ZERO
```

True

```
(Leib[X₁[e₁], E₁[1], X₁[e₄]] // Expand) == ZERO
```

True

```
T[E₂[1], X₁[e₁], X₁[e₄]] == ZERO
```

True

```
(Leib[E₂[1], X₁[e₁], X₁[e₄]] // Expand) == ZERO
```

True

```
T[X₁[e₄], X₁[e₁], X₁[e₁]] == ZERO
```

True

```
(Leib[X₁[e₄], X₁[e₁], X₁[e₁]] // Expand) == ZERO
```

True

```
T[X₁[e₁], X₁[e₄], X₁[e₁]] == ZERO
```

True

```
(Leib[X₁[e₁], X₁[e₄], X₁[e₁]] // Expand) == ZERO
```

True

```
T[X₁[e₄], X₁[e₃], X₁[e₁]] == ZERO
```

True

```
(Leib[X₁[e₄], X₁[e₃], X₁[e₁]] // Expand) == ZERO
```

True

```
T[X₁[e₃], X₁[e₄], X₁[e₄]] == ZERO
```

True

```
(Leib[X₁[e₃], X₁[e₄], X₁[e₄]] // Expand) == ZERO
```

True

```
T[X₁[e₄], X₁[e₄], X₁[e₄]] == ZERO
```

True

```
(Leib[X₁[e₄], X₁[e₄], X₁[e₄]] // Expand) == ZERO
```

True

```
T[X₁[e₄], X₁[e₁], X₁[e₃]] == ZERO
```

True

```
(Leib[X₁[e₄], X₁[e₁], X₁[e₃]] // Expand) == ZERO
```

True

$$\Delta[X_1[e_4]] = \Delta[X_1[e_4]] \;//.\; \{\lambda_6 \to \epsilon_1,\; \lambda_7 \to \epsilon_2,\; \lambda_8 \to \epsilon_3,\; \lambda_{11} \to \epsilon_4\}$$

$$\begin{pmatrix}
\begin{pmatrix} \alpha_7 & \{0,0,0\} \\ \{0,0,0\} & \alpha_7 \end{pmatrix} & \begin{pmatrix} \rho_5 & \{\epsilon_1,\epsilon_2,\epsilon_3\} \\ \{-\eta_4,0,\epsilon_4\} & \delta_6 \end{pmatrix} & \begin{pmatrix} \gamma_7 & \{0,\gamma_2,0\} \\ \{\gamma_5,0,-\gamma_3\} & 0 \end{pmatrix} \\
\begin{pmatrix} \delta_6 & \{-\epsilon_1,-\epsilon_2,-\epsilon_3\} \\ \{\eta_4,0,-\epsilon_4\} & \rho_5 \end{pmatrix} & \begin{pmatrix} -\alpha_7 & \{0,0,0\} \\ \{0,0,0\} & -\alpha_7 \end{pmatrix} & \begin{pmatrix} -\beta_7 & \{0,-\beta_2,0\} \\ \{-\beta_5,0,\beta_3\} & 0 \end{pmatrix} \\
\begin{pmatrix} 0 & \{0,-\gamma_2,0\} \\ \{-\gamma_5,0,\gamma_3\} & \gamma_7 \end{pmatrix} & \begin{pmatrix} 0 & \{0,\beta_2,0\} \\ \{\beta_5,0,-\beta_3\} & -\beta_7 \end{pmatrix} & \begin{pmatrix} 0 & \{0,0,0\} \\ \{0,0,0\} & 0 \end{pmatrix}
\end{pmatrix}$$

- $X_1[e_5]$

$\Delta[X_1[e_5]] = \text{generic};$

$U_{X_1[e_5]}[E_1[1]] == \text{ZERO}$

True

$\text{defin}[X_1[e_5],\, E_1[1]] \;//\; \text{Expand}$

$$\begin{pmatrix}
\begin{pmatrix} 0 & \{0,0,0\} \\ \{0,0,0\} & 0 \end{pmatrix} & \begin{pmatrix} 0 & \{0,0,\lambda_1-\alpha_8\} \\ \{0,0,0\} & 0 \end{pmatrix} & \begin{pmatrix} 0 & \{0,0,0\} \\ \{0,0,0\} & 0 \end{pmatrix} \\
\begin{pmatrix} 0 & \{0,0,\alpha_8-\lambda_1\} \\ \{0,0,0\} & 0 \end{pmatrix} & \begin{pmatrix} -\lambda_{11} & \{0,0,0\} \\ \{0,0,0\} & -\lambda_{11} \end{pmatrix} & \begin{pmatrix} -\lambda_{27} & \{0,0,-\lambda_{21}\} \\ \{\lambda_{23},-\lambda_{22},0\} & 0 \end{pmatrix} \\
\begin{pmatrix} 0 & \{0,0,0\} \\ \{0,0,0\} & 0 \end{pmatrix} & \begin{pmatrix} 0 & \{0,0,\lambda_{21}\} \\ \{-\lambda_{23},\lambda_{22},0\} & -\lambda_{27} \end{pmatrix} & \begin{pmatrix} 0 & \{0,0,0\} \\ \{0,0,0\} & 0 \end{pmatrix}
\end{pmatrix}$$

$\Delta[X_1[e_5]] =$
$\Delta[X_1[e_5]] \;//.\; \{\lambda_1 \to \alpha_8,\; \lambda_{11} \to 0,\; \lambda_{21} \to 0,\; \lambda_{22} \to 0,\; \lambda_{23} \to 0,\; \lambda_{27} \to 0\}$

$$\begin{pmatrix}
\begin{pmatrix} \alpha_8 & \{0,0,0\} \\ \{0,0,0\} & \alpha_8 \end{pmatrix} & \begin{pmatrix} \lambda_4 & \{\lambda_6,\lambda_7,\lambda_8\} \\ \{\lambda_9,\lambda_{10},0\} & \lambda_5 \end{pmatrix} & \begin{pmatrix} \lambda_{20} & \{0,0,\lambda\ldots\} \\ \{\lambda_{25},\lambda_{26},0\} & 0 \end{pmatrix} \\
\begin{pmatrix} \lambda_5 & \{-\lambda_6,-\lambda_7,-\lambda_8\} \\ \{-\lambda_9,-\lambda_{10},0\} & \lambda_4 \end{pmatrix} & \begin{pmatrix} \lambda_2 & \{0,0,0\} \\ \{0,0,0\} & \lambda_2 \end{pmatrix} & \begin{pmatrix} \lambda_{12} & \{\lambda_{14},\lambda_{15}\ldots\} \\ \{\lambda_{17},\lambda_{18},\lambda_{19}\} & \lambda_{13} \end{pmatrix} \\
\begin{pmatrix} 0 & \{0,0,-\lambda_{24}\} \\ \{-\lambda_{25},-\lambda_{26},0\} & \lambda_{20} \end{pmatrix} & \begin{pmatrix} \lambda_{13} & \{-\lambda_{14},-\lambda_{15},-\lambda_{16}\} \\ \{-\lambda_{17},-\lambda_{18},-\lambda_{19}\} & \lambda_{12} \end{pmatrix} & \begin{pmatrix} \lambda_3 & \{0,0,0\} \\ \{0,0,0\} & \lambda_3 \end{pmatrix}
\end{pmatrix}$$

$U_{X_1[e_5]}[E_2[1]] == \text{ZERO}$

True

$\text{defin}[X_1[e_5],\, E_2[1]] \;//\; \text{Expand}$

$$\begin{pmatrix}
\begin{pmatrix} 0 & \{0,0,0\} \\ \{0,0,0\} & 0 \end{pmatrix} & \begin{pmatrix} 0 & \{0,0,\alpha_8+\lambda_2\} \\ \{0,0,0\} & 0 \end{pmatrix} & \begin{pmatrix} \lambda_{19} & \{0,0,\lambda_{13}\} \\ \{-\lambda_{15},\lambda_{14},0\} & 0 \end{pmatrix} \\
\begin{pmatrix} 0 & \{0,0,-\alpha_8-\lambda_2\} \\ \{0,0,0\} & 0 \end{pmatrix} & \begin{pmatrix} 0 & \{0,0,0\} \\ \{0,0,0\} & 0 \end{pmatrix} & \begin{pmatrix} 0 & \{0,0,0\} \\ \{0,0,0\} & 0 \end{pmatrix} \\
\begin{pmatrix} 0 & \{0,0,-\lambda_{13}\} \\ \{\lambda_{15},-\lambda_{14},0\} & \lambda_{19} \end{pmatrix} & \begin{pmatrix} 0 & \{0,0,0\} \\ \{0,0,0\} & 0 \end{pmatrix} & \begin{pmatrix} 0 & \{0,0,0\} \\ \{0,0,0\} & 0 \end{pmatrix}
\end{pmatrix}$$

$\Delta[X_1[e_5]] = \Delta[X_1[e_5]] \,//.\, \{\lambda_2 \to -\alpha_8,\, \lambda_{13} \to 0,\, \lambda_{14} \to 0,\, \lambda_{15} \to 0,\, \lambda_{19} \to 0\}$

$$\begin{pmatrix} \begin{pmatrix} \alpha_8 & \{0,0,0\} \\ \{0,0,0\} & \alpha_8 \end{pmatrix} & \begin{pmatrix} \lambda_4 & \{\lambda_6,\lambda_7,\lambda_8\} \\ \{\lambda_9,\lambda_{10},0\} & \lambda_5 \end{pmatrix} & \begin{pmatrix} \lambda_{20} & \{0,0,\lambda_{24}\} \\ \{\lambda_{25},\lambda_{26},0\} & 0 \end{pmatrix} \\ \begin{pmatrix} \lambda_5 & \{-\lambda_6,-\lambda_7,-\lambda_8\} \\ \{-\lambda_9,-\lambda_{10},0\} & \lambda_4 \end{pmatrix} & \begin{pmatrix} -\alpha_8 & \{0,0,0\} \\ \{0,0,0\} & -\alpha_8 \end{pmatrix} & \begin{pmatrix} \lambda_{12} & \{0,0,\lambda_{16}\} \\ \{\lambda_{17},\lambda_{18},0\} & 0 \end{pmatrix} \\ \begin{pmatrix} 0 & \{0,0,-\lambda_{24}\} \\ \{-\lambda_{25},-\lambda_{26},0\} & \lambda_{20} \end{pmatrix} & \begin{pmatrix} 0 & \{0,0,-\lambda_{16}\} \\ \{-\lambda_{17},-\lambda_{18},0\} & \lambda_{12} \end{pmatrix} & \begin{pmatrix} \lambda_3 & \{0,0,0\} \\ \{0,0,0\} & \lambda_3 \end{pmatrix} \end{pmatrix}$$

$U_{X_1[e_5]}[E_3[1]] == \text{ZERO}$

True

$(\text{defin}[X_1[e_5], E_3[1]] \,//\, \text{Expand}) == \text{ZERO}$

True

$U_{X_1[e_5]}[X_1[e_5]] == \text{ZERO}$

True

$(\text{defin}[X_1[e_5], X_1[e_5]] \,//\, \text{Expand}) == \text{ZERO}$

True

$U_{X_1[e_5]}[X_1[e_4]] == \text{ZERO}$

True

$\text{defin}[X_1[e_5], X_1[e_4]] \,//\, \text{Expand}$

$$\begin{pmatrix} \begin{pmatrix} 0 & \{0,0,0\} \\ \{0,0,0\} & 0 \end{pmatrix} & \begin{pmatrix} 0 & \{0,0,-\epsilon_4-\lambda_{10}\} \\ \{0,0,0\} & 0 \end{pmatrix} & \begin{pmatrix} 0 & \{0,0,0\} \\ \{0,0,0\} & 0 \end{pmatrix} \\ \begin{pmatrix} 0 & \{0,0,\epsilon_4+\lambda_{10}\} \\ \{0,0,0\} & 0 \end{pmatrix} & \begin{pmatrix} 0 & \{0,0,0\} \\ \{0,0,0\} & 0 \end{pmatrix} & \begin{pmatrix} 0 & \{0,0,0\} \\ \{0,0,0\} & 0 \end{pmatrix} \\ \begin{pmatrix} 0 & \{0,0,0\} \\ \{0,0,0\} & 0 \end{pmatrix} & \begin{pmatrix} 0 & \{0,0,0\} \\ \{0,0,0\} & 0 \end{pmatrix} & \begin{pmatrix} 0 & \{0,0,0\} \\ \{0,0,0\} & 0 \end{pmatrix} \end{pmatrix}$$

$\Delta[X_1[e_5]] = \Delta[X_1[e_5]] \,//.\, \{\lambda_{10} \to -\epsilon_4\}$

$$\begin{pmatrix} \begin{pmatrix} \alpha_8 & \{0,0,0\} \\ \{0,0,0\} & \alpha_8 \end{pmatrix} & \begin{pmatrix} \lambda_4 & \{\lambda_6,\lambda_7,\lambda_8\} \\ \{\lambda_9,-\epsilon_4,0\} & \lambda_5 \end{pmatrix} & \begin{pmatrix} \lambda_{20} & \{0,0,\lambda_{24}\} \\ \{\lambda_{25},\lambda_{26},0\} & 0 \end{pmatrix} \\ \begin{pmatrix} \lambda_5 & \{-\lambda_6,-\lambda_7,-\lambda_8\} \\ \{-\lambda_9,\epsilon_4,0\} & \lambda_4 \end{pmatrix} & \begin{pmatrix} -\alpha_8 & \{0,0,0\} \\ \{0,0,0\} & -\alpha_8 \end{pmatrix} & \begin{pmatrix} \lambda_{12} & \{0,0,\lambda_{16}\} \\ \{\lambda_{17},\lambda_{18},0\} & 0 \end{pmatrix} \\ \begin{pmatrix} 0 & \{0,0,-\lambda_{24}\} \\ \{-\lambda_{25},-\lambda_{26},0\} & \lambda_{20} \end{pmatrix} & \begin{pmatrix} 0 & \{0,0,-\lambda_{16}\} \\ \{-\lambda_{17},-\lambda_{18},0\} & \lambda_{12} \end{pmatrix} & \begin{pmatrix} \lambda_3 & \{0,0,0\} \\ \{0,0,0\} & \lambda_3 \end{pmatrix} \end{pmatrix}$$

$U_{X_1[e_5]}[X_1[e_3]] == \text{ZERO}$

True

**defin[X$_1$[e$_5$], X$_1$[e$_3$]] // Expand**

$$\begin{pmatrix} \begin{pmatrix} 0 & \{0,0,0\} \\ \{0,0,0\} & 0 \end{pmatrix} & \begin{pmatrix} 0 & \{0,0,-\eta_5-\lambda_9\} \\ \{0,0,0\} & 0 \end{pmatrix} & \begin{pmatrix} 0 & \{0,0,0\} \\ \{0,0,0\} & 0 \end{pmatrix} \\ \begin{pmatrix} 0 & \{0,0,\eta_5+\lambda_9\} \\ \{0,0,0\} & 0 \end{pmatrix} & \begin{pmatrix} 0 & \{0,0,0\} \\ \{0,0,0\} & 0 \end{pmatrix} & \begin{pmatrix} 0 & \{0,0,0\} \\ \{0,0,0\} & 0 \end{pmatrix} \\ \begin{pmatrix} 0 & \{0,0,0\} \\ \{0,0,0\} & 0 \end{pmatrix} & \begin{pmatrix} 0 & \{0,0,0\} \\ \{0,0,0\} & 0 \end{pmatrix} & \begin{pmatrix} 0 & \{0,0,0\} \\ \{0,0,0\} & 0 \end{pmatrix} \end{pmatrix}$$

**Δ[X$_1$[e$_5$]] = Δ[X$_1$[e$_5$]] //. {λ$_9$ → -η$_5$}**

$$\begin{pmatrix} \begin{pmatrix} \alpha_8 & \{0,0,0\} \\ \{0,0,0\} & \alpha_8 \end{pmatrix} & \begin{pmatrix} \lambda_4 & \{\lambda_6,\lambda_7,\lambda_8\} \\ \{-\eta_5,-\epsilon_4,0\} & \lambda_5 \end{pmatrix} & \begin{pmatrix} \lambda_{20} & \{0,0,\lambda_{24}\} \\ \{\lambda_{25},\lambda_{26},0\} & 0 \end{pmatrix} \\ \begin{pmatrix} \lambda_5 & \{-\lambda_6,-\lambda_7,-\lambda_8\} \\ \{\eta_5,\epsilon_4,0\} & \lambda_4 \end{pmatrix} & \begin{pmatrix} -\alpha_8 & \{0,0,0\} \\ \{0,0,0\} & -\alpha_8 \end{pmatrix} & \begin{pmatrix} \lambda_{12} & \{0,0,\lambda_{16}\} \\ \{\lambda_{17},\lambda_{18},0\} & 0 \end{pmatrix} \\ \begin{pmatrix} 0 & \{0,0,-\lambda_{24}\} \\ \{-\lambda_{25},-\lambda_{26},0\} & \lambda_{20} \end{pmatrix} & \begin{pmatrix} 0 & \{0,0,-\lambda_{16}\} \\ \{-\lambda_{17},-\lambda_{18},0\} & \lambda_{12} \end{pmatrix} & \begin{pmatrix} \lambda_3 & \{0,0,0\} \\ \{0,0,0\} & \lambda_3 \end{pmatrix} \end{pmatrix}$$

**U$_{X_1[e_5]}$[X$_1$[e$_2$]] == ZERO**

True

**defin[X$_1$[e$_5$], X$_1$[e$_2$]] // Expand**

$$\begin{pmatrix} \begin{pmatrix} 0 & \{0,0,0\} \\ \{0,0,0\} & 0 \end{pmatrix} & \begin{pmatrix} 0 & \{0,0,\lambda_4-\rho_6\} \\ \{0,0,0\} & 0 \end{pmatrix} & \begin{pmatrix} 0 & \{0,0,0\} \\ \{0,0,0\} & 0 \end{pmatrix} \\ \begin{pmatrix} 0 & \{0,0,\rho_6-\lambda_4\} \\ \{0,0,0\} & 0 \end{pmatrix} & \begin{pmatrix} 0 & \{0,0,0\} \\ \{0,0,0\} & 0 \end{pmatrix} & \begin{pmatrix} 0 & \{0,0,0\} \\ \{0,0,0\} & 0 \end{pmatrix} \\ \begin{pmatrix} 0 & \{0,0,0\} \\ \{0,0,0\} & 0 \end{pmatrix} & \begin{pmatrix} 0 & \{0,0,0\} \\ \{0,0,0\} & 0 \end{pmatrix} & \begin{pmatrix} 0 & \{0,0,0\} \\ \{0,0,0\} & 0 \end{pmatrix} \end{pmatrix}$$

**Δ[X$_1$[e$_5$]] = Δ[X$_1$[e$_5$]] //. {λ$_4$ → ρ$_6$}**

$$\begin{pmatrix} \begin{pmatrix} \alpha_8 & \{0,0,0\} \\ \{0,0,0\} & \alpha_8 \end{pmatrix} & \begin{pmatrix} \rho_6 & \{\lambda_6,\lambda_7,\lambda_8\} \\ \{-\eta_5,-\epsilon_4,0\} & \lambda_5 \end{pmatrix} & \begin{pmatrix} \lambda_{20} & \{0,0,\lambda_{24}\} \\ \{\lambda_{25},\lambda_{26},0\} & 0 \end{pmatrix} \\ \begin{pmatrix} \lambda_5 & \{-\lambda_6,-\lambda_7,-\lambda_8\} \\ \{\eta_5,\epsilon_4,0\} & \rho_6 \end{pmatrix} & \begin{pmatrix} -\alpha_8 & \{0,0,0\} \\ \{0,0,0\} & -\alpha_8 \end{pmatrix} & \begin{pmatrix} \lambda_{12} & \{0,0,\lambda_{16}\} \\ \{\lambda_{17},\lambda_{18},0\} & 0 \end{pmatrix} \\ \begin{pmatrix} 0 & \{0,0,-\lambda_{24}\} \\ \{-\lambda_{25},-\lambda_{26},0\} & \lambda_{20} \end{pmatrix} & \begin{pmatrix} 0 & \{0,0,-\lambda_{16}\} \\ \{-\lambda_{17},-\lambda_{18},0\} & \lambda_{12} \end{pmatrix} & \begin{pmatrix} \lambda_3 & \{0,0,0\} \\ \{0,0,0\} & \lambda_3 \end{pmatrix} \end{pmatrix}$$

**U$_{X_1[e_5]}$[X$_1$[e$_1$]] == ZERO**

True

**defin[X$_1$[e$_5$], X$_1$[e$_1$]] // Expand**

$$\begin{pmatrix} \begin{pmatrix} 0 & \{0,0,0\} \\ \{0,0,0\} & 0 \end{pmatrix} & \begin{pmatrix} 0 & \{0,0,\lambda_5-\delta_7\} \\ \{0,0,0\} & 0 \end{pmatrix} & \begin{pmatrix} 0 & \{0,0,0\} \\ \{0,0,0\} & 0 \end{pmatrix} \\ \begin{pmatrix} 0 & \{0,0,\delta_7-\lambda_5\} \\ \{0,0,0\} & 0 \end{pmatrix} & \begin{pmatrix} 0 & \{0,0,0\} \\ \{0,0,0\} & 0 \end{pmatrix} & \begin{pmatrix} 0 & \{0,0,0\} \\ \{0,0,0\} & 0 \end{pmatrix} \\ \begin{pmatrix} 0 & \{0,0,0\} \\ \{0,0,0\} & 0 \end{pmatrix} & \begin{pmatrix} 0 & \{0,0,0\} \\ \{0,0,0\} & 0 \end{pmatrix} & \begin{pmatrix} 0 & \{0,0,0\} \\ \{0,0,0\} & 0 \end{pmatrix} \end{pmatrix}$$

**Δ[X$_1$[e$_5$]] = Δ[X$_1$[e$_5$]] //. {λ$_5$ → δ$_7$}**

$$\begin{pmatrix} \begin{pmatrix} \alpha_8 & \{0,0,0\} \\ \{0,0,0\} & \alpha_8 \end{pmatrix} & \begin{pmatrix} \rho_6 & \{\lambda_6,\lambda_7,\lambda_8\} \\ \{-\eta_5,-\epsilon_4,0\} & \delta_7 \end{pmatrix} & \begin{pmatrix} \lambda_{20} & \{0,0,\lambda_{24}\} \\ \{\lambda_{25},\lambda_{26},0\} & 0 \end{pmatrix} \\ \begin{pmatrix} \delta_7 & \{-\lambda_6,-\lambda_7,-\lambda_8\} \\ \{\eta_5,\epsilon_4,0\} & \rho_6 \end{pmatrix} & \begin{pmatrix} -\alpha_8 & \{0,0,0\} \\ \{0,0,0\} & -\alpha_8 \end{pmatrix} & \begin{pmatrix} \lambda_{12} & \{0,0,\lambda_{16}\} \\ \{\lambda_{17},\lambda_{18},0\} & 0 \end{pmatrix} \\ \begin{pmatrix} 0 & \{0,0,-\lambda_{24}\} \\ \{-\lambda_{25},-\lambda_{26},0\} & \lambda_{20} \end{pmatrix} & \begin{pmatrix} 0 & \{0,0,-\lambda_{16}\} \\ \{-\lambda_{17},-\lambda_{18},0\} & \lambda_{12} \end{pmatrix} & \begin{pmatrix} \lambda_3 & \{0,0,0\} \\ \{0,0,0\} & \lambda_3 \end{pmatrix} \end{pmatrix}$$

**T[E$_1$[1], E$_1$[1], X$_1$[e$_5$]] == X$_1$[e$_5$]**

True

**Δ[X$_1$[e$_5$]] - Leib[E$_1$[1], E$_1$[1], X$_1$[e$_5$]] // Expand**

$$\begin{pmatrix} \begin{pmatrix} 0 & \{0,0,0\} \\ \{0,0,0\} & 0 \end{pmatrix} & \begin{pmatrix} 0 & \{0,0,0\} \\ \{0,0,0\} & 0 \end{pmatrix} & \begin{pmatrix} 0 & \{0,0,0\} \\ \{0,0,0\} & 0 \end{pmatrix} \\ \begin{pmatrix} 0 & \{0,0,0\} \\ \{0,0,0\} & 0 \end{pmatrix} & \begin{pmatrix} 0 & \{0,0,0\} \\ \{0,0,0\} & 0 \end{pmatrix} & \begin{pmatrix} \beta_8+\lambda_{12} & \{0,0,\beta_2+\lambda_{16}\} \\ \{\lambda_{17}-\beta_4,\beta_3+\lambda_{18},0\} & 0 \end{pmatrix} \\ \begin{pmatrix} 0 & \{0,0,0\} \\ \{0,0,0\} & 0 \end{pmatrix} & \begin{pmatrix} 0 & \{0,0,-\beta_2-\lambda_{16}\} \\ \{\beta_4-\lambda_{17},-\beta_3-\lambda_{18},0\} & \beta_8+\lambda_{12} \end{pmatrix} & \begin{pmatrix} \lambda_3 & \{0,0,0\} \\ \{0,0,0\} & \lambda_3 \end{pmatrix} \end{pmatrix}$$

**Δ[X$_1$[e$_5$]] =**
**Δ[X$_1$[e$_5$]] //. {λ$_3$ → 0, λ$_{12}$ → -β$_8$, λ$_{16}$ → -β$_2$, λ$_{17}$ → β$_4$, λ$_{18}$ → -β$_3$}**

$$\begin{pmatrix} \begin{pmatrix} \alpha_8 & \{0,0,0\} \\ \{0,0,0\} & \alpha_8 \end{pmatrix} & \begin{pmatrix} \rho_6 & \{\lambda_6,\lambda_7,\lambda_8\} \\ \{-\eta_5,-\epsilon_4,0\} & \delta_7 \end{pmatrix} & \begin{pmatrix} \lambda_{20} & \{0,0,\lambda_{24}\} \\ \{\lambda_{25},\lambda_{26},0\} & 0 \end{pmatrix} \\ \begin{pmatrix} \delta_7 & \{-\lambda_6,-\lambda_7,-\lambda_8\} \\ \{\eta_5,\epsilon_4,0\} & \rho_6 \end{pmatrix} & \begin{pmatrix} -\alpha_8 & \{0,0,0\} \\ \{0,0,0\} & -\alpha_8 \end{pmatrix} & \begin{pmatrix} -\beta_8 & \{0,0,-\beta_2\} \\ \{\beta_4,-\beta_3,0\} & 0 \end{pmatrix} \\ \begin{pmatrix} 0 & \{0,0,-\lambda_{24}\} \\ \{-\lambda_{25},-\lambda_{26},0\} & \lambda_{20} \end{pmatrix} & \begin{pmatrix} 0 & \{0,0,\beta_2\} \\ \{-\beta_4,\beta_3,0\} & -\beta_8 \end{pmatrix} & \begin{pmatrix} 0 & \{0,0,0\} \\ \{0,0,0\} & 0 \end{pmatrix} \end{pmatrix}$$

**T[E$_1$[1], X$_1$[e$_5$], E$_1$[1]] == ZERO**

True

**Leib[E$_1$[1], X$_1$[e$_5$], E$_1$[1]] == ZERO**

True

**T[E₂[1], E₂[1], X₁[e₅]] == X₁[e₅]**

True

**Δ[X₁[e₅]] - Leib[E₂[1], E₂[1], X₁[e₅]] // Expand**

$$\begin{pmatrix} \begin{pmatrix} 0 & \{0,0,0\} \\ \{0,0,0\} & 0 \end{pmatrix} & \begin{pmatrix} 0 & \{0,0,0\} \\ \{0,0,0\} & 0 \end{pmatrix} & \begin{pmatrix} \lambda_{20}-\gamma_8 & \{0,0,\lambda_{24}-\gamma_2\} \\ \{\gamma_4+\lambda_{25},\lambda_{26}-\gamma_3,0\} & 0 \end{pmatrix} \\ \begin{pmatrix} 0 & \{0,0,0\} \\ \{0,0,0\} & 0 \end{pmatrix} & \begin{pmatrix} 0 & \{0,0,0\} \\ \{0,0,0\} & 0 \end{pmatrix} & \begin{pmatrix} 0 & \{0,0,0\} \\ \{0,0,0\} & 0 \end{pmatrix} \\ \begin{pmatrix} 0 & \{0,0,\gamma_2-\lambda_{24}\} \\ \{-\gamma_4-\lambda_{25},\gamma_3-\lambda_{26},0\} & \lambda_{20}-\gamma_8 \end{pmatrix} & \begin{pmatrix} 0 & \{0,0,0\} \\ \{0,0,0\} & 0 \end{pmatrix} & \begin{pmatrix} 0 & \{0,0,0\} \\ \{0,0,0\} & 0 \end{pmatrix} \end{pmatrix}$$

**Δ[X₁[e₅]] = Δ[X₁[e₅]] //. {λ₂₀ → γ₈, λ₂₄ → γ₂, λ₂₅ → -γ₄, λ₂₆ → γ₃}**

$$\begin{pmatrix} \begin{pmatrix} \alpha_8 & \{0,0,0\} \\ \{0,0,0\} & \alpha_8 \end{pmatrix} & \begin{pmatrix} \rho_6 & \{\lambda_6,\lambda_7,\lambda_8\} \\ \{-\eta_5,-\epsilon_4,0\} & \delta_7 \end{pmatrix} & \begin{pmatrix} \gamma_8 & \{0,0,\gamma_2\} \\ \{-\gamma_4,\gamma_3,0\} & 0 \end{pmatrix} \\ \begin{pmatrix} \delta_7 & \{-\lambda_6,-\lambda_7,-\lambda_8\} \\ \{\eta_5,\epsilon_4,0\} & \rho_6 \end{pmatrix} & \begin{pmatrix} -\alpha_8 & \{0,0,0\} \\ \{0,0,0\} & -\alpha_8 \end{pmatrix} & \begin{pmatrix} -\beta_8 & \{0,0,-\beta_2\} \\ \{\beta_4,-\beta_3,0\} & 0 \end{pmatrix} \\ \begin{pmatrix} 0 & \{0,0,-\gamma_2\} \\ \{\gamma_4,-\gamma_3,0\} & \gamma_8 \end{pmatrix} & \begin{pmatrix} 0 & \{0,0,\beta_2\} \\ \{-\beta_4,\beta_3,0\} & -\beta_8 \end{pmatrix} & \begin{pmatrix} 0 & \{0,0,0\} \\ \{0,0,0\} & 0 \end{pmatrix} \end{pmatrix}$$

**T[E₂[1], X₁[e₅], E₂[1]] == ZERO**

True

**Leib[E₂[1], X₁[e₅], E₂[1]] == ZERO**

True

**T[E₁[1], E₂[1], X₁[e₅]] == ZERO**

True

**(Leib[E₁[1], E₂[1], X₁[e₅]] // Expand) == ZERO**

True

**T[E₂[1], E₁[1], X₁[e₅]] == ZERO**

True

**(Leib[E₂[1], E₁[1], X₁[e₅]] // Expand) == ZERO**

True

**T[E₁[1], X₁[e₅], E₂[1]] == X₁[e₅]**

True

```
Δ[X₁[e₅]] - Leib[E₁[1], X₁[e₅], E₂[1]] == ZERO
```
True

```
T[E₁[1], X₁[e₁], X₁[e₅]] == ZERO
```
True

```
(Leib[E₁[1], X₁[e₁], X₁[e₅]] // Expand) == ZERO
```
True

```
T[X₁[e₁], E₁[1], X₁[e₅]] == ZERO
```
True

```
(Leib[X₁[e₁], E₁[1], X₁[e₅]] // Expand) == ZERO
```
True

```
T[E₂[1], X₁[e₁], X₁[e₅]] == ZERO
```
True

```
(Leib[E₂[1], X₁[e₁], X₁[e₅]] // Expand) == ZERO
```
True

```
T[X₁[e₅], X₁[e₁], X₁[e₁]] == ZERO
```
True

```
(Leib[X₁[e₅], X₁[e₁], X₁[e₁]] // Expand) == ZERO
```
True

```
T[X₁[e₅], X₁[e₄], X₁[e₁]] == ZERO
```
True

```
(Leib[X₁[e₅], X₁[e₄], X₁[e₁]] // Expand) == ZERO
```
True

```
T[X₁[e₄], X₁[e₅], X₁[e₂]] == ZERO
```
True

```
(Leib[X₁[e₄], X₁[e₅], X₁[e₂]] // Expand) == ZERO
```

True

```
T[X₁[e₃], X₁[e₅], X₁[e₅]] == ZERO
```

True

```
(Leib[X₁[e₃], X₁[e₅], X₁[e₅]] // Expand) == ZERO
```

True

```
T[X₁[e₅], X₁[e₅], X₁[e₅]] == ZERO
```

True

```
(Leib[X₁[e₅], X₁[e₅], X₁[e₅]] // Expand) == ZERO
```

True

```
T[X₁[e₄], X₁[e₅], X₁[e₃]] == ZERO
```

True

```
(Leib[X₁[e₄], X₁[e₅], X₁[e₃]] // Expand) == ZERO
```

True

$$\Delta[X_1[e_5]] = \Delta[X_1[e_5]] //. \{\lambda_6 \to \phi_1, \lambda_7 \to \phi_2, \lambda_8 \to \phi_3\}$$

$$\begin{pmatrix} \begin{pmatrix} \alpha_8 & \{0,0,0\} \\ \{0,0,0\} & \alpha_8 \end{pmatrix} & \begin{pmatrix} \rho_6 & \{\phi_1,\phi_2,\phi_3\} \\ \{-\eta_5,-\epsilon_4,0\} & \delta_7 \end{pmatrix} & \begin{pmatrix} \gamma_8 & \{0,0,\gamma_2\} \\ \{-\gamma_4,\gamma_3,0\} & 0 \end{pmatrix} \\ \begin{pmatrix} \delta_7 & \{-\phi_1,-\phi_2,-\phi_3\} \\ \{\eta_5,\epsilon_4,0\} & \rho_6 \end{pmatrix} & \begin{pmatrix} -\alpha_8 & \{0,0,0\} \\ \{0,0,0\} & -\alpha_8 \end{pmatrix} & \begin{pmatrix} -\beta_8 & \{0,0,-\beta_2\} \\ \{\beta_4,-\beta_3,0\} & 0 \end{pmatrix} \\ \begin{pmatrix} 0 & \{0,0,-\gamma_2\} \\ \{\gamma_4,-\gamma_3,0\} & \gamma_8 \end{pmatrix} & \begin{pmatrix} 0 & \{0,0,\beta_2\} \\ \{-\beta_4,\beta_3,0\} & -\beta_8 \end{pmatrix} & \begin{pmatrix} 0 & \{0,0,0\} \\ \{0,0,0\} & 0 \end{pmatrix} \end{pmatrix}$$

- $X_1[e_6]$

  ```
  Δ[X₁[e₆]] = generic;
  U_{X₁[e₆]}[E₁[1]] == ZERO
  ```

  True

**defin[X$_1$[e$_6$], E$_1$[1]] // Expand**

$$\begin{pmatrix} \begin{pmatrix} 0 & \{0,0,0\} \\ \{0,0,0\} & 0 \end{pmatrix} & \begin{pmatrix} 0 & \{0,0,0\} \\ \{\lambda_1 - \alpha_3, 0, 0\} & 0 \end{pmatrix} & \begin{pmatrix} 0 & \{0,0,0\} \\ \{0,0,0\} & 0 \end{pmatrix} \\ \begin{pmatrix} 0 & \{0,0,0\} \\ \{\alpha_3 - \lambda_1, 0, 0\} & 0 \end{pmatrix} & \begin{pmatrix} -\lambda_6 & \{0,0,0\} \\ \{0,0,0\} & -\lambda_6 \end{pmatrix} & \begin{pmatrix} 0 & \{0, -\lambda_{27}, \lambda_{26}\} \\ \{-\lambda_{20}, 0, 0\} & -\lambda_{22} \end{pmatrix} \\ \begin{pmatrix} 0 & \{0,0,0\} \\ \{0,0,0\} & 0 \end{pmatrix} & \begin{pmatrix} -\lambda_{22} & \{0, \lambda_{27}, -\lambda_{26}\} \\ \{\lambda_{20}, 0, 0\} & 0 \end{pmatrix} & \begin{pmatrix} 0 & \{0,0,0\} \\ \{0,0,0\} & 0 \end{pmatrix} \end{pmatrix}$$

**Δ[X$_1$[e$_6$]] =**
 **Δ[X$_1$[e$_6$]] //. {λ$_1$ → α$_3$, λ$_6$ → 0, λ$_{20}$ → 0, λ$_{22}$ → 0, λ$_{26}$ → 0, λ$_{27}$ → 0}**

$$\begin{pmatrix} \begin{pmatrix} \alpha_3 & \{0,0,0\} \\ \{0,0,0\} & \alpha_3 \end{pmatrix} & \begin{pmatrix} \lambda_4 & \{0, \lambda_7, \lambda_8\} \\ \{\lambda_9, \lambda_{10}, \lambda_{11}\} & \lambda_5 \end{pmatrix} & \begin{pmatrix} 0 & \{0, \lambda_{23}, \lambda_{24}\} \\ \{\lambda_{25}, 0, 0\} & \lambda_{21} \end{pmatrix} \\ \begin{pmatrix} \lambda_5 & \{0, -\lambda_7, -\lambda_8\} \\ \{-\lambda_9, -\lambda_{10}, -\lambda_{11}\} & \lambda_4 \end{pmatrix} & \begin{pmatrix} \lambda_2 & \{0,0,0\} \\ \{0,0,0\} & \lambda_2 \end{pmatrix} & \begin{pmatrix} \lambda_{12} & \{\lambda_{14}, \lambda_{15}, \lambda_{16}\} \\ \{\lambda_{17}, \lambda_{18}, \lambda_{19}\} & \lambda_{13} \end{pmatrix} \\ \begin{pmatrix} \lambda_{21} & \{0, -\lambda_{23}, -\lambda_{24}\} \\ \{-\lambda_{25}, 0, 0\} & 0 \end{pmatrix} & \begin{pmatrix} \lambda_{13} & \{-\lambda_{14}, -\lambda_{15}, -\lambda_{16}\} \\ \{-\lambda_{17}, -\lambda_{18}, -\lambda_{19}\} & \lambda_{12} \end{pmatrix} & \begin{pmatrix} \lambda_3 & \{0,0,0\} \\ \{0,0,0\} & \lambda_3 \end{pmatrix} \end{pmatrix}$$

**U$_{X_1[e_6]}$[E$_2$[1]] == ZERO**

True

**defin[X$_1$[e$_6$], E$_2$[1]] // Expand**

$$\begin{pmatrix} \begin{pmatrix} 0 & \{0,0,0\} \\ \{0,0,0\} & 0 \end{pmatrix} & \begin{pmatrix} 0 & \{0,0,0\} \\ \{\alpha_3 + \lambda_2, 0, 0\} & 0 \end{pmatrix} & \begin{pmatrix} 0 & \{0, \lambda_{19}, -\lambda_{18}\} \\ \{\lambda_{12}, 0, 0\} & \lambda_{14} \end{pmatrix} \\ \begin{pmatrix} 0 & \{0,0,0\} \\ \{-\alpha_3 - \lambda_2, 0, 0\} & 0 \end{pmatrix} & \begin{pmatrix} 0 & \{0,0,0\} \\ \{0,0,0\} & 0 \end{pmatrix} & \begin{pmatrix} 0 & \{0,0,0\} \\ \{0,0,0\} & 0 \end{pmatrix} \\ \begin{pmatrix} \lambda_{14} & \{0, -\lambda_{19}, \lambda_{18}\} \\ \{-\lambda_{12}, 0, 0\} & 0 \end{pmatrix} & \begin{pmatrix} 0 & \{0,0,0\} \\ \{0,0,0\} & 0 \end{pmatrix} & \begin{pmatrix} 0 & \{0,0,0\} \\ \{0,0,0\} & 0 \end{pmatrix} \end{pmatrix}$$

**Δ[X$_1$[e$_6$]] = Δ[X$_1$[e$_6$]] //. {λ$_2$ → −α$_3$, λ$_{12}$ → 0, λ$_{14}$ → 0, λ$_{18}$ → 0, λ$_{19}$ → 0}**

$$\begin{pmatrix} \begin{pmatrix} \alpha_3 & \{0,0,0\} \\ \{0,0,0\} & \alpha_3 \end{pmatrix} & \begin{pmatrix} \lambda_4 & \{0, \lambda_7, \lambda_8\} \\ \{\lambda_9, \lambda_{10}, \lambda_{11}\} & \lambda_5 \end{pmatrix} & \begin{pmatrix} 0 & \{0, \lambda_{23}, \lambda_{24}\} \\ \{\lambda_{25}, 0, 0\} & \lambda_{21} \end{pmatrix} \\ \begin{pmatrix} \lambda_5 & \{0, -\lambda_7, -\lambda_8\} \\ \{-\lambda_9, -\lambda_{10}, -\lambda_{11}\} & \lambda_4 \end{pmatrix} & \begin{pmatrix} -\alpha_3 & \{0,0,0\} \\ \{0,0,0\} & -\alpha_3 \end{pmatrix} & \begin{pmatrix} 0 & \{0, \lambda_{15}, \lambda_{16}\} \\ \{\lambda_{17}, 0, 0\} & \lambda_{13} \end{pmatrix} \\ \begin{pmatrix} \lambda_{21} & \{0, -\lambda_{23}, -\lambda_{24}\} \\ \{-\lambda_{25}, 0, 0\} & 0 \end{pmatrix} & \begin{pmatrix} \lambda_{13} & \{0, -\lambda_{15}, -\lambda_{16}\} \\ \{-\lambda_{17}, 0, 0\} & 0 \end{pmatrix} & \begin{pmatrix} \lambda_3 & \{0,0,0\} \\ \{0,0,0\} & \lambda_3 \end{pmatrix} \end{pmatrix}$$

**U$_{X_1[e_6]}$[X$_1$[e$_6$]] == ZERO**

True

**(defin[X$_1$[e$_6$], X$_1$[e$_6$]] // Expand) == ZERO**

True

**U$_{X_1[e_6]}$[X$_1$[e$_5$]] == ZERO**

True

**defin[X$_1$[e$_6$], X$_1$[e$_5$]] // Expand**

$$\left( \begin{array}{ccc} \begin{pmatrix} 0 & \{0,0,0\} \\ \{0,0,0\} & 0 \end{pmatrix} & \begin{pmatrix} 0 & \{0,0,0\} \\ \{-\lambda_{11}-\phi_1,0,0\} & 0 \end{pmatrix} & \begin{pmatrix} 0 & \{0,0,0\} \\ \{0,0,0\} & 0 \end{pmatrix} \\ \begin{pmatrix} 0 & \{0,0,0\} \\ \{\lambda_{11}+\phi_1,0,0\} & 0 \end{pmatrix} & \begin{pmatrix} 0 & \{0,0,0\} \\ \{0,0,0\} & 0 \end{pmatrix} & \begin{pmatrix} 0 & \{0,0,0\} \\ \{0,0,0\} & 0 \end{pmatrix} \\ \begin{pmatrix} 0 & \{0,0,0\} \\ \{0,0,0\} & 0 \end{pmatrix} & \begin{pmatrix} 0 & \{0,0,0\} \\ \{0,0,0\} & 0 \end{pmatrix} & \begin{pmatrix} 0 & \{0,0,0\} \\ \{0,0,0\} & 0 \end{pmatrix} \end{array} \right)$$

**Δ[X$_1$[e$_6$]] = Δ[X$_1$[e$_6$]] //. {λ$_{11}$ → -φ$_1$}**

$$\left( \begin{array}{ccc} \begin{pmatrix} \alpha_3 & \{0,0,0\} \\ \{0,0,0\} & \alpha_3 \end{pmatrix} & \begin{pmatrix} \lambda_4 & \{0,\lambda_7,\lambda_8\} \\ \{\lambda_9,\lambda_{10},-\phi_1\} & \lambda_5 \end{pmatrix} & \begin{pmatrix} 0 & \{0,\lambda_{23},\lambda_{24}\} \\ \{\lambda_{25},0,0\} & \lambda_{21} \end{pmatrix} \\ \begin{pmatrix} \lambda_5 & \{0,-\lambda_7,-\lambda_8\} \\ \{-\lambda_9,-\lambda_{10},\phi_1\} & \lambda_4 \end{pmatrix} & \begin{pmatrix} -\alpha_3 & \{0,0,0\} \\ \{0,0,0\} & -\alpha_3 \end{pmatrix} & \begin{pmatrix} 0 & \{0,\lambda_{15},\lambda_{16}\} \\ \{\lambda_{17},0,0\} & \lambda_{13} \end{pmatrix} \\ \begin{pmatrix} \lambda_{21} & \{0,-\lambda_{23},-\lambda_{24}\} \\ \{-\lambda_{25},0,0\} & 0 \end{pmatrix} & \begin{pmatrix} \lambda_{13} & \{0,-\lambda_{15},-\lambda_{16}\} \\ \{-\lambda_{17},0,0\} & 0 \end{pmatrix} & \begin{pmatrix} \lambda_3 & \{0,0,0\} \\ \{0,0,0\} & \lambda_3 \end{pmatrix} \end{array} \right)$$

**U$_{X_1[e_6]}$[X$_1$[e$_4$]] == ZERO**

True

**defin[X$_1$[e$_6$], X$_1$[e$_4$]] // Expand**

$$\left( \begin{array}{ccc} \begin{pmatrix} 0 & \{0,0,0\} \\ \{0,0,0\} & 0 \end{pmatrix} & \begin{pmatrix} 0 & \{0,0,0\} \\ \{-\epsilon_1-\lambda_{10},0,0\} & 0 \end{pmatrix} & \begin{pmatrix} 0 & \{0,0,0\} \\ \{0,0,0\} & 0 \end{pmatrix} \\ \begin{pmatrix} 0 & \{0,0,0\} \\ \{\epsilon_1+\lambda_{10},0,0\} & 0 \end{pmatrix} & \begin{pmatrix} 0 & \{0,0,0\} \\ \{0,0,0\} & 0 \end{pmatrix} & \begin{pmatrix} 0 & \{0,0,0\} \\ \{0,0,0\} & 0 \end{pmatrix} \\ \begin{pmatrix} 0 & \{0,0,0\} \\ \{0,0,0\} & 0 \end{pmatrix} & \begin{pmatrix} 0 & \{0,0,0\} \\ \{0,0,0\} & 0 \end{pmatrix} & \begin{pmatrix} 0 & \{0,0,0\} \\ \{0,0,0\} & 0 \end{pmatrix} \end{array} \right)$$

**Δ[X$_1$[e$_6$]] = Δ[X$_1$[e$_6$]] //. {λ$_{10}$ → -ε$_1$}**

$$\left( \begin{array}{ccc} \begin{pmatrix} \alpha_3 & \{0,0,0\} \\ \{0,0,0\} & \alpha_3 \end{pmatrix} & \begin{pmatrix} \lambda_4 & \{0,\lambda_7,\lambda_8\} \\ \{\lambda_9,-\epsilon_1,-\phi_1\} & \lambda_5 \end{pmatrix} & \begin{pmatrix} 0 & \{0,\lambda_{23},\lambda_{24}\} \\ \{\lambda_{25},0,0\} & \lambda_{21} \end{pmatrix} \\ \begin{pmatrix} \lambda_5 & \{0,-\lambda_7,-\lambda_8\} \\ \{-\lambda_9,\epsilon_1,\phi_1\} & \lambda_4 \end{pmatrix} & \begin{pmatrix} -\alpha_3 & \{0,0,0\} \\ \{0,0,0\} & -\alpha_3 \end{pmatrix} & \begin{pmatrix} 0 & \{0,\lambda_{15},\lambda_{16}\} \\ \{\lambda_{17},0,0\} & \lambda_{13} \end{pmatrix} \\ \begin{pmatrix} \lambda_{21} & \{0,-\lambda_{23},-\lambda_{24}\} \\ \{-\lambda_{25},0,0\} & 0 \end{pmatrix} & \begin{pmatrix} \lambda_{13} & \{0,-\lambda_{15},-\lambda_{16}\} \\ \{-\lambda_{17},0,0\} & 0 \end{pmatrix} & \begin{pmatrix} \lambda_3 & \{0,0,0\} \\ \{0,0,0\} & \lambda_3 \end{pmatrix} \end{array} \right)$$

**U$_{X_1[e_6]}$[X$_1$[e$_3$]] == -X$_1$[e$_6$]**

True

**Δ[X₁[e₆]] + defin[X₁[e₆], X₁[e₃]] // Expand**

$$\left( \begin{pmatrix} 0 & \{0,0,0\} \\ \{0,0,0\} & 0 \end{pmatrix} \quad \begin{pmatrix} 0 & \{0,0,0\} \\ \{-\eta_1-\lambda_9,0,0\} & 0 \end{pmatrix} \quad \begin{pmatrix} 0 & \{0,0,0\} \\ \{0,0,0\} & 0 \end{pmatrix} \right.$$
$$\begin{pmatrix} 0 & \{0,0,0\} \\ \{\eta_1+\lambda_9,0,0\} & 0 \end{pmatrix} \quad \begin{pmatrix} 0 & \{0,0,0\} \\ \{0,0,0\} & 0 \end{pmatrix} \quad \begin{pmatrix} 0 & \{0,0,0\} \\ \{0,0,0\} & 0 \end{pmatrix}$$
$$\left. \begin{pmatrix} 0 & \{0,0,0\} \\ \{0,0,0\} & 0 \end{pmatrix} \quad \begin{pmatrix} 0 & \{0,0,0\} \\ \{0,0,0\} & 0 \end{pmatrix} \quad \begin{pmatrix} \lambda_3 & \{0,0,0\} \\ \{0,0,0\} & \lambda_3 \end{pmatrix} \right)$$

**Δ[X₁[e₆]] = Δ[X₁[e₆]] //. {λ₉ → -η₁}**

$$\left( \begin{pmatrix} \alpha_3 & \{0,0,0\} \\ \{0,0,0\} & \alpha_3 \end{pmatrix} \quad \begin{pmatrix} \lambda_4 & \{0,\lambda_7,\lambda_8\} \\ \{-\eta_1,-\epsilon_1,-\phi_1\} & \lambda_5 \end{pmatrix} \quad \begin{pmatrix} 0 & \{0,\lambda_{23},\lambda_{24}\} \\ \{\lambda_{25},0,0\} & \lambda_{21} \end{pmatrix} \right.$$
$$\begin{pmatrix} \lambda_5 & \{0,-\lambda_7,-\lambda_8\} \\ \{\eta_1,\epsilon_1,\phi_1\} & \lambda_4 \end{pmatrix} \quad \begin{pmatrix} -\alpha_3 & \{0,0,0\} \\ \{0,0,0\} & -\alpha_3 \end{pmatrix} \quad \begin{pmatrix} 0 & \{0,\lambda_{15},\lambda_{16}\} \\ \{\lambda_{17},0,0\} & \lambda_{13} \end{pmatrix}$$
$$\left. \begin{pmatrix} \lambda_{21} & \{0,-\lambda_{23},-\lambda_{24}\} \\ \{-\lambda_{25},0,0\} & 0 \end{pmatrix} \quad \begin{pmatrix} \lambda_{13} & \{0,-\lambda_{15},-\lambda_{16}\} \\ \{-\lambda_{17},0,0\} & 0 \end{pmatrix} \quad \begin{pmatrix} \lambda_3 & \{0,0,0\} \\ \{0,0,0\} & \lambda_3 \end{pmatrix} \right)$$

**U_{X₁[e₆]} [X₁[e₂]] == ZERO**

True

**defin[X₁[e₆], X₁[e₂]] // Expand**

$$\left( \begin{pmatrix} 0 & \{0,0,0\} \\ \{0,0,0\} & 0 \end{pmatrix} \quad \begin{pmatrix} 0 & \{0,0,0\} \\ \{\lambda_4-\rho_1,0,0\} & 0 \end{pmatrix} \quad \begin{pmatrix} 0 & \{0,0,0\} \\ \{0,0,0\} & 0 \end{pmatrix} \right.$$
$$\begin{pmatrix} 0 & \{0,0,0\} \\ \{\rho_1-\lambda_4,0,0\} & 0 \end{pmatrix} \quad \begin{pmatrix} 0 & \{0,0,0\} \\ \{0,0,0\} & 0 \end{pmatrix} \quad \begin{pmatrix} 0 & \{0,0,0\} \\ \{0,0,0\} & 0 \end{pmatrix}$$
$$\left. \begin{pmatrix} 0 & \{0,0,0\} \\ \{0,0,0\} & 0 \end{pmatrix} \quad \begin{pmatrix} 0 & \{0,0,0\} \\ \{0,0,0\} & 0 \end{pmatrix} \quad \begin{pmatrix} 0 & \{0,0,0\} \\ \{0,0,0\} & 0 \end{pmatrix} \right)$$

**Δ[X₁[e₆]] = Δ[X₁[e₆]] //. {λ₄ → ρ₁}**

$$\left( \begin{pmatrix} \alpha_3 & \{0,0,0\} \\ \{0,0,0\} & \alpha_3 \end{pmatrix} \quad \begin{pmatrix} \rho_1 & \{0,\lambda_7,\lambda_8\} \\ \{-\eta_1,-\epsilon_1,-\phi_1\} & \lambda_5 \end{pmatrix} \quad \begin{pmatrix} 0 & \{0,\lambda_{23},\lambda_{24}\} \\ \{\lambda_{25},0,0\} & \lambda_{21} \end{pmatrix} \right.$$
$$\begin{pmatrix} \lambda_5 & \{0,-\lambda_7,-\lambda_8\} \\ \{\eta_1,\epsilon_1,\phi_1\} & \rho_1 \end{pmatrix} \quad \begin{pmatrix} -\alpha_3 & \{0,0,0\} \\ \{0,0,0\} & -\alpha_3 \end{pmatrix} \quad \begin{pmatrix} 0 & \{0,\lambda_{15},\lambda_{16}\} \\ \{\lambda_{17},0,0\} & \lambda_{13} \end{pmatrix}$$
$$\left. \begin{pmatrix} \lambda_{21} & \{0,-\lambda_{23},-\lambda_{24}\} \\ \{-\lambda_{25},0,0\} & 0 \end{pmatrix} \quad \begin{pmatrix} \lambda_{13} & \{0,-\lambda_{15},-\lambda_{16}\} \\ \{-\lambda_{17},0,0\} & 0 \end{pmatrix} \quad \begin{pmatrix} \lambda_3 & \{0,0,0\} \\ \{0,0,0\} & \lambda_3 \end{pmatrix} \right)$$

**U_{X₁[e₆]} [X₁[e₁]] == ZERO**

True

**defin[X₁[e₆], X₁[e₁]] // Expand**

$$\begin{pmatrix} \begin{pmatrix} 0 & \{0,0,0\} \\ \{0,0,0\} & 0 \end{pmatrix} & \begin{pmatrix} 0 & \{0,0,0\} \\ \{\lambda_5 - \delta_2, 0, 0\} & 0 \end{pmatrix} & \begin{pmatrix} 0 & \{0,0,0\} \\ \{0,0,0\} & 0 \end{pmatrix} \\ \begin{pmatrix} 0 & \{0,0,0\} \\ \{\delta_2 - \lambda_5, 0, 0\} & 0 \end{pmatrix} & \begin{pmatrix} 0 & \{0,0,0\} \\ \{0,0,0\} & 0 \end{pmatrix} & \begin{pmatrix} 0 & \{0,0,0\} \\ \{0,0,0\} & 0 \end{pmatrix} \\ \begin{pmatrix} 0 & \{0,0,0\} \\ \{0,0,0\} & 0 \end{pmatrix} & \begin{pmatrix} 0 & \{0,0,0\} \\ \{0,0,0\} & 0 \end{pmatrix} & \begin{pmatrix} 0 & \{0,0,0\} \\ \{0,0,0\} & 0 \end{pmatrix} \end{pmatrix}$$

**Δ[X₁[e₆]] = Δ[X₁[e₆]] //. {λ₅ → δ₂}**

$$\begin{pmatrix} \begin{pmatrix} \alpha_3 & \{0,0,0\} \\ \{0,0,0\} & \alpha_3 \end{pmatrix} & \begin{pmatrix} \rho_1 & \{0,\lambda_7,\lambda_8\} \\ \{-\eta_1,-\epsilon_1,-\phi_1\} & \delta_2 \end{pmatrix} & \begin{pmatrix} 0 & \{0,\lambda_{23},\lambda_{24}\} \\ \{\lambda_{25},0,0\} & \lambda_{21} \end{pmatrix} \\ \begin{pmatrix} \delta_2 & \{0,-\lambda_7,-\lambda_8\} \\ \{\eta_1,\epsilon_1,\phi_1\} & \rho_1 \end{pmatrix} & \begin{pmatrix} -\alpha_3 & \{0,0,0\} \\ \{0,0,0\} & -\alpha_3 \end{pmatrix} & \begin{pmatrix} 0 & \{0,\lambda_{15},\lambda_{16}\} \\ \{\lambda_{17},0,0\} & \lambda_{13} \end{pmatrix} \\ \begin{pmatrix} \lambda_{21} & \{0,-\lambda_{23},-\lambda_{24}\} \\ \{-\lambda_{25},0,0\} & 0 \end{pmatrix} & \begin{pmatrix} \lambda_{13} & \{0,-\lambda_{15},-\lambda_{16}\} \\ \{-\lambda_{17},0,0\} & 0 \end{pmatrix} & \begin{pmatrix} \lambda_3 & \{0,0,0\} \\ \{0,0,0\} & \lambda_3 \end{pmatrix} \end{pmatrix}$$

**T[E₁[1], E₁[1], X₁[e₆]] == X₁[e₆]**

True

**Δ[X₁[e₆]] - Leib[E₁[1], E₁[1], X₁[e₆]] // Expand**

$$\begin{pmatrix} \begin{pmatrix} 0 & \{0,0,0\} \\ \{0,0,0\} & 0 \end{pmatrix} & \begin{pmatrix} 0 & \{0,0,0\} \\ \{0,0,0\} & 0 \end{pmatrix} & \begin{pmatrix} 0 & \{0,0,0\} \\ \{0,0,0\} & 0 \end{pmatrix} \\ \begin{pmatrix} 0 & \{0,0,0\} \\ \{0,0,0\} & 0 \end{pmatrix} & \begin{pmatrix} 0 & \{0,0,0\} \\ \{0,0,0\} & 0 \end{pmatrix} & \begin{pmatrix} 0 & \{0,\beta_8+\lambda_{15},\lambda_{16}-\beta_7\} \\ \{\beta_1+\lambda_{17},0,0\} & \beta_3+\lambda_{13} \end{pmatrix} \\ \begin{pmatrix} 0 & \{0,0,0\} \\ \{0,0,0\} & 0 \end{pmatrix} & \begin{pmatrix} \beta_3+\lambda_{13} & \{0,-\beta_8-\lambda_{15},\beta_7-\lambda_{16}\} \\ \{-\beta_1-\lambda_{17},0,0\} & 0 \end{pmatrix} & \begin{pmatrix} \lambda_3 & \{0,0,0\} \\ \{0,0,0\} & \lambda_3 \end{pmatrix} \end{pmatrix}$$

**Δ[X₁[e₆]] =
  Δ[X₁[e₆]] //. {λ₃ → 0, λ₁₃ → -β₃, λ₁₅ → -β₈, λ₁₆ → β₇, λ₁₇ → -β₁}**

$$\begin{pmatrix} \begin{pmatrix} \alpha_3 & \{0,0,0\} \\ \{0,0,0\} & \alpha_3 \end{pmatrix} & \begin{pmatrix} \rho_1 & \{0,\lambda_7,\lambda_8\} \\ \{-\eta_1,-\epsilon_1,-\phi_1\} & \delta_2 \end{pmatrix} & \begin{pmatrix} 0 & \{0,\lambda_{23},\lambda_{24}\} \\ \{\lambda_{25},0,0\} & \lambda_{21} \end{pmatrix} \\ \begin{pmatrix} \delta_2 & \{0,-\lambda_7,-\lambda_8\} \\ \{\eta_1,\epsilon_1,\phi_1\} & \rho_1 \end{pmatrix} & \begin{pmatrix} -\alpha_3 & \{0,0,0\} \\ \{0,0,0\} & -\alpha_3 \end{pmatrix} & \begin{pmatrix} 0 & \{0,-\beta_8,\beta_7\} \\ \{-\beta_1,0,0\} & -\beta_3 \end{pmatrix} \\ \begin{pmatrix} \lambda_{21} & \{0,-\lambda_{23},-\lambda_{24}\} \\ \{-\lambda_{25},0,0\} & 0 \end{pmatrix} & \begin{pmatrix} -\beta_3 & \{0,\beta_8,-\beta_7\} \\ \{\beta_1,0,0\} & 0 \end{pmatrix} & \begin{pmatrix} 0 & \{0,0,0\} \\ \{0,0,0\} & 0 \end{pmatrix} \end{pmatrix}$$

**T[E₁[1], X₁[e₆], E₁[1]] == ZERO**

True

**Leib[E₁[1], X₁[e₆], E₁[1]] == ZERO**

True

**T[E₂[1], E₂[1], X₁[e₆]] == X₁[e₆]**

True

**Δ[X₁[e₆]] - Leib[E₂[1], E₂[1], X₁[e₆]] // Expand**

$$\begin{pmatrix} \begin{pmatrix} 0 & \{0,0,0\} \\ \{0,0,0\} & 0 \end{pmatrix} & \begin{pmatrix} 0 & \{0,0,0\} \\ \{0,0,0\} & 0 \end{pmatrix} & \begin{pmatrix} 0 & \{0, \lambda_{23} - \gamma_8, \gamma_7 + \lambda_{24}\} \\ \{\lambda_{25} - \gamma_1, 0, 0\} & \lambda_{21} - \gamma_3 \end{pmatrix} \\ \begin{pmatrix} 0 & \{0,0,0\} \\ \{0,0,0\} & 0 \end{pmatrix} & \begin{pmatrix} 0 & \{0,0,0\} \\ \{0,0,0\} & 0 \end{pmatrix} & \begin{pmatrix} 0 & \{0,0,0\} \\ \{0,0,0\} & 0 \end{pmatrix} \\ \begin{pmatrix} \lambda_{21} - \gamma_3 & \{0, \gamma_8 - \lambda_{23}, -\gamma_7 - \lambda_{24}\} \\ \{\gamma_1 - \lambda_{25}, 0, 0\} & 0 \end{pmatrix} & \begin{pmatrix} 0 & \{0,0,0\} \\ \{0,0,0\} & 0 \end{pmatrix} & \begin{pmatrix} 0 & \{0,0,0\} \\ \{0,0,0\} & 0 \end{pmatrix} \end{pmatrix}$$

**Δ[X₁[e₆]] = Δ[X₁[e₆]] //. {λ₂₁ → γ₃, λ₂₃ → γ₈, λ₂₄ → -γ₇, λ₂₅ → γ₁}**

$$\begin{pmatrix} \begin{pmatrix} \alpha_3 & \{0,0,0\} \\ \{0,0,0\} & \alpha_3 \end{pmatrix} & \begin{pmatrix} \rho_1 & \{0, \lambda_7, \lambda_8\} \\ \{-\eta_1, -\epsilon_1, -\phi_1\} & \delta_2 \end{pmatrix} & \begin{pmatrix} 0 & \{0, \gamma_8, -\gamma_7\} \\ \{\gamma_1, 0, 0\} & \gamma_3 \end{pmatrix} \\ \begin{pmatrix} \delta_2 & \{0, -\lambda_7, -\lambda_8\} \\ \{\eta_1, \epsilon_1, \phi_1\} & \rho_1 \end{pmatrix} & \begin{pmatrix} -\alpha_3 & \{0,0,0\} \\ \{0,0,0\} & -\alpha_3 \end{pmatrix} & \begin{pmatrix} 0 & \{0, -\beta_8, \beta_7\} \\ \{-\beta_1, 0, 0\} & -\beta_3 \end{pmatrix} \\ \begin{pmatrix} \gamma_3 & \{0, -\gamma_8, \gamma_7\} \\ \{-\gamma_1, 0, 0\} & 0 \end{pmatrix} & \begin{pmatrix} -\beta_3 & \{0, \beta_8, -\beta_7\} \\ \{\beta_1, 0, 0\} & 0 \end{pmatrix} & \begin{pmatrix} 0 & \{0,0,0\} \\ \{0,0,0\} & 0 \end{pmatrix} \end{pmatrix}$$

**T[E₂[1], X₁[e₆], E₂[1]] == ZERO**

True

**Leib[E₂[1], X₁[e₆], E₂[1]] == ZERO**

True

**T[E₁[1], E₂[1], X₁[e₆]] == ZERO**

True

**(Leib[E₁[1], E₂[1], X₁[e₆]] // Expand) == ZERO**

True

**T[E₂[1], E₁[1], X₁[e₆]] == ZERO**

True

**(Leib[E₂[1], E₁[1], X₁[e₆]] // Expand) == ZERO**

True

**T[E₁[1], X₁[e₆], E₂[1]] == X₁[e₆]**

True

```
Δ[X₁[e₆]] - Leib[E₁[1], X₁[e₆], E₂[1]] == ZERO
```

True

```
T[E₁[1], X₁[e₁], X₁[e₆]] == ZERO
```

True

```
(Leib[E₁[1], X₁[e₁], X₁[e₆]] // Expand) == ZERO
```

True

```
T[X₁[e₁], E₁[1], X₁[e₆]] == ZERO
```

True

```
(Leib[X₁[e₁], E₁[1], X₁[e₆]] // Expand) == ZERO
```

True

```
T[E₂[1], X₁[e₁], X₁[e₆]] == ZERO
```

True

```
(Leib[E₂[1], X₁[e₁], X₁[e₆]] // Expand) == ZERO
```

True

```
T[X₁[e₆], X₁[e₁], X₁[e₁]] == ZERO
```

True

```
(Leib[X₁[e₆], X₁[e₁], X₁[e₁]] // Expand) == ZERO
```

True

```
T[X₁[e₅], X₁[e₄], X₁[e₆]] == ZERO
```

True

```
(Leib[X₁[e₅], X₁[e₄], X₁[e₆]] // Expand) == ZERO
```

True

```
T[X₁[e₄], X₁[e₆], X₁[e₂]] == ZERO
```

True

**(Leib[X₁[e₄], X₁[e₆], X₁[e₂]] // Expand) == ZERO**

True

**T[X₁[e₃], X₁[e₅], X₁[e₆]] == X₁[e₅]**

True

**(Δ[X₁[e₅]] - Leib[X₁[e₃], X₁[e₅], X₁[e₆]] // Expand) == ZERO**

True

**T[X₁[e₅], X₁[e₆], X₁[e₅]] == ZERO**

True

**(Leib[X₁[e₅], X₁[e₆], X₁[e₅]] // Expand) == ZERO**

True

**T[X₁[e₆], X₁[e₆], X₁[e₆]] == ZERO**

True

**(Leib[X₁[e₆], X₁[e₆], X₁[e₆]] // Expand) == ZERO**

True

**Δ[X₁[e₆]] = Δ[X₁[e₆]] //. {λ₇ → ψ₁, λ₈ → ψ₂}**

$$\begin{pmatrix}
\begin{pmatrix} \alpha_3 & \{0,0,0\} \\ \{0,0,0\} & \alpha_3 \end{pmatrix} & \begin{pmatrix} \rho_1 & \{0,\psi_1,\psi_2\} \\ \{-\eta_1,-\epsilon_1,-\phi_1\} & \delta_2 \end{pmatrix} & \begin{pmatrix} 0 & \{0,\gamma_8,-\gamma_7\} \\ \{\gamma_1,0,0\} & \gamma_3 \end{pmatrix} \\
\begin{pmatrix} \delta_2 & \{0,-\psi_1,-\psi_2\} \\ \{\eta_1,\epsilon_1,\phi_1\} & \rho_1 \end{pmatrix} & \begin{pmatrix} -\alpha_3 & \{0,0,0\} \\ \{0,0,0\} & -\alpha_3 \end{pmatrix} & \begin{pmatrix} 0 & \{0,-\beta_8,\beta_7\} \\ \{-\beta_1,0,0\} & -\beta_3 \end{pmatrix} \\
\begin{pmatrix} \gamma_3 & \{0,-\gamma_8,\gamma_7\} \\ \{-\gamma_1,0,0\} & 0 \end{pmatrix} & \begin{pmatrix} -\beta_3 & \{0,\beta_8,-\beta_7\} \\ \{\beta_1,0,0\} & 0 \end{pmatrix} & \begin{pmatrix} 0 & \{0,0,0\} \\ \{0,0,0\} & 0 \end{pmatrix}
\end{pmatrix}$$

- **X₁[e₇]**

    **Δ[X₁[e₇]] = generic;**

    **U_{X₁[e₇]}[E₁[1]] == ZERO**

    True

**defin[X₁[e₇], E₁[1]] // Expand**

$$\begin{pmatrix} \begin{pmatrix} 0 & \{0,0,0\} \\ \{0,0,0\} & 0 \end{pmatrix} & \begin{pmatrix} 0 & \{0,0,0\} \\ \{0, \lambda_1 - \alpha_4, 0\} & 0 \end{pmatrix} & \begin{pmatrix} 0 & \{0,0,0\} \\ \{0,0,0\} & 0 \end{pmatrix} \\ \begin{pmatrix} 0 & \{0,0,0\} \\ \{0, \alpha_4 - \lambda_1, 0\} & 0 \end{pmatrix} & \begin{pmatrix} -\lambda_7 & \{0,0,0\} \\ \{0,0,0\} & -\lambda_7 \end{pmatrix} & \begin{pmatrix} 0 & \{\lambda_{27}, 0, -\lambda_{25}\} \\ \{0, -\lambda_{20}, 0\} & -\lambda_{23} \end{pmatrix} \\ \begin{pmatrix} 0 & \{0,0,0\} \\ \{0,0,0\} & 0 \end{pmatrix} & \begin{pmatrix} -\lambda_{23} & \{-\lambda_{27}, 0, \lambda_{25}\} \\ \{0, \lambda_{20}, 0\} & 0 \end{pmatrix} & \begin{pmatrix} 0 & \{0,0,0\} \\ \{0,0,0\} & 0 \end{pmatrix} \end{pmatrix}$$

**Δ[X₁[e₇]] =**
 **Δ[X₁[e₇]] //. {λ₁ → α₄, λ₇ → 0, λ₂₀ → 0, λ₂₃ → 0, λ₂₅ → 0, λ₂₇ → 0}**

$$\begin{pmatrix} \begin{pmatrix} \alpha_4 & \{0,0,0\} \\ \{0,0,0\} & \alpha_4 \end{pmatrix} & \begin{pmatrix} \lambda_4 & \{\lambda_6, 0, \lambda_8\} \\ \{\lambda_9, \lambda_{10}, \lambda_{11}\} & \lambda_5 \end{pmatrix} & \begin{pmatrix} 0 & \{\lambda_{22}, 0, \lambda_{24}\} \\ \{0, \lambda_{26}, 0\} & \lambda_{21} \end{pmatrix} \\ \begin{pmatrix} \lambda_5 & \{-\lambda_6, 0, -\lambda_8\} \\ \{-\lambda_9, -\lambda_{10}, -\lambda_{11}\} & \lambda_4 \end{pmatrix} & \begin{pmatrix} \lambda_2 & \{0,0,0\} \\ \{0,0,0\} & \lambda_2 \end{pmatrix} & \begin{pmatrix} \lambda_{12} & \{\lambda_{14}, \lambda_{15}, \lambda_{16}\} \\ \{\lambda_{17}, \lambda_{18}, \lambda_{19}\} & \lambda_{13} \end{pmatrix} \\ \begin{pmatrix} \lambda_{21} & \{-\lambda_{22}, 0, -\lambda_{24}\} \\ \{0, -\lambda_{26}, 0\} & 0 \end{pmatrix} & \begin{pmatrix} \lambda_{13} & \{-\lambda_{14}, -\lambda_{15}, -\lambda_{16}\} \\ \{-\lambda_{17}, -\lambda_{18}, -\lambda_{19}\} & \lambda_{12} \end{pmatrix} & \begin{pmatrix} \lambda_3 & \{0,0,0\} \\ \{0,0,0\} & \lambda_3 \end{pmatrix} \end{pmatrix}$$

**U_{X₁[e₇]}[E₂[1]] == ZERO**

True

**defin[X₁[e₇], E₂[1]] // Expand**

$$\begin{pmatrix} \begin{pmatrix} 0 & \{0,0,0\} \\ \{0,0,0\} & 0 \end{pmatrix} & \begin{pmatrix} 0 & \{0,0,0\} \\ \{0, \alpha_4 + \lambda_2, 0\} & 0 \end{pmatrix} & \begin{pmatrix} 0 & \{-\lambda_{19}, 0, \lambda_{17}\} \\ \{0, \lambda_{12}, 0\} & \lambda_{15} \end{pmatrix} \\ \begin{pmatrix} 0 & \{0,0,0\} \\ \{0, -\alpha_4 - \lambda_2, 0\} & 0 \end{pmatrix} & \begin{pmatrix} 0 & \{0,0,0\} \\ \{0,0,0\} & 0 \end{pmatrix} & \begin{pmatrix} 0 & \{0,0,0\} \\ \{0,0,0\} & 0 \end{pmatrix} \\ \begin{pmatrix} \lambda_{15} & \{\lambda_{19}, 0, -\lambda_{17}\} \\ \{0, -\lambda_{12}, 0\} & 0 \end{pmatrix} & \begin{pmatrix} 0 & \{0,0,0\} \\ \{0,0,0\} & 0 \end{pmatrix} & \begin{pmatrix} 0 & \{0,0,0\} \\ \{0,0,0\} & 0 \end{pmatrix} \end{pmatrix}$$

**Δ[X₁[e₇]] = Δ[X₁[e₇]] //. {λ₂ → -α₄, λ₁₂ → 0, λ₁₅ → 0, λ₁₇ → 0, λ₁₉ → 0}**

$$\begin{pmatrix} \begin{pmatrix} \alpha_4 & \{0,0,0\} \\ \{0,0,0\} & \alpha_4 \end{pmatrix} & \begin{pmatrix} \lambda_4 & \{\lambda_6, 0, \lambda_8\} \\ \{\lambda_9, \lambda_{10}, \lambda_{11}\} & \lambda_5 \end{pmatrix} & \begin{pmatrix} 0 & \{\lambda_{22}, 0, \lambda_{24}\} \\ \{0, \lambda_{26}, 0\} & \lambda_{21} \end{pmatrix} \\ \begin{pmatrix} \lambda_5 & \{-\lambda_6, 0, -\lambda_8\} \\ \{-\lambda_9, -\lambda_{10}, -\lambda_{11}\} & \lambda_4 \end{pmatrix} & \begin{pmatrix} -\alpha_4 & \{0,0,0\} \\ \{0,0,0\} & -\alpha_4 \end{pmatrix} & \begin{pmatrix} 0 & \{\lambda_{14}, 0, \lambda_{16}\} \\ \{0, \lambda_{18}, 0\} & \lambda_{13} \end{pmatrix} \\ \begin{pmatrix} \lambda_{21} & \{-\lambda_{22}, 0, -\lambda_{24}\} \\ \{0, -\lambda_{26}, 0\} & 0 \end{pmatrix} & \begin{pmatrix} \lambda_{13} & \{-\lambda_{14}, 0, -\lambda_{16}\} \\ \{0, -\lambda_{18}, 0\} & 0 \end{pmatrix} & \begin{pmatrix} \lambda_3 & \{0,0,0\} \\ \{0,0,0\} & \lambda_3 \end{pmatrix} \end{pmatrix}$$

**U_{X₁[e₇]}[X₁[e₇]] == ZERO**

True

**(defin[X₁[e₇], X₁[e₇]] // Expand) == ZERO**

True

**U$_{X_1[e_7]}$ [X$_1$ [e$_6$]] == ZERO**

True

**defin[X$_1$[e$_7$], X$_1$[e$_6$]] // Expand**

$$\left(\begin{array}{ccc} \begin{pmatrix} 0 & \{0,0,0\} \\ \{0,0,0\} & 0 \end{pmatrix} & \begin{pmatrix} 0 & \{0,0,0\} \\ \{0,-\lambda_6-\psi_1,0\} & 0 \end{pmatrix} & \begin{pmatrix} 0 & \{0,0,0\} \\ \{0,0,0\} & 0 \end{pmatrix} \\ \begin{pmatrix} 0 & \{0,0,0\} \\ \{0,\lambda_6+\psi_1,0\} & 0 \end{pmatrix} & \begin{pmatrix} 0 & \{0,0,0\} \\ \{0,0,0\} & 0 \end{pmatrix} & \begin{pmatrix} 0 & \{0,0,0\} \\ \{0,0,0\} & 0 \end{pmatrix} \\ \begin{pmatrix} 0 & \{0,0,0\} \\ \{0,0,0\} & 0 \end{pmatrix} & \begin{pmatrix} 0 & \{0,0,0\} \\ \{0,0,0\} & 0 \end{pmatrix} & \begin{pmatrix} 0 & \{0,0,0\} \\ \{0,0,0\} & 0 \end{pmatrix} \end{array}\right)$$

**Δ[X$_1$[e$_7$]] = Δ[X$_1$[e$_7$]] //. {λ$_6$ → −ψ$_1$}**

$$\left(\begin{array}{ccc} \begin{pmatrix} \alpha_4 & \{0,0,0\} \\ \{0,0,0\} & \alpha_4 \end{pmatrix} & \begin{pmatrix} \lambda_4 & \{-\psi_1,0,\lambda_8\} \\ \{\lambda_9,\lambda_{10},\lambda_{11}\} & \lambda_5 \end{pmatrix} & \begin{pmatrix} 0 & \{\lambda_{22},0,\lambda_{24}\} \\ \{0,\lambda_{26},0\} & \lambda_{21} \end{pmatrix} \\ \begin{pmatrix} \lambda_5 & \{\psi_1,0,-\lambda_8\} \\ \{-\lambda_9,-\lambda_{10},-\lambda_{11}\} & \lambda_4 \end{pmatrix} & \begin{pmatrix} -\alpha_4 & \{0,0,0\} \\ \{0,0,0\} & -\alpha_4 \end{pmatrix} & \begin{pmatrix} 0 & \{\lambda_{14},0,\lambda_{16}\} \\ \{0,\lambda_{18},0\} & \lambda_{13} \end{pmatrix} \\ \begin{pmatrix} \lambda_{21} & \{-\lambda_{22},0,-\lambda_{24}\} \\ \{0,-\lambda_{26},0\} & 0 \end{pmatrix} & \begin{pmatrix} \lambda_{13} & \{-\lambda_{14},0,-\lambda_{16}\} \\ \{0,-\lambda_{18},0\} & 0 \end{pmatrix} & \begin{pmatrix} \lambda_3 & \{0,0,0\} \\ \{0,0,0\} & \lambda_3 \end{pmatrix} \end{array}\right)$$

**U$_{X_1[e_7]}$ [X$_1$ [e$_5$]] == ZERO**

True

**defin[X$_1$[e$_7$], X$_1$[e$_5$]] // Expand**

$$\left(\begin{array}{ccc} \begin{pmatrix} 0 & \{0,0,0\} \\ \{0,0,0\} & 0 \end{pmatrix} & \begin{pmatrix} 0 & \{0,0,0\} \\ \{0,-\lambda_{11}-\phi_2,0\} & 0 \end{pmatrix} & \begin{pmatrix} 0 & \{0,0,0\} \\ \{0,0,0\} & 0 \end{pmatrix} \\ \begin{pmatrix} 0 & \{0,0,0\} \\ \{0,\lambda_{11}+\phi_2,0\} & 0 \end{pmatrix} & \begin{pmatrix} 0 & \{0,0,0\} \\ \{0,0,0\} & 0 \end{pmatrix} & \begin{pmatrix} 0 & \{0,0,0\} \\ \{0,0,0\} & 0 \end{pmatrix} \\ \begin{pmatrix} 0 & \{0,0,0\} \\ \{0,0,0\} & 0 \end{pmatrix} & \begin{pmatrix} 0 & \{0,0,0\} \\ \{0,0,0\} & 0 \end{pmatrix} & \begin{pmatrix} 0 & \{0,0,0\} \\ \{0,0,0\} & 0 \end{pmatrix} \end{array}\right)$$

**Δ[X$_1$[e$_7$]] = Δ[X$_1$[e$_7$]] //. {λ$_{11}$ → −φ$_2$}**

$$\left(\begin{array}{ccc} \begin{pmatrix} \alpha_4 & \{0,0,0\} \\ \{0,0,0\} & \alpha_4 \end{pmatrix} & \begin{pmatrix} \lambda_4 & \{-\psi_1,0,\lambda_8\} \\ \{\lambda_9,\lambda_{10},-\phi_2\} & \lambda_5 \end{pmatrix} & \begin{pmatrix} 0 & \{\lambda_{22},0,\lambda_{24}\} \\ \{0,\lambda_{26},0\} & \lambda_{21} \end{pmatrix} \\ \begin{pmatrix} \lambda_5 & \{\psi_1,0,-\lambda_8\} \\ \{-\lambda_9,-\lambda_{10},\phi_2\} & \lambda_4 \end{pmatrix} & \begin{pmatrix} -\alpha_4 & \{0,0,0\} \\ \{0,0,0\} & -\alpha_4 \end{pmatrix} & \begin{pmatrix} 0 & \{\lambda_{14},0,\lambda_{16}\} \\ \{0,\lambda_{18},0\} & \lambda_{13} \end{pmatrix} \\ \begin{pmatrix} \lambda_{21} & \{-\lambda_{22},0,-\lambda_{24}\} \\ \{0,-\lambda_{26},0\} & 0 \end{pmatrix} & \begin{pmatrix} \lambda_{13} & \{-\lambda_{14},0,-\lambda_{16}\} \\ \{0,-\lambda_{18},0\} & 0 \end{pmatrix} & \begin{pmatrix} \lambda_3 & \{0,0,0\} \\ \{0,0,0\} & \lambda_3 \end{pmatrix} \end{array}\right)$$

**U$_{X_1[e_7]}$ [X$_1$ [e$_4$]] == −X$_1$ [e$_7$]**

True

**Δ[X₁[e₇]] + defin[X₁[e₇], X₁[e₄]] // Expand**

$$\begin{pmatrix} \begin{pmatrix} 0 & \{0,0,0\} \\ \{0,0,0\} & 0 \end{pmatrix} & \begin{pmatrix} 0 & \{0,0,0\} \\ \{0,-\epsilon_2-\lambda_{10},0\} & 0 \end{pmatrix} & \begin{pmatrix} 0 & \{0,0,0\} \\ \{0,0,0\} & 0 \end{pmatrix} \\ \begin{pmatrix} 0 & \{0,0,0\} \\ \{0,\epsilon_2+\lambda_{10},0\} & 0 \end{pmatrix} & \begin{pmatrix} 0 & \{0,0,0\} \\ \{0,0,0\} & 0 \end{pmatrix} & \begin{pmatrix} 0 & \{0,0,0\} \\ \{0,0,0\} & 0 \end{pmatrix} \\ \begin{pmatrix} 0 & \{0,0,0\} \\ \{0,0,0\} & 0 \end{pmatrix} & \begin{pmatrix} 0 & \{0,0,0\} \\ \{0,0,0\} & 0 \end{pmatrix} & \begin{pmatrix} \lambda_3 & \{0,0,0\} \\ \{0,0,0\} & \lambda_3 \end{pmatrix} \end{pmatrix}$$

**Δ[X₁[e₇]] = Δ[X₁[e₇]] //. {λ₁₀ → -ϵ₂}**

$$\begin{pmatrix} \begin{pmatrix} \alpha_4 & \{0,0,0\} \\ \{0,0,0\} & \alpha_4 \end{pmatrix} & \begin{pmatrix} \lambda_4 & \{-\psi_1,0,\lambda_8\} \\ \{\lambda_9,-\epsilon_2,-\phi_2\} & \lambda_5 \end{pmatrix} & \begin{pmatrix} 0 & \{\lambda_{22},0,\lambda_{24}\} \\ \{0,\lambda_{26},0\} & \lambda_{21} \end{pmatrix} \\ \begin{pmatrix} \lambda_5 & \{\psi_1,0,-\lambda_8\} \\ \{-\lambda_9,\epsilon_2,\phi_2\} & \lambda_4 \end{pmatrix} & \begin{pmatrix} -\alpha_4 & \{0,0,0\} \\ \{0,0,0\} & -\alpha_4 \end{pmatrix} & \begin{pmatrix} 0 & \{\lambda_{14},0,\lambda_{16}\} \\ \{0,\lambda_{18},0\} & \lambda_{13} \end{pmatrix} \\ \begin{pmatrix} \lambda_{21} & \{-\lambda_{22},0,-\lambda_{24}\} \\ \{0,-\lambda_{26},0\} & 0 \end{pmatrix} & \begin{pmatrix} \lambda_{13} & \{-\lambda_{14},0,-\lambda_{16}\} \\ \{0,-\lambda_{18},0\} & 0 \end{pmatrix} & \begin{pmatrix} \lambda_3 & \{0,0,0\} \\ \{0,0,0\} & \lambda_3 \end{pmatrix} \end{pmatrix}$$

**U_{X₁[e₇]}[X₁[e₃]] == ZERO**

True

**defin[X₁[e₇], X₁[e₃]] // Expand**

$$\begin{pmatrix} \begin{pmatrix} 0 & \{0,0,0\} \\ \{0,0,0\} & 0 \end{pmatrix} & \begin{pmatrix} 0 & \{0,0,0\} \\ \{0,-\eta_2-\lambda_9,0\} & 0 \end{pmatrix} & \begin{pmatrix} 0 & \{0,0,0\} \\ \{0,0,0\} & 0 \end{pmatrix} \\ \begin{pmatrix} 0 & \{0,0,0\} \\ \{0,\eta_2+\lambda_9,0\} & 0 \end{pmatrix} & \begin{pmatrix} 0 & \{0,0,0\} \\ \{0,0,0\} & 0 \end{pmatrix} & \begin{pmatrix} 0 & \{0,0,0\} \\ \{0,0,0\} & 0 \end{pmatrix} \\ \begin{pmatrix} 0 & \{0,0,0\} \\ \{0,0,0\} & 0 \end{pmatrix} & \begin{pmatrix} 0 & \{0,0,0\} \\ \{0,0,0\} & 0 \end{pmatrix} & \begin{pmatrix} 0 & \{0,0,0\} \\ \{0,0,0\} & 0 \end{pmatrix} \end{pmatrix}$$

**Δ[X₁[e₇]] = Δ[X₁[e₇]] //. {λ₉ → -η₂}**

$$\begin{pmatrix} \begin{pmatrix} \alpha_4 & \{0,0,0\} \\ \{0,0,0\} & \alpha_4 \end{pmatrix} & \begin{pmatrix} \lambda_4 & \{-\psi_1,0,\lambda_8\} \\ \{-\eta_2,-\epsilon_2,-\phi_2\} & \lambda_5 \end{pmatrix} & \begin{pmatrix} 0 & \{\lambda_{22},0,\lambda_{24}\} \\ \{0,\lambda_{26},0\} & \lambda_{21} \end{pmatrix} \\ \begin{pmatrix} \lambda_5 & \{\psi_1,0,-\lambda_8\} \\ \{\eta_2,\epsilon_2,\phi_2\} & \lambda_4 \end{pmatrix} & \begin{pmatrix} -\alpha_4 & \{0,0,0\} \\ \{0,0,0\} & -\alpha_4 \end{pmatrix} & \begin{pmatrix} 0 & \{\lambda_{14},0,\lambda_{16}\} \\ \{0,\lambda_{18},0\} & \lambda_{13} \end{pmatrix} \\ \begin{pmatrix} \lambda_{21} & \{-\lambda_{22},0,-\lambda_{24}\} \\ \{0,-\lambda_{26},0\} & 0 \end{pmatrix} & \begin{pmatrix} \lambda_{13} & \{-\lambda_{14},0,-\lambda_{16}\} \\ \{0,-\lambda_{18},0\} & 0 \end{pmatrix} & \begin{pmatrix} \lambda_3 & \{0,0,0\} \\ \{0,0,0\} & \lambda_3 \end{pmatrix} \end{pmatrix}$$

**U_{X₁[e₇]}[X₁[e₂]] == ZERO**

True

**defin[X₁[e₇], X₁[e₂]] // Expand**

$$\begin{pmatrix} \begin{pmatrix} 0 & \{0,0,0\} \\ \{0,0,0\} & 0 \end{pmatrix} & \begin{pmatrix} 0 & \{0,0,0\} \\ \{0, \lambda_4 - \rho_2, 0\} & 0 \end{pmatrix} & \begin{pmatrix} 0 & \{0,0,0\} \\ \{0,0,0\} & 0 \end{pmatrix} \\ \begin{pmatrix} 0 & \{0,0,0\} \\ \{0, \rho_2 - \lambda_4, 0\} & 0 \end{pmatrix} & \begin{pmatrix} 0 & \{0,0,0\} \\ \{0,0,0\} & 0 \end{pmatrix} & \begin{pmatrix} 0 & \{0,0,0\} \\ \{0,0,0\} & 0 \end{pmatrix} \\ \begin{pmatrix} 0 & \{0,0,0\} \\ \{0,0,0\} & 0 \end{pmatrix} & \begin{pmatrix} 0 & \{0,0,0\} \\ \{0,0,0\} & 0 \end{pmatrix} & \begin{pmatrix} 0 & \{0,0,0\} \\ \{0,0,0\} & 0 \end{pmatrix} \end{pmatrix}$$

**Δ[X₁[e₇]] = Δ[X₁[e₇]] //. {λ₄ → ρ₂}**

$$\begin{pmatrix} \begin{pmatrix} \alpha_4 & \{0,0,0\} \\ \{0,0,0\} & \alpha_4 \end{pmatrix} & \begin{pmatrix} \rho_2 & \{-\psi_1, 0, \lambda_8\} \\ \{-\eta_2, -\epsilon_2, -\phi_2\} & \lambda_5 \end{pmatrix} & \begin{pmatrix} 0 & \{\lambda_{22}, 0, \lambda_{24}\} \\ \{0, \lambda_{26}, 0\} & \lambda_{21} \end{pmatrix} \\ \begin{pmatrix} \lambda_5 & \{\psi_1, 0, -\lambda_8\} \\ \{\eta_2, \epsilon_2, \phi_2\} & \rho_2 \end{pmatrix} & \begin{pmatrix} -\alpha_4 & \{0,0,0\} \\ \{0,0,0\} & -\alpha_4 \end{pmatrix} & \begin{pmatrix} 0 & \{\lambda_{14}, 0, \lambda_{16}\} \\ \{0, \lambda_{18}, 0\} & \lambda_{13} \end{pmatrix} \\ \begin{pmatrix} \lambda_{21} & \{-\lambda_{22}, 0, -\lambda_{24}\} \\ \{0, -\lambda_{26}, 0\} & 0 \end{pmatrix} & \begin{pmatrix} \lambda_{13} & \{-\lambda_{14}, 0, -\lambda_{16}\} \\ \{0, -\lambda_{18}, 0\} & 0 \end{pmatrix} & \begin{pmatrix} \lambda_3 & \{0,0,0\} \\ \{0,0,0\} & \lambda_3 \end{pmatrix} \end{pmatrix}$$

**U_{X₁[e₇]}[X₁[e₁]] == ZERO**

True

**defin[X₁[e₇], X₁[e₁]] // Expand**

$$\begin{pmatrix} \begin{pmatrix} 0 & \{0,0,0\} \\ \{0,0,0\} & 0 \end{pmatrix} & \begin{pmatrix} 0 & \{0,0,0\} \\ \{0, \lambda_5 - \delta_3, 0\} & 0 \end{pmatrix} & \begin{pmatrix} 0 & \{0,0,0\} \\ \{0,0,0\} & 0 \end{pmatrix} \\ \begin{pmatrix} 0 & \{0,0,0\} \\ \{0, \delta_3 - \lambda_5, 0\} & 0 \end{pmatrix} & \begin{pmatrix} 0 & \{0,0,0\} \\ \{0,0,0\} & 0 \end{pmatrix} & \begin{pmatrix} 0 & \{0,0,0\} \\ \{0,0,0\} & 0 \end{pmatrix} \\ \begin{pmatrix} 0 & \{0,0,0\} \\ \{0,0,0\} & 0 \end{pmatrix} & \begin{pmatrix} 0 & \{0,0,0\} \\ \{0,0,0\} & 0 \end{pmatrix} & \begin{pmatrix} 0 & \{0,0,0\} \\ \{0,0,0\} & 0 \end{pmatrix} \end{pmatrix}$$

**Δ[X₁[e₇]] = Δ[X₁[e₇]] //. {λ₅ → δ₃}**

$$\begin{pmatrix} \begin{pmatrix} \alpha_4 & \{0,0,0\} \\ \{0,0,0\} & \alpha_4 \end{pmatrix} & \begin{pmatrix} \rho_2 & \{-\psi_1, 0, \lambda_8\} \\ \{-\eta_2, -\epsilon_2, -\phi_2\} & \delta_3 \end{pmatrix} & \begin{pmatrix} 0 & \{\lambda_{22}, 0, \lambda_{24}\} \\ \{0, \lambda_{26}, 0\} & \lambda_{21} \end{pmatrix} \\ \begin{pmatrix} \delta_3 & \{\psi_1, 0, -\lambda_8\} \\ \{\eta_2, \epsilon_2, \phi_2\} & \rho_2 \end{pmatrix} & \begin{pmatrix} -\alpha_4 & \{0,0,0\} \\ \{0,0,0\} & -\alpha_4 \end{pmatrix} & \begin{pmatrix} 0 & \{\lambda_{14}, 0, \lambda_{16}\} \\ \{0, \lambda_{18}, 0\} & \lambda_{13} \end{pmatrix} \\ \begin{pmatrix} \lambda_{21} & \{-\lambda_{22}, 0, -\lambda_{24}\} \\ \{0, -\lambda_{26}, 0\} & 0 \end{pmatrix} & \begin{pmatrix} \lambda_{13} & \{-\lambda_{14}, 0, -\lambda_{16}\} \\ \{0, -\lambda_{18}, 0\} & 0 \end{pmatrix} & \begin{pmatrix} \lambda_3 & \{0,0,0\} \\ \{0,0,0\} & \lambda_3 \end{pmatrix} \end{pmatrix}$$

**T[E₁[1], E₁[1], X₁[e₇]] == X₁[e₇]**

True

**Δ[X₁[e₇]] - Leib[E₁[1], E₁[1], X₁[e₇]] // Expand**

$$\left(\begin{pmatrix} 0 & \{0,0,0\} \\ \{0,0,0\} & 0 \end{pmatrix} \quad \begin{pmatrix} 0 & \{0,0,0\} \\ \{0,0,0\} & 0 \end{pmatrix} \quad \begin{pmatrix} 0 & \{0,0,0\} \\ \{0,0,0\} & 0 \end{pmatrix}\right.$$

$$\begin{pmatrix} 0 & \{0,0,0\} \\ \{0,0,0\} & 0 \end{pmatrix} \quad \begin{pmatrix} 0 & \{0,0,0\} \\ \{0,0,0\} & 0 \end{pmatrix} \quad \begin{pmatrix} 0 & \{\lambda_{14}-\beta_8, 0, \beta_6+\lambda_{16}\} \\ \{0, \beta_1+\lambda_{18}, 0\} & \beta_4+\lambda_{13} \end{pmatrix}$$

$$\left.\begin{pmatrix} 0 & \{0,0,0\} \\ \{0,0,0\} & 0 \end{pmatrix} \quad \begin{pmatrix} \beta_4+\lambda_{13} & \{\beta_8-\lambda_{14}, 0, -\beta_6-\lambda_{16}\} \\ \{0, -\beta_1-\lambda_{18}, 0\} & 0 \end{pmatrix} \quad \begin{pmatrix} \lambda_3 & \{0,0,0\} \\ \{0,0,0\} & \lambda_3 \end{pmatrix}\right)$$

**Δ[X₁[e₇]] =
Δ[X₁[e₇]] //. {λ₃ → 0, λ₁₃ → -β₄, λ₁₄ → β₈, λ₁₆ → -β₆, λ₁₈ → -β₁}**

$$\left(\begin{pmatrix} \alpha_4 & \{0,0,0\} \\ \{0,0,0\} & \alpha_4 \end{pmatrix} \quad \begin{pmatrix} \rho_2 & \{-\psi_1, 0, \lambda_8\} \\ \{-\eta_2, -\epsilon_2, -\phi_2\} & \delta_3 \end{pmatrix} \quad \begin{pmatrix} 0 & \{\lambda_{22}, 0, \lambda_{24}\} \\ \{0, \lambda_{26}, 0\} & \lambda_{21} \end{pmatrix}\right.$$

$$\begin{pmatrix} \delta_3 & \{\psi_1, 0, -\lambda_8\} \\ \{\eta_2, \epsilon_2, \phi_2\} & \rho_2 \end{pmatrix} \quad \begin{pmatrix} -\alpha_4 & \{0,0,0\} \\ \{0,0,0\} & -\alpha_4 \end{pmatrix} \quad \begin{pmatrix} 0 & \{\beta_8, 0, -\beta_6\} \\ \{0, -\beta_1, 0\} & -\beta_4 \end{pmatrix}$$

$$\left.\begin{pmatrix} \lambda_{21} & \{-\lambda_{22}, 0, -\lambda_{24}\} \\ \{0, -\lambda_{26}, 0\} & 0 \end{pmatrix} \quad \begin{pmatrix} -\beta_4 & \{-\beta_8, 0, \beta_6\} \\ \{0, \beta_1, 0\} & 0 \end{pmatrix} \quad \begin{pmatrix} 0 & \{0,0,0\} \\ \{0,0,0\} & 0 \end{pmatrix}\right)$$

**T[E₁[1], X₁[e₇], E₁[1]] == ZERO**

True

**Leib[E₁[1], X₁[e₇], E₁[1]] == ZERO**

True

**T[E₂[1], E₂[1], X₁[e₇]] == X₁[e₇]**

True

**Δ[X₁[e₇]] - Leib[E₂[1], E₂[1], X₁[e₇]] // Expand**

$$\left(\begin{pmatrix} 0 & \{0,0,0\} \\ \{0,0,0\} & 0 \end{pmatrix} \quad \begin{pmatrix} 0 & \{0,0,0\} \\ \{0,0,0\} & 0 \end{pmatrix} \quad \begin{pmatrix} 0 & \{\gamma_8+\lambda_{22}, 0, \lambda_{24}-\gamma_6\} \\ \{0, \lambda_{26}-\gamma_1, 0\} & \lambda_{21}-\gamma_4 \end{pmatrix}\right.$$

$$\begin{pmatrix} 0 & \{0,0,0\} \\ \{0,0,0\} & 0 \end{pmatrix} \quad \begin{pmatrix} 0 & \{0,0,0\} \\ \{0,0,0\} & 0 \end{pmatrix} \quad \begin{pmatrix} 0 & \{0,0,0\} \\ \{0,0,0\} & 0 \end{pmatrix}$$

$$\left.\begin{pmatrix} \lambda_{21}-\gamma_4 & \{-\gamma_8-\lambda_{22}, 0, \gamma_6-\lambda_{24}\} \\ \{0, \gamma_1-\lambda_{26}, 0\} & 0 \end{pmatrix} \quad \begin{pmatrix} 0 & \{0,0,0\} \\ \{0,0,0\} & 0 \end{pmatrix} \quad \begin{pmatrix} 0 & \{0,0,0\} \\ \{0,0,0\} & 0 \end{pmatrix}\right)$$

**Δ[X₁[e₇]] = Δ[X₁[e₇]] //. {λ₂₁ → γ₄, λ₂₂ → -γ₈, λ₂₄ → γ₆, λ₂₆ → γ₁}**

$$\left(\begin{pmatrix} \alpha_4 & \{0,0,0\} \\ \{0,0,0\} & \alpha_4 \end{pmatrix} \quad \begin{pmatrix} \rho_2 & \{-\psi_1, 0, \lambda_8\} \\ \{-\eta_2, -\epsilon_2, -\phi_2\} & \delta_3 \end{pmatrix} \quad \begin{pmatrix} 0 & \{-\gamma_8, 0, \gamma_6\} \\ \{0, \gamma_1, 0\} & \gamma_4 \end{pmatrix}\right.$$

$$\begin{pmatrix} \delta_3 & \{\psi_1, 0, -\lambda_8\} \\ \{\eta_2, \epsilon_2, \phi_2\} & \rho_2 \end{pmatrix} \quad \begin{pmatrix} -\alpha_4 & \{0,0,0\} \\ \{0,0,0\} & -\alpha_4 \end{pmatrix} \quad \begin{pmatrix} 0 & \{\beta_8, 0, -\beta_6\} \\ \{0, -\beta_1, 0\} & -\beta_4 \end{pmatrix}$$

$$\left.\begin{pmatrix} \gamma_4 & \{\gamma_8, 0, -\gamma_6\} \\ \{0, -\gamma_1, 0\} & 0 \end{pmatrix} \quad \begin{pmatrix} -\beta_4 & \{-\beta_8, 0, \beta_6\} \\ \{0, \beta_1, 0\} & 0 \end{pmatrix} \quad \begin{pmatrix} 0 & \{0,0,0\} \\ \{0,0,0\} & 0 \end{pmatrix}\right)$$

**T[E₂[1], X₁[e₇], E₂[1]] == ZERO**

True

**Leib[E₂[1], X₁[e₇], E₂[1]] == ZERO**

True

**T[E₁[1], E₂[1], X₁[e₇]] == ZERO**

True

**(Leib[E₁[1], E₂[1], X₁[e₇]] // Expand) == ZERO**

True

**T[E₂[1], E₁[1], X₁[e₇]] == ZERO**

True

**(Leib[E₂[1], E₁[1], X₁[e₇]] // Expand) == ZERO**

True

**T[E₁[1], X₁[e₇], E₂[1]] == X₁[e₇]**

True

**Δ[X₁[e₇]] − Leib[E₁[1], X₁[e₇], E₂[1]] == ZERO**

True

**T[E₁[1], X₁[e₁], X₁[e₇]] == ZERO**

True

**(Leib[E₁[1], X₁[e₁], X₁[e₇]] // Expand) == ZERO**

True

**T[X₁[e₁], E₁[1], X₁[e₇]] == ZERO**

True

**(Leib[X₁[e₁], E₁[1], X₁[e₇]] // Expand) == ZERO**

True

```
T[E₂[1], X₁[e₁], X₁[e₇]] == ZERO
```

True

```
(Leib[E₂[1], X₁[e₁], X₁[e₇]] // Expand) == ZERO
```

True

```
T[X₁[e₆], X₁[e₇], X₁[e₁]] == ZERO
```

True

```
(Leib[X₁[e₆], X₁[e₇], X₁[e₁]] // Expand) == ZERO
```

True

```
T[X₁[e₅], X₁[e₇], X₁[e₆]] == ZERO
```

True

```
(Leib[X₁[e₅], X₁[e₇], X₁[e₆]] // Expand) == ZERO
```

True

```
T[X₁[e₄], X₁[e₆], X₁[e₇]] == X₁[e₆]
```

True

```
(Δ[X₁[e₆]] - Leib[X₁[e₄], X₁[e₆], X₁[e₇]] // Expand) == ZERO
```

True

```
T[X₁[e₃], X₁[e₁], X₁[e₇]] == ZERO
```

True

```
(Leib[X₁[e₃], X₁[e₁], X₁[e₇]] // Expand) == ZERO
```

True

```
T[X₁[e₇], X₁[e₆], X₁[e₅]] == ZERO
```

True

```
(Leib[X₁[e₇], X₁[e₆], X₁[e₅]] // Expand) == ZERO
```

True

```
T[X₁[e₇], X₁[e₇], X₁[e₇]] == ZERO
```

True

```
(Leib[X₁[e₇], X₁[e₇], X₁[e₇]] // Expand) == ZERO
```

True

```
Δ[X₁[e₇]] = Δ[X₁[e₇]] //. {λ₈ → χ₁}
```

$$\begin{pmatrix} \begin{pmatrix} \alpha_4 & \{0,0,0\} \\ \{0,0,0\} & \alpha_4 \end{pmatrix} & \begin{pmatrix} \rho_2 & \{-\psi_1, 0, \chi_1\} \\ \{-\eta_2, -\epsilon_2, -\phi_2\} & \delta_3 \end{pmatrix} & \begin{pmatrix} 0 & \{-\gamma_8, 0, \gamma_6\} \\ \{0, \gamma_1, 0\} & \gamma_4 \end{pmatrix} \\ \begin{pmatrix} \delta_3 & \{\psi_1, 0, -\chi_1\} \\ \{\eta_2, \epsilon_2, \phi_2\} & \rho_2 \end{pmatrix} & \begin{pmatrix} -\alpha_4 & \{0,0,0\} \\ \{0,0,0\} & -\alpha_4 \end{pmatrix} & \begin{pmatrix} 0 & \{\beta_8, 0, -\beta_6\} \\ \{0, -\beta_1, 0\} & -\beta_4 \end{pmatrix} \\ \begin{pmatrix} \gamma_4 & \{\gamma_8, 0, -\gamma_6\} \\ \{0, -\gamma_1, 0\} & 0 \end{pmatrix} & \begin{pmatrix} -\beta_4 & \{-\beta_8, 0, \beta_6\} \\ \{0, \beta_1, 0\} & 0 \end{pmatrix} & \begin{pmatrix} 0 & \{0,0,0\} \\ \{0,0,0\} & 0 \end{pmatrix} \end{pmatrix}$$

- **$X_1[e_8]$**

  ```
  Δ[X₁[e₈]] = generic;
  ```

  ```
  U_{X₁[e₈]}[E₁[1]] == ZERO
  ```

  True

  ```
  defin[X₁[e₈], E₁[1]] // Expand
  ```

  $$\begin{pmatrix} \begin{pmatrix} 0 & \{0,0,0\} \\ \{0,0,0\} & 0 \end{pmatrix} & \begin{pmatrix} 0 & \{0,0,0\} \\ \{0,0,\lambda_1-\alpha_5\} & 0 \end{pmatrix} & \begin{pmatrix} 0 & \{0,0,0\} \\ \{0,0,0\} & 0 \end{pmatrix} \\ \begin{pmatrix} 0 & \{0,0,0\} \\ \{0,0,\alpha_5-\lambda_1\} & 0 \end{pmatrix} & \begin{pmatrix} -\lambda_8 & \{0,0,0\} \\ \{0,0,0\} & -\lambda_8 \end{pmatrix} & \begin{pmatrix} 0 & \{-\lambda_{26}, \lambda_{25}, 0\} \\ \{0,0,-\lambda_{20}\} & -\lambda_{24} \end{pmatrix} \\ \begin{pmatrix} 0 & \{0,0,0\} \\ \{0,0,0\} & 0 \end{pmatrix} & \begin{pmatrix} -\lambda_{24} & \{\lambda_{26}, -\lambda_{25}, 0\} \\ \{0,0,\lambda_{20}\} & 0 \end{pmatrix} & \begin{pmatrix} 0 & \{0,0,0\} \\ \{0,0,0\} & 0 \end{pmatrix} \end{pmatrix}$$

  ```
  Δ[X₁[e₈]] =
    Δ[X₁[e₈]] //. {λ₁ → α₅, λ₈ → 0, λ₂₀ → 0, λ₂₄ → 0, λ₂₅ → 0, λ₂₆ → 0}
  ```

  $$\begin{pmatrix} \begin{pmatrix} \alpha_5 & \{0,0,0\} \\ \{0,0,0\} & \alpha_5 \end{pmatrix} & \begin{pmatrix} \lambda_4 & \{\lambda_6, \lambda_7, 0\} \\ \{\lambda_9, \lambda_{10}, \lambda_{11}\} & \lambda_5 \end{pmatrix} & \begin{pmatrix} 0 & \{\lambda_{22}, \lambda_{23} \\ \{0,0,\lambda_{27}\} & \lambda_{21} \end{pmatrix} \\ \begin{pmatrix} \lambda_5 & \{-\lambda_6, -\lambda_7, 0\} \\ \{-\lambda_9, -\lambda_{10}, -\lambda_{11}\} & \lambda_4 \end{pmatrix} & \begin{pmatrix} \lambda_2 & \{0,0,0\} \\ \{0,0,0\} & \lambda_2 \end{pmatrix} & \begin{pmatrix} \lambda_{12} & \{\lambda_{14}, \lambda_1 \\ \{\lambda_{17}, \lambda_{18}, \lambda_{19}\} & \lambda_1 \end{pmatrix} \\ \begin{pmatrix} \lambda_{21} & \{-\lambda_{22}, -\lambda_{23}, 0\} \\ \{0,0,-\lambda_{27}\} & 0 \end{pmatrix} & \begin{pmatrix} \lambda_{13} & \{-\lambda_{14}, -\lambda_{15}, -\lambda_{16}\} \\ \{-\lambda_{17}, -\lambda_{18}, -\lambda_{19}\} & \lambda_{12} \end{pmatrix} & \begin{pmatrix} \lambda_3 & \{0,0,0 \\ \{0,0,0\} & \lambda_3 \end{pmatrix} \end{pmatrix}$$

  ```
  U_{X₁[e₈]}[E₂[1]] == ZERO
  ```

  True

```
defin[X₁[e₈], E₂[1]] // Expand
```

$$\left(\begin{array}{ccc} \begin{pmatrix} 0 & \{0,0,0\} \\ \{0,0,0\} & 0 \end{pmatrix} & \begin{pmatrix} 0 & \{0,0,0\} \\ \{0,0,\alpha_5+\lambda_2\} & 0 \end{pmatrix} & \begin{pmatrix} 0 & \{\lambda_{18},-\lambda_{17},0\} \\ \{0,0,\lambda_{12}\} & \lambda_{16} \end{pmatrix} \\ \begin{pmatrix} 0 & \{0,0,0\} \\ \{0,0,-\alpha_5-\lambda_2\} & 0 \end{pmatrix} & \begin{pmatrix} 0 & \{0,0,0\} \\ \{0,0,0\} & 0 \end{pmatrix} & \begin{pmatrix} 0 & \{0,0,0\} \\ \{0,0,0\} & 0 \end{pmatrix} \\ \begin{pmatrix} \lambda_{16} & \{-\lambda_{18},\lambda_{17},0\} \\ \{0,0,-\lambda_{12}\} & 0 \end{pmatrix} & \begin{pmatrix} 0 & \{0,0,0\} \\ \{0,0,0\} & 0 \end{pmatrix} & \begin{pmatrix} 0 & \{0,0,0\} \\ \{0,0,0\} & 0 \end{pmatrix} \end{array}\right)$$

```
Δ[X₁[e₈]] = Δ[X₁[e₈]] //. {λ₂ → -α₅, λ₁₂ → 0, λ₁₆ → 0, λ₁₇ → 0, λ₁₈ → 0}
```

$$\left(\begin{array}{ccc} \begin{pmatrix} \alpha_5 & \{0,0,0\} \\ \{0,0,0\} & \alpha_5 \end{pmatrix} & \begin{pmatrix} \lambda_4 & \{\lambda_6,\lambda_7,0\} \\ \{\lambda_9,\lambda_{10},\lambda_{11}\} & \lambda_5 \end{pmatrix} & \begin{pmatrix} 0 & \{\lambda_{22},\lambda_{23},0\} \\ \{0,0,\lambda_{27}\} & \lambda_{21} \end{pmatrix} \\ \begin{pmatrix} \lambda_5 & \{-\lambda_6,-\lambda_7,0\} \\ \{-\lambda_9,-\lambda_{10},-\lambda_{11}\} & \lambda_4 \end{pmatrix} & \begin{pmatrix} -\alpha_5 & \{0,0,0\} \\ \{0,0,0\} & -\alpha_5 \end{pmatrix} & \begin{pmatrix} 0 & \{\lambda_{14},\lambda_{15},0\} \\ \{0,0,\lambda_{19}\} & \lambda_{13} \end{pmatrix} \\ \begin{pmatrix} \lambda_{21} & \{-\lambda_{22},-\lambda_{23},0\} \\ \{0,0,-\lambda_{27}\} & 0 \end{pmatrix} & \begin{pmatrix} \lambda_{13} & \{-\lambda_{14},-\lambda_{15},0\} \\ \{0,0,-\lambda_{19}\} & 0 \end{pmatrix} & \begin{pmatrix} \lambda_3 & \{0,0,0\} \\ \{0,0,0\} & \lambda_3 \end{pmatrix} \end{array}\right)$$

```
U_{X₁[e₈]}[X₁[e₈]] == ZERO
```

True

```
(defin[X₁[e₈], X₁[e₈]] // Expand) == ZERO
```

True

```
U_{X₁[e₈]}[X₁[e₇]] == ZERO
```

True

```
defin[X₁[e₈], X₁[e₇]] // Expand
```

$$\left(\begin{array}{ccc} \begin{pmatrix} 0 & \{0,0,0\} \\ \{0,0,0\} & 0 \end{pmatrix} & \begin{pmatrix} 0 & \{0,0,0\} \\ \{0,0,-\lambda_7-\chi_1\} & 0 \end{pmatrix} & \begin{pmatrix} 0 & \{0,0,0\} \\ \{0,0,0\} & 0 \end{pmatrix} \\ \begin{pmatrix} 0 & \{0,0,0\} \\ \{0,0,\lambda_7+\chi_1\} & 0 \end{pmatrix} & \begin{pmatrix} 0 & \{0,0,0\} \\ \{0,0,0\} & 0 \end{pmatrix} & \begin{pmatrix} 0 & \{0,0,0\} \\ \{0,0,0\} & 0 \end{pmatrix} \\ \begin{pmatrix} 0 & \{0,0,0\} \\ \{0,0,0\} & 0 \end{pmatrix} & \begin{pmatrix} 0 & \{0,0,0\} \\ \{0,0,0\} & 0 \end{pmatrix} & \begin{pmatrix} 0 & \{0,0,0\} \\ \{0,0,0\} & 0 \end{pmatrix} \end{array}\right)$$

```
Δ[X₁[e₈]] = Δ[X₁[e₈]] //. {λ₇ → -χ₁}
```

$$\left(\begin{array}{ccc} \begin{pmatrix} \alpha_5 & \{0,0,0\} \\ \{0,0,0\} & \alpha_5 \end{pmatrix} & \begin{pmatrix} \lambda_4 & \{\lambda_6,-\chi_1,0\} \\ \{\lambda_9,\lambda_{10},\lambda_{11}\} & \lambda_5 \end{pmatrix} & \begin{pmatrix} 0 & \{\lambda_{22},\lambda_{23},0\} \\ \{0,0,\lambda_{27}\} & \lambda_{21} \end{pmatrix} \\ \begin{pmatrix} \lambda_5 & \{-\lambda_6,\chi_1,0\} \\ \{-\lambda_9,-\lambda_{10},-\lambda_{11}\} & \lambda_4 \end{pmatrix} & \begin{pmatrix} -\alpha_5 & \{0,0,0\} \\ \{0,0,0\} & -\alpha_5 \end{pmatrix} & \begin{pmatrix} 0 & \{\lambda_{14},\lambda_{15},0\} \\ \{0,0,\lambda_{19}\} & \lambda_{13} \end{pmatrix} \\ \begin{pmatrix} \lambda_{21} & \{-\lambda_{22},-\lambda_{23},0\} \\ \{0,0,-\lambda_{27}\} & 0 \end{pmatrix} & \begin{pmatrix} \lambda_{13} & \{-\lambda_{14},-\lambda_{15},0\} \\ \{0,0,-\lambda_{19}\} & 0 \end{pmatrix} & \begin{pmatrix} \lambda_3 & \{0,0,0\} \\ \{0,0,0\} & \lambda_3 \end{pmatrix} \end{array}\right)$$

$U_{X_1[e_8]}[X_1[e_6]] == \text{ZERO}$

True

$\text{defin}[X_1[e_8], X_1[e_6]] \,/\!/\, \text{Expand}$

$$\begin{pmatrix} \begin{pmatrix} 0 & \{0,0,0\} \\ \{0,0,0\} & 0 \end{pmatrix} & \begin{pmatrix} 0 & \{0,0,0\} \\ \{0,0,-\lambda_6-\psi_2\} & 0 \end{pmatrix} & \begin{pmatrix} 0 & \{0,0,0\} \\ \{0,0,0\} & 0 \end{pmatrix} \\ \begin{pmatrix} 0 & \{0,0,0\} \\ \{0,0,\lambda_6+\psi_2\} & 0 \end{pmatrix} & \begin{pmatrix} 0 & \{0,0,0\} \\ \{0,0,0\} & 0 \end{pmatrix} & \begin{pmatrix} 0 & \{0,0,0\} \\ \{0,0,0\} & 0 \end{pmatrix} \\ \begin{pmatrix} 0 & \{0,0,0\} \\ \{0,0,0\} & 0 \end{pmatrix} & \begin{pmatrix} 0 & \{0,0,0\} \\ \{0,0,0\} & 0 \end{pmatrix} & \begin{pmatrix} 0 & \{0,0,0\} \\ \{0,0,0\} & 0 \end{pmatrix} \end{pmatrix}$$

$\Delta[X_1[e_8]] = \Delta[X_1[e_8]] \,/\!/.\, \{\lambda_6 \to -\psi_2\}$

$$\begin{pmatrix} \begin{pmatrix} \alpha_5 & \{0,0,0\} \\ \{0,0,0\} & \alpha_5 \end{pmatrix} & \begin{pmatrix} \lambda_4 & \{-\psi_2,-\chi_1,0\} \\ \{\lambda_9,\lambda_{10},\lambda_{11}\} & \lambda_5 \end{pmatrix} & \begin{pmatrix} 0 & \{\lambda_{22},\lambda_{23},0\} \\ \{0,0,\lambda_{27}\} & \lambda_{21} \end{pmatrix} \\ \begin{pmatrix} \lambda_5 & \{\psi_2,\chi_1,0\} \\ \{-\lambda_9,-\lambda_{10},-\lambda_{11}\} & \lambda_4 \end{pmatrix} & \begin{pmatrix} -\alpha_5 & \{0,0,0\} \\ \{0,0,0\} & -\alpha_5 \end{pmatrix} & \begin{pmatrix} 0 & \{\lambda_{14},\lambda_{15},0\} \\ \{0,0,\lambda_{19}\} & \lambda_{13} \end{pmatrix} \\ \begin{pmatrix} \lambda_{21} & \{-\lambda_{22},-\lambda_{23},0\} \\ \{0,0,-\lambda_{27}\} & 0 \end{pmatrix} & \begin{pmatrix} \lambda_{13} & \{-\lambda_{14},-\lambda_{15},0\} \\ \{0,0,-\lambda_{19}\} & 0 \end{pmatrix} & \begin{pmatrix} \lambda_3 & \{0,0,0\} \\ \{0,0,0\} & \lambda_3 \end{pmatrix} \end{pmatrix}$$

$U_{X_1[e_8]}[X_1[e_5]] == -X_1[e_8]$

True

$\Delta[X_1[e_8]] + \text{defin}[X_1[e_8], X_1[e_5]] \,/\!/\, \text{Expand}$

$$\begin{pmatrix} \begin{pmatrix} 0 & \{0,0,0\} \\ \{0,0,0\} & 0 \end{pmatrix} & \begin{pmatrix} 0 & \{0,0,0\} \\ \{0,0,-\lambda_{11}-\phi_3\} & 0 \end{pmatrix} & \begin{pmatrix} 0 & \{0,0,0\} \\ \{0,0,0\} & 0 \end{pmatrix} \\ \begin{pmatrix} 0 & \{0,0,0\} \\ \{0,0,\lambda_{11}+\phi_3\} & 0 \end{pmatrix} & \begin{pmatrix} 0 & \{0,0,0\} \\ \{0,0,0\} & 0 \end{pmatrix} & \begin{pmatrix} 0 & \{0,0,0\} \\ \{0,0,0\} & 0 \end{pmatrix} \\ \begin{pmatrix} 0 & \{0,0,0\} \\ \{0,0,0\} & 0 \end{pmatrix} & \begin{pmatrix} 0 & \{0,0,0\} \\ \{0,0,0\} & 0 \end{pmatrix} & \begin{pmatrix} \lambda_3 & \{0,0,0\} \\ \{0,0,0\} & \lambda_3 \end{pmatrix} \end{pmatrix}$$

$\Delta[X_1[e_8]] = \Delta[X_1[e_8]] \,/\!/.\, \{\lambda_3 \to 0,\, \lambda_{11} \to -\phi_3\}$

$$\begin{pmatrix} \begin{pmatrix} \alpha_5 & \{0,0,0\} \\ \{0,0,0\} & \alpha_5 \end{pmatrix} & \begin{pmatrix} \lambda_4 & \{-\psi_2,-\chi_1,0\} \\ \{\lambda_9,\lambda_{10},-\phi_3\} & \lambda_5 \end{pmatrix} & \begin{pmatrix} 0 & \{\lambda_{22},\lambda_{23},0\} \\ \{0,0,\lambda_{27}\} & \lambda_{21} \end{pmatrix} \\ \begin{pmatrix} \lambda_5 & \{\psi_2,\chi_1,0\} \\ \{-\lambda_9,-\lambda_{10},\phi_3\} & \lambda_4 \end{pmatrix} & \begin{pmatrix} -\alpha_5 & \{0,0,0\} \\ \{0,0,0\} & -\alpha_5 \end{pmatrix} & \begin{pmatrix} 0 & \{\lambda_{14},\lambda_{15},0\} \\ \{0,0,\lambda_{19}\} & \lambda_{13} \end{pmatrix} \\ \begin{pmatrix} \lambda_{21} & \{-\lambda_{22},-\lambda_{23},0\} \\ \{0,0,-\lambda_{27}\} & 0 \end{pmatrix} & \begin{pmatrix} \lambda_{13} & \{-\lambda_{14},-\lambda_{15},0\} \\ \{0,0,-\lambda_{19}\} & 0 \end{pmatrix} & \begin{pmatrix} 0 & \{0,0,0\} \\ \{0,0,0\} & 0 \end{pmatrix} \end{pmatrix}$$

$U_{X_1[e_8]}[X_1[e_4]] == \text{ZERO}$

True

**defin[X₁[e₈], X₁[e₄]] // Expand**

$$\begin{pmatrix} \begin{pmatrix} 0 & \{0,0,0\} \\ \{0,0,0\} & 0 \end{pmatrix} & \begin{pmatrix} 0 & \{0,0,0\} \\ \{0,0,-\epsilon_3-\lambda_{10}\} & 0 \end{pmatrix} & \begin{pmatrix} 0 & \{0,0,0\} \\ \{0,0,0\} & 0 \end{pmatrix} \\ \begin{pmatrix} 0 & \{0,0,0\} \\ \{0,0,\epsilon_3+\lambda_{10}\} & 0 \end{pmatrix} & \begin{pmatrix} 0 & \{0,0,0\} \\ \{0,0,0\} & 0 \end{pmatrix} & \begin{pmatrix} 0 & \{0,0,0\} \\ \{0,0,0\} & 0 \end{pmatrix} \\ \begin{pmatrix} 0 & \{0,0,0\} \\ \{0,0,0\} & 0 \end{pmatrix} & \begin{pmatrix} 0 & \{0,0,0\} \\ \{0,0,0\} & 0 \end{pmatrix} & \begin{pmatrix} 0 & \{0,0,0\} \\ \{0,0,0\} & 0 \end{pmatrix} \end{pmatrix}$$

**Δ[X₁[e₈]] = Δ[X₁[e₈]] //. {λ₁₀ → -ε₃}**

$$\begin{pmatrix} \begin{pmatrix} \alpha_5 & \{0,0,0\} \\ \{0,0,0\} & \alpha_5 \end{pmatrix} & \begin{pmatrix} \lambda_4 & \{-\psi_2,-\chi_1,0\} \\ \{\lambda_9,-\epsilon_3,-\phi_3\} & \lambda_5 \end{pmatrix} & \begin{pmatrix} 0 & \{\lambda_{22},\lambda_{23},0\} \\ \{0,0,\lambda_{27}\} & \lambda_{21} \end{pmatrix} \\ \begin{pmatrix} \lambda_5 & \{\psi_2,\chi_1,0\} \\ \{-\lambda_9,\epsilon_3,\phi_3\} & \lambda_4 \end{pmatrix} & \begin{pmatrix} -\alpha_5 & \{0,0,0\} \\ \{0,0,0\} & -\alpha_5 \end{pmatrix} & \begin{pmatrix} 0 & \{\lambda_{14},\lambda_{15},0\} \\ \{0,0,\lambda_{19}\} & \lambda_{13} \end{pmatrix} \\ \begin{pmatrix} \lambda_{21} & \{-\lambda_{22},-\lambda_{23},0\} \\ \{0,0,-\lambda_{27}\} & 0 \end{pmatrix} & \begin{pmatrix} \lambda_{13} & \{-\lambda_{14},-\lambda_{15},0\} \\ \{0,0,-\lambda_{19}\} & 0 \end{pmatrix} & \begin{pmatrix} 0 & \{0,0,0\} \\ \{0,0,0\} & 0 \end{pmatrix} \end{pmatrix}$$

**U_{X₁[e₈]}[X₁[e₃]] == ZERO**

True

**defin[X₁[e₈], X₁[e₃]] // Expand**

$$\begin{pmatrix} \begin{pmatrix} 0 & \{0,0,0\} \\ \{0,0,0\} & 0 \end{pmatrix} & \begin{pmatrix} 0 & \{0,0,0\} \\ \{0,0,-\eta_3-\lambda_9\} & 0 \end{pmatrix} & \begin{pmatrix} 0 & \{0,0,0\} \\ \{0,0,0\} & 0 \end{pmatrix} \\ \begin{pmatrix} 0 & \{0,0,0\} \\ \{0,0,\eta_3+\lambda_9\} & 0 \end{pmatrix} & \begin{pmatrix} 0 & \{0,0,0\} \\ \{0,0,0\} & 0 \end{pmatrix} & \begin{pmatrix} 0 & \{0,0,0\} \\ \{0,0,0\} & 0 \end{pmatrix} \\ \begin{pmatrix} 0 & \{0,0,0\} \\ \{0,0,0\} & 0 \end{pmatrix} & \begin{pmatrix} 0 & \{0,0,0\} \\ \{0,0,0\} & 0 \end{pmatrix} & \begin{pmatrix} 0 & \{0,0,0\} \\ \{0,0,0\} & 0 \end{pmatrix} \end{pmatrix}$$

**Δ[X₁[e₈]] = Δ[X₁[e₈]] //. {λ₉ → -η₃}**

$$\begin{pmatrix} \begin{pmatrix} \alpha_5 & \{0,0,0\} \\ \{0,0,0\} & \alpha_5 \end{pmatrix} & \begin{pmatrix} \lambda_4 & \{-\psi_2,-\chi_1,0\} \\ \{-\eta_3,-\epsilon_3,-\phi_3\} & \lambda_5 \end{pmatrix} & \begin{pmatrix} 0 & \{\lambda_{22},\lambda_{23},0\} \\ \{0,0,\lambda_{27}\} & \lambda_{21} \end{pmatrix} \\ \begin{pmatrix} \lambda_5 & \{\psi_2,\chi_1,0\} \\ \{\eta_3,\epsilon_3,\phi_3\} & \lambda_4 \end{pmatrix} & \begin{pmatrix} -\alpha_5 & \{0,0,0\} \\ \{0,0,0\} & -\alpha_5 \end{pmatrix} & \begin{pmatrix} 0 & \{\lambda_{14},\lambda_{15},0\} \\ \{0,0,\lambda_{19}\} & \lambda_{13} \end{pmatrix} \\ \begin{pmatrix} \lambda_{21} & \{-\lambda_{22},-\lambda_{23},0\} \\ \{0,0,-\lambda_{27}\} & 0 \end{pmatrix} & \begin{pmatrix} \lambda_{13} & \{-\lambda_{14},-\lambda_{15},0\} \\ \{0,0,-\lambda_{19}\} & 0 \end{pmatrix} & \begin{pmatrix} 0 & \{0,0,0\} \\ \{0,0,0\} & 0 \end{pmatrix} \end{pmatrix}$$

**U_{X₁[e₈]}[X₁[e₂]] == ZERO**

True

**defin[X$_1$[e$_8$], X$_1$[e$_2$]] // Expand**

$$\begin{pmatrix} \begin{pmatrix} 0 & \{0,0,0\} \\ \{0,0,0\} & 0 \end{pmatrix} & \begin{pmatrix} 0 & \{0,0,0\} \\ \{0,0,\lambda_4-\rho_3\} & 0 \end{pmatrix} & \begin{pmatrix} 0 & \{0,0,0\} \\ \{0,0,0\} & 0 \end{pmatrix} \\ \begin{pmatrix} 0 & \{0,0,0\} \\ \{0,0,\rho_3-\lambda_4\} & 0 \end{pmatrix} & \begin{pmatrix} 0 & \{0,0,0\} \\ \{0,0,0\} & 0 \end{pmatrix} & \begin{pmatrix} 0 & \{0,0,0\} \\ \{0,0,0\} & 0 \end{pmatrix} \\ \begin{pmatrix} 0 & \{0,0,0\} \\ \{0,0,0\} & 0 \end{pmatrix} & \begin{pmatrix} 0 & \{0,0,0\} \\ \{0,0,0\} & 0 \end{pmatrix} & \begin{pmatrix} 0 & \{0,0,0\} \\ \{0,0,0\} & 0 \end{pmatrix} \end{pmatrix}$$

**Δ[X$_1$[e$_8$]] = Δ[X$_1$[e$_8$]] //. {λ$_4$ → ρ$_3$}**

$$\begin{pmatrix} \begin{pmatrix} \alpha_5 & \{0,0,0\} \\ \{0,0,0\} & \alpha_5 \end{pmatrix} & \begin{pmatrix} \rho_3 & \{-\psi_2,-\chi_1,0\} \\ \{-\eta_3,-\epsilon_3,-\phi_3\} & \lambda_5 \end{pmatrix} & \begin{pmatrix} 0 & \{\lambda_{22},\lambda_{23},0\} \\ \{0,0,\lambda_{27}\} & \lambda_{21} \end{pmatrix} \\ \begin{pmatrix} \lambda_5 & \{\psi_2,\chi_1,0\} \\ \{\eta_3,\epsilon_3,\phi_3\} & \rho_3 \end{pmatrix} & \begin{pmatrix} -\alpha_5 & \{0,0,0\} \\ \{0,0,0\} & -\alpha_5 \end{pmatrix} & \begin{pmatrix} 0 & \{\lambda_{14},\lambda_{15},0\} \\ \{0,0,\lambda_{19}\} & \lambda_{13} \end{pmatrix} \\ \begin{pmatrix} \lambda_{21} & \{-\lambda_{22},-\lambda_{23},0\} \\ \{0,0,-\lambda_{27}\} & 0 \end{pmatrix} & \begin{pmatrix} \lambda_{13} & \{-\lambda_{14},-\lambda_{15},0\} \\ \{0,0,-\lambda_{19}\} & 0 \end{pmatrix} & \begin{pmatrix} 0 & \{0,0,0\} \\ \{0,0,0\} & 0 \end{pmatrix} \end{pmatrix}$$

**U$_{X_1[e_8]}$[X$_1$[e$_1$]] == ZERO**

True

**defin[X$_1$[e$_8$], X$_1$[e$_1$]] // Expand**

$$\begin{pmatrix} \begin{pmatrix} 0 & \{0,0,0\} \\ \{0,0,0\} & 0 \end{pmatrix} & \begin{pmatrix} 0 & \{0,0,0\} \\ \{0,0,\lambda_5-\delta_4\} & 0 \end{pmatrix} & \begin{pmatrix} 0 & \{0,0,0\} \\ \{0,0,0\} & 0 \end{pmatrix} \\ \begin{pmatrix} 0 & \{0,0,0\} \\ \{0,0,\delta_4-\lambda_5\} & 0 \end{pmatrix} & \begin{pmatrix} 0 & \{0,0,0\} \\ \{0,0,0\} & 0 \end{pmatrix} & \begin{pmatrix} 0 & \{0,0,0\} \\ \{0,0,0\} & 0 \end{pmatrix} \\ \begin{pmatrix} 0 & \{0,0,0\} \\ \{0,0,0\} & 0 \end{pmatrix} & \begin{pmatrix} 0 & \{0,0,0\} \\ \{0,0,0\} & 0 \end{pmatrix} & \begin{pmatrix} 0 & \{0,0,0\} \\ \{0,0,0\} & 0 \end{pmatrix} \end{pmatrix}$$

**Δ[X$_1$[e$_8$]] = Δ[X$_1$[e$_8$]] //. {λ$_5$ → δ$_4$}**

$$\begin{pmatrix} \begin{pmatrix} \alpha_5 & \{0,0,0\} \\ \{0,0,0\} & \alpha_5 \end{pmatrix} & \begin{pmatrix} \rho_3 & \{-\psi_2,-\chi_1,0\} \\ \{-\eta_3,-\epsilon_3,-\phi_3\} & \delta_4 \end{pmatrix} & \begin{pmatrix} 0 & \{\lambda_{22},\lambda_{23},0\} \\ \{0,0,\lambda_{27}\} & \lambda_{21} \end{pmatrix} \\ \begin{pmatrix} \delta_4 & \{\psi_2,\chi_1,0\} \\ \{\eta_3,\epsilon_3,\phi_3\} & \rho_3 \end{pmatrix} & \begin{pmatrix} -\alpha_5 & \{0,0,0\} \\ \{0,0,0\} & -\alpha_5 \end{pmatrix} & \begin{pmatrix} 0 & \{\lambda_{14},\lambda_{15},0\} \\ \{0,0,\lambda_{19}\} & \lambda_{13} \end{pmatrix} \\ \begin{pmatrix} \lambda_{21} & \{-\lambda_{22},-\lambda_{23},0\} \\ \{0,0,-\lambda_{27}\} & 0 \end{pmatrix} & \begin{pmatrix} \lambda_{13} & \{-\lambda_{14},-\lambda_{15},0\} \\ \{0,0,-\lambda_{19}\} & 0 \end{pmatrix} & \begin{pmatrix} 0 & \{0,0,0\} \\ \{0,0,0\} & 0 \end{pmatrix} \end{pmatrix}$$

**T[E$_1$[1], E$_1$[1], X$_1$[e$_8$]] == X$_1$[e$_8$]**

True

```
Δ[X₁[e₈]] - Leib[E₁[1], E₁[1], X₁[e₈]] // Expand
```

$$\left(\begin{pmatrix} 0 & \{0,0,0\} \\ \{0,0,0\} & 0 \end{pmatrix} \quad \begin{pmatrix} 0 & \{0,0,0\} \\ \{0,0,0\} & 0 \end{pmatrix} \quad \begin{pmatrix} 0 & \{0,0,0\} \\ \{0,0,0\} & 0 \end{pmatrix}\right.$$
$$\begin{pmatrix} 0 & \{0,0,0\} \\ \{0,0,0\} & 0 \end{pmatrix} \quad \begin{pmatrix} 0 & \{0,0,0\} \\ \{0,0,0\} & 0 \end{pmatrix} \quad \begin{pmatrix} 0 & \{\beta_7+\lambda_{14}, \lambda_{15}-\beta_6, \\ \{0,0,\beta_1+\lambda_{19}\} & \beta_5+\lambda_{13} \end{pmatrix}$$
$$\left.\begin{pmatrix} 0 & \{0,0,0\} \\ \{0,0,0\} & 0 \end{pmatrix} \quad \begin{pmatrix} \beta_5+\lambda_{13} & \{-\beta_7-\lambda_{14},\beta_6-\lambda_{15},0\} \\ \{0,0,-\beta_1-\lambda_{19}\} & 0 \end{pmatrix} \quad \begin{pmatrix} 0 & \{0,0,0\} \\ \{0,0,0\} & 0 \end{pmatrix}\right)$$

```
Δ[X₁[e₈]] = Δ[X₁[e₈]] //. {λ₁₃ → -β₅, λ₁₄ → -β₇, λ₁₅ → β₆, λ₁₉ → -β₁}
```

$$\left(\begin{pmatrix} \alpha_5 & \{0,0,0\} \\ \{0,0,0\} & \alpha_5 \end{pmatrix} \quad \begin{pmatrix} \rho_3 & \{-\psi_2,-\chi_1,0\} \\ \{-\eta_3,-\epsilon_3,-\phi_3\} & \delta_4 \end{pmatrix} \quad \begin{pmatrix} 0 & \{\lambda_{22},\lambda_{23},0\} \\ \{0,0,\lambda_{27}\} & \lambda_{21} \end{pmatrix}\right.$$
$$\begin{pmatrix} \delta_4 & \{\psi_2,\chi_1,0\} \\ \{\eta_3,\epsilon_3,\phi_3\} & \rho_3 \end{pmatrix} \quad \begin{pmatrix} -\alpha_5 & \{0,0,0\} \\ \{0,0,0\} & -\alpha_5 \end{pmatrix} \quad \begin{pmatrix} 0 & \{-\beta_7,\beta_6,0\} \\ \{0,0,-\beta_1\} & -\beta_5 \end{pmatrix}$$
$$\left.\begin{pmatrix} \lambda_{21} & \{-\lambda_{22},-\lambda_{23},0\} \\ \{0,0,-\lambda_{27}\} & 0 \end{pmatrix} \quad \begin{pmatrix} -\beta_5 & \{\beta_7,-\beta_6,0\} \\ \{0,0,\beta_1\} & 0 \end{pmatrix} \quad \begin{pmatrix} 0 & \{0,0,0\} \\ \{0,0,0\} & 0 \end{pmatrix}\right)$$

```
T[E₂[1], E₂[1], X₁[e₈]] == X₁[e₈]
```

True

```
Δ[X₁[e₈]] - Leib[E₂[1], E₂[1], X₁[e₈]] // Expand
```

$$\left(\begin{pmatrix} 0 & \{0,0,0\} \\ \{0,0,0\} & 0 \end{pmatrix} \quad \begin{pmatrix} 0 & \{0,0,0\} \\ \{0,0,0\} & 0 \end{pmatrix} \quad \begin{pmatrix} 0 & \{\lambda_{22}-\gamma_7, \gamma_6+\lambda_{23}, 0\} \\ \{0,0,\lambda_{27}-\gamma_1\} & \lambda_{21}-\gamma_5 \end{pmatrix}\right.$$
$$\begin{pmatrix} 0 & \{0,0,0\} \\ \{0,0,0\} & 0 \end{pmatrix} \quad \begin{pmatrix} 0 & \{0,0,0\} \\ \{0,0,0\} & 0 \end{pmatrix} \quad \begin{pmatrix} 0 & \{0,0,0\} \\ \{0,0,0\} & 0 \end{pmatrix}$$
$$\left.\begin{pmatrix} \lambda_{21}-\gamma_5 & \{\gamma_7-\lambda_{22}, -\gamma_6-\lambda_{23}, 0\} \\ \{0,0,\gamma_1-\lambda_{27}\} & 0 \end{pmatrix} \quad \begin{pmatrix} 0 & \{0,0,0\} \\ \{0,0,0\} & 0 \end{pmatrix} \quad \begin{pmatrix} 0 & \{0,0,0\} \\ \{0,0,0\} & 0 \end{pmatrix}\right)$$

```
Δ[X₁[e₈]] = Δ[X₁[e₈]] //. {λ₂₁ → γ₅, λ₂₂ → γ₇, λ₂₃ → -γ₆, λ₂₇ → γ₁}
```

$$\left(\begin{pmatrix} \alpha_5 & \{0,0,0\} \\ \{0,0,0\} & \alpha_5 \end{pmatrix} \quad \begin{pmatrix} \rho_3 & \{-\psi_2,-\chi_1,0\} \\ \{-\eta_3,-\epsilon_3,-\phi_3\} & \delta_4 \end{pmatrix} \quad \begin{pmatrix} 0 & \{\gamma_7,-\gamma_6,0\} \\ \{0,0,\gamma_1\} & \gamma_5 \end{pmatrix}\right.$$
$$\begin{pmatrix} \delta_4 & \{\psi_2,\chi_1,0\} \\ \{\eta_3,\epsilon_3,\phi_3\} & \rho_3 \end{pmatrix} \quad \begin{pmatrix} -\alpha_5 & \{0,0,0\} \\ \{0,0,0\} & -\alpha_5 \end{pmatrix} \quad \begin{pmatrix} 0 & \{-\beta_7,\beta_6,0\} \\ \{0,0,-\beta_1\} & -\beta_5 \end{pmatrix}$$
$$\left.\begin{pmatrix} \gamma_5 & \{-\gamma_7,\gamma_6,0\} \\ \{0,0,-\gamma_1\} & 0 \end{pmatrix} \quad \begin{pmatrix} -\beta_5 & \{\beta_7,-\beta_6,0\} \\ \{0,0,\beta_1\} & 0 \end{pmatrix} \quad \begin{pmatrix} 0 & \{0,0,0\} \\ \{0,0,0\} & 0 \end{pmatrix}\right)$$

```
Variables[Δ[X₁[e₈]]]
```

$\{\alpha_5, \beta_1, \beta_5, \beta_6, \beta_7, \gamma_1, \gamma_5, \gamma_6, \gamma_7, \delta_4, \epsilon_3, \eta_3, \rho_3, \phi_3, \chi_1, \psi_2\}$

- $X_2[e_1]$

```
Δ[X₂[e₁]] = generic;
```

**U$_{X_2[e_1]}$[E$_1$[1]] == ZERO**

True

**(defin[X$_2$[e$_1$], E$_1$[1]] // Expand) == ZERO**

True

**U$_{E_1[1]}$[X$_2$[e$_1$]] == ZERO**

True

**(defin[E$_1$[1], X$_2$[e$_1$]] // Expand)**

$$\begin{pmatrix} \begin{pmatrix} \lambda_1 & \{0,0,0\} \\ \{0,0,0\} & \lambda_1 \end{pmatrix} & \begin{pmatrix} 0 & \{0,0,0\} \\ \{0,0,0\} & 0 \end{pmatrix} & \begin{pmatrix} 0 & \{0,0,0\} \\ \{0,0,0\} & 0 \end{pmatrix} \\ \begin{pmatrix} 0 & \{0,0,0\} \\ \{0,0,0\} & 0 \end{pmatrix} & \begin{pmatrix} 0 & \{0,0,0\} \\ \{0,0,0\} & 0 \end{pmatrix} & \begin{pmatrix} 0 & \{0,0,0\} \\ \{0,0,0\} & 0 \end{pmatrix} \\ \begin{pmatrix} 0 & \{0,0,0\} \\ \{0,0,0\} & 0 \end{pmatrix} & \begin{pmatrix} 0 & \{0,0,0\} \\ \{0,0,0\} & 0 \end{pmatrix} & \begin{pmatrix} 0 & \{0,0,0\} \\ \{0,0,0\} & 0 \end{pmatrix} \end{pmatrix}$$

**Δ[X$_2$[e$_1$]] = Δ[X$_2$[e$_1$]] //. {λ$_1$ → 0}**

$$\begin{pmatrix} \begin{pmatrix} 0 & \{0,0,0\} \\ \{0,0,0\} & 0 \end{pmatrix} & \begin{pmatrix} \lambda_4 & \{\lambda_6,\lambda_7,\lambda_8\} \\ \{\lambda_9,\lambda_{10},\lambda_{11}\} & \lambda_5 \end{pmatrix} & \begin{pmatrix} \lambda_{20} & \{\lambda_{22},\lambda_{23},\lambda_{24}\} \\ \{\lambda_{25},\lambda_{26},\lambda_{27}\} & \lambda_{21} \end{pmatrix} \\ \begin{pmatrix} \lambda_5 & \{-\lambda_6,-\lambda_7,-\lambda_8\} \\ \{-\lambda_9,-\lambda_{10},-\lambda_{11}\} & \lambda_4 \end{pmatrix} & \begin{pmatrix} \lambda_2 & \{0,0,0\} \\ \{0,0,0\} & \lambda_2 \end{pmatrix} & \begin{pmatrix} \lambda_{12} & \{\lambda_{14},\lambda_{15},\lambda_{16}\} \\ \{\lambda_{17},\lambda_{18},\lambda_{19}\} & \lambda_{13} \end{pmatrix} \\ \begin{pmatrix} \lambda_{21} & \{-\lambda_{22},-\lambda_{23},-\lambda_{24}\} \\ \{-\lambda_{25},-\lambda_{26},-\lambda_{27}\} & \lambda_{20} \end{pmatrix} & \begin{pmatrix} \lambda_{13} & \{-\lambda_{14},-\lambda_{15},-\lambda_{16}\} \\ \{-\lambda_{17},-\lambda_{18},-\lambda_{19}\} & \lambda_{12} \end{pmatrix} & \begin{pmatrix} \lambda_3 & \{0,0,0\} \\ \{0,0,0\} & \lambda_3 \end{pmatrix} \end{pmatrix}$$

**U$_{X_2[e_1]}$[E$_2$[1]] == ZERO**

True

**defin[X$_2$[e$_1$], E$_2$[1]] // Expand**

$$\begin{pmatrix} \begin{pmatrix} 0 & \{0,0,0\} \\ \{0,0,0\} & 0 \end{pmatrix} & \begin{pmatrix} 0 & \{0,0,0\} \\ \{0,0,0\} & 0 \end{pmatrix} & \begin{pmatrix} \lambda_4 & \{0,0,0\} \\ \{\lambda_9,\lambda_{10},\lambda_{11}\} & 0 \end{pmatrix} \\ \begin{pmatrix} 0 & \{0,0,0\} \\ \{0,0,0\} & 0 \end{pmatrix} & \begin{pmatrix} 0 & \{0,0,0\} \\ \{0,0,0\} & 0 \end{pmatrix} & \begin{pmatrix} \gamma_2+\lambda_2 & \{0,0,0\} \\ \{0,0,0\} & 0 \end{pmatrix} \\ \begin{pmatrix} 0 & \{0,0,0\} \\ \{-\lambda_9,-\lambda_{10},-\lambda_{11}\} & \lambda_4 \end{pmatrix} & \begin{pmatrix} 0 & \{0,0,0\} \\ \{0,0,0\} & \gamma_2+\lambda_2 \end{pmatrix} & \begin{pmatrix} \lambda_{13} & \{0,0,0\} \\ \{0,0,0\} & \lambda_{13} \end{pmatrix} \end{pmatrix}$$

**Δ[X₂[e₁]] =**
  **Δ[X₂[e₁]] //. {λ₂ → -γ₂, λ₄ → 0, λ₉ → 0, λ₁₀ → 0, λ₁₁ → 0, λ₁₃ → 0}**

$$\left(\begin{array}{ccc}
\begin{pmatrix} 0 & \{0,0,0\} \\ \{0,0,0\} & 0 \end{pmatrix} & \begin{pmatrix} 0 & \{\lambda_6,\lambda_7,\lambda_8\} \\ \{0,0,0\} & \lambda_5 \end{pmatrix} & \begin{pmatrix} \lambda_{20} & \{\lambda_{22},\lambda_{23},\lambda_{24}\} \\ \{\lambda_{25},\lambda_{26},\lambda_{27}\} & \lambda_{21} \end{pmatrix} \\
\begin{pmatrix} \lambda_5 & \{-\lambda_6,-\lambda_7,-\lambda_8\} \\ \{0,0,0\} & 0 \end{pmatrix} & \begin{pmatrix} -\gamma_2 & \{0,0,0\} \\ \{0,0,0\} & -\gamma_2 \end{pmatrix} & \begin{pmatrix} \lambda_{12} & \{\lambda_{14},\lambda_{15},\lambda_{16}\} \\ \{\lambda_{17},\lambda_{18},\lambda_{19}\} & 0 \end{pmatrix} \\
\begin{pmatrix} \lambda_{21} & \{-\lambda_{22},-\lambda_{23},-\lambda_{24}\} \\ \{-\lambda_{25},-\lambda_{26},-\lambda_{27}\} & \lambda_{20} \end{pmatrix} & \begin{pmatrix} 0 & \{-\lambda_{14},-\lambda_{15},-\lambda_{16}\} \\ \{-\lambda_{17},-\lambda_{18},-\lambda_{19}\} & \lambda_{12} \end{pmatrix} & \begin{pmatrix} \lambda_3 & \{0,0,0\} \\ \{0,0,0\} & \lambda_3 \end{pmatrix}
\end{array}\right)$$

**U_{E₂[1]}[X₂[e₁]] == ZERO**

True

**(defin[E₂[1], X₂[e₁]] // Expand) == ZERO**

True

**U_{X₂[e₁]}[E₃[1]] == ZERO**

True

**defin[X₂[e₁], E₃[1]] // Expand**

$$\left(\begin{array}{ccc}
\begin{pmatrix} 0 & \{0,0,0\} \\ \{0,0,0\} & 0 \end{pmatrix} & \begin{pmatrix} 0 & \{\lambda_{22},\lambda_{23},\lambda_{24}\} \\ \{0,0,0\} & \lambda_{21} \end{pmatrix} & \begin{pmatrix} 0 & \{0,0,0\} \\ \{0,0,0\} & 0 \end{pmatrix} \\
\begin{pmatrix} \lambda_{21} & \{-\lambda_{22},-\lambda_{23},-\lambda_{24}\} \\ \{0,0,0\} & 0 \end{pmatrix} & \begin{pmatrix} 0 & \{0,0,0\} \\ \{0,0,0\} & 0 \end{pmatrix} & \begin{pmatrix} \lambda_3-\gamma_2 & \{0,0,0\} \\ \{0,0,0\} & 0 \end{pmatrix} \\
\begin{pmatrix} 0 & \{0,0,0\} \\ \{0,0,0\} & 0 \end{pmatrix} & \begin{pmatrix} 0 & \{0,0,0\} \\ \{0,0,0\} & \lambda_3-\gamma_2 \end{pmatrix} & \begin{pmatrix} 0 & \{0,0,0\} \\ \{0,0,0\} & 0 \end{pmatrix}
\end{array}\right)$$

**Δ[X₂[e₁]] = Δ[X₂[e₁]] //. {λ₃ → γ₂, λ₂₁ → 0, λ₂₂ → 0, λ₂₃ → 0, λ₂₄ → 0}**

$$\left(\begin{array}{ccc}
\begin{pmatrix} 0 & \{0,0,0\} \\ \{0,0,0\} & 0 \end{pmatrix} & \begin{pmatrix} 0 & \{\lambda_6,\lambda_7,\lambda_8\} \\ \{0,0,0\} & \lambda_5 \end{pmatrix} & \begin{pmatrix} \lambda_{20} & \{0,0,0\} \\ \{\lambda_{25},\lambda_{26},\lambda_{27}\} & 0 \end{pmatrix} \\
\begin{pmatrix} \lambda_5 & \{-\lambda_6,-\lambda_7,-\lambda_8\} \\ \{0,0,0\} & 0 \end{pmatrix} & \begin{pmatrix} -\gamma_2 & \{0,0,0\} \\ \{0,0,0\} & -\gamma_2 \end{pmatrix} & \begin{pmatrix} \lambda_{12} & \{\lambda_{14},\lambda_{15},\lambda_{16}\} \\ \{\lambda_{17},\lambda_{18},\lambda_{19}\} & 0 \end{pmatrix} \\
\begin{pmatrix} 0 & \{0,0,0\} \\ \{-\lambda_{25},-\lambda_{26},-\lambda_{27}\} & \lambda_{20} \end{pmatrix} & \begin{pmatrix} 0 & \{-\lambda_{14},-\lambda_{15},-\lambda_{16}\} \\ \{-\lambda_{17},-\lambda_{18},-\lambda_{19}\} & \lambda_{12} \end{pmatrix} & \begin{pmatrix} \gamma_2 & \{0,0,0\} \\ \{0,0,0\} & \gamma_2 \end{pmatrix}
\end{array}\right)$$

**U_{E₃[1]}[X₂[e₁]] == ZERO**

True

**(defin[E₃[1], X₂[e₁]] // Expand) == ZERO**

True

**U_{X₂[e₁]}[X₁[e₈]] == ZERO**

True

**defin[X₂[e₁], X₁[e₈]] // Expand**

$$\left( \begin{pmatrix} 0 & \{0,0,0\} \\ \{0,0,0\} & 0 \end{pmatrix} \begin{pmatrix} 0 & \{0,0,0\} \\ \{0,0,0\} & 0 \end{pmatrix} \begin{pmatrix} 0 & \{0,0,0\} \\ \{0,0,0\} & 0 \end{pmatrix} \right.$$
$$\begin{pmatrix} 0 & \{0,0,0\} \\ \{0,0,0\} & 0 \end{pmatrix} \begin{pmatrix} 0 & \{0,0,0\} \\ \{0,0,0\} & 0 \end{pmatrix} \begin{pmatrix} -\beta_5 - \lambda_8 & \{0,0,0\} \\ \{0,0,0\} & 0 \end{pmatrix}$$
$$\left. \begin{pmatrix} 0 & \{0,0,0\} \\ \{0,0,0\} & 0 \end{pmatrix} \begin{pmatrix} 0 & \{0,0,0\} \\ \{0,0,0\} & -\beta_5 - \lambda_8 \end{pmatrix} \begin{pmatrix} 0 & \{0,0,0\} \\ \{0,0,0\} & 0 \end{pmatrix} \right)$$

**Δ[X₂[e₁]] = Δ[X₂[e₁]] //. {λ₈ → −β₅}**

$$\left( \begin{pmatrix} 0 & \{0,0,0\} \\ \{0,0,0\} & 0 \end{pmatrix} \begin{pmatrix} 0 & \{\lambda_6, \lambda_7, -\beta_5\} \\ \{0,0,0\} & \lambda_5 \end{pmatrix} \begin{pmatrix} \lambda_{20} & \{0,0,0\} \\ \{\lambda_{25}, \lambda_{26}, \lambda_{27}\} & 0 \end{pmatrix} \right.$$
$$\begin{pmatrix} \lambda_5 & \{-\lambda_6, -\lambda_7, \beta_5\} \\ \{0,0,0\} & 0 \end{pmatrix} \begin{pmatrix} -\gamma_2 & \{0,0,0\} \\ \{0,0,0\} & -\gamma_2 \end{pmatrix} \begin{pmatrix} \lambda_{12} & \{\lambda_{14}, \lambda_{15}, \lambda_{16}\} \\ \{\lambda_{17}, \lambda_{18}, \lambda_{19}\} & 0 \end{pmatrix}$$
$$\left. \begin{pmatrix} 0 & \{0,0,0\} \\ \{-\lambda_{25}, -\lambda_{26}, -\lambda_{27}\} & \lambda_{20} \end{pmatrix} \begin{pmatrix} 0 & \{-\lambda_{14}, -\lambda_{15}, -\lambda_{16}\} \\ \{-\lambda_{17}, -\lambda_{18}, -\lambda_{19}\} & \lambda_{12} \end{pmatrix} \begin{pmatrix} \gamma_2 & \{0,0,0\} \\ \{0,0,0\} & \gamma_2 \end{pmatrix} \right)$$

**U_{X₁[e₈]}[X₂[e₁]] == ZERO**

True

**(defin[X₁[e₈], X₂[e₁]] // Expand) == ZERO**

True

**U_{X₂[e₁]}[X₁[e₇]] == ZERO**

True

**defin[X₂[e₁], X₁[e₇]] // Expand**

$$\left( \begin{pmatrix} 0 & \{0,0,0\} \\ \{0,0,0\} & 0 \end{pmatrix} \begin{pmatrix} 0 & \{0,0,0\} \\ \{0,0,0\} & 0 \end{pmatrix} \begin{pmatrix} 0 & \{0,0,0\} \\ \{0,0,0\} & 0 \end{pmatrix} \right.$$
$$\begin{pmatrix} 0 & \{0,0,0\} \\ \{0,0,0\} & 0 \end{pmatrix} \begin{pmatrix} 0 & \{0,0,0\} \\ \{0,0,0\} & 0 \end{pmatrix} \begin{pmatrix} -\beta_4 - \lambda_7 & \{0,0,0\} \\ \{0,0,0\} & 0 \end{pmatrix}$$
$$\left. \begin{pmatrix} 0 & \{0,0,0\} \\ \{0,0,0\} & 0 \end{pmatrix} \begin{pmatrix} 0 & \{0,0,0\} \\ \{0,0,0\} & -\beta_4 - \lambda_7 \end{pmatrix} \begin{pmatrix} 0 & \{0,0,0\} \\ \{0,0,0\} & 0 \end{pmatrix} \right)$$

**Δ[X₂[e₁]] = Δ[X₂[e₁]] //. {λ₇ → −β₄}**

$$\left( \begin{pmatrix} 0 & \{0,0,0\} \\ \{0,0,0\} & 0 \end{pmatrix} \begin{pmatrix} 0 & \{\lambda_6, -\beta_4, -\beta_5\} \\ \{0,0,0\} & \lambda_5 \end{pmatrix} \begin{pmatrix} \lambda_{20} & \{0,0,0\} \\ \{\lambda_{25}, \lambda_{26}, \lambda_{27}\} & 0 \end{pmatrix} \right.$$
$$\begin{pmatrix} \lambda_5 & \{-\lambda_6, \beta_4, \beta_5\} \\ \{0,0,0\} & 0 \end{pmatrix} \begin{pmatrix} -\gamma_2 & \{0,0,0\} \\ \{0,0,0\} & -\gamma_2 \end{pmatrix} \begin{pmatrix} \lambda_{12} & \{\lambda_{14}, \lambda_{15}, \lambda_{16}\} \\ \{\lambda_{17}, \lambda_{18}, \lambda_{19}\} & 0 \end{pmatrix}$$
$$\left. \begin{pmatrix} 0 & \{0,0,0\} \\ \{-\lambda_{25}, -\lambda_{26}, -\lambda_{27}\} & \lambda_{20} \end{pmatrix} \begin{pmatrix} 0 & \{-\lambda_{14}, -\lambda_{15}, -\lambda_{16}\} \\ \{-\lambda_{17}, -\lambda_{18}, -\lambda_{19}\} & \lambda_{12} \end{pmatrix} \begin{pmatrix} \gamma_2 & \{0,0,0\} \\ \{0,0,0\} & \gamma_2 \end{pmatrix} \right)$$

**U$_{X_2[e_1]}$ [X$_1$[e$_6$]] == ZERO**

True

**defin[X$_2$[e$_1$], X$_1$[e$_6$]] // Expand**

$$\begin{pmatrix} \begin{pmatrix} 0 & \{0,0,0\} \\ \{0,0,0\} & 0 \end{pmatrix} & \begin{pmatrix} 0 & \{0,0,0\} \\ \{0,0,0\} & 0 \end{pmatrix} & \begin{pmatrix} 0 & \{0,0,0\} \\ \{0,0,0\} & 0 \end{pmatrix} \\ \begin{pmatrix} 0 & \{0,0,0\} \\ \{0,0,0\} & 0 \end{pmatrix} & \begin{pmatrix} 0 & \{0,0,0\} \\ \{0,0,0\} & 0 \end{pmatrix} & \begin{pmatrix} -\beta_3-\lambda_6 & \{0,0,0\} \\ \{0,0,0\} & 0 \end{pmatrix} \\ \begin{pmatrix} 0 & \{0,0,0\} \\ \{0,0,0\} & 0 \end{pmatrix} & \begin{pmatrix} 0 & \{0,0,0\} \\ \{0,0,0\} & -\beta_3-\lambda_6 \end{pmatrix} & \begin{pmatrix} 0 & \{0,0,0\} \\ \{0,0,0\} & 0 \end{pmatrix} \end{pmatrix}$$

**Δ[X$_2$[e$_1$]] = Δ[X$_2$[e$_1$]] //. {λ$_6$ → -β$_3$}**

$$\begin{pmatrix} \begin{pmatrix} 0 & \{0,0,0\} \\ \{0,0,0\} & 0 \end{pmatrix} & \begin{pmatrix} 0 & \{-\beta_3,-\beta_4,-\beta_5\} \\ \{0,0,0\} & \lambda_5 \end{pmatrix} & \begin{pmatrix} \lambda_{20} & \{0,0,0\} \\ \{\lambda_{25},\lambda_{26},\lambda_{27}\} & 0 \end{pmatrix} \\ \begin{pmatrix} \lambda_5 & \{\beta_3,\beta_4,\beta_5\} \\ \{0,0,0\} & 0 \end{pmatrix} & \begin{pmatrix} -\gamma_2 & \{0,0,0\} \\ \{0,0,0\} & -\gamma_2 \end{pmatrix} & \begin{pmatrix} \lambda_{12} & \{\lambda_{14},\lambda_{15},\lambda_{16}\} \\ \{\lambda_{17},\lambda_{18},\lambda_{19}\} & 0 \end{pmatrix} \\ \begin{pmatrix} 0 & \{0,0,0\} \\ \{-\lambda_{25},-\lambda_{26},-\lambda_{27}\} & \lambda_{20} \end{pmatrix} & \begin{pmatrix} 0 & \{-\lambda_{14},-\lambda_{15},-\lambda_{16}\} \\ \{-\lambda_{17},-\lambda_{18},-\lambda_{19}\} & \lambda_{12} \end{pmatrix} & \begin{pmatrix} \gamma_2 & \{0,0,0\} \\ \{0,0,0\} & \gamma_2 \end{pmatrix} \end{pmatrix}$$

**U$_{X_1[e_6]}$ [X$_2$[e$_1$]] == ZERO**

True

**(defin[X$_1$[e$_6$], X$_2$[e$_1$]] // Expand) == ZERO**

True

**U$_{X_2[e_1]}$ [X$_1$[e$_5$]] == ZERO**

True

**(defin[X$_2$[e$_1$], X$_1$[e$_5$]] // Expand) == ZERO**

True

**U$_{X_1[e_5]}$ [X$_2$[e$_1$]] == ZERO**

True

**(defin[X$_1$[e$_5$], X$_2$[e$_1$]] // Expand) == ZERO**

True

**U$_{X_2[e_1]}$ [X$_1$[e$_4$]] == ZERO**

True

```
(defin[X₂[e₁], X₁[e₄]] // Expand) == ZERO
```
True

```
U_{X₁[e₄]}[X₂[e₁]] == ZERO
```
True

```
(defin[X₁[e₄], X₂[e₁]] // Expand) == ZERO
```
True

```
U_{X₂[e₁]}[X₁[e₃]] == ZERO
```
True

```
(defin[X₂[e₁], X₁[e₃]] // Expand) == ZERO
```
True

```
U_{X₁[e₃]}[X₂[e₁]] == ZERO
```
True

```
(defin[X₁[e₃], X₂[e₁]] // Expand) == ZERO
```
True

```
U_{X₂[e₁]}[X₁[e₂]] == ZERO
```
True

```
(defin[X₂[e₁], X₁[e₂]] // Expand) == ZERO
```
True

```
U_{X₁[e₂]}[X₂[e₁]] == ZERO
```
True

```
(defin[X₁[e₂], X₂[e₁]] // Expand) == ZERO
```
True

```
U_{X₂[e₁]}[X₁[e₁]] == ZERO
```
True

```
defin[X₂[e₁], X₁[e₁]] // Expand
```

$$\begin{pmatrix}
\begin{pmatrix} 0 & \{0,0,0\} \\ \{0,0,0\} & 0 \end{pmatrix} & \begin{pmatrix} 0 & \{0,0,0\} \\ \{0,0,0\} & 0 \end{pmatrix} & \begin{pmatrix} 0 & \{0,0,0\} \\ \{0,0,0\} & 0 \end{pmatrix} \\
\begin{pmatrix} 0 & \{0,0,0\} \\ \{0,0,0\} & 0 \end{pmatrix} & \begin{pmatrix} 0 & \{0,0,0\} \\ \{0,0,0\} & 0 \end{pmatrix} & \begin{pmatrix} \beta_2+\lambda_5 & \{0,0,0\} \\ \{0,0,0\} & 0 \end{pmatrix} \\
\begin{pmatrix} 0 & \{0,0,0\} \\ \{0,0,0\} & 0 \end{pmatrix} & \begin{pmatrix} 0 & \{0,0,0\} \\ \{0,0,0\} & \beta_2+\lambda_5 \end{pmatrix} & \begin{pmatrix} 0 & \{0,0,0\} \\ \{0,0,0\} & 0 \end{pmatrix}
\end{pmatrix}$$

```
Δ[X₂[e₁]] = Δ[X₂[e₁]] //. {λ₅ → -β₂}
```

$$\begin{pmatrix}
\begin{pmatrix} 0 & \{0,0,0\} \\ \{0,0,0\} & 0 \end{pmatrix} & \begin{pmatrix} 0 & \{-\beta_3,-\beta_4,-\beta_5\} \\ \{0,0,0\} & -\beta_2 \end{pmatrix} & \begin{pmatrix} \lambda_{20} & \{0,0,0\} \\ \{\lambda_{25},\lambda_{26},\lambda_{27}\} & 0 \end{pmatrix} \\
\begin{pmatrix} -\beta_2 & \{\beta_3,\beta_4,\beta_5\} \\ \{0,0,0\} & 0 \end{pmatrix} & \begin{pmatrix} -\gamma_2 & \{0,0,0\} \\ \{0,0,0\} & -\gamma_2 \end{pmatrix} & \begin{pmatrix} \lambda_{12} & \{\lambda_{14},\lambda_{15},\lambda_{16}\} \\ \{\lambda_{17},\lambda_{18},\lambda_{19}\} & 0 \end{pmatrix} \\
\begin{pmatrix} 0 & \{0,0,0\} \\ \{-\lambda_{25},-\lambda_{26},-\lambda_{27}\} & \lambda_{20} \end{pmatrix} & \begin{pmatrix} 0 & \{-\lambda_{14},-\lambda_{15},-\lambda_{16}\} \\ \{-\lambda_{17},-\lambda_{18},-\lambda_{19}\} & \lambda_{12} \end{pmatrix} & \begin{pmatrix} \gamma_2 & \{0,0,0\} \\ \{0,0,0\} & \gamma_2 \end{pmatrix}
\end{pmatrix}$$

```
U_{X₁[e₁]}[X₂[e₁]] == ZERO
```

True

```
(defin[X₁[e₁], X₂[e₁]] // Expand) == ZERO
```

True

```
T[E₁[1], E₁[1], X₂[e₁]] == ZERO
```

True

```
Leib[E₁[1], E₁[1], X₂[e₁]] // Expand
```

$$\begin{pmatrix}
\begin{pmatrix} 0 & \{0,0,0\} \\ \{0,0,0\} & 0 \end{pmatrix} & \begin{pmatrix} 0 & \{0,0,0\} \\ \{0,0,0\} & 0 \end{pmatrix} & \begin{pmatrix} \alpha_1+\lambda_{20} & \{0,0,0\} \\ \{\alpha_6+\lambda_{25},\alpha_7+\lambda_{26},\alpha_8+\lambda_{27}\} & 0 \end{pmatrix} \\
\begin{pmatrix} 0 & \{0,0,0\} \\ \{0,0,0\} & 0 \end{pmatrix} & \begin{pmatrix} 0 & \{0,0,0\} \\ \{0,0,0\} & 0 \end{pmatrix} & \begin{pmatrix} 0 & \{0,0,0\} \\ \{0,0,0\} & 0 \end{pmatrix} \\
\begin{pmatrix} 0 & \{0,0,0\} \\ \{-\alpha_6-\lambda_{25},-\alpha_7-\lambda_{26},-\alpha_8-\lambda_{27}\} & \alpha_1+\lambda_{20} \end{pmatrix} & \begin{pmatrix} 0 & \{0,0,0\} \\ \{0,0,0\} & 0 \end{pmatrix} & \begin{pmatrix} 0 & \{0,0,0\} \\ \{0,0,0\} & 0 \end{pmatrix}
\end{pmatrix}$$

```
Δ[X₂[e₁]] = Δ[X₂[e₁]] //. {λ₂₀ → -α₁, λ₂₅ → -α₆, λ₂₆ → -α₇, λ₂₇ → -α₈}
```

$$\begin{pmatrix}
\begin{pmatrix} 0 & \{0,0,0\} \\ \{0,0,0\} & 0 \end{pmatrix} & \begin{pmatrix} 0 & \{-\beta_3,-\beta_4,-\beta_5\} \\ \{0,0,0\} & -\beta_2 \end{pmatrix} & \begin{pmatrix} -\alpha_1 & \{0,0,0\} \\ \{-\alpha_6,-\alpha_7,-\alpha_8\} & 0 \end{pmatrix} \\
\begin{pmatrix} -\beta_2 & \{\beta_3,\beta_4,\beta_5\} \\ \{0,0,0\} & 0 \end{pmatrix} & \begin{pmatrix} -\gamma_2 & \{0,0,0\} \\ \{0,0,0\} & -\gamma_2 \end{pmatrix} & \begin{pmatrix} \lambda_{12} & \{\lambda_{14},\lambda_{15},\lambda_{16}\} \\ \{\lambda_{17},\lambda_{18},\lambda_{19}\} & 0 \end{pmatrix} \\
\begin{pmatrix} 0 & \{0,0,0\} \\ \{\alpha_6,\alpha_7,\alpha_8\} & -\alpha_1 \end{pmatrix} & \begin{pmatrix} 0 & \{-\lambda_{14},-\lambda_{15},-\lambda_{16}\} \\ \{-\lambda_{17},-\lambda_{18},-\lambda_{19}\} & \lambda_{12} \end{pmatrix} & \begin{pmatrix} \gamma_2 & \{0,0,0\} \\ \{0,0,0\} & \gamma_2 \end{pmatrix}
\end{pmatrix}$$

```
Δ[ONE] = ZERO;
```

```
T[X₂[e₁], E₂[1], E₂[1]] == X₂[e₁]
```

True

```
(Δ[X₂[e₁]] - Leib[X₂[e₁], E₂[1], E₂[1]] // Expand) == ZERO
```

True

```
T[X₁[e₁], X₁[e₂], X₂[e₁]] == ZERO
```

True

```
Leib[X₁[e₁], X₁[e₂], X₂[e₁]] // Expand
```

$$\begin{pmatrix} \begin{pmatrix} 0 & \{0,0,0\} \\ \{0,0,0\} & 0 \end{pmatrix} & \begin{pmatrix} 0 & \{0,0,0\} \\ \{0,0,0\} & 0 \end{pmatrix} & \begin{pmatrix} 0 & \{0,0,0\} \\ \{0,0,0\} & 0 \end{pmatrix} \\ \begin{pmatrix} 0 & \{0,0,0\} \\ \{0,0,0\} & 0 \end{pmatrix} & \begin{pmatrix} 0 & \{0,0,0\} \\ \{0,0,0\} & 0 \end{pmatrix} & \begin{pmatrix} 0 & \{0,0,0\} \\ \{\lambda_{17}+\rho_4, \lambda_{18}+\rho_5, \lambda_{19}+\rho_6\} & 0 \end{pmatrix} \\ \begin{pmatrix} 0 & \{0,0,0\} \\ \{0,0,0\} & 0 \end{pmatrix} & \begin{pmatrix} 0 & \{0,0,0\} \\ \{-\lambda_{17}-\rho_4, -\lambda_{18}-\rho_5, -\lambda_{19}-\rho_6\} & 0 \end{pmatrix} & \begin{pmatrix} 0 & \{0,0,0\} \\ \{0,0,0\} & 0 \end{pmatrix} \end{pmatrix}$$

```
Δ[X₂[e₁]] = Δ[X₂[e₁]] //. {λ₁₇ → -ρ₄, λ₁₈ → -ρ₅, λ₁₉ → -ρ₆}
```

$$\begin{pmatrix} \begin{pmatrix} 0 & \{0,0,0\} \\ \{0,0,0\} & 0 \end{pmatrix} & \begin{pmatrix} 0 & \{-\beta_3,-\beta_4,-\beta_5\} \\ \{0,0,0\} & -\beta_2 \end{pmatrix} & \begin{pmatrix} -\alpha_1 & \{0,0,0\} \\ \{-\alpha_6,-\alpha_7,-\alpha_8\} & 0 \end{pmatrix} \\ \begin{pmatrix} -\beta_2 & \{\beta_3,\beta_4,\beta_5\} \\ \{0,0,0\} & 0 \end{pmatrix} & \begin{pmatrix} -\gamma_2 & \{0,0,0\} \\ \{0,0,0\} & -\gamma_2 \end{pmatrix} & \begin{pmatrix} \lambda_{12} & \{\lambda_{14},\lambda_{15},\lambda_{16}\} \\ \{-\rho_4,-\rho_5,-\rho_6\} & 0 \end{pmatrix} \\ \begin{pmatrix} 0 & \{0,0,0\} \\ \{\alpha_6,\alpha_7,\alpha_8\} & -\alpha_1 \end{pmatrix} & \begin{pmatrix} 0 & \{-\lambda_{14},-\lambda_{15},-\lambda_{16}\} \\ \{\rho_4,\rho_5,\rho_6\} & \lambda_{12} \end{pmatrix} & \begin{pmatrix} \gamma_2 & \{0,0,0\} \\ \{0,0,0\} & \gamma_2 \end{pmatrix} \end{pmatrix}$$

```
T[X₁[e₁], X₁[e₃], X₂[e₁]] == ZERO
```

True

```
Leib[X₁[e₁], X₁[e₃], X₂[e₁]] // Expand
```

$$\begin{pmatrix} \begin{pmatrix} 0 & \{0,0,0\} \\ \{0,0,0\} & 0 \end{pmatrix} & \begin{pmatrix} 0 & \{0,0,0\} \\ \{0,0,0\} & 0 \end{pmatrix} & \begin{pmatrix} 0 & \{0,0,0\} \\ \{0,0,0\} & 0 \end{pmatrix} \\ \begin{pmatrix} 0 & \{0,0,0\} \\ \{0,0,0\} & 0 \end{pmatrix} & \begin{pmatrix} 0 & \{0,0,0\} \\ \{0,0,0\} & 0 \end{pmatrix} & \begin{pmatrix} 0 & \{0,0,0\} \\ \{0, \eta_4-\lambda_{16}, \eta_5+\lambda_{15}\} & 0 \end{pmatrix} \\ \begin{pmatrix} 0 & \{0,0,0\} \\ \{0,0,0\} & 0 \end{pmatrix} & \begin{pmatrix} 0 & \{0,0,0\} \\ \{0, \lambda_{16}-\eta_4, -\eta_5-\lambda_{15}\} & 0 \end{pmatrix} & \begin{pmatrix} 0 & \{0,0,0\} \\ \{0,0,0\} & 0 \end{pmatrix} \end{pmatrix}$$

```
Δ[X₂[e₁]] = Δ[X₂[e₁]] //. {λ₁₅ → -η₅, λ₁₆ → η₄}
```

$$\begin{pmatrix} \begin{pmatrix} 0 & \{0,0,0\} \\ \{0,0,0\} & 0 \end{pmatrix} & \begin{pmatrix} 0 & \{-\beta_3,-\beta_4,-\beta_5\} \\ \{0,0,0\} & -\beta_2 \end{pmatrix} & \begin{pmatrix} -\alpha_1 & \{0,0,0\} \\ \{-\alpha_6,-\alpha_7,-\alpha_8\} & 0 \end{pmatrix} \\ \begin{pmatrix} -\beta_2 & \{\beta_3,\beta_4,\beta_5\} \\ \{0,0,0\} & 0 \end{pmatrix} & \begin{pmatrix} -\gamma_2 & \{0,0,0\} \\ \{0,0,0\} & -\gamma_2 \end{pmatrix} & \begin{pmatrix} \lambda_{12} & \{\lambda_{14},-\eta_5,\eta_4\} \\ \{-\rho_4,-\rho_5,-\rho_6\} & 0 \end{pmatrix} \\ \begin{pmatrix} 0 & \{0,0,0\} \\ \{\alpha_6,\alpha_7,\alpha_8\} & -\alpha_1 \end{pmatrix} & \begin{pmatrix} 0 & \{-\lambda_{14},\eta_5,-\eta_4\} \\ \{\rho_4,\rho_5,\rho_6\} & \lambda_{12} \end{pmatrix} & \begin{pmatrix} \gamma_2 & \{0,0,0\} \\ \{0,0,0\} & \gamma_2 \end{pmatrix} \end{pmatrix}$$

**T[X₁[e₁], X₁[e₄], X₂[e₁]] == ZERO**

True

**Leib[X₁[e₁], X₁[e₄], X₂[e₁]] // Expand**

$$\left( \begin{array}{ccc} \begin{pmatrix} 0 & \{0,0,0\} \\ \{0,0,0\} & 0 \end{pmatrix} & \begin{pmatrix} 0 & \{0,0,0\} \\ \{0,0,0\} & 0 \end{pmatrix} & \begin{pmatrix} 0 & \{0,0,0\} \\ \{0,0,0\} & 0 \end{pmatrix} \\ \begin{pmatrix} 0 & \{0,0,0\} \\ \{0,0,0\} & 0 \end{pmatrix} & \begin{pmatrix} 0 & \{0,0,0\} \\ \{0,0,0\} & 0 \end{pmatrix} & \begin{pmatrix} 0 & \{0,0,0\} \\ \{0,0,\epsilon_4-\lambda_{14}\} & 0 \end{pmatrix} \\ \begin{pmatrix} 0 & \{0,0,0\} \\ \{0,0,0\} & 0 \end{pmatrix} & \begin{pmatrix} 0 & \{0,0,0\} \\ \{0,0,\lambda_{14}-\epsilon_4\} & 0 \end{pmatrix} & \begin{pmatrix} 0 & \{0,0,0\} \\ \{0,0,0\} & 0 \end{pmatrix} \end{array} \right)$$

**Δ[X₂[e₁]] = Δ[X₂[e₁]] //. {λ₁₄ → ε₄}**

$$\left( \begin{array}{ccc} \begin{pmatrix} 0 & \{0,0,0\} \\ \{0,0,0\} & 0 \end{pmatrix} & \begin{pmatrix} 0 & \{-\beta_3,-\beta_4,-\beta_5\} \\ \{0,0,0\} & -\beta_2 \end{pmatrix} & \begin{pmatrix} -\alpha_1 & \{0,0,0\} \\ \{-\alpha_6,-\alpha_7,-\alpha_8\} & 0 \end{pmatrix} \\ \begin{pmatrix} -\beta_2 & \{\beta_3,\beta_4,\beta_5\} \\ \{0,0,0\} & 0 \end{pmatrix} & \begin{pmatrix} -\gamma_2 & \{0,0,0\} \\ \{0,0,0\} & -\gamma_2 \end{pmatrix} & \begin{pmatrix} \lambda_{12} & \{\epsilon_4,-\eta_5,\eta_4\} \\ \{-\rho_4,-\rho_5,-\rho_6\} & 0 \end{pmatrix} \\ \begin{pmatrix} 0 & \{0,0,0\} \\ \{\alpha_6,\alpha_7,\alpha_8\} & -\alpha_1 \end{pmatrix} & \begin{pmatrix} 0 & \{-\epsilon_4,\eta_5,-\eta_4\} \\ \{\rho_4,\rho_5,\rho_6\} & \lambda_{12} \end{pmatrix} & \begin{pmatrix} \gamma_2 & \{0,0,0\} \\ \{0,0,0\} & \gamma_2 \end{pmatrix} \end{array} \right)$$

**Variables[Δ[X₂[e₁]]]**

$\{\alpha_1, \alpha_6, \alpha_7, \alpha_8, \beta_2, \beta_3, \beta_4, \beta_5, \gamma_2, \epsilon_4, \eta_4, \eta_5, \lambda_{12}, \rho_4, \rho_5, \rho_6\}$

We will determine the parameter $\xi$ in a moment

**Δ[X₂[e₁]] = Δ[X₂[e₁]] //. {λ₁₂ → ξ}**

$$\left( \begin{array}{ccc} \begin{pmatrix} 0 & \{0,0,0\} \\ \{0,0,0\} & 0 \end{pmatrix} & \begin{pmatrix} 0 & \{-\beta_3,-\beta_4,-\beta_5\} \\ \{0,0,0\} & -\beta_2 \end{pmatrix} & \begin{pmatrix} -\alpha_1 & \{0,0,0\} \\ \{-\alpha_6,-\alpha_7,-\alpha_8\} & 0 \end{pmatrix} \\ \begin{pmatrix} -\beta_2 & \{\beta_3,\beta_4,\beta_5\} \\ \{0,0,0\} & 0 \end{pmatrix} & \begin{pmatrix} -\gamma_2 & \{0,0,0\} \\ \{0,0,0\} & -\gamma_2 \end{pmatrix} & \begin{pmatrix} \xi & \{\epsilon_4,-\eta_5,\eta_4\} \\ \{-\rho_4,-\rho_5,-\rho_6\} & 0 \end{pmatrix} \\ \begin{pmatrix} 0 & \{0,0,0\} \\ \{\alpha_6,\alpha_7,\alpha_8\} & -\alpha_1 \end{pmatrix} & \begin{pmatrix} 0 & \{-\epsilon_4,\eta_5,-\eta_4\} \\ \{\rho_4,\rho_5,\rho_6\} & \xi \end{pmatrix} & \begin{pmatrix} \gamma_2 & \{0,0,0\} \\ \{0,0,0\} & \gamma_2 \end{pmatrix} \end{array} \right)$$

- **X₂[e₂]**

  **Δ[X₂[e₂]] = generic;**

  **U_{E₁[1]}[X₂[e₂]] == ZERO**

  True

**(defin[E$_1$[1], X$_2$[e$_2$]] // Expand)**

$$\begin{pmatrix} \begin{pmatrix} \lambda_1 & \{0,0,0\} \\ \{0,0,0\} & \lambda_1 \end{pmatrix} & \begin{pmatrix} 0 & \{0,0,0\} \\ \{0,0,0\} & 0 \end{pmatrix} & \begin{pmatrix} 0 & \{0,0,0\} \\ \{0,0,0\} & 0 \end{pmatrix} \\ \begin{pmatrix} 0 & \{0,0,0\} \\ \{0,0,0\} & 0 \end{pmatrix} & \begin{pmatrix} 0 & \{0,0,0\} \\ \{0,0,0\} & 0 \end{pmatrix} & \begin{pmatrix} 0 & \{0,0,0\} \\ \{0,0,0\} & 0 \end{pmatrix} \\ \begin{pmatrix} 0 & \{0,0,0\} \\ \{0,0,0\} & 0 \end{pmatrix} & \begin{pmatrix} 0 & \{0,0,0\} \\ \{0,0,0\} & 0 \end{pmatrix} & \begin{pmatrix} 0 & \{0,0,0\} \\ \{0,0,0\} & 0 \end{pmatrix} \end{pmatrix}$$

**Δ[X$_2$[e$_2$]] = Δ[X$_2$[e$_2$]] //. {λ$_1$ → 0}**

$$\begin{pmatrix} \begin{pmatrix} 0 & \{0,0,0\} \\ \{0,0,0\} & 0 \end{pmatrix} & \begin{pmatrix} \lambda_4 & \{\lambda_6,\lambda_7,\lambda_8\} \\ \{\lambda_9,\lambda_{10},\lambda_{11}\} & \lambda_5 \end{pmatrix} & \begin{pmatrix} \lambda_{20} & \{\lambda_{22},\lambda_{23},\lambda_{24}\} \\ \{\lambda_{25},\lambda_{26},\lambda_{27}\} & \lambda_{21} \end{pmatrix} \\ \begin{pmatrix} \lambda_5 & \{-\lambda_6,-\lambda_7,-\lambda_8\} \\ \{-\lambda_9,-\lambda_{10},-\lambda_{11}\} & \lambda_4 \end{pmatrix} & \begin{pmatrix} \lambda_2 & \{0,0,0\} \\ \{0,0,0\} & \lambda_2 \end{pmatrix} & \begin{pmatrix} \lambda_{12} & \{\lambda_{14},\lambda_{15},\lambda_{16}\} \\ \{\lambda_{17},\lambda_{18},\lambda_{19}\} & \lambda_{13} \end{pmatrix} \\ \begin{pmatrix} \lambda_{21} & \{-\lambda_{22},-\lambda_{23},-\lambda_{24}\} \\ \{-\lambda_{25},-\lambda_{26},-\lambda_{27}\} & \lambda_{20} \end{pmatrix} & \begin{pmatrix} \lambda_{13} & \{-\lambda_{14},-\lambda_{15},-\lambda_{16}\} \\ \{-\lambda_{17},-\lambda_{18},-\lambda_{19}\} & \lambda_{12} \end{pmatrix} & \begin{pmatrix} \lambda_3 & \{0,0,0\} \\ \{0,0,0\} & \lambda_3 \end{pmatrix} \end{pmatrix}$$

**U$_{X_2[e_2]}$[E$_2$[1]] == ZERO**

True

**defin[X$_2$[e$_2$], E$_2$[1]] // Expand**

$$\begin{pmatrix} \begin{pmatrix} 0 & \{0,0,0\} \\ \{0,0,0\} & 0 \end{pmatrix} & \begin{pmatrix} 0 & \{0,0,0\} \\ \{0,0,0\} & 0 \end{pmatrix} & \begin{pmatrix} 0 & \{\lambda_6,\lambda_7,\lambda_8\} \\ \{0,0,0\} & \lambda_5 \end{pmatrix} \\ \begin{pmatrix} 0 & \{0,0,0\} \\ \{0,0,0\} & 0 \end{pmatrix} & \begin{pmatrix} 0 & \{0,0,0\} \\ \{0,0,0\} & 0 \end{pmatrix} & \begin{pmatrix} 0 & \{0,0,0\} \\ \{0,0,0\} & \gamma_1+\lambda_2 \end{pmatrix} \\ \begin{pmatrix} \lambda_5 & \{-\lambda_6,-\lambda_7,-\lambda_8\} \\ \{0,0,0\} & 0 \end{pmatrix} & \begin{pmatrix} \gamma_1+\lambda_2 & \{0,0,0\} \\ \{0,0,0\} & 0 \end{pmatrix} & \begin{pmatrix} \lambda_{12} & \{0,0,0\} \\ \{0,0,0\} & \lambda_{12} \end{pmatrix} \end{pmatrix}$$

**Δ[X$_2$[e$_2$]] =
Δ[X$_2$[e$_2$]] //. {λ$_2$ → −γ$_1$, λ$_5$ → 0, λ$_6$ → 0, λ$_7$ → 0, λ$_8$ → 0, λ$_{12}$ → 0}**

$$\begin{pmatrix} \begin{pmatrix} 0 & \{0,0,0\} \\ \{0,0,0\} & 0 \end{pmatrix} & \begin{pmatrix} \lambda_4 & \{0,0,0\} \\ \{\lambda_9,\lambda_{10},\lambda_{11}\} & 0 \end{pmatrix} & \begin{pmatrix} \lambda_{20} & \{\lambda_{22},\lambda_{23},\lambda_{24}\} \\ \{\lambda_{25},\lambda_{26},\lambda_{27}\} & \lambda_{21} \end{pmatrix} \\ \begin{pmatrix} 0 & \{0,0,0\} \\ \{-\lambda_9,-\lambda_{10},-\lambda_{11}\} & \lambda_4 \end{pmatrix} & \begin{pmatrix} -\gamma_1 & \{0,0,0\} \\ \{0,0,0\} & -\gamma_1 \end{pmatrix} & \begin{pmatrix} 0 & \{\lambda_{14},\lambda_{15},\lambda_{16}\} \\ \{\lambda_{17},\lambda_{18},\lambda_{19}\} & \lambda_{13} \end{pmatrix} \\ \begin{pmatrix} \lambda_{21} & \{-\lambda_{22},-\lambda_{23},-\lambda_{24}\} \\ \{-\lambda_{25},-\lambda_{26},-\lambda_{27}\} & \lambda_{20} \end{pmatrix} & \begin{pmatrix} \lambda_{13} & \{-\lambda_{14},-\lambda_{15},-\lambda_{16}\} \\ \{-\lambda_{17},-\lambda_{18},-\lambda_{19}\} & 0 \end{pmatrix} & \begin{pmatrix} \lambda_3 & \{0,0,0\} \\ \{0,0,0\} & \lambda_3 \end{pmatrix} \end{pmatrix}$$

**U$_{X_2[e_2]}$[E$_3$[1]] == ZERO**

True

**defin[X$_2$[e$_2$], E$_3$[1]] // Expand**

$$\begin{pmatrix} \begin{pmatrix} 0 & \{0,0,0\} \\ \{0,0,0\} & 0 \end{pmatrix} & \begin{pmatrix} \lambda_{20} & \{0,0,0\} \\ \{\lambda_{25},\lambda_{26},\lambda_{27}\} & 0 \end{pmatrix} & \begin{pmatrix} 0 & \{0,0,0\} \\ \{0,0,0\} & 0 \end{pmatrix} \\ \begin{pmatrix} 0 & \{0,0,0\} \\ \{-\lambda_{25},-\lambda_{26},-\lambda_{27}\} & \lambda_{20} \end{pmatrix} & \begin{pmatrix} 0 & \{0,0,0\} \\ \{0,0,0\} & 0 \end{pmatrix} & \begin{pmatrix} 0 & \{0,0,0\} \\ \{0,0,0\} & \lambda_3-\gamma_1 \end{pmatrix} \\ \begin{pmatrix} 0 & \{0,0,0\} \\ \{0,0,0\} & 0 \end{pmatrix} & \begin{pmatrix} \lambda_3-\gamma_1 & \{0,0,0\} \\ \{0,0,0\} & 0 \end{pmatrix} & \begin{pmatrix} 0 & \{0,0,0\} \\ \{0,0,0\} & 0 \end{pmatrix} \end{pmatrix}$$

**Δ[X$_2$[e$_2$]] = Δ[X$_2$[e$_2$]] //. {λ$_3$ → γ$_1$, λ$_{20}$ → 0, λ$_{25}$ → 0, λ$_{26}$ → 0, λ$_{27}$ → 0}**

$$\begin{pmatrix} \begin{pmatrix} 0 & \{0,0,0\} \\ \{0,0,0\} & 0 \end{pmatrix} & \begin{pmatrix} \lambda_4 & \{0,0,0\} \\ \{\lambda_9,\lambda_{10},\lambda_{11}\} & 0 \end{pmatrix} & \begin{pmatrix} 0 & \{\lambda_{22},\lambda_{23},\lambda_{24}\} \\ \{0,0,0\} & \lambda_{21} \end{pmatrix} \\ \begin{pmatrix} 0 & \{0,0,0\} \\ \{-\lambda_9,-\lambda_{10},-\lambda_{11}\} & \lambda_4 \end{pmatrix} & \begin{pmatrix} -\gamma_1 & \{0,0,0\} \\ \{0,0,0\} & -\gamma_1 \end{pmatrix} & \begin{pmatrix} 0 & \{\lambda_{14},\lambda_{15},\lambda_{16}\} \\ \{\lambda_{17},\lambda_{18},\lambda_{19}\} & \lambda_{13} \end{pmatrix} \\ \begin{pmatrix} \lambda_{21} & \{-\lambda_{22},-\lambda_{23},-\lambda_{24}\} \\ \{0,0,0\} & 0 \end{pmatrix} & \begin{pmatrix} \lambda_{13} & \{-\lambda_{14},-\lambda_{15},-\lambda_{16}\} \\ \{-\lambda_{17},-\lambda_{18},-\lambda_{19}\} & 0 \end{pmatrix} & \begin{pmatrix} \gamma_1 & \{0,0,0\} \\ \{0,0,0\} & \gamma_1 \end{pmatrix} \end{pmatrix}$$

**U$_{X_2[e_2]}$[X$_1$[e$_5$]] == ZERO**

True

**defin[X$_2$[e$_2$], X$_1$[e$_5$]] // Expand**

$$\begin{pmatrix} \begin{pmatrix} 0 & \{0,0,0\} \\ \{0,0,0\} & 0 \end{pmatrix} & \begin{pmatrix} 0 & \{0,0,0\} \\ \{0,0,0\} & 0 \end{pmatrix} & \begin{pmatrix} 0 & \{0,0,0\} \\ \{0,0,0\} & 0 \end{pmatrix} \\ \begin{pmatrix} 0 & \{0,0,0\} \\ \{0,0,0\} & 0 \end{pmatrix} & \begin{pmatrix} 0 & \{0,0,0\} \\ \{0,0,0\} & 0 \end{pmatrix} & \begin{pmatrix} 0 & \{0,0,0\} \\ \{0,0,0\} & -\beta_8-\lambda_{11} \end{pmatrix} \\ \begin{pmatrix} 0 & \{0,0,0\} \\ \{0,0,0\} & 0 \end{pmatrix} & \begin{pmatrix} -\beta_8-\lambda_{11} & \{0,0,0\} \\ \{0,0,0\} & 0 \end{pmatrix} & \begin{pmatrix} 0 & \{0,0,0\} \\ \{0,0,0\} & 0 \end{pmatrix} \end{pmatrix}$$

**Δ[X$_2$[e$_2$]] = Δ[X$_2$[e$_2$]] //. {λ$_{11}$ → -β$_8$}**

$$\begin{pmatrix} \begin{pmatrix} 0 & \{0,0,0\} \\ \{0,0,0\} & 0 \end{pmatrix} & \begin{pmatrix} \lambda_4 & \{0,0,0\} \\ \{\lambda_9,\lambda_{10},-\beta_8\} & 0 \end{pmatrix} & \begin{pmatrix} 0 & \{\lambda_{22},\lambda_{23},\lambda_{24}\} \\ \{0,0,0\} & \lambda_{21} \end{pmatrix} \\ \begin{pmatrix} 0 & \{0,0,0\} \\ \{-\lambda_9,-\lambda_{10},\beta_8\} & \lambda_4 \end{pmatrix} & \begin{pmatrix} -\gamma_1 & \{0,0,0\} \\ \{0,0,0\} & -\gamma_1 \end{pmatrix} & \begin{pmatrix} 0 & \{\lambda_{14},\lambda_{15},\lambda_{16}\} \\ \{\lambda_{17},\lambda_{18},\lambda_{19}\} & \lambda_{13} \end{pmatrix} \\ \begin{pmatrix} \lambda_{21} & \{-\lambda_{22},-\lambda_{23},-\lambda_{24}\} \\ \{0,0,0\} & 0 \end{pmatrix} & \begin{pmatrix} \lambda_{13} & \{-\lambda_{14},-\lambda_{15},-\lambda_{16}\} \\ \{-\lambda_{17},-\lambda_{18},-\lambda_{19}\} & 0 \end{pmatrix} & \begin{pmatrix} \gamma_1 & \{0,0,0\} \\ \{0,0,0\} & \gamma_1 \end{pmatrix} \end{pmatrix}$$

**U$_{X_2[e_2]}$[X$_1$[e$_4$]] == ZERO**

True

**defin[X$_2$[e$_2$], X$_1$[e$_4$]] // Expand**

$$\begin{pmatrix} \begin{pmatrix} 0 & \{0,0,0\} \\ \{0,0,0\} & 0 \end{pmatrix} & \begin{pmatrix} 0 & \{0,0,0\} \\ \{0,0,0\} & 0 \end{pmatrix} & \begin{pmatrix} 0 & \{0,0,0\} \\ \{0,0,0\} & 0 \end{pmatrix} \\ \begin{pmatrix} 0 & \{0,0,0\} \\ \{0,0,0\} & 0 \end{pmatrix} & \begin{pmatrix} 0 & \{0,0,0\} \\ \{0,0,0\} & 0 \end{pmatrix} & \begin{pmatrix} 0 & \{0,0,0\} \\ \{0,0,0\} & -\beta_7-\lambda_{10} \end{pmatrix} \\ \begin{pmatrix} 0 & \{0,0,0\} \\ \{0,0,0\} & 0 \end{pmatrix} & \begin{pmatrix} -\beta_7-\lambda_{10} & \{0,0,0\} \\ \{0,0,0\} & 0 \end{pmatrix} & \begin{pmatrix} 0 & \{0,0,0\} \\ \{0,0,0\} & 0 \end{pmatrix} \end{pmatrix}$$

**Δ[X$_2$[e$_2$]] = Δ[X$_2$[e$_2$]] //. {λ$_{10}$ → -β$_7$}**

$$\begin{pmatrix} \begin{pmatrix} 0 & \{0,0,0\} \\ \{0,0,0\} & 0 \end{pmatrix} & \begin{pmatrix} \lambda_4 & \{0,0,0\} \\ \{\lambda_9,-\beta_7,-\beta_8\} & 0 \end{pmatrix} & \begin{pmatrix} 0 & \{\lambda_{22},\lambda_{23},\lambda_{24}\} \\ \{0,0,0\} & \lambda_{21} \end{pmatrix} \\ \begin{pmatrix} 0 & \{0,0,0\} \\ \{-\lambda_9,\beta_7,\beta_8\} & \lambda_4 \end{pmatrix} & \begin{pmatrix} -\gamma_1 & \{0,0,0\} \\ \{0,0,0\} & -\gamma_1 \end{pmatrix} & \begin{pmatrix} 0 & \{\lambda_{14},\lambda_{15},\lambda_{16}\} \\ \{\lambda_{17},\lambda_{18},\lambda_{19}\} & \lambda_{13} \end{pmatrix} \\ \begin{pmatrix} \lambda_{21} & \{-\lambda_{22},-\lambda_{23},-\lambda_{24}\} \\ \{0,0,0\} & 0 \end{pmatrix} & \begin{pmatrix} \lambda_{13} & \{-\lambda_{14},-\lambda_{15},-\lambda_{16}\} \\ \{-\lambda_{17},-\lambda_{18},-\lambda_{19}\} & 0 \end{pmatrix} & \begin{pmatrix} \gamma_1 & \{0,0,0\} \\ \{0,0,0\} & \gamma_1 \end{pmatrix} \end{pmatrix}$$

**U$_{X_2[e_2]}$[X$_1$[e$_3$]] == ZERO**

True

**defin[X$_2$[e$_2$], X$_1$[e$_3$]] // Expand**

$$\begin{pmatrix} \begin{pmatrix} 0 & \{0,0,0\} \\ \{0,0,0\} & 0 \end{pmatrix} & \begin{pmatrix} 0 & \{0,0,0\} \\ \{0,0,0\} & 0 \end{pmatrix} & \begin{pmatrix} 0 & \{0,0,0\} \\ \{0,0,0\} & 0 \end{pmatrix} \\ \begin{pmatrix} 0 & \{0,0,0\} \\ \{0,0,0\} & 0 \end{pmatrix} & \begin{pmatrix} 0 & \{0,0,0\} \\ \{0,0,0\} & 0 \end{pmatrix} & \begin{pmatrix} 0 & \{0,0,0\} \\ \{0,0,0\} & -\beta_6-\lambda_9 \end{pmatrix} \\ \begin{pmatrix} 0 & \{0,0,0\} \\ \{0,0,0\} & 0 \end{pmatrix} & \begin{pmatrix} -\beta_6-\lambda_9 & \{0,0,0\} \\ \{0,0,0\} & 0 \end{pmatrix} & \begin{pmatrix} 0 & \{0,0,0\} \\ \{0,0,0\} & 0 \end{pmatrix} \end{pmatrix}$$

**Δ[X$_2$[e$_2$]] = Δ[X$_2$[e$_2$]] //. {λ$_9$ → -β$_6$}**

$$\begin{pmatrix} \begin{pmatrix} 0 & \{0,0,0\} \\ \{0,0,0\} & 0 \end{pmatrix} & \begin{pmatrix} \lambda_4 & \{0,0,0\} \\ \{-\beta_6,-\beta_7,-\beta_8\} & 0 \end{pmatrix} & \begin{pmatrix} 0 & \{\lambda_{22},\lambda_{23},\lambda_{24}\} \\ \{0,0,0\} & \lambda_{21} \end{pmatrix} \\ \begin{pmatrix} 0 & \{0,0,0\} \\ \{\beta_6,\beta_7,\beta_8\} & \lambda_4 \end{pmatrix} & \begin{pmatrix} -\gamma_1 & \{0,0,0\} \\ \{0,0,0\} & -\gamma_1 \end{pmatrix} & \begin{pmatrix} 0 & \{\lambda_{14},\lambda_{15},\lambda_{16}\} \\ \{\lambda_{17},\lambda_{18},\lambda_{19}\} & \lambda_{13} \end{pmatrix} \\ \begin{pmatrix} \lambda_{21} & \{-\lambda_{22},-\lambda_{23},-\lambda_{24}\} \\ \{0,0,0\} & 0 \end{pmatrix} & \begin{pmatrix} \lambda_{13} & \{-\lambda_{14},-\lambda_{15},-\lambda_{16}\} \\ \{-\lambda_{17},-\lambda_{18},-\lambda_{19}\} & 0 \end{pmatrix} & \begin{pmatrix} \gamma_1 & \{0,0,0\} \\ \{0,0,0\} & \gamma_1 \end{pmatrix} \end{pmatrix}$$

**U$_{X_1[e_2]}$[X$_2$[e$_2$]] == ZERO**

True

**(defin[X$_1$[e$_2$], X$_2$[e$_2$]] // Expand)**

$$\left( \begin{pmatrix} 0 & \{0,0,0\} \\ \{0,0,0\} & 0 \end{pmatrix} \begin{pmatrix} 0 & \{0,0,0\} \\ \{0,0,0\} & \beta_1+\lambda_4 \end{pmatrix} \begin{pmatrix} 0 & \{0,0,0\} \\ \{0,0,0\} & 0 \end{pmatrix} \right.$$
$$\begin{pmatrix} \beta_1+\lambda_4 & \{0,0,0\} \\ \{0,0,0\} & 0 \end{pmatrix} \begin{pmatrix} 0 & \{0,0,0\} \\ \{0,0,0\} & 0 \end{pmatrix} \begin{pmatrix} 0 & \{0,0,0\} \\ \{0,0,0\} & 0 \end{pmatrix}$$
$$\left. \begin{pmatrix} 0 & \{0,0,0\} \\ \{0,0,0\} & 0 \end{pmatrix} \begin{pmatrix} 0 & \{0,0,0\} \\ \{0,0,0\} & 0 \end{pmatrix} \begin{pmatrix} 0 & \{0,0,0\} \\ \{0,0,0\} & 0 \end{pmatrix} \right)$$

**Δ[X$_2$[e$_2$]] = Δ[X$_2$[e$_2$]] //. {λ$_4$ → -β$_1$}**

$$\left( \begin{pmatrix} 0 & \{0,0,0\} \\ \{0,0,0\} & 0 \end{pmatrix} \begin{pmatrix} -\beta_1 & \{0,0,0\} \\ \{-\beta_6,-\beta_7,-\beta_8\} & 0 \end{pmatrix} \begin{pmatrix} 0 & \{\lambda_{22},\lambda_{23},\lambda_{24}\} \\ \{0,0,0\} & \lambda_{21} \end{pmatrix} \right.$$
$$\begin{pmatrix} 0 & \{0,0,0\} \\ \{\beta_6,\beta_7,\beta_8\} & -\beta_1 \end{pmatrix} \begin{pmatrix} -\gamma_1 & \{0,0,0\} \\ \{0,0,0\} & -\gamma_1 \end{pmatrix} \begin{pmatrix} 0 & \{\lambda_{14},\lambda_{15},\lambda_{16}\} \\ \{\lambda_{17},\lambda_{18},\lambda_{19}\} & \lambda_{13} \end{pmatrix}$$
$$\left. \begin{pmatrix} \lambda_{21} & \{-\lambda_{22},-\lambda_{23},-\lambda_{24}\} \\ \{0,0,0\} & 0 \end{pmatrix} \begin{pmatrix} \lambda_{13} & \{-\lambda_{14},-\lambda_{15},-\lambda_{16}\} \\ \{-\lambda_{17},-\lambda_{18},-\lambda_{19}\} & 0 \end{pmatrix} \begin{pmatrix} \gamma_1 & \{0,0,0\} \\ \{0,0,0\} & \gamma_1 \end{pmatrix} \right)$$

**T[E$_1$[1], E$_1$[1], X$_2$[e$_2$]] == ZERO**

True

**Leib[E$_1$[1], E$_1$[1], X$_2$[e$_2$]] // Expand**

$$\left( \begin{pmatrix} 0 & \{0,0,0\} \\ \{0,0,0\} & 0 \end{pmatrix} \begin{pmatrix} 0 & \{0,0,0\} \\ \{0,0,0\} & 0 \end{pmatrix} \begin{pmatrix} 0 & \{\alpha_3+\lambda_{22},\alpha_4+\lambda_{23},\alpha_5+\lambda_{24}\} \\ \{0,0,0\} & \alpha_2+\lambda_{21} \end{pmatrix} \right.$$
$$\begin{pmatrix} 0 & \{0,0,0\} \\ \{0,0,0\} & 0 \end{pmatrix} \begin{pmatrix} 0 & \{0,0,0\} \\ \{0,0,0\} & 0 \end{pmatrix} \begin{pmatrix} 0 & \{0,0,0\} \\ \{0,0,0\} & 0 \end{pmatrix}$$
$$\left. \begin{pmatrix} \alpha_2+\lambda_{21} & \{-\alpha_3-\lambda_{22},-\alpha_4-\lambda_{23},-\alpha_5-\lambda_{24}\} \\ \{0,0,0\} & 0 \end{pmatrix} \begin{pmatrix} 0 & \{0,0,0\} \\ \{0,0,0\} & 0 \end{pmatrix} \begin{pmatrix} 0 & \{0,0,0\} \\ \{0,0,0\} & 0 \end{pmatrix} \right)$$

**Δ[X$_2$[e$_2$]] = Δ[X$_2$[e$_2$]] //. {λ$_{21}$ → -α$_2$, λ$_{22}$ → -α$_3$, λ$_{23}$ → -α$_4$, λ$_{24}$ → -α$_5$}**

$$\left( \begin{pmatrix} 0 & \{0,0,0\} \\ \{0,0,0\} & 0 \end{pmatrix} \begin{pmatrix} -\beta_1 & \{0,0,0\} \\ \{-\beta_6,-\beta_7,-\beta_8\} & 0 \end{pmatrix} \begin{pmatrix} 0 & \{-\alpha_3,-\alpha_4,-\alpha_5\} \\ \{0,0,0\} & -\alpha_2 \end{pmatrix} \right.$$
$$\begin{pmatrix} 0 & \{0,0,0\} \\ \{\beta_6,\beta_7,\beta_8\} & -\beta_1 \end{pmatrix} \begin{pmatrix} -\gamma_1 & \{0,0,0\} \\ \{0,0,0\} & -\gamma_1 \end{pmatrix} \begin{pmatrix} 0 & \{\lambda_{14},\lambda_{15},\lambda_{16}\} \\ \{\lambda_{17},\lambda_{18},\lambda_{19}\} & \lambda_{13} \end{pmatrix}$$
$$\left. \begin{pmatrix} -\alpha_2 & \{\alpha_3,\alpha_4,\alpha_5\} \\ \{0,0,0\} & 0 \end{pmatrix} \begin{pmatrix} \lambda_{13} & \{-\lambda_{14},-\lambda_{15},-\lambda_{16}\} \\ \{-\lambda_{17},-\lambda_{18},-\lambda_{19}\} & 0 \end{pmatrix} \begin{pmatrix} \gamma_1 & \{0,0,0\} \\ \{0,0,0\} & \gamma_1 \end{pmatrix} \right)$$

**T[X$_2$[e$_2$], X$_1$[e$_1$], X$_1$[e$_2$]]**

$$\left( \begin{pmatrix} 0 & \{0,0,0\} \\ \{0,0,0\} & 0 \end{pmatrix} \begin{pmatrix} 0 & \{0,0,0\} \\ \{0,0,0\} & 0 \end{pmatrix} \begin{pmatrix} 0 & \{0,0,0\} \\ \{0,0,0\} & 0 \end{pmatrix} \right.$$
$$\begin{pmatrix} 0 & \{0,0,0\} \\ \{0,0,0\} & 0 \end{pmatrix} \begin{pmatrix} 0 & \{0,0,0\} \\ \{0,0,0\} & 0 \end{pmatrix} \begin{pmatrix} 0 & \{0,0,0\} \\ \{0,0,0\} & 0 \end{pmatrix}$$
$$\left. \begin{pmatrix} 0 & \{0,0,0\} \\ \{0,0,0\} & 0 \end{pmatrix} \begin{pmatrix} 0 & \{0,0,0\} \\ \{0,0,0\} & 0 \end{pmatrix} \begin{pmatrix} 0 & \{0,0,0\} \\ \{0,0,0\} & 0 \end{pmatrix} \right)$$

**Leib[X$_2$[e$_2$], X$_1$[e$_1$], X$_1$[e$_2$]] // Expand**

$$\left(\begin{array}{ccc} \begin{pmatrix} 0 & \{0,0,0\} \\ \{0,0,0\} & 0 \end{pmatrix} & \begin{pmatrix} 0 & \{0,0,0\} \\ \{0,0,0\} & 0 \end{pmatrix} & \begin{pmatrix} 0 & \{0,0,0\} \\ \{0,0,0\} & 0 \end{pmatrix} \\ \begin{pmatrix} 0 & \{0,0,0\} \\ \{0,0,0\} & 0 \end{pmatrix} & \begin{pmatrix} 0 & \{0,0,0\} \\ \{0,0,0\} & 0 \end{pmatrix} & \begin{pmatrix} 0 & \{\delta_2+\lambda_{14},\delta_3+\lambda_{15},\delta_4+\lambda_{16}\} \\ \{0,0,0\} & 0 \end{pmatrix} \\ \begin{pmatrix} 0 & \{0,0,0\} \\ \{0,0,0\} & 0 \end{pmatrix} & \begin{pmatrix} 0 & \{-\delta_2-\lambda_{14},-\delta_3-\lambda_{15},-\delta_4-\lambda_{16}\} \\ \{0,0,0\} & 0 \end{pmatrix} & \begin{pmatrix} 0 & \{0,0,0\} \\ \{0,0,0\} & 0 \end{pmatrix} \end{array}\right)$$

**Δ[X$_2$[e$_2$]] = Δ[X$_2$[e$_2$]] //. {λ$_{14}$ → -δ$_2$, λ$_{15}$ → -δ$_3$, λ$_{16}$ → -δ$_4$}**

$$\left(\begin{array}{ccc} \begin{pmatrix} 0 & \{0,0,0\} \\ \{0,0,0\} & 0 \end{pmatrix} & \begin{pmatrix} -\beta_1 & \{0,0,0\} \\ \{-\beta_6,-\beta_7,-\beta_8\} & 0 \end{pmatrix} & \begin{pmatrix} 0 & \{-\alpha_3,-\alpha_4,-\alpha_5\} \\ \{0,0,0\} & -\alpha_2 \end{pmatrix} \\ \begin{pmatrix} 0 & \{0,0,0\} \\ \{\beta_6,\beta_7,\beta_8\} & -\beta_1 \end{pmatrix} & \begin{pmatrix} -\gamma_1 & \{0,0,0\} \\ \{0,0,0\} & -\gamma_1 \end{pmatrix} & \begin{pmatrix} 0 & \{-\delta_2,-\delta_3,-\delta_4\} \\ \{\lambda_{17},\lambda_{18},\lambda_{19}\} & \lambda_{13} \end{pmatrix} \\ \begin{pmatrix} -\alpha_2 & \{\alpha_3,\alpha_4,\alpha_5\} \\ \{0,0,0\} & 0 \end{pmatrix} & \begin{pmatrix} \lambda_{13} & \{\delta_2,\delta_3,\delta_4\} \\ \{-\lambda_{17},-\lambda_{18},-\lambda_{19}\} & 0 \end{pmatrix} & \begin{pmatrix} \gamma_1 & \{0,0,0\} \\ \{0,0,0\} & \gamma_1 \end{pmatrix} \end{array}\right)$$

**T[X$_2$[e$_2$], X$_1$[e$_2$], X$_1$[e$_8$]] == ZERO**

True

**Leib[X$_2$[e$_2$], X$_1$[e$_2$], X$_1$[e$_8$]] // Expand**

$$\left(\begin{array}{ccc} \begin{pmatrix} 0 & \{0,0,0\} \\ \{0,0,0\} & 0 \end{pmatrix} & \begin{pmatrix} 0 & \{0,0,0\} \\ \{0,0,0\} & 0 \end{pmatrix} & \begin{pmatrix} 0 & \{0,0,0\} \\ \{0,0,0\} & 0 \end{pmatrix} \\ \begin{pmatrix} 0 & \{0,0,0\} \\ \{0,0,0\} & 0 \end{pmatrix} & \begin{pmatrix} 0 & \{0,0,0\} \\ \{0,0,0\} & 0 \end{pmatrix} & \begin{pmatrix} 0 & \{\psi_2-\lambda_{18},\lambda_{17}+\chi_1,0\} \\ \{0,0,0\} & 0 \end{pmatrix} \\ \begin{pmatrix} 0 & \{0,0,0\} \\ \{0,0,0\} & 0 \end{pmatrix} & \begin{pmatrix} 0 & \{\lambda_{18}-\psi_2,-\lambda_{17}-\chi_1,0\} \\ \{0,0,0\} & 0 \end{pmatrix} & \begin{pmatrix} 0 & \{0,0,0\} \\ \{0,0,0\} & 0 \end{pmatrix} \end{array}\right)$$

**Δ[X$_2$[e$_2$]] = Δ[X$_2$[e$_2$]] //. {λ$_{17}$ → -χ$_1$, λ$_{18}$ → ψ$_2$}**

$$\left(\begin{array}{ccc} \begin{pmatrix} 0 & \{0,0,0\} \\ \{0,0,0\} & 0 \end{pmatrix} & \begin{pmatrix} -\beta_1 & \{0,0,0\} \\ \{-\beta_6,-\beta_7,-\beta_8\} & 0 \end{pmatrix} & \begin{pmatrix} 0 & \{-\alpha_3,-\alpha_4,-\alpha_5\} \\ \{0,0,0\} & -\alpha_2 \end{pmatrix} \\ \begin{pmatrix} 0 & \{0,0,0\} \\ \{\beta_6,\beta_7,\beta_8\} & -\beta_1 \end{pmatrix} & \begin{pmatrix} -\gamma_1 & \{0,0,0\} \\ \{0,0,0\} & -\gamma_1 \end{pmatrix} & \begin{pmatrix} 0 & \{-\delta_2,-\delta_3,-\delta_4\} \\ \{-\chi_1,\psi_2,\lambda_{19}\} & \lambda_{13} \end{pmatrix} \\ \begin{pmatrix} -\alpha_2 & \{\alpha_3,\alpha_4,\alpha_5\} \\ \{0,0,0\} & 0 \end{pmatrix} & \begin{pmatrix} \lambda_{13} & \{\delta_2,\delta_3,\delta_4\} \\ \{\chi_1,-\psi_2,-\lambda_{19}\} & 0 \end{pmatrix} & \begin{pmatrix} \gamma_1 & \{0,0,0\} \\ \{0,0,0\} & \gamma_1 \end{pmatrix} \end{array}\right)$$

**T[X$_2$[e$_2$], X$_1$[e$_2$], X$_1$[e$_7$]]**

$$\left(\begin{array}{ccc} \begin{pmatrix} 0 & \{0,0,0\} \\ \{0,0,0\} & 0 \end{pmatrix} & \begin{pmatrix} 0 & \{0,0,0\} \\ \{0,0,0\} & 0 \end{pmatrix} & \begin{pmatrix} 0 & \{0,0,0\} \\ \{0,0,0\} & 0 \end{pmatrix} \\ \begin{pmatrix} 0 & \{0,0,0\} \\ \{0,0,0\} & 0 \end{pmatrix} & \begin{pmatrix} 0 & \{0,0,0\} \\ \{0,0,0\} & 0 \end{pmatrix} & \begin{pmatrix} 0 & \{0,0,0\} \\ \{0,0,0\} & 0 \end{pmatrix} \\ \begin{pmatrix} 0 & \{0,0,0\} \\ \{0,0,0\} & 0 \end{pmatrix} & \begin{pmatrix} 0 & \{0,0,0\} \\ \{0,0,0\} & 0 \end{pmatrix} & \begin{pmatrix} 0 & \{0,0,0\} \\ \{0,0,0\} & 0 \end{pmatrix} \end{array}\right)$$

**Leib[X₂[e₂], X₁[e₂], X₁[e₇]] // Expand**

$$\left( \begin{pmatrix} 0 & \{0,0,0\} \\ \{0,0,0\} & 0 \end{pmatrix} \begin{pmatrix} 0 & \{0,0,0\} \\ \{0,0,0\} & 0 \end{pmatrix} \begin{pmatrix} 0 & \{0,0,0\} \\ \{0,0,0\} & 0 \end{pmatrix} \\ \begin{pmatrix} 0 & \{0,0,0\} \\ \{0,0,0\} & 0 \end{pmatrix} \begin{pmatrix} 0 & \{0,0,0\} \\ \{0,0,0\} & 0 \end{pmatrix} \begin{pmatrix} 0 & \{\lambda_{19}+\psi_1,0,0\} \\ \{0,0,0\} & 0 \end{pmatrix} \\ \begin{pmatrix} 0 & \{0,0,0\} \\ \{0,0,0\} & 0 \end{pmatrix} \begin{pmatrix} 0 & \{-\lambda_{19}-\psi_1,0,0\} \\ \{0,0,0\} & 0 \end{pmatrix} \begin{pmatrix} 0 & \{0,0,0\} \\ \{0,0,0\} & 0 \end{pmatrix} \right)$$

**Δ[X₂[e₂]] = Δ[X₂[e₂]] //. {λ₁₉ → -ψ₁}**

$$\left( \begin{pmatrix} 0 & \{0,0,0\} \\ \{0,0,0\} & 0 \end{pmatrix} \begin{pmatrix} -\beta_1 & \{0,0,0\} \\ \{-\beta_6,-\beta_7,-\beta_8\} & 0 \end{pmatrix} \begin{pmatrix} 0 & \{-\alpha_3,-\alpha_4,-\alpha_5\} \\ \{0,0,0\} & -\alpha_2 \end{pmatrix} \\ \begin{pmatrix} 0 & \{0,0,0\} \\ \{\beta_6,\beta_7,\beta_8\} & -\beta_1 \end{pmatrix} \begin{pmatrix} -\gamma_1 & \{0,0,0\} \\ \{0,0,0\} & -\gamma_1 \end{pmatrix} \begin{pmatrix} 0 & \{-\delta_2,-\delta_3,-\delta_4\} \\ \{-\chi_1,\psi_2,-\psi_1\} & \lambda_{13} \end{pmatrix} \\ \begin{pmatrix} -\alpha_2 & \{\alpha_3,\alpha_4,\alpha_5\} \\ \{0,0,0\} & 0 \end{pmatrix} \begin{pmatrix} \lambda_{13} & \{\delta_2,\delta_3,\delta_4\} \\ \{\chi_1,-\psi_2,\psi_1\} & 0 \end{pmatrix} \begin{pmatrix} \gamma_1 & \{0,0,0\} \\ \{0,0,0\} & \gamma_1 \end{pmatrix} \right)$$

**T[X₂[e₂], X₂[e₁], X₁[e₁]] == X₁[e₁]**

True

**Δ[X₁[e₁]] - Leib[X₂[e₂], X₂[e₁], X₁[e₁]] // Expand**

$$\left( \begin{pmatrix} 0 & \{0,0,0\} \\ \{0,0,0\} & 0 \end{pmatrix} \begin{pmatrix} -\xi-\lambda_{13} & \{0,0,0\} \\ \{0,0,0\} & 0 \end{pmatrix} \begin{pmatrix} 0 & \{0,0,0\} \\ \{0,0,0\} & 0 \end{pmatrix} \\ \begin{pmatrix} 0 & \{0,0,0\} \\ \{0,0,0\} & -\xi-\lambda_{13} \end{pmatrix} \begin{pmatrix} 0 & \{0,0,0\} \\ \{0,0,0\} & 0 \end{pmatrix} \begin{pmatrix} 0 & \{0,0,0\} \\ \{0,0,0\} & 0 \end{pmatrix} \\ \begin{pmatrix} 0 & \{0,0,0\} \\ \{0,0,0\} & 0 \end{pmatrix} \begin{pmatrix} 0 & \{0,0,0\} \\ \{0,0,0\} & 0 \end{pmatrix} \begin{pmatrix} 0 & \{0,0,0\} \\ \{0,0,0\} & 0 \end{pmatrix} \right)$$

**Δ[X₂[e₂]] = Δ[X₂[e₂]] //. {λ₁₃ → -ξ}**

$$\left( \begin{pmatrix} 0 & \{0,0,0\} \\ \{0,0,0\} & 0 \end{pmatrix} \begin{pmatrix} -\beta_1 & \{0,0,0\} \\ \{-\beta_6,-\beta_7,-\beta_8\} & 0 \end{pmatrix} \begin{pmatrix} 0 & \{-\alpha_3,-\alpha_4,-\alpha_5\} \\ \{0,0,0\} & -\alpha_2 \end{pmatrix} \\ \begin{pmatrix} 0 & \{0,0,0\} \\ \{\beta_6,\beta_7,\beta_8\} & -\beta_1 \end{pmatrix} \begin{pmatrix} -\gamma_1 & \{0,0,0\} \\ \{0,0,0\} & -\gamma_1 \end{pmatrix} \begin{pmatrix} 0 & \{-\delta_2,-\delta_3,-\delta_4\} \\ \{-\chi_1,\psi_2,-\psi_1\} & -\xi \end{pmatrix} \\ \begin{pmatrix} -\alpha_2 & \{\alpha_3,\alpha_4,\alpha_5\} \\ \{0,0,0\} & 0 \end{pmatrix} \begin{pmatrix} -\xi & \{\delta_2,\delta_3,\delta_4\} \\ \{\chi_1,-\psi_2,\psi_1\} & 0 \end{pmatrix} \begin{pmatrix} \gamma_1 & \{0,0,0\} \\ \{0,0,0\} & \gamma_1 \end{pmatrix} \right)$$

**Variables[Δ[X₂[e₂]]]**

$\{\xi, \alpha_2, \alpha_3, \alpha_4, \alpha_5, \beta_1, \beta_6, \beta_7, \beta_8, \gamma_1, \delta_2, \delta_3, \delta_4, \chi_1, \psi_1, \psi_2\}$

- **X₂[e₃]**

  **Δ[X₂[e₃]] = generic;**

**U$_{E_1[1]}$[X$_2$[e$_3$]] == ZERO**

True

**(defin[E$_1$[1], X$_2$[e$_3$]] // Expand)**

$$\begin{pmatrix} \begin{pmatrix} \lambda_1 & \{0,0,0\} \\ \{0,0,0\} & \lambda_1 \end{pmatrix} & \begin{pmatrix} 0 & \{0,0,0\} \\ \{0,0,0\} & 0 \end{pmatrix} & \begin{pmatrix} 0 & \{0,0,0\} \\ \{0,0,0\} & 0 \end{pmatrix} \\ \begin{pmatrix} 0 & \{0,0,0\} \\ \{0,0,0\} & 0 \end{pmatrix} & \begin{pmatrix} 0 & \{0,0,0\} \\ \{0,0,0\} & 0 \end{pmatrix} & \begin{pmatrix} 0 & \{0,0,0\} \\ \{0,0,0\} & 0 \end{pmatrix} \\ \begin{pmatrix} 0 & \{0,0,0\} \\ \{0,0,0\} & 0 \end{pmatrix} & \begin{pmatrix} 0 & \{0,0,0\} \\ \{0,0,0\} & 0 \end{pmatrix} & \begin{pmatrix} 0 & \{0,0,0\} \\ \{0,0,0\} & 0 \end{pmatrix} \end{pmatrix}$$

**Δ[X$_2$[e$_3$]] = Δ[X$_2$[e$_3$]] //. {λ$_1$ → 0}**

$$\begin{pmatrix} \begin{pmatrix} 0 & \{0,0,0\} \\ \{0,0,0\} & 0 \end{pmatrix} & \begin{pmatrix} \lambda_4 & \{\lambda_6, \lambda_7, \lambda_8\} \\ \{\lambda_9, \lambda_{10}, \lambda_{11}\} & \lambda_5 \end{pmatrix} & \begin{pmatrix} \lambda_{20} & \{\lambda_{22}, \lambda_{23}, \lambda_{24}\} \\ \{\lambda_{25}, \lambda_{26}, \lambda_{27}\} & \lambda_{21} \end{pmatrix} \\ \begin{pmatrix} \lambda_5 & \{-\lambda_6, -\lambda_7, -\lambda_8\} \\ \{-\lambda_9, -\lambda_{10}, -\lambda_{11}\} & \lambda_4 \end{pmatrix} & \begin{pmatrix} \lambda_2 & \{0,0,0\} \\ \{0,0,0\} & \lambda_2 \end{pmatrix} & \begin{pmatrix} \lambda_{12} & \{\lambda_{14}, \lambda_{15}, \lambda_{16}\} \\ \{\lambda_{17}, \lambda_{18}, \lambda_{19}\} & \lambda_{13} \end{pmatrix} \\ \begin{pmatrix} \lambda_{21} & \{-\lambda_{22}, -\lambda_{23}, -\lambda_{24}\} \\ \{-\lambda_{25}, -\lambda_{26}, -\lambda_{27}\} & \lambda_{20} \end{pmatrix} & \begin{pmatrix} \lambda_{13} & \{-\lambda_{14}, -\lambda_{15}, -\lambda_{16}\} \\ \{-\lambda_{17}, -\lambda_{18}, -\lambda_{19}\} & \lambda_{12} \end{pmatrix} & \begin{pmatrix} \lambda_3 & \{0,0,0\} \\ \{0,0,0\} & \lambda_3 \end{pmatrix} \end{pmatrix}$$

**U$_{X_2[e_3]}$[E$_2$[1]] == ZERO**

True

**defin[X$_2$[e$_3$], E$_2$[1]] // Expand**

$$\begin{pmatrix} \begin{pmatrix} 0 & \{0,0,0\} \\ \{0,0,0\} & 0 \end{pmatrix} & \begin{pmatrix} 0 & \{0,0,0\} \\ \{0,0,0\} & 0 \end{pmatrix} & \begin{pmatrix} 0 & \{\lambda_4, 0, 0\} \\ \{0, \lambda_8, -\lambda_7\} & \lambda_9 \end{pmatrix} \\ \begin{pmatrix} 0 & \{0,0,0\} \\ \{0,0,0\} & 0 \end{pmatrix} & \begin{pmatrix} 0 & \{0,0,0\} \\ \{0,0,0\} & 0 \end{pmatrix} & \begin{pmatrix} 0 & \{\lambda_2 - \gamma_6, 0, 0\} \\ \{0,0,0\} & 0 \end{pmatrix} \\ \begin{pmatrix} \lambda_9 & \{-\lambda_4, 0, 0\} \\ \{0, -\lambda_8, \lambda_7\} & 0 \end{pmatrix} & \begin{pmatrix} 0 & \{\gamma_6 - \lambda_2, 0, 0\} \\ \{0,0,0\} & 0 \end{pmatrix} & \begin{pmatrix} -\lambda_{17} & \{0,0,0\} \\ \{0,0,0\} & -\lambda_{17} \end{pmatrix} \end{pmatrix}$$

**Δ[X$_2$[e$_3$]] =**
 **Δ[X$_2$[e$_3$]] //. {λ$_2$ → γ$_6$, λ$_4$ → 0, λ$_7$ → 0, λ$_8$ → 0, λ$_9$ → 0, λ$_{17}$ → 0}**

$$\begin{pmatrix} \begin{pmatrix} 0 & \{0,0,0\} \\ \{0,0,0\} & 0 \end{pmatrix} & \begin{pmatrix} 0 & \{\lambda_6, 0, 0\} \\ \{0, \lambda_{10}, \lambda_{11}\} & \lambda_5 \end{pmatrix} & \begin{pmatrix} \lambda_{20} & \{\lambda_{22}, \ldots\} \\ \{\lambda_{25}, \lambda_{26}, \lambda_{27}\} & \ldots \end{pmatrix} \\ \begin{pmatrix} \lambda_5 & \{-\lambda_6, 0, 0\} \\ \{0, -\lambda_{10}, -\lambda_{11}\} & 0 \end{pmatrix} & \begin{pmatrix} \gamma_6 & \{0,0,0\} \\ \{0,0,0\} & \gamma_6 \end{pmatrix} & \begin{pmatrix} \lambda_{12} & \{\lambda_{14}, \ldots \\ \{0, \lambda_{18}, \lambda_{19}\} & \ldots \end{pmatrix} \\ \begin{pmatrix} \lambda_{21} & \{-\lambda_{22}, -\lambda_{23}, -\lambda_{24}\} \\ \{-\lambda_{25}, -\lambda_{26}, -\lambda_{27}\} & \lambda_{20} \end{pmatrix} & \begin{pmatrix} \lambda_{13} & \{-\lambda_{14}, -\lambda_{15}, -\lambda_{16}\} \\ \{0, -\lambda_{18}, -\lambda_{19}\} & \lambda_{12} \end{pmatrix} & \begin{pmatrix} \lambda_3 & \{0, 0, \ldots \\ \{0,0,0\} & \lambda_3 \end{pmatrix} \end{pmatrix}$$

**U$_{X_2[e_3]}$[E$_3$[1]] == ZERO**

True

**defin[X₂[e₃], E₃[1]] // Expand**

$$\begin{pmatrix} \begin{pmatrix} 0 & \{0,0,0\} \\ \{0,0,0\} & 0 \end{pmatrix} & \begin{pmatrix} 0 & \{-\lambda_{20},0,0\} \\ \{0,-\lambda_{24},\lambda_{23}\} & -\lambda_{25} \end{pmatrix} & \begin{pmatrix} 0 & \{0,0,0\} \\ \{0,0,0\} & 0 \end{pmatrix} \\ \begin{pmatrix} -\lambda_{25} & \{\lambda_{20},0,0\} \\ \{0,\lambda_{24},-\lambda_{23}\} & 0 \end{pmatrix} & \begin{pmatrix} 0 & \{0,0,0\} \\ \{0,0,0\} & 0 \end{pmatrix} & \begin{pmatrix} 0 & \{\gamma_6+\lambda_3,0,0\} \\ \{0,0,0\} & 0 \end{pmatrix} \\ \begin{pmatrix} 0 & \{0,0,0\} \\ \{0,0,0\} & 0 \end{pmatrix} & \begin{pmatrix} 0 & \{-\gamma_6-\lambda_3,0,0\} \\ \{0,0,0\} & 0 \end{pmatrix} & \begin{pmatrix} 0 & \{0,0,0\} \\ \{0,0,0\} & 0 \end{pmatrix} \end{pmatrix}$$

**Δ[X₂[e₃]] = Δ[X₂[e₃]] //. {λ₃ → -γ₆, λ₂₀ → 0, λ₂₃ → 0, λ₂₄ → 0, λ₂₅ → 0}**

$$\begin{pmatrix} \begin{pmatrix} 0 & \{0,0,0\} \\ \{0,0,0\} & 0 \end{pmatrix} & \begin{pmatrix} 0 & \{\lambda_6,0,0\} \\ \{0,\lambda_{10},\lambda_{11}\} & \lambda_5 \end{pmatrix} & \begin{pmatrix} 0 & \{\lambda_{22},0,0\} \\ \{0,\lambda_{26},\lambda_{27}\} & \lambda_{21} \end{pmatrix} \\ \begin{pmatrix} \lambda_5 & \{-\lambda_6,0,0\} \\ \{0,-\lambda_{10},-\lambda_{11}\} & 0 \end{pmatrix} & \begin{pmatrix} \gamma_6 & \{0,0,0\} \\ \{0,0,0\} & \gamma_6 \end{pmatrix} & \begin{pmatrix} \lambda_{12} & \{\lambda_{14},\lambda_{15},\lambda_{16}\} \\ \{0,\lambda_{18},\lambda_{19}\} & \lambda_{13} \end{pmatrix} \\ \begin{pmatrix} \lambda_{21} & \{-\lambda_{22},0,0\} \\ \{0,-\lambda_{26},-\lambda_{27}\} & 0 \end{pmatrix} & \begin{pmatrix} \lambda_{13} & \{-\lambda_{14},-\lambda_{15},-\lambda_{16}\} \\ \{0,-\lambda_{18},-\lambda_{19}\} & \lambda_{12} \end{pmatrix} & \begin{pmatrix} -\gamma_6 & \{0,0,0\} \\ \{0,0,0\} & -\gamma_6 \end{pmatrix} \end{pmatrix}$$

**U_{X₂[e₃]}[X₁[e₁]] == ZERO**

True

**defin[X₂[e₃], X₁[e₁]] // Expand**

$$\begin{pmatrix} \begin{pmatrix} 0 & \{0,0,0\} \\ \{0,0,0\} & 0 \end{pmatrix} & \begin{pmatrix} 0 & \{0,0,0\} \\ \{0,0,0\} & 0 \end{pmatrix} & \begin{pmatrix} 0 & \{0,0,0\} \\ \{0,0,0\} & 0 \end{pmatrix} \\ \begin{pmatrix} 0 & \{0,0,0\} \\ \{0,0,0\} & 0 \end{pmatrix} & \begin{pmatrix} 0 & \{0,0,0\} \\ \{0,0,0\} & 0 \end{pmatrix} & \begin{pmatrix} 0 & \{\lambda_5-\beta_6,0,0\} \\ \{0,0,0\} & 0 \end{pmatrix} \\ \begin{pmatrix} 0 & \{0,0,0\} \\ \{0,0,0\} & 0 \end{pmatrix} & \begin{pmatrix} 0 & \{\beta_6-\lambda_5,0,0\} \\ \{0,0,0\} & 0 \end{pmatrix} & \begin{pmatrix} 0 & \{0,0,0\} \\ \{0,0,0\} & 0 \end{pmatrix} \end{pmatrix}$$

**Δ[X₂[e₃]] = Δ[X₂[e₃]] //. {λ₅ → β₆}**

$$\begin{pmatrix} \begin{pmatrix} 0 & \{0,0,0\} \\ \{0,0,0\} & 0 \end{pmatrix} & \begin{pmatrix} 0 & \{\lambda_6,0,0\} \\ \{0,\lambda_{10},\lambda_{11}\} & \beta_6 \end{pmatrix} & \begin{pmatrix} 0 & \{\lambda_{22},0,0\} \\ \{0,\lambda_{26},\lambda_{27}\} & \lambda_{21} \end{pmatrix} \\ \begin{pmatrix} \beta_6 & \{-\lambda_6,0,0\} \\ \{0,-\lambda_{10},-\lambda_{11}\} & 0 \end{pmatrix} & \begin{pmatrix} \gamma_6 & \{0,0,0\} \\ \{0,0,0\} & \gamma_6 \end{pmatrix} & \begin{pmatrix} \lambda_{12} & \{\lambda_{14},\lambda_{15},\lambda_{16}\} \\ \{0,\lambda_{18},\lambda_{19}\} & \lambda_{13} \end{pmatrix} \\ \begin{pmatrix} \lambda_{21} & \{-\lambda_{22},0,0\} \\ \{0,-\lambda_{26},-\lambda_{27}\} & 0 \end{pmatrix} & \begin{pmatrix} \lambda_{13} & \{-\lambda_{14},-\lambda_{15},-\lambda_{16}\} \\ \{0,-\lambda_{18},-\lambda_{19}\} & \lambda_{12} \end{pmatrix} & \begin{pmatrix} -\gamma_6 & \{0,0,0\} \\ \{0,0,0\} & -\gamma_6 \end{pmatrix} \end{pmatrix}$$

**U_{X₂[e₃]}[X₁[e₄]] == ZERO**

True

**defin[X$_2$[e$_3$], X$_1$[e$_4$]] // Expand**

$$\left(\begin{array}{ccc} \begin{pmatrix} 0 & \{0,0,0\} \\ \{0,0,0\} & 0 \end{pmatrix} & \begin{pmatrix} 0 & \{0,0,0\} \\ \{0,0,0\} & 0 \end{pmatrix} & \begin{pmatrix} 0 & \{0,0,0\} \\ \{0,0,0\} & 0 \end{pmatrix} \\ \begin{pmatrix} 0 & \{0,0,0\} \\ \{0,0,0\} & 0 \end{pmatrix} & \begin{pmatrix} 0 & \{0,0,0\} \\ \{0,0,0\} & 0 \end{pmatrix} & \begin{pmatrix} 0 & \{\beta_5-\lambda_{10},0,0\} \\ \{0,0,0\} & 0 \end{pmatrix} \\ \begin{pmatrix} 0 & \{0,0,0\} \\ \{0,0,0\} & 0 \end{pmatrix} & \begin{pmatrix} 0 & \{\lambda_{10}-\beta_5,0,0\} \\ \{0,0,0\} & 0 \end{pmatrix} & \begin{pmatrix} 0 & \{0,0,0\} \\ \{0,0,0\} & 0 \end{pmatrix} \end{array}\right)$$

**Δ[X$_2$[e$_3$]] = Δ[X$_2$[e$_3$]] //. {λ$_{10}$ → β$_5$}**

$$\left(\begin{array}{ccc} \begin{pmatrix} 0 & \{0,0,0\} \\ \{0,0,0\} & 0 \end{pmatrix} & \begin{pmatrix} 0 & \{\lambda_6,0,0\} \\ \{0,\beta_5,\lambda_{11}\} & \beta_6 \end{pmatrix} & \begin{pmatrix} 0 & \{\lambda_{22},0,0\} \\ \{0,\lambda_{26},\lambda_{27}\} & \lambda_{21} \end{pmatrix} \\ \begin{pmatrix} \beta_6 & \{-\lambda_6,0,0\} \\ \{0,-\beta_5,-\lambda_{11}\} & 0 \end{pmatrix} & \begin{pmatrix} \gamma_6 & \{0,0,0\} \\ \{0,0,0\} & \gamma_6 \end{pmatrix} & \begin{pmatrix} \lambda_{12} & \{\lambda_{14},\lambda_{15},\lambda_{16}\} \\ \{0,\lambda_{18},\lambda_{19}\} & \lambda_{13} \end{pmatrix} \\ \begin{pmatrix} \lambda_{21} & \{-\lambda_{22},0,0\} \\ \{0,-\lambda_{26},-\lambda_{27}\} & 0 \end{pmatrix} & \begin{pmatrix} \lambda_{13} & \{-\lambda_{14},-\lambda_{15},-\lambda_{16}\} \\ \{0,-\lambda_{18},-\lambda_{19}\} & \lambda_{12} \end{pmatrix} & \begin{pmatrix} -\gamma_6 & \{0,0,0\} \\ \{0,0,0\} & -\gamma_6 \end{pmatrix} \end{array}\right)$$

**U$_{X_2[e_3]}$[X$_1$[e$_5$]] == ZERO**

True

**defin[X$_2$[e$_3$], X$_1$[e$_5$]] // Expand**

$$\left(\begin{array}{ccc} \begin{pmatrix} 0 & \{0,0,0\} \\ \{0,0,0\} & 0 \end{pmatrix} & \begin{pmatrix} 0 & \{0,0,0\} \\ \{0,0,0\} & 0 \end{pmatrix} & \begin{pmatrix} 0 & \{0,0,0\} \\ \{0,0,0\} & 0 \end{pmatrix} \\ \begin{pmatrix} 0 & \{0,0,0\} \\ \{0,0,0\} & 0 \end{pmatrix} & \begin{pmatrix} 0 & \{0,0,0\} \\ \{0,0,0\} & 0 \end{pmatrix} & \begin{pmatrix} 0 & \{-\beta_4-\lambda_{11},0,0\} \\ \{0,0,0\} & 0 \end{pmatrix} \\ \begin{pmatrix} 0 & \{0,0,0\} \\ \{0,0,0\} & 0 \end{pmatrix} & \begin{pmatrix} 0 & \{\beta_4+\lambda_{11},0,0\} \\ \{0,0,0\} & 0 \end{pmatrix} & \begin{pmatrix} 0 & \{0,0,0\} \\ \{0,0,0\} & 0 \end{pmatrix} \end{array}\right)$$

**Δ[X$_2$[e$_3$]] = Δ[X$_2$[e$_3$]] //. {λ$_{11}$ → −β$_4$}**

$$\left(\begin{array}{ccc} \begin{pmatrix} 0 & \{0,0,0\} \\ \{0,0,0\} & 0 \end{pmatrix} & \begin{pmatrix} 0 & \{\lambda_6,0,0\} \\ \{0,\beta_5,-\beta_4\} & \beta_6 \end{pmatrix} & \begin{pmatrix} 0 & \{\lambda_{22},0,0\} \\ \{0,\lambda_{26},\lambda_{27}\} & \lambda_{21} \end{pmatrix} \\ \begin{pmatrix} \beta_6 & \{-\lambda_6,0,0\} \\ \{0,-\beta_5,\beta_4\} & 0 \end{pmatrix} & \begin{pmatrix} \gamma_6 & \{0,0,0\} \\ \{0,0,0\} & \gamma_6 \end{pmatrix} & \begin{pmatrix} \lambda_{12} & \{\lambda_{14},\lambda_{15},\lambda_{16}\} \\ \{0,\lambda_{18},\lambda_{19}\} & \lambda_{13} \end{pmatrix} \\ \begin{pmatrix} \lambda_{21} & \{-\lambda_{22},0,0\} \\ \{0,-\lambda_{26},-\lambda_{27}\} & 0 \end{pmatrix} & \begin{pmatrix} \lambda_{13} & \{-\lambda_{14},-\lambda_{15},-\lambda_{16}\} \\ \{0,-\lambda_{18},-\lambda_{19}\} & \lambda_{12} \end{pmatrix} & \begin{pmatrix} -\gamma_6 & \{0,0,0\} \\ \{0,0,0\} & -\gamma_6 \end{pmatrix} \end{array}\right)$$

**U$_{X_2[e_3]}$[X$_1$[e$_6$]] == ZERO**

True

**(defin[X₂[e₃], X₁[e₆]] // Expand)**

$$\begin{pmatrix} \begin{pmatrix} 0 & \{0,0,0\} \\ \{0,0,0\} & 0 \end{pmatrix} & \begin{pmatrix} 0 & \{0,0,0\} \\ \{0,0,0\} & 0 \end{pmatrix} & \begin{pmatrix} 0 & \{0,0,0\} \\ \{0,0,0\} & 0 \end{pmatrix} \\ \begin{pmatrix} 0 & \{0,0,0\} \\ \{0,0,0\} & 0 \end{pmatrix} & \begin{pmatrix} 0 & \{0,0,0\} \\ \{0,0,0\} & 0 \end{pmatrix} & \begin{pmatrix} 0 & \{\beta_1-\lambda_6,0,0\} \\ \{0,0,0\} & 0 \end{pmatrix} \\ \begin{pmatrix} 0 & \{0,0,0\} \\ \{0,0,0\} & 0 \end{pmatrix} & \begin{pmatrix} 0 & \{\lambda_6-\beta_1,0,0\} \\ \{0,0,0\} & 0 \end{pmatrix} & \begin{pmatrix} 0 & \{0,0,0\} \\ \{0,0,0\} & 0 \end{pmatrix} \end{pmatrix}$$

**Δ[X₂[e₃]] = Δ[X₂[e₃]] //. {λ₆ → β₁}**

$$\begin{pmatrix} \begin{pmatrix} 0 & \{0,0,0\} \\ \{0,0,0\} & 0 \end{pmatrix} & \begin{pmatrix} 0 & \{\beta_1,0,0\} \\ \{0,\beta_5,-\beta_4\} & \beta_6 \end{pmatrix} & \begin{pmatrix} 0 & \{\lambda_{22},0,0\} \\ \{0,\lambda_{26},\lambda_{27}\} & \lambda_{21} \end{pmatrix} \\ \begin{pmatrix} \beta_6 & \{-\beta_1,0,0\} \\ \{0,-\beta_5,\beta_4\} & 0 \end{pmatrix} & \begin{pmatrix} \gamma_6 & \{0,0,0\} \\ \{0,0,0\} & \gamma_6 \end{pmatrix} & \begin{pmatrix} \lambda_{12} & \{\lambda_{14},\lambda_{15},\lambda_{16}\} \\ \{0,\lambda_{18},\lambda_{19}\} & \lambda_{13} \end{pmatrix} \\ \begin{pmatrix} \lambda_{21} & \{-\lambda_{22},0,0\} \\ \{0,-\lambda_{26},-\lambda_{27}\} & 0 \end{pmatrix} & \begin{pmatrix} \lambda_{13} & \{-\lambda_{14},-\lambda_{15},-\lambda_{16}\} \\ \{0,-\lambda_{18},-\lambda_{19}\} & \lambda_{12} \end{pmatrix} & \begin{pmatrix} -\gamma_6 & \{0,0,0\} \\ \{0,0,0\} & -\gamma_6 \end{pmatrix} \end{pmatrix}$$

**U_{X₂[e₁]}[X₂[e₃]] == ZERO**

True

**defin[X₂[e₁], X₂[e₃]] // Expand**

$$\begin{pmatrix} \begin{pmatrix} 0 & \{0,0,0\} \\ \{0,0,0\} & 0 \end{pmatrix} & \begin{pmatrix} 0 & \{0,0,0\} \\ \{0,0,0\} & 0 \end{pmatrix} & \begin{pmatrix} 0 & \{0,0,0\} \\ \{0,0,0\} & 0 \end{pmatrix} \\ \begin{pmatrix} 0 & \{0,0,0\} \\ \{0,0,0\} & 0 \end{pmatrix} & \begin{pmatrix} 0 & \{0,0,0\} \\ \{0,0,0\} & 0 \end{pmatrix} & \begin{pmatrix} \lambda_{13}+\rho_4 & \{0,0,0\} \\ \{0,0,0\} & 0 \end{pmatrix} \\ \begin{pmatrix} 0 & \{0,0,0\} \\ \{0,0,0\} & 0 \end{pmatrix} & \begin{pmatrix} 0 & \{0,0,0\} \\ \{0,0,0\} & \lambda_{13}+\rho_4 \end{pmatrix} & \begin{pmatrix} 0 & \{0,0,0\} \\ \{0,0,0\} & 0 \end{pmatrix} \end{pmatrix}$$

**Δ[X₂[e₃]] = Δ[X₂[e₃]] //. {λ₁₃ → -ρ₄}**

$$\begin{pmatrix} \begin{pmatrix} 0 & \{0,0,0\} \\ \{0,0,0\} & 0 \end{pmatrix} & \begin{pmatrix} 0 & \{\beta_1,0,0\} \\ \{0,\beta_5,-\beta_4\} & \beta_6 \end{pmatrix} & \begin{pmatrix} 0 & \{\lambda_{22},0,0\} \\ \{0,\lambda_{26},\lambda_{27}\} & \lambda_{21} \end{pmatrix} \\ \begin{pmatrix} \beta_6 & \{-\beta_1,0,0\} \\ \{0,-\beta_5,\beta_4\} & 0 \end{pmatrix} & \begin{pmatrix} \gamma_6 & \{0,0,0\} \\ \{0,0,0\} & \gamma_6 \end{pmatrix} & \begin{pmatrix} \lambda_{12} & \{\lambda_{14},\lambda_{15},\lambda_{16}\} \\ \{0,\lambda_{18},\lambda_{19}\} & -\rho_4 \end{pmatrix} \\ \begin{pmatrix} \lambda_{21} & \{-\lambda_{22},0,0\} \\ \{0,-\lambda_{26},-\lambda_{27}\} & 0 \end{pmatrix} & \begin{pmatrix} -\rho_4 & \{-\lambda_{14},-\lambda_{15},-\lambda_{16}\} \\ \{0,-\lambda_{18},-\lambda_{19}\} & \lambda_{12} \end{pmatrix} & \begin{pmatrix} -\gamma_6 & \{0,0,0\} \\ \{0,0,0\} & -\gamma_6 \end{pmatrix} \end{pmatrix}$$

**U_{X₂[e₂]}[X₂[e₃]] == ZERO**

True

**defin[X₂[e₂], X₂[e₃]] // Expand**

$$\begin{pmatrix} \begin{pmatrix} 0 & \{0,0,0\} \\ \{0,0,0\} & 0 \end{pmatrix} & \begin{pmatrix} 0 & \{0,0,0\} \\ \{0,0,0\} & 0 \end{pmatrix} & \begin{pmatrix} 0 & \{0,0,0\} \\ \{0,0,0\} & 0 \end{pmatrix} \\ \begin{pmatrix} 0 & \{0,0,0\} \\ \{0,0,0\} & 0 \end{pmatrix} & \begin{pmatrix} 0 & \{0,0,0\} \\ \{0,0,0\} & 0 \end{pmatrix} & \begin{pmatrix} 0 & \{0,0,0\} \\ \{0,0,0\} & \lambda_{12}+\chi_1 \end{pmatrix} \\ \begin{pmatrix} 0 & \{0,0,0\} \\ \{0,0,0\} & 0 \end{pmatrix} & \begin{pmatrix} \lambda_{12}+\chi_1 & \{0,0,0\} \\ \{0,0,0\} & 0 \end{pmatrix} & \begin{pmatrix} 0 & \{0,0,0\} \\ \{0,0,0\} & 0 \end{pmatrix} \end{pmatrix}$$

**Δ[X₂[e₃]] = Δ[X₂[e₃]] //. {λ₁₂ → −χ₁}**

$$\begin{pmatrix} \begin{pmatrix} 0 & \{0,0,0\} \\ \{0,0,0\} & 0 \end{pmatrix} & \begin{pmatrix} 0 & \{\beta_1,0,0\} \\ \{0,\beta_5,-\beta_4\} & \beta_6 \end{pmatrix} & \begin{pmatrix} 0 & \{\lambda_{22},0,0\} \\ \{0,\lambda_{26},\lambda_{27}\} & \lambda_{21} \end{pmatrix} \\ \begin{pmatrix} \beta_6 & \{-\beta_1,0,0\} \\ \{0,-\beta_5,\beta_4\} & 0 \end{pmatrix} & \begin{pmatrix} \gamma_6 & \{0,0,0\} \\ \{0,0,0\} & \gamma_6 \end{pmatrix} & \begin{pmatrix} -\chi_1 & \{\lambda_{14},\lambda_{15},\lambda_{16}\} \\ \{0,\lambda_{18},\lambda_{19}\} & -\rho_4 \end{pmatrix} \\ \begin{pmatrix} \lambda_{21} & \{-\lambda_{22},0,0\} \\ \{0,-\lambda_{26},-\lambda_{27}\} & 0 \end{pmatrix} & \begin{pmatrix} -\rho_4 & \{-\lambda_{14},-\lambda_{15},-\lambda_{16}\} \\ \{0,-\lambda_{18},-\lambda_{19}\} & -\chi_1 \end{pmatrix} & \begin{pmatrix} -\gamma_6 & \{0,0,0\} \\ \{0,0,0\} & -\gamma_6 \end{pmatrix} \end{pmatrix}$$

**T[E₁[1], E₁[1], X₂[e₃]] == ZERO**

True

**Leib[E₁[1], E₁[1], X₂[e₃]] // Expand**

$$\begin{pmatrix} \begin{pmatrix} 0 & \{0,0,0\} \\ \{0,0,0\} & 0 \end{pmatrix} & \begin{pmatrix} 0 & \{0,0,0\} \\ \{0,0,0\} & 0 \end{pmatrix} & \begin{pmatrix} 0 & \{\alpha_1+\lambda_{22},0,0\} \\ \{0,\alpha_5+\lambda_{26},\lambda_{27}-\alpha_4\} & \alpha_6+\lambda_{21} \end{pmatrix} \\ \begin{pmatrix} 0 & \{0,0,0\} \\ \{0,0,0\} & 0 \end{pmatrix} & \begin{pmatrix} 0 & \{0,0,0\} \\ \{0,0,0\} & 0 \end{pmatrix} & \begin{pmatrix} 0 & \{0,0,0\} \\ \{0,0,0\} & 0 \end{pmatrix} \\ \begin{pmatrix} \alpha_6+\lambda_{21} & \{-\alpha_1-\lambda_{22},0,0\} \\ \{0,-\alpha_5-\lambda_{26},\alpha_4-\lambda_{27}\} & 0 \end{pmatrix} & \begin{pmatrix} 0 & \{0,0,0\} \\ \{0,0,0\} & 0 \end{pmatrix} & \begin{pmatrix} 0 & \{0,0,0\} \\ \{0,0,0\} & 0 \end{pmatrix} \end{pmatrix}$$

**Δ[X₂[e₃]] = Δ[X₂[e₃]] //. {λ₂₁ → −α₆, λ₂₂ → −α₁, λ₂₆ → −α₅, λ₂₇ → α₄}**

$$\begin{pmatrix} \begin{pmatrix} 0 & \{0,0,0\} \\ \{0,0,0\} & 0 \end{pmatrix} & \begin{pmatrix} 0 & \{\beta_1,0,0\} \\ \{0,\beta_5,-\beta_4\} & \beta_6 \end{pmatrix} & \begin{pmatrix} 0 & \{-\alpha_1,0,0\} \\ \{0,-\alpha_5,\alpha_4\} & -\alpha_6 \end{pmatrix} \\ \begin{pmatrix} \beta_6 & \{-\beta_1,0,0\} \\ \{0,-\beta_5,\beta_4\} & 0 \end{pmatrix} & \begin{pmatrix} \gamma_6 & \{0,0,0\} \\ \{0,0,0\} & \gamma_6 \end{pmatrix} & \begin{pmatrix} -\chi_1 & \{\lambda_{14},\lambda_{15},\lambda_{16}\} \\ \{0,\lambda_{18},\lambda_{19}\} & -\rho_4 \end{pmatrix} \\ \begin{pmatrix} -\alpha_6 & \{\alpha_1,0,0\} \\ \{0,\alpha_5,-\alpha_4\} & 0 \end{pmatrix} & \begin{pmatrix} -\rho_4 & \{-\lambda_{14},-\lambda_{15},-\lambda_{16}\} \\ \{0,-\lambda_{18},-\lambda_{19}\} & -\chi_1 \end{pmatrix} & \begin{pmatrix} -\gamma_6 & \{0,0,0\} \\ \{0,0,0\} & -\gamma_6 \end{pmatrix} \end{pmatrix}$$

**T[X₁[e₁], X₁[e₂], X₂[e₃]] == ZERO**

True

**Leib[X₁[e₁], X₁[e₂], X₂[e₃]] // Expand**

$$\begin{pmatrix} \begin{pmatrix} 0 & \{0,0,0\} \\ \{0,0,0\} & 0 \end{pmatrix} & \begin{pmatrix} 0 & \{0,0,0\} \\ \{0,0,0\} & 0 \end{pmatrix} & \begin{pmatrix} 0 & \{0,0,0\} \\ \{0,0,0\} & 0 \end{pmatrix} \\ \begin{pmatrix} 0 & \{0,0,0\} \\ \{0,0,0\} & 0 \end{pmatrix} & \begin{pmatrix} 0 & \{0,0,0\} \\ \{0,0,0\} & 0 \end{pmatrix} & \begin{pmatrix} 0 & \{0,0,0\} \\ \{0,\lambda_{18}+\rho_3,\lambda_{19}-\rho_2\} & 0 \end{pmatrix} \\ \begin{pmatrix} 0 & \{0,0,0\} \\ \{0,0,0\} & 0 \end{pmatrix} & \begin{pmatrix} 0 & \{0,0,0\} \\ \{0,-\lambda_{18}-\rho_3,\rho_2-\lambda_{19}\} & 0 \end{pmatrix} & \begin{pmatrix} 0 & \{0,0,0\} \\ \{0,0,0\} & 0 \end{pmatrix} \end{pmatrix}$$

$\Delta[X_2[e_3]] = \Delta[X_2[e_3]] //. \{\lambda_{18} \to -\rho_3, \lambda_{19} \to \rho_2\}$

$$\begin{pmatrix} \begin{pmatrix} 0 & \{0,0,0\} \\ \{0,0,0\} & 0 \end{pmatrix} & \begin{pmatrix} 0 & \{\beta_1,0,0\} \\ \{0,\beta_5,-\beta_4\} & \beta_6 \end{pmatrix} & \begin{pmatrix} 0 & \{-\alpha_1,0,0\} \\ \{0,-\alpha_5,\alpha_4\} & -\alpha_6 \end{pmatrix} \\ \begin{pmatrix} \beta_6 & \{-\beta_1,0,0\} \\ \{0,-\beta_5,\beta_4\} & 0 \end{pmatrix} & \begin{pmatrix} \gamma_6 & \{0,0,0\} \\ \{0,0,0\} & \gamma_6 \end{pmatrix} & \begin{pmatrix} -\chi_1 & \{\lambda_{14},\lambda_{15},\lambda_{16}\} \\ \{0,-\rho_3,\rho_2\} & -\rho_4 \end{pmatrix} \\ \begin{pmatrix} -\alpha_6 & \{\alpha_1,0,0\} \\ \{0,\alpha_5,-\alpha_4\} & 0 \end{pmatrix} & \begin{pmatrix} -\rho_4 & \{-\lambda_{14},-\lambda_{15},-\lambda_{16}\} \\ \{0,\rho_3,-\rho_2\} & -\chi_1 \end{pmatrix} & \begin{pmatrix} -\gamma_6 & \{0,0,0\} \\ \{0,0,0\} & -\gamma_6 \end{pmatrix} \end{pmatrix}$$

**T[X$_2$[e$_1$], X$_2$[e$_3$], X$_1$[e$_7$]] == ZERO**

True

**Leib[X$_2$[e$_1$], X$_2$[e$_3$], X$_1$[e$_7$]] // Expand**

$$\begin{pmatrix} \begin{pmatrix} 0 & \{0,0,0\} \\ \{0,0,0\} & 0 \end{pmatrix} & \begin{pmatrix} 0 & \{0,0,0\} \\ \{0,0,0\} & \lambda_{15}-\eta_2 \end{pmatrix} & \begin{pmatrix} 0 & \{0,0,0\} \\ \{0,0,0\} & 0 \end{pmatrix} \\ \begin{pmatrix} \lambda_{15}-\eta_2 & \{0,0,0\} \\ \{0,0,0\} & 0 \end{pmatrix} & \begin{pmatrix} 0 & \{0,0,0\} \\ \{0,0,0\} & 0 \end{pmatrix} & \begin{pmatrix} 0 & \{0,0,0\} \\ \{0,0,0\} & 0 \end{pmatrix} \\ \begin{pmatrix} 0 & \{0,0,0\} \\ \{0,0,0\} & 0 \end{pmatrix} & \begin{pmatrix} 0 & \{0,0,0\} \\ \{0,0,0\} & 0 \end{pmatrix} & \begin{pmatrix} 0 & \{0,0,0\} \\ \{0,0,0\} & 0 \end{pmatrix} \end{pmatrix}$$

$\Delta[X_2[e_3]] = \Delta[X_2[e_3]] //. \{\lambda_{15} \to \eta_2\}$

$$\begin{pmatrix} \begin{pmatrix} 0 & \{0,0,0\} \\ \{0,0,0\} & 0 \end{pmatrix} & \begin{pmatrix} 0 & \{\beta_1,0,0\} \\ \{0,\beta_5,-\beta_4\} & \beta_6 \end{pmatrix} & \begin{pmatrix} 0 & \{-\alpha_1,0,0\} \\ \{0,-\alpha_5,\alpha_4\} & -\alpha_6 \end{pmatrix} \\ \begin{pmatrix} \beta_6 & \{-\beta_1,0,0\} \\ \{0,-\beta_5,\beta_4\} & 0 \end{pmatrix} & \begin{pmatrix} \gamma_6 & \{0,0,0\} \\ \{0,0,0\} & \gamma_6 \end{pmatrix} & \begin{pmatrix} -\chi_1 & \{\lambda_{14},\eta_2,\lambda_{16}\} \\ \{0,-\rho_3,\rho_2\} & -\rho_4 \end{pmatrix} \\ \begin{pmatrix} -\alpha_6 & \{\alpha_1,0,0\} \\ \{0,\alpha_5,-\alpha_4\} & 0 \end{pmatrix} & \begin{pmatrix} -\rho_4 & \{-\lambda_{14},-\eta_2,-\lambda_{16}\} \\ \{0,\rho_3,-\rho_2\} & -\chi_1 \end{pmatrix} & \begin{pmatrix} -\gamma_6 & \{0,0,0\} \\ \{0,0,0\} & -\gamma_6 \end{pmatrix} \end{pmatrix}$$

**T[X$_2$[e$_2$], X$_2$[e$_3$], X$_1$[e$_3$]] == ZERO**

True

**Leib[X$_2$[e$_2$], X$_2$[e$_3$], X$_1$[e$_3$]] // Expand**

$$\begin{pmatrix} \begin{pmatrix} 0 & \{0,0,0\} \\ \{0,0,0\} & 0 \end{pmatrix} & \begin{pmatrix} 0 & \{0,0,0\} \\ \{0,\eta_3-\lambda_{16},0\} & 0 \end{pmatrix} & \begin{pmatrix} 0 & \{0,0,0\} \\ \{0,0,0\} & 0 \end{pmatrix} \\ \begin{pmatrix} 0 & \{0,0,0\} \\ \{0,\lambda_{16}-\eta_3,0\} & 0 \end{pmatrix} & \begin{pmatrix} 0 & \{0,0,0\} \\ \{0,0,0\} & 0 \end{pmatrix} & \begin{pmatrix} 0 & \{0,0,0\} \\ \{0,0,0\} & 0 \end{pmatrix} \\ \begin{pmatrix} 0 & \{0,0,0\} \\ \{0,0,0\} & 0 \end{pmatrix} & \begin{pmatrix} 0 & \{0,0,0\} \\ \{0,0,0\} & 0 \end{pmatrix} & \begin{pmatrix} 0 & \{0,0,0\} \\ \{0,0,0\} & 0 \end{pmatrix} \end{pmatrix}$$

$\Delta[X_2[e_3]] = \Delta[X_2[e_3]] \,//\, \{\lambda_{16} \to \eta_3\}$

$$\left(\begin{pmatrix} 0 & \{0,0,0\} \\ \{0,0,0\} & 0 \end{pmatrix} \begin{pmatrix} 0 & \{\beta_1,0,0\} \\ \{0,\beta_5,-\beta_4\} & \beta_6 \end{pmatrix} \begin{pmatrix} 0 & \{-\alpha_1,0,0\} \\ \{0,-\alpha_5,\alpha_4\} & -\alpha_6 \end{pmatrix} \\ \begin{pmatrix} \beta_6 & \{-\beta_1,0,0\} \\ \{0,-\beta_5,\beta_4\} & 0 \end{pmatrix} \begin{pmatrix} \gamma_6 & \{0,0,0\} \\ \{0,0,0\} & \gamma_6 \end{pmatrix} \begin{pmatrix} -\chi_1 & \{\lambda_{14},\eta_2,\eta_3\} \\ \{0,-\rho_3,\rho_2\} & -\rho_4 \end{pmatrix} \\ \begin{pmatrix} -\alpha_6 & \{\alpha_1,0,0\} \\ \{0,\alpha_5,-\alpha_4\} & 0 \end{pmatrix} \begin{pmatrix} -\rho_4 & \{-\lambda_{14},-\eta_2,-\eta_3\} \\ \{0,\rho_3,-\rho_2\} & -\chi_1 \end{pmatrix} \begin{pmatrix} -\gamma_6 & \{0,0,0\} \\ \{0,0,0\} & -\gamma_6 \end{pmatrix} \right)$$

$T[X_2[e_3], X_2[e_1], X_1[e_1]] == -X_1[e_3]$

True

$\Delta[X_1[e_3]] + \text{Leib}[X_2[e_3], X_2[e_1], X_1[e_1]] \,//\, \text{Expand}$

$$\left(\begin{pmatrix} 0 & \{0,0,0\} \\ \{0,0,0\} & 0 \end{pmatrix} \begin{pmatrix} 0 & \{-\xi-\delta_1+\eta_1-\lambda_{14},0,0\} \\ \{0,0,0\} & 0 \end{pmatrix} \begin{pmatrix} 0 & \{0,0,0\} \\ \{0,0,0\} & 0 \end{pmatrix} \\ \begin{pmatrix} 0 & \{\xi+\delta_1-\eta_1+\lambda_{14},0,0\} \\ \{0,0,0\} & 0 \end{pmatrix} \begin{pmatrix} 0 & \{0,0,0\} \\ \{0,0,0\} & 0 \end{pmatrix} \begin{pmatrix} 0 & \{0,0,0\} \\ \{0,0,0\} & 0 \end{pmatrix} \\ \begin{pmatrix} 0 & \{0,0,0\} \\ \{0,0,0\} & 0 \end{pmatrix} \begin{pmatrix} 0 & \{0,0,0\} \\ \{0,0,0\} & 0 \end{pmatrix} \begin{pmatrix} 0 & \{0,0,0\} \\ \{0,0,0\} & 0 \end{pmatrix} \right)$$

$\Delta[X_2[e_3]] = \Delta[X_2[e_3]] \,//\, \{\lambda_{14} \to -\xi - \delta_1 + \eta_1\}$

$$\left(\begin{pmatrix} 0 & \{0,0,0\} \\ \{0,0,0\} & 0 \end{pmatrix} \begin{pmatrix} 0 & \{\beta_1,0,0\} \\ \{0,\beta_5,-\beta_4\} & \beta_6 \end{pmatrix} \begin{pmatrix} 0 & \{-\alpha_1,0,0\} \\ \{0,-\alpha_5,\alpha_4\} & -\alpha_6 \end{pmatrix} \\ \begin{pmatrix} \beta_6 & \{-\beta_1,0,0\} \\ \{0,-\beta_5,\beta_4\} & 0 \end{pmatrix} \begin{pmatrix} \gamma_6 & \{0,0,0\} \\ \{0,0,0\} & \gamma_6 \end{pmatrix} \begin{pmatrix} -\chi_1 & \{-\xi-\delta_1+\eta_1,\eta_2,\eta_3\} \\ \{0,-\rho_3,\rho_2\} & -\rho_4 \end{pmatrix} \\ \begin{pmatrix} -\alpha_6 & \{\alpha_1,0,0\} \\ \{0,\alpha_5,-\alpha_4\} & 0 \end{pmatrix} \begin{pmatrix} -\rho_4 & \{\xi+\delta_1-\eta_1,-\eta_2,-\eta_3\} \\ \{0,\rho_3,-\rho_2\} & -\chi_1 \end{pmatrix} \begin{pmatrix} -\gamma_6 & \{0,0,0\} \\ \{0,0,0\} & -\gamma_6 \end{pmatrix} \right)$$

$\text{Variables}[\Delta[X_2[e_3]]]$

$\{\xi, \alpha_1, \alpha_4, \alpha_5, \alpha_6, \beta_1, \beta_4, \beta_5, \beta_6, \gamma_6, \delta_1, \eta_1, \eta_2, \eta_3, \rho_2, \rho_3, \rho_4, \chi_1\}$

- **Eliminating the extra parameter**

   $T[X_1[e_5], X_1[e_4], X_2[e_3]] == X_2[e_1]$

   True

   $\Delta[X_2[e_1]] - \text{Leib}[X_1[e_5], X_1[e_4], X_2[e_3]] \,//\, \text{Expand}$

$$\left(\begin{pmatrix} 0 & \{0,0,0\} \\ \{0,0,0\} & 0 \end{pmatrix} \begin{pmatrix} 0 & \{0,0,0\} \\ \{0,0,0\} & 0 \end{pmatrix} \begin{pmatrix} 0 & \{0,0,0\} \\ \{0,0,0\} & 0 \end{pmatrix} \\ \begin{pmatrix} 0 & \{0,0,0\} \\ \{0,0,0\} & 0 \end{pmatrix} \begin{pmatrix} 0 & \{0,0,0\} \\ \{0,0,0\} & 0 \end{pmatrix} \begin{pmatrix} 2\xi+\delta_1-\epsilon_2-\eta_1-\phi_3 & \{0,0,0\} \\ \{0,0,0\} & 0 \end{pmatrix} \\ \begin{pmatrix} 0 & \{0,0,0\} \\ \{0,0,0\} & 0 \end{pmatrix} \begin{pmatrix} 0 & \{0,0,0\} \\ \{0,0,0\} & 2\xi+\delta_1-\epsilon_2-\eta_1-\phi_3 \end{pmatrix} \begin{pmatrix} 0 & \{0,0,0\} \\ \{0,0,0\} & 0 \end{pmatrix} \right)$$

```
Do[
    Δ[X₁[eᵢ]] = Δ[X₁[eᵢ]] //. {ϕ₃ → 2 ξ + δ₁ - ϵ₂ - η₁}, {i, 1, 8}]

Δ[X₂[e₁]] = Δ[X₂[e₁]] //. {ϕ₃ → 2 ξ + δ₁ - ϵ₂ - η₁};
Δ[X₂[e₂]] = Δ[X₂[e₂]] //. {ϕ₃ → 2 ξ + δ₁ - ϵ₂ - η₁};
Δ[X₂[e₃]] = Δ[X₂[e₃]] //. {ϕ₃ → 2 ξ + δ₁ - ϵ₂ - η₁};

Union[Variables[Δ[X₁[e₁]]], Variables[Δ[X₁[e₂]]],
 Variables[Δ[X₁[e₃]]], Variables[Δ[X₁[e₄]]],
 Variables[Δ[X₁[e₅]]], Variables[Δ[X₁[e₆]]],
 Variables[Δ[X₁[e₇]]], Variables[Δ[X₁[e₈]]],
 Variables[Δ[E₁[1]]], Variables[Δ[E₂[1]]],
 Variables[Δ[E₃[1]]], Variables[Δ[X₂[e₁]]],
 Variables[Δ[X₂[e₂]]], Variables[Δ[X₂[e₃]]]]
```

$\{\xi, \alpha_1, \alpha_2, \alpha_3, \alpha_4, \alpha_5, \alpha_6, \alpha_7, \alpha_8, \beta_1, \beta_2, \beta_3, \beta_4, \beta_5, \beta_6, \beta_7,$
$\beta_8, \gamma_1, \gamma_2, \gamma_3, \gamma_4, \gamma_5, \gamma_6, \gamma_7, \gamma_8, \delta_1, \delta_2, \delta_3, \delta_4, \delta_5, \delta_6, \delta_7, \epsilon_1, \epsilon_2,$
$\epsilon_3, \epsilon_4, \eta_1, \eta_2, \eta_3, \eta_4, \eta_5, \rho_1, \rho_2, \rho_3, \rho_4, \rho_5, \rho_6, \phi_1, \phi_2, \chi_1, \psi_1, \psi_2\}$

```
Length[%]
```

- $X_2[e_4]$

```
Δ[X₂[e₄]] = generic;

U_{E₁[1]}[X₂[e₄]] == ZERO
```

True

```
(defin[E₁[1], X₂[e₄]] // Expand)
```

$$\begin{pmatrix} \begin{pmatrix} \lambda_1 & \{0,0,0\} \\ \{0,0,0\} & \lambda_1 \end{pmatrix} & \begin{pmatrix} 0 & \{0,0,0\} \\ \{0,0,0\} & 0 \end{pmatrix} & \begin{pmatrix} 0 & \{0,0,0\} \\ \{0,0,0\} & 0 \end{pmatrix} \\ \begin{pmatrix} 0 & \{0,0,0\} \\ \{0,0,0\} & 0 \end{pmatrix} & \begin{pmatrix} 0 & \{0,0,0\} \\ \{0,0,0\} & 0 \end{pmatrix} & \begin{pmatrix} 0 & \{0,0,0\} \\ \{0,0,0\} & 0 \end{pmatrix} \\ \begin{pmatrix} 0 & \{0,0,0\} \\ \{0,0,0\} & 0 \end{pmatrix} & \begin{pmatrix} 0 & \{0,0,0\} \\ \{0,0,0\} & 0 \end{pmatrix} & \begin{pmatrix} 0 & \{0,0,0\} \\ \{0,0,0\} & 0 \end{pmatrix} \end{pmatrix}$$

```
Δ[X₂[e₄]] = Δ[X₂[e₄]] //. {λ₁ → 0}
```

$$\begin{pmatrix} \begin{pmatrix} 0 & \{0,0,0\} \\ \{0,0,0\} & 0 \end{pmatrix} & \begin{pmatrix} \lambda_4 & \{\lambda_6, \lambda_7, \lambda_8\} \\ \{\lambda_9, \lambda_{10}, \lambda_{11}\} & \lambda_5 \end{pmatrix} & \begin{pmatrix} \lambda_{20} & \{\lambda_{22}, \lambda_{23}, \lambda_{24}\} \\ \{\lambda_{25}, \lambda_{26}, \lambda_{27}\} & \lambda_{21} \end{pmatrix} \\ \begin{pmatrix} \lambda_5 & \{-\lambda_6, -\lambda_7, -\lambda_8\} \\ \{-\lambda_9, -\lambda_{10}, -\lambda_{11}\} & \lambda_4 \end{pmatrix} & \begin{pmatrix} \lambda_2 & \{0,0,0\} \\ \{0,0,0\} & \lambda_2 \end{pmatrix} & \begin{pmatrix} \lambda_{12} & \{\lambda_{14}, \lambda_{15}, \lambda_{16}\} \\ \{\lambda_{17}, \lambda_{18}, \lambda_{19}\} & \lambda_{13} \end{pmatrix} \\ \begin{pmatrix} \lambda_{21} & \{-\lambda_{22}, -\lambda_{23}, -\lambda_{24}\} \\ \{-\lambda_{25}, -\lambda_{26}, -\lambda_{27}\} & \lambda_{20} \end{pmatrix} & \begin{pmatrix} \lambda_{13} & \{-\lambda_{14}, -\lambda_{15}, -\lambda_{16}\} \\ \{-\lambda_{17}, -\lambda_{18}, -\lambda_{19}\} & \lambda_{12} \end{pmatrix} & \begin{pmatrix} \lambda_3 & \{0,0,0\} \\ \{0,0,0\} & \lambda_3 \end{pmatrix} \end{pmatrix}$$

```
U_{X₂[e₄]}[E₂[1]] == ZERO
```

True

**defin[X₂[e₄], E₂[1]] // Expand**

$$\begin{pmatrix} \begin{pmatrix} 0 & \{0,0,0\} \\ \{0,0,0\} & 0 \end{pmatrix} & \begin{pmatrix} 0 & \{0,0,0\} \\ \{0,0,0\} & 0 \end{pmatrix} & \begin{pmatrix} 0 & \{0,\lambda_4,0\} \\ \{-\lambda_8,0,\lambda_6\} & \lambda_{10} \end{pmatrix} \\ \begin{pmatrix} 0 & \{0,0,0\} \\ \{0,0,0\} & 0 \end{pmatrix} & \begin{pmatrix} 0 & \{0,0,0\} \\ \{0,0,0\} & 0 \end{pmatrix} & \begin{pmatrix} 0 & \{0,\lambda_2-\gamma_7,0\} \\ \{0,0,0\} & 0 \end{pmatrix} \\ \begin{pmatrix} \lambda_{10} & \{0,-\lambda_4,0\} \\ \{\lambda_8,0,-\lambda_6\} & 0 \end{pmatrix} & \begin{pmatrix} 0 & \{0,\gamma_7-\lambda_2,0\} \\ \{0,0,0\} & 0 \end{pmatrix} & \begin{pmatrix} -\lambda_{18} & \{0,0,0\} \\ \{0,0,0\} & -\lambda_{18} \end{pmatrix} \end{pmatrix}$$

**Δ[X₂[e₄]] =
  Δ[X₂[e₄]] //. {λ₂ → γ₇, λ₄ → 0, λ₆ → 0, λ₈ → 0, λ₁₀ → 0, λ₁₈ → 0}**

$$\begin{pmatrix} \begin{pmatrix} 0 & \{0,0,0\} \\ \{0,0,0\} & 0 \end{pmatrix} & \begin{pmatrix} 0 & \{0,\lambda_7,0\} \\ \{\lambda_9,0,\lambda_{11}\} & \lambda_5 \end{pmatrix} & \begin{pmatrix} \lambda_{20} & \{\lambda_{22},\lambda_{23},\lambda_{24}\} \\ \{\lambda_{25},\lambda_{26},\lambda_{27}\} & \lambda_{21} \end{pmatrix} \\ \begin{pmatrix} \lambda_5 & \{0,-\lambda_7,0\} \\ \{-\lambda_9,0,-\lambda_{11}\} & 0 \end{pmatrix} & \begin{pmatrix} \gamma_7 & \{0,0,0\} \\ \{0,0,0\} & \gamma_7 \end{pmatrix} & \begin{pmatrix} \lambda_{12} & \{\lambda_{14},\lambda_{15},\lambda_{16}\} \\ \{\lambda_{17},0,\lambda_{19}\} & \lambda_{13} \end{pmatrix} \\ \begin{pmatrix} \lambda_{21} & \{-\lambda_{22},-\lambda_{23},-\lambda_{24}\} \\ \{-\lambda_{25},-\lambda_{26},-\lambda_{27}\} & \lambda_{20} \end{pmatrix} & \begin{pmatrix} \lambda_{13} & \{-\lambda_{14},-\lambda_{15},-\lambda_{16}\} \\ \{-\lambda_{17},0,-\lambda_{19}\} & \lambda_{12} \end{pmatrix} & \begin{pmatrix} \lambda_3 & \{0,0,0\} \\ \{0,0,0\} & \lambda_3 \end{pmatrix} \end{pmatrix}$$

**U_{X₂[e₄]} [E₃[1]] == ZERO**

True

**defin[X₂[e₄], E₃[1]] // Expand**

$$\begin{pmatrix} \begin{pmatrix} 0 & \{0,0,0\} \\ \{0,0,0\} & 0 \end{pmatrix} & \begin{pmatrix} 0 & \{0,-\lambda_{20},0\} \\ \{\lambda_{24},0,-\lambda_{22}\} & -\lambda_{26} \end{pmatrix} & \begin{pmatrix} 0 & \{0,0,0\} \\ \{0,0,0\} & 0 \end{pmatrix} \\ \begin{pmatrix} -\lambda_{26} & \{0,\lambda_{20},0\} \\ \{-\lambda_{24},0,\lambda_{22}\} & 0 \end{pmatrix} & \begin{pmatrix} 0 & \{0,0,0\} \\ \{0,0,0\} & 0 \end{pmatrix} & \begin{pmatrix} 0 & \{0,\gamma_7+\lambda_3,0\} \\ \{0,0,0\} & 0 \end{pmatrix} \\ \begin{pmatrix} 0 & \{0,0,0\} \\ \{0,0,0\} & 0 \end{pmatrix} & \begin{pmatrix} 0 & \{0,-\gamma_7-\lambda_3,0\} \\ \{0,0,0\} & 0 \end{pmatrix} & \begin{pmatrix} 0 & \{0,0,0\} \\ \{0,0,0\} & 0 \end{pmatrix} \end{pmatrix}$$

**Δ[X₂[e₄]] = Δ[X₂[e₄]] //. {λ₃ → −γ₇, λ₂₀ → 0, λ₂₂ → 0, λ₂₄ → 0, λ₂₆ → 0}**

$$\begin{pmatrix} \begin{pmatrix} 0 & \{0,0,0\} \\ \{0,0,0\} & 0 \end{pmatrix} & \begin{pmatrix} 0 & \{0,\lambda_7,0\} \\ \{\lambda_9,0,\lambda_{11}\} & \lambda_5 \end{pmatrix} & \begin{pmatrix} 0 & \{0,\lambda_{23},0\} \\ \{\lambda_{25},0,\lambda_{27}\} & \lambda_{21} \end{pmatrix} \\ \begin{pmatrix} \lambda_5 & \{0,-\lambda_7,0\} \\ \{-\lambda_9,0,-\lambda_{11}\} & 0 \end{pmatrix} & \begin{pmatrix} \gamma_7 & \{0,0,0\} \\ \{0,0,0\} & \gamma_7 \end{pmatrix} & \begin{pmatrix} \lambda_{12} & \{\lambda_{14},\lambda_{15},\lambda_{16}\} \\ \{\lambda_{17},0,\lambda_{19}\} & \lambda_{13} \end{pmatrix} \\ \begin{pmatrix} \lambda_{21} & \{0,-\lambda_{23},0\} \\ \{-\lambda_{25},0,-\lambda_{27}\} & 0 \end{pmatrix} & \begin{pmatrix} \lambda_{13} & \{-\lambda_{14},-\lambda_{15},-\lambda_{16}\} \\ \{-\lambda_{17},0,-\lambda_{19}\} & \lambda_{12} \end{pmatrix} & \begin{pmatrix} -\gamma_7 & \{0,0,0\} \\ \{0,0,0\} & -\gamma_7 \end{pmatrix} \end{pmatrix}$$

**U_{X₂[e₄]} [X₁[e₁]] == ZERO**

True

**defin[X$_2$[e$_4$], X$_1$[e$_1$]] // Expand**

$$\begin{pmatrix} \begin{pmatrix} 0 & \{0,0,0\} \\ \{0,0,0\} & 0 \end{pmatrix} & \begin{pmatrix} 0 & \{0,0,0\} \\ \{0,0,0\} & 0 \end{pmatrix} & \begin{pmatrix} 0 & \{0,0,0\} \\ \{0,0,0\} & 0 \end{pmatrix} \\ \begin{pmatrix} 0 & \{0,0,0\} \\ \{0,0,0\} & 0 \end{pmatrix} & \begin{pmatrix} 0 & \{0,0,0\} \\ \{0,0,0\} & 0 \end{pmatrix} & \begin{pmatrix} 0 & \{0, \lambda_5 - \beta_7, 0\} \\ \{0,0,0\} & 0 \end{pmatrix} \\ \begin{pmatrix} 0 & \{0,0,0\} \\ \{0,0,0\} & 0 \end{pmatrix} & \begin{pmatrix} 0 & \{0, \beta_7 - \lambda_5, 0\} \\ \{0,0,0\} & 0 \end{pmatrix} & \begin{pmatrix} 0 & \{0,0,0\} \\ \{0,0,0\} & 0 \end{pmatrix} \end{pmatrix}$$

**Δ[X$_2$[e$_4$]] = Δ[X$_2$[e$_4$]] //. {λ$_5$ → β$_7$}**

$$\begin{pmatrix} \begin{pmatrix} 0 & \{0,0,0\} \\ \{0,0,0\} & 0 \end{pmatrix} & \begin{pmatrix} 0 & \{0, \lambda_7, 0\} \\ \{\lambda_9, 0, \lambda_{11}\} & \beta_7 \end{pmatrix} & \begin{pmatrix} 0 & \{0, \lambda_{23}, 0\} \\ \{\lambda_{25}, 0, \lambda_{27}\} & \lambda_{21} \end{pmatrix} \\ \begin{pmatrix} \beta_7 & \{0, -\lambda_7, 0\} \\ \{-\lambda_9, 0, -\lambda_{11}\} & 0 \end{pmatrix} & \begin{pmatrix} \gamma_7 & \{0,0,0\} \\ \{0,0,0\} & \gamma_7 \end{pmatrix} & \begin{pmatrix} \lambda_{12} & \{\lambda_{14}, \lambda_{15}, \lambda_{16}\} \\ \{\lambda_{17}, 0, \lambda_{19}\} & \lambda_{13} \end{pmatrix} \\ \begin{pmatrix} \lambda_{21} & \{0, -\lambda_{23}, 0\} \\ \{-\lambda_{25}, 0, -\lambda_{27}\} & 0 \end{pmatrix} & \begin{pmatrix} \lambda_{13} & \{-\lambda_{14}, -\lambda_{15}, -\lambda_{16}\} \\ \{-\lambda_{17}, 0, -\lambda_{19}\} & \lambda_{12} \end{pmatrix} & \begin{pmatrix} -\gamma_7 & \{0,0,0\} \\ \{0,0,0\} & -\gamma_7 \end{pmatrix} \end{pmatrix}$$

**U$_{X_2[e_4]}$[X$_1$[e$_3$]] == ZERO**

True

**defin[X$_2$[e$_4$], X$_1$[e$_3$]] // Expand**

$$\begin{pmatrix} \begin{pmatrix} 0 & \{0,0,0\} \\ \{0,0,0\} & 0 \end{pmatrix} & \begin{pmatrix} 0 & \{0,0,0\} \\ \{0,0,0\} & 0 \end{pmatrix} & \begin{pmatrix} 0 & \{0,0,0\} \\ \{0,0,0\} & 0 \end{pmatrix} \\ \begin{pmatrix} 0 & \{0,0,0\} \\ \{0,0,0\} & 0 \end{pmatrix} & \begin{pmatrix} 0 & \{0,0,0\} \\ \{0,0,0\} & 0 \end{pmatrix} & \begin{pmatrix} 0 & \{0, -\beta_5 - \lambda_9, 0\} \\ \{0,0,0\} & 0 \end{pmatrix} \\ \begin{pmatrix} 0 & \{0,0,0\} \\ \{0,0,0\} & 0 \end{pmatrix} & \begin{pmatrix} 0 & \{0, \beta_5 + \lambda_9, 0\} \\ \{0,0,0\} & 0 \end{pmatrix} & \begin{pmatrix} 0 & \{0,0,0\} \\ \{0,0,0\} & 0 \end{pmatrix} \end{pmatrix}$$

**Δ[X$_2$[e$_4$]] = Δ[X$_2$[e$_4$]] //. {λ$_9$ → -β$_5$}**

$$\begin{pmatrix} \begin{pmatrix} 0 & \{0,0,0\} \\ \{0,0,0\} & 0 \end{pmatrix} & \begin{pmatrix} 0 & \{0, \lambda_7, 0\} \\ \{-\beta_5, 0, \lambda_{11}\} & \beta_7 \end{pmatrix} & \begin{pmatrix} 0 & \{0, \lambda_{23}, 0\} \\ \{\lambda_{25}, 0, \lambda_{27}\} & \lambda_{21} \end{pmatrix} \\ \begin{pmatrix} \beta_7 & \{0, -\lambda_7, 0\} \\ \{\beta_5, 0, -\lambda_{11}\} & 0 \end{pmatrix} & \begin{pmatrix} \gamma_7 & \{0,0,0\} \\ \{0,0,0\} & \gamma_7 \end{pmatrix} & \begin{pmatrix} \lambda_{12} & \{\lambda_{14}, \lambda_{15}, \lambda_{16}\} \\ \{\lambda_{17}, 0, \lambda_{19}\} & \lambda_{13} \end{pmatrix} \\ \begin{pmatrix} \lambda_{21} & \{0, -\lambda_{23}, 0\} \\ \{-\lambda_{25}, 0, -\lambda_{27}\} & 0 \end{pmatrix} & \begin{pmatrix} \lambda_{13} & \{-\lambda_{14}, -\lambda_{15}, -\lambda_{16}\} \\ \{-\lambda_{17}, 0, -\lambda_{19}\} & \lambda_{12} \end{pmatrix} & \begin{pmatrix} -\gamma_7 & \{0,0,0\} \\ \{0,0,0\} & -\gamma_7 \end{pmatrix} \end{pmatrix}$$

**U$_{X_2[e_4]}$[X$_1$[e$_5$]] == ZERO**

True

**defin[X₂[e₄], X₁[e₅]] // Expand**

$$\begin{pmatrix} \begin{pmatrix} 0 & \{0,0,0\} \\ \{0,0,0\} & 0 \end{pmatrix} & \begin{pmatrix} 0 & \{0,0,0\} \\ \{0,0,0\} & 0 \end{pmatrix} & \begin{pmatrix} 0 & \{0,0,0\} \\ \{0,0,0\} & 0 \end{pmatrix} \\ \begin{pmatrix} 0 & \{0,0,0\} \\ \{0,0,0\} & 0 \end{pmatrix} & \begin{pmatrix} 0 & \{0,0,0\} \\ \{0,0,0\} & 0 \end{pmatrix} & \begin{pmatrix} 0 & \{0,\beta_3-\lambda_{11},0\} \\ \{0,0,0\} & 0 \end{pmatrix} \\ \begin{pmatrix} 0 & \{0,0,0\} \\ \{0,0,0\} & 0 \end{pmatrix} & \begin{pmatrix} 0 & \{0,\lambda_{11}-\beta_3,0\} \\ \{0,0,0\} & 0 \end{pmatrix} & \begin{pmatrix} 0 & \{0,0,0\} \\ \{0,0,0\} & 0 \end{pmatrix} \end{pmatrix}$$

**Δ[X₂[e₄]] = Δ[X₂[e₄]] //. {λ₁₁ → β₃}**

$$\begin{pmatrix} \begin{pmatrix} 0 & \{0,0,0\} \\ \{0,0,0\} & 0 \end{pmatrix} & \begin{pmatrix} 0 & \{0,\lambda_7,0\} \\ \{-\beta_5,0,\beta_3\} & \beta_7 \end{pmatrix} & \begin{pmatrix} 0 & \{0,\lambda_{23},0\} \\ \{\lambda_{25},0,\lambda_{27}\} & \lambda_{21} \end{pmatrix} \\ \begin{pmatrix} \beta_7 & \{0,-\lambda_7,0\} \\ \{\beta_5,0,-\beta_3\} & 0 \end{pmatrix} & \begin{pmatrix} \gamma_7 & \{0,0,0\} \\ \{0,0,0\} & \gamma_7 \end{pmatrix} & \begin{pmatrix} \lambda_{12} & \{\lambda_{14},\lambda_{15},\lambda_{16}\} \\ \{\lambda_{17},0,\lambda_{19}\} & \lambda_{13} \end{pmatrix} \\ \begin{pmatrix} \lambda_{21} & \{0,-\lambda_{23},0\} \\ \{-\lambda_{25},0,-\lambda_{27}\} & 0 \end{pmatrix} & \begin{pmatrix} \lambda_{13} & \{-\lambda_{14},-\lambda_{15},-\lambda_{16}\} \\ \{-\lambda_{17},0,-\lambda_{19}\} & \lambda_{12} \end{pmatrix} & \begin{pmatrix} -\gamma_7 & \{0,0,0\} \\ \{0,0,0\} & -\gamma_7 \end{pmatrix} \end{pmatrix}$$

**U_{X₁[e₇]}[X₂[e₄]] == ZERO**

True

**defin[X₁[e₇], X₂[e₄]] // Expand**

$$\begin{pmatrix} \begin{pmatrix} 0 & \{0,0,0\} \\ \{0,0,0\} & 0 \end{pmatrix} & \begin{pmatrix} 0 & \{0,0,0\} \\ \{0,\beta_1-\lambda_7,0\} & 0 \end{pmatrix} & \begin{pmatrix} 0 & \{0,0,0\} \\ \{0,0,0\} & 0 \end{pmatrix} \\ \begin{pmatrix} 0 & \{0,0,0\} \\ \{0,\lambda_7-\beta_1,0\} & 0 \end{pmatrix} & \begin{pmatrix} 0 & \{0,0,0\} \\ \{0,0,0\} & 0 \end{pmatrix} & \begin{pmatrix} 0 & \{0,0,0\} \\ \{0,0,0\} & 0 \end{pmatrix} \\ \begin{pmatrix} 0 & \{0,0,0\} \\ \{0,0,0\} & 0 \end{pmatrix} & \begin{pmatrix} 0 & \{0,0,0\} \\ \{0,0,0\} & 0 \end{pmatrix} & \begin{pmatrix} 0 & \{0,0,0\} \\ \{0,0,0\} & 0 \end{pmatrix} \end{pmatrix}$$

**Δ[X₂[e₄]] = Δ[X₂[e₄]] //. {λ₇ → β₁}**

$$\begin{pmatrix} \begin{pmatrix} 0 & \{0,0,0\} \\ \{0,0,0\} & 0 \end{pmatrix} & \begin{pmatrix} 0 & \{0,\beta_1,0\} \\ \{-\beta_5,0,\beta_3\} & \beta_7 \end{pmatrix} & \begin{pmatrix} 0 & \{0,\lambda_{23},0\} \\ \{\lambda_{25},0,\lambda_{27}\} & \lambda_{21} \end{pmatrix} \\ \begin{pmatrix} \beta_7 & \{0,-\beta_1,0\} \\ \{\beta_5,0,-\beta_3\} & 0 \end{pmatrix} & \begin{pmatrix} \gamma_7 & \{0,0,0\} \\ \{0,0,0\} & \gamma_7 \end{pmatrix} & \begin{pmatrix} \lambda_{12} & \{\lambda_{14},\lambda_{15},\lambda_{16}\} \\ \{\lambda_{17},0,\lambda_{19}\} & \lambda_{13} \end{pmatrix} \\ \begin{pmatrix} \lambda_{21} & \{0,-\lambda_{23},0\} \\ \{-\lambda_{25},0,-\lambda_{27}\} & 0 \end{pmatrix} & \begin{pmatrix} \lambda_{13} & \{-\lambda_{14},-\lambda_{15},-\lambda_{16}\} \\ \{-\lambda_{17},0,-\lambda_{19}\} & \lambda_{12} \end{pmatrix} & \begin{pmatrix} -\gamma_7 & \{0,0,0\} \\ \{0,0,0\} & -\gamma_7 \end{pmatrix} \end{pmatrix}$$

**T[E₁[1], E₁[1], X₂[e₄]] == ZERO**

True

**Leib[E₁[1], E₁[1], X₂[e₄]] // Expand**

$$\begin{pmatrix} \begin{pmatrix} 0 & \{0,0,0\} \\ \{0,0,0\} & 0 \end{pmatrix} & \begin{pmatrix} 0 & \{0,0,0\} \\ \{0,0,0\} & 0 \end{pmatrix} & \begin{pmatrix} 0 & \{0,\alpha_1+\lambda_{23},0\} \\ \{\lambda_{25}-\alpha_5,0,\alpha_3+\lambda_{27}\} & \alpha_7+\lambda_{21} \end{pmatrix} \\ \begin{pmatrix} 0 & \{0,0,0\} \\ \{0,0,0\} & 0 \end{pmatrix} & \begin{pmatrix} 0 & \{0,0,0\} \\ \{0,0,0\} & 0 \end{pmatrix} & \begin{pmatrix} 0 & \{0,0,0\} \\ \{0,0,0\} & 0 \end{pmatrix} \\ \begin{pmatrix} \alpha_7+\lambda_{21} & \{0,-\alpha_1-\lambda_{23},0\} \\ \{\alpha_5-\lambda_{25},0,-\alpha_3-\lambda_{27}\} & 0 \end{pmatrix} & \begin{pmatrix} 0 & \{0,0,0\} \\ \{0,0,0\} & 0 \end{pmatrix} & \begin{pmatrix} 0 & \{0,0,0\} \\ \{0,0,0\} & 0 \end{pmatrix} \end{pmatrix}$$

$\Delta[X_2[e_4]] = \Delta[X_2[e_4]] \text{ //. } \{\lambda_{21} \to -\alpha_7, \lambda_{23} \to -\alpha_1, \lambda_{25} \to \alpha_5, \lambda_{27} \to -\alpha_3\}$

$$\begin{pmatrix} \begin{pmatrix} 0 & \{0,0,0\} \\ \{0,0,0\} & 0 \end{pmatrix} & \begin{pmatrix} 0 & \{0,\beta_1,0\} \\ \{-\beta_5,0,\beta_3\} & \beta_7 \end{pmatrix} & \begin{pmatrix} 0 & \{0,-\alpha_1,0\} \\ \{\alpha_5,0,-\alpha_3\} & -\alpha_7 \end{pmatrix} \\ \begin{pmatrix} \beta_7 & \{0,-\beta_1,0\} \\ \{\beta_5,0,-\beta_3\} & 0 \end{pmatrix} & \begin{pmatrix} \gamma_7 & \{0,0,0\} \\ \{0,0,0\} & \gamma_7 \end{pmatrix} & \begin{pmatrix} \lambda_{12} & \{\lambda_{14},\lambda_{15},\lambda_{16}\} \\ \{\lambda_{17},0,\lambda_{19}\} & \lambda_{13} \end{pmatrix} \\ \begin{pmatrix} -\alpha_7 & \{0,\alpha_1,0\} \\ \{-\alpha_5,0,\alpha_3\} & 0 \end{pmatrix} & \begin{pmatrix} \lambda_{13} & \{-\lambda_{14},-\lambda_{15},-\lambda_{16}\} \\ \{-\lambda_{17},0,-\lambda_{19}\} & \lambda_{12} \end{pmatrix} & \begin{pmatrix} -\gamma_7 & \{0,0,0\} \\ \{0,0,0\} & -\gamma_7 \end{pmatrix} \end{pmatrix}$$

**T[X$_1$[e$_1$], X$_1$[e$_2$], X$_2$[e$_4$]] == ZERO**

True

**Leib[X$_1$[e$_1$], X$_1$[e$_2$], X$_2$[e$_4$]] // Expand**

$$\begin{pmatrix} \begin{pmatrix} 0 & \{0,0,0\} \\ \{0,0,0\} & 0 \end{pmatrix} & \begin{pmatrix} 0 & \{0,0,0\} \\ \{0,0,0\} & 0 \end{pmatrix} & \begin{pmatrix} 0 & \{0,0,0\} \\ \{0,0,0\} & 0 \end{pmatrix} \\ \begin{pmatrix} 0 & \{0,0,0\} \\ \{0,0,0\} & 0 \end{pmatrix} & \begin{pmatrix} 0 & \{0,0,0\} \\ \{0,0,0\} & 0 \end{pmatrix} & \begin{pmatrix} 0 & \{0,0,0\} \\ \{\lambda_{17}-\rho_3,0,\lambda_{19}+\rho_1\} & \lambda_{13}+\rho_5 \end{pmatrix} \\ \begin{pmatrix} 0 & \{0,0,0\} \\ \{0,0,0\} & 0 \end{pmatrix} & \begin{pmatrix} \lambda_{13}+\rho_5 & \{0,0,0\} \\ \{\rho_3-\lambda_{17},0,-\lambda_{19}-\rho_1\} & 0 \end{pmatrix} & \begin{pmatrix} 0 & \{0,0,0\} \\ \{0,0,0\} & 0 \end{pmatrix} \end{pmatrix}$$

$\Delta[X_2[e_4]] = \Delta[X_2[e_4]] \text{ //. } \{\lambda_{13} \to -\rho_5, \lambda_{17} \to \rho_3, \lambda_{19} \to -\rho_1\}$

$$\begin{pmatrix} \begin{pmatrix} 0 & \{0,0,0\} \\ \{0,0,0\} & 0 \end{pmatrix} & \begin{pmatrix} 0 & \{0,\beta_1,0\} \\ \{-\beta_5,0,\beta_3\} & \beta_7 \end{pmatrix} & \begin{pmatrix} 0 & \{0,-\alpha_1,0\} \\ \{\alpha_5,0,-\alpha_3\} & -\alpha_7 \end{pmatrix} \\ \begin{pmatrix} \beta_7 & \{0,-\beta_1,0\} \\ \{\beta_5,0,-\beta_3\} & 0 \end{pmatrix} & \begin{pmatrix} \gamma_7 & \{0,0,0\} \\ \{0,0,0\} & \gamma_7 \end{pmatrix} & \begin{pmatrix} \lambda_{12} & \{\lambda_{14},\lambda_{15},\lambda_{16}\} \\ \{\rho_3,0,-\rho_1\} & -\rho_5 \end{pmatrix} \\ \begin{pmatrix} -\alpha_7 & \{0,\alpha_1,0\} \\ \{-\alpha_5,0,\alpha_3\} & 0 \end{pmatrix} & \begin{pmatrix} -\rho_5 & \{-\lambda_{14},-\lambda_{15},-\lambda_{16}\} \\ \{-\rho_3,0,\rho_1\} & \lambda_{12} \end{pmatrix} & \begin{pmatrix} -\gamma_7 & \{0,0,0\} \\ \{0,0,0\} & -\gamma_7 \end{pmatrix} \end{pmatrix}$$

**T[X$_1$[e$_1$], X$_1$[e$_4$], X$_2$[e$_4$]] == ZERO**

True

**Leib[X$_1$[e$_1$], X$_1$[e$_4$], X$_2$[e$_4$]] // Expand**

$$\begin{pmatrix} \begin{pmatrix} 0 & \{0,0,0\} \\ \{0,0,0\} & 0 \end{pmatrix} & \begin{pmatrix} 0 & \{0,0,0\} \\ \{0,0,0\} & 0 \end{pmatrix} & \begin{pmatrix} 0 & \{0,0,0\} \\ \{0,0,0\} & 0 \end{pmatrix} \\ \begin{pmatrix} 0 & \{0,0,0\} \\ \{0,0,0\} & 0 \end{pmatrix} & \begin{pmatrix} 0 & \{0,0,0\} \\ \{0,0,0\} & 0 \end{pmatrix} & \begin{pmatrix} 0 & \{0,0,0\} \\ \{\lambda_{16}-\epsilon_3,0,\epsilon_1-\lambda_{14}\} & 0 \end{pmatrix} \\ \begin{pmatrix} 0 & \{0,0,0\} \\ \{0,0,0\} & 0 \end{pmatrix} & \begin{pmatrix} 0 & \{0,0,0\} \\ \{\epsilon_3-\lambda_{16},0,\lambda_{14}-\epsilon_1\} & 0 \end{pmatrix} & \begin{pmatrix} 0 & \{0,0,0\} \\ \{0,0,0\} & 0 \end{pmatrix} \end{pmatrix}$$

**Δ[X₂[e₄]] = Δ[X₂[e₄]] //. {λ₁₄ → ϵ₁, λ₁₆ → ϵ₃}**

$$\begin{pmatrix} \begin{pmatrix} 0 & \{0,0,0\} \\ \{0,0,0\} & 0 \end{pmatrix} & \begin{pmatrix} 0 & \{0,\beta_1,0\} \\ \{-\beta_5,0,\beta_3\} & \beta_7 \end{pmatrix} & \begin{pmatrix} 0 & \{0,-\alpha_1,0\} \\ \{\alpha_5,0,-\alpha_3\} & -\alpha_7 \end{pmatrix} \\ \begin{pmatrix} \beta_7 & \{0,-\beta_1,0\} \\ \{\beta_5,0,-\beta_3\} & 0 \end{pmatrix} & \begin{pmatrix} \gamma_7 & \{0,0,0\} \\ \{0,0,0\} & \gamma_7 \end{pmatrix} & \begin{pmatrix} \lambda_{12} & \{\epsilon_1,\lambda_{15},\epsilon_3\} \\ \{\rho_3,0,-\rho_1\} & -\rho_5 \end{pmatrix} \\ \begin{pmatrix} -\alpha_7 & \{0,\alpha_1,0\} \\ \{-\alpha_5,0,\alpha_3\} & 0 \end{pmatrix} & \begin{pmatrix} -\rho_5 & \{-\epsilon_1,-\lambda_{15},-\epsilon_3\} \\ \{-\rho_3,0,\rho_1\} & \lambda_{12} \end{pmatrix} & \begin{pmatrix} -\gamma_7 & \{0,0,0\} \\ \{0,0,0\} & -\gamma_7 \end{pmatrix} \end{pmatrix}$$

**T[X₁[e₁], X₁[e₆], X₂[e₄]] == ZERO**

True

**Leib[X₁[e₁], X₁[e₆], X₂[e₄]] // Expand**

$$\begin{pmatrix} \begin{pmatrix} 0 & \{0,0,0\} \\ \{0,0,0\} & 0 \end{pmatrix} & \begin{pmatrix} 0 & \{0,0,0\} \\ \{0,0,0\} & 0 \end{pmatrix} & \begin{pmatrix} 0 & \{0,0,0\} \\ \{0,0,0\} & 0 \end{pmatrix} \\ \begin{pmatrix} 0 & \{0,0,0\} \\ \{0,0,0\} & 0 \end{pmatrix} & \begin{pmatrix} 0 & \{0,0,0\} \\ \{0,0,0\} & 0 \end{pmatrix} & \begin{pmatrix} 0 & \{0,0,0\} \\ \{\lambda_{12}-\psi_2,0,0\} & 0 \end{pmatrix} \\ \begin{pmatrix} 0 & \{0,0,0\} \\ \{0,0,0\} & 0 \end{pmatrix} & \begin{pmatrix} 0 & \{0,0,0\} \\ \{\psi_2-\lambda_{12},0,0\} & 0 \end{pmatrix} & \begin{pmatrix} 0 & \{0,0,0\} \\ \{0,0,0\} & 0 \end{pmatrix} \end{pmatrix}$$

**Δ[X₂[e₄]] = Δ[X₂[e₄]] //. {λ₁₂ → ψ₂}**

$$\begin{pmatrix} \begin{pmatrix} 0 & \{0,0,0\} \\ \{0,0,0\} & 0 \end{pmatrix} & \begin{pmatrix} 0 & \{0,\beta_1,0\} \\ \{-\beta_5,0,\beta_3\} & \beta_7 \end{pmatrix} & \begin{pmatrix} 0 & \{0,-\alpha_1,0\} \\ \{\alpha_5,0,-\alpha_3\} & -\alpha_7 \end{pmatrix} \\ \begin{pmatrix} \beta_7 & \{0,-\beta_1,0\} \\ \{\beta_5,0,-\beta_3\} & 0 \end{pmatrix} & \begin{pmatrix} \gamma_7 & \{0,0,0\} \\ \{0,0,0\} & \gamma_7 \end{pmatrix} & \begin{pmatrix} \psi_2 & \{\epsilon_1,\lambda_{15},\epsilon_3\} \\ \{\rho_3,0,-\rho_1\} & -\rho_5 \end{pmatrix} \\ \begin{pmatrix} -\alpha_7 & \{0,\alpha_1,0\} \\ \{-\alpha_5,0,\alpha_3\} & 0 \end{pmatrix} & \begin{pmatrix} -\rho_5 & \{-\epsilon_1,-\lambda_{15},-\epsilon_3\} \\ \{-\rho_3,0,\rho_1\} & \psi_2 \end{pmatrix} & \begin{pmatrix} -\gamma_7 & \{0,0,0\} \\ \{0,0,0\} & -\gamma_7 \end{pmatrix} \end{pmatrix}$$

**T[X₁[e₁], X₁[e₇], X₂[e₄]] == X₂[e₂]**

True

**Δ[X₂[e₂]] − Leib[X₁[e₁], X₁[e₇], X₂[e₄]] // Expand**

$$\begin{pmatrix} \begin{pmatrix} 0 & \{0,0,0\} \\ \{0,0,0\} & 0 \end{pmatrix} & \begin{pmatrix} 0 & \{0,0,0\} \\ \{0,0,0\} & 0 \end{pmatrix} & \begin{pmatrix} 0 & \{0,0,0\} \\ \{0,0,0\} & 0 \end{pmatrix} \\ \begin{pmatrix} 0 & \{0,0,0\} \\ \{0,0,0\} & 0 \end{pmatrix} & \begin{pmatrix} 0 & \{0,0,0\} \\ \{0,0,0\} & 0 \end{pmatrix} & \begin{pmatrix} 0 & \{0,0,0\} \\ \{0,0,0\} & -\xi-\delta_1+\epsilon_2-\lambda_{15} \end{pmatrix} \\ \begin{pmatrix} 0 & \{0,0,0\} \\ \{0,0,0\} & 0 \end{pmatrix} & \begin{pmatrix} -\xi-\delta_1+\epsilon_2-\lambda_{15} & \{0,0,0\} \\ \{0,0,0\} & 0 \end{pmatrix} & \begin{pmatrix} 0 & \{0,0,0\} \\ \{0,0,0\} & 0 \end{pmatrix} \end{pmatrix}$$

$\Delta[X_2[e_4]] = \Delta[X_2[e_4]] \,//.\, \{\lambda_{15} \to -\xi - \delta_1 + \epsilon_2\}$

$$\begin{pmatrix} \begin{pmatrix} 0 & \{0,0,0\} \\ \{0,0,0\} & 0 \end{pmatrix} & \begin{pmatrix} 0 & \{0,\beta_1,0\} \\ \{-\beta_5,0,\beta_3\} & \beta_7 \end{pmatrix} & \begin{pmatrix} 0 & \{0,-\alpha_1,0\} \\ \{\alpha_5,0,-\alpha_3\} & -\alpha_7 \end{pmatrix} \\ \begin{pmatrix} \beta_7 & \{0,-\beta_1,0\} \\ \{\beta_5,0,-\beta_3\} & 0 \end{pmatrix} & \begin{pmatrix} \gamma_7 & \{0,0,0\} \\ \{0,0,0\} & \gamma_7 \end{pmatrix} & \begin{pmatrix} \psi_2 & \{\epsilon_1,-\xi-\delta_1+\epsilon_2,\epsilon_3\} \\ \{\rho_3,0,-\rho_1\} & -\rho_5 \end{pmatrix} \\ \begin{pmatrix} -\alpha_7 & \{0,\alpha_1,0\} \\ \{-\alpha_5,0,\alpha_3\} & 0 \end{pmatrix} & \begin{pmatrix} -\rho_5 & \{-\epsilon_1,\xi+\delta_1-\epsilon_2,-\epsilon_3\} \\ \{-\rho_3,0,\rho_1\} & \psi_2 \end{pmatrix} & \begin{pmatrix} -\gamma_7 & \{0,0,0\} \\ \{0,0,0\} & -\gamma_7 \end{pmatrix} \end{pmatrix}$$

**Variables[Δ[X₂[e₄]]]**

$\{\xi, \alpha_1, \alpha_3, \alpha_5, \alpha_7, \beta_1, \beta_3, \beta_5, \beta_7, \gamma_7, \delta_1, \epsilon_1, \epsilon_2, \epsilon_3, \rho_1, \rho_3, \rho_5, \psi_2\}$

- $X_2[e_5]$

  $\Delta[X_2[e_5]] = \text{generic};$

  $U_{E_1[1]}[X_2[e_5]] == \text{ZERO}$

  True

  $(\text{defin}[E_1[1], X_2[e_5]] \,//\, \text{Expand})$

  $$\begin{pmatrix} \begin{pmatrix} \lambda_1 & \{0,0,0\} \\ \{0,0,0\} & \lambda_1 \end{pmatrix} & \begin{pmatrix} 0 & \{0,0,0\} \\ \{0,0,0\} & 0 \end{pmatrix} & \begin{pmatrix} 0 & \{0,0,0\} \\ \{0,0,0\} & 0 \end{pmatrix} \\ \begin{pmatrix} 0 & \{0,0,0\} \\ \{0,0,0\} & 0 \end{pmatrix} & \begin{pmatrix} 0 & \{0,0,0\} \\ \{0,0,0\} & 0 \end{pmatrix} & \begin{pmatrix} 0 & \{0,0,0\} \\ \{0,0,0\} & 0 \end{pmatrix} \\ \begin{pmatrix} 0 & \{0,0,0\} \\ \{0,0,0\} & 0 \end{pmatrix} & \begin{pmatrix} 0 & \{0,0,0\} \\ \{0,0,0\} & 0 \end{pmatrix} & \begin{pmatrix} 0 & \{0,0,0\} \\ \{0,0,0\} & 0 \end{pmatrix} \end{pmatrix}$$

  $\Delta[X_2[e_5]] = \Delta[X_2[e_5]] \,//.\, \{\lambda_1 \to 0\}$

  $$\begin{pmatrix} \begin{pmatrix} 0 & \{0,0,0\} \\ \{0,0,0\} & 0 \end{pmatrix} & \begin{pmatrix} \lambda_4 & \{\lambda_6,\lambda_7,\lambda_8\} \\ \{\lambda_9,\lambda_{10},\lambda_{11}\} & \lambda_5 \end{pmatrix} & \begin{pmatrix} \lambda_{20} & \{\lambda_{22},\lambda_{23},\lambda_{24}\} \\ \{\lambda_{25},\lambda_{26},\lambda_{27}\} & \lambda_{21} \end{pmatrix} \\ \begin{pmatrix} \lambda_5 & \{-\lambda_6,-\lambda_7,-\lambda_8\} \\ \{-\lambda_9,-\lambda_{10},-\lambda_{11}\} & \lambda_4 \end{pmatrix} & \begin{pmatrix} \lambda_2 & \{0,0,0\} \\ \{0,0,0\} & \lambda_2 \end{pmatrix} & \begin{pmatrix} \lambda_{12} & \{\lambda_{14},\lambda_{15},\lambda_{16}\} \\ \{\lambda_{17},\lambda_{18},\lambda_{19}\} & \lambda_{13} \end{pmatrix} \\ \begin{pmatrix} \lambda_{21} & \{-\lambda_{22},-\lambda_{23},-\lambda_{24}\} \\ \{-\lambda_{25},-\lambda_{26},-\lambda_{27}\} & \lambda_{20} \end{pmatrix} & \begin{pmatrix} \lambda_{13} & \{-\lambda_{14},-\lambda_{15},-\lambda_{16}\} \\ \{-\lambda_{17},-\lambda_{18},-\lambda_{19}\} & \lambda_{12} \end{pmatrix} & \begin{pmatrix} \lambda_3 & \{0,0,0\} \\ \{0,0,0\} & \lambda_3 \end{pmatrix} \end{pmatrix}$$

  $U_{X_2[e_5]}[E_2[1]] == \text{ZERO}$

  True

  $\text{defin}[X_2[e_5], E_2[1]] \,//\, \text{Expand}$

  $$\begin{pmatrix} \begin{pmatrix} 0 & \{0,0,0\} \\ \{0,0,0\} & 0 \end{pmatrix} & \begin{pmatrix} 0 & \{0,0,0\} \\ \{0,0,0\} & 0 \end{pmatrix} & \begin{pmatrix} 0 & \{0,0,\lambda_4\} \\ \{\lambda_7,-\lambda_6,0\} & \lambda_{11} \end{pmatrix} \\ \begin{pmatrix} 0 & \{0,0,0\} \\ \{0,0,0\} & 0 \end{pmatrix} & \begin{pmatrix} 0 & \{0,0,0\} \\ \{0,0,0\} & 0 \end{pmatrix} & \begin{pmatrix} 0 & \{0,0,\lambda_2-\gamma_8\} \\ \{0,0,0\} & 0 \end{pmatrix} \\ \begin{pmatrix} \lambda_{11} & \{0,0,-\lambda_4\} \\ \{-\lambda_7,\lambda_6,0\} & 0 \end{pmatrix} & \begin{pmatrix} 0 & \{0,0,\gamma_8-\lambda_2\} \\ \{0,0,0\} & 0 \end{pmatrix} & \begin{pmatrix} -\lambda_{19} & \{0,0,0\} \\ \{0,0,0\} & -\lambda_{19} \end{pmatrix} \end{pmatrix}$$

```
Δ[X₂[e₅]] =
 Δ[X₂[e₅]] //. {λ₂ → γ₈, λ₄ → 0, λ₆ → 0, λ₇ → 0, λ₁₁ → 0, λ₁₉ → 0}
```

$$\left(\begin{array}{ccc}\begin{pmatrix}0 & \{0,0,0\}\\ \{0,0,0\} & 0\end{pmatrix} & \begin{pmatrix}0 & \{0,0,\lambda_8\}\\ \{\lambda_9,\lambda_{10},0\} & \lambda_5\end{pmatrix} & \begin{pmatrix}\lambda_{20} & \{\lambda_{22},\lambda_{23},\lambda_{24}\}\\ \{\lambda_{25},\lambda_{26},\lambda_{27}\} & \lambda_{21}\end{pmatrix}\\ \begin{pmatrix}\lambda_5 & \{0,0,-\lambda_8\}\\ \{-\lambda_9,-\lambda_{10},0\} & 0\end{pmatrix} & \begin{pmatrix}\gamma_8 & \{0,0,0\}\\ \{0,0,0\} & \gamma_8\end{pmatrix} & \begin{pmatrix}\lambda_{12} & \{\lambda_{14},\lambda_{15},\lambda_{16}\}\\ \{\lambda_{17},\lambda_{18},0\} & \lambda_{13}\end{pmatrix}\\ \begin{pmatrix}\lambda_{21} & \{-\lambda_{22},-\lambda_{23},-\lambda_{24}\}\\ \{-\lambda_{25},-\lambda_{26},-\lambda_{27}\} & \lambda_{20}\end{pmatrix} & \begin{pmatrix}\lambda_{13} & \{-\lambda_{14},-\lambda_{15},-\lambda_{16}\}\\ \{-\lambda_{17},-\lambda_{18},0\} & \lambda_{12}\end{pmatrix} & \begin{pmatrix}\lambda_3 & \{0,0,0\}\\ \{0,0,0\} & \lambda_3\end{pmatrix}\end{array}\right)$$

```
U_{E₃[1]}[X₂[e₅]] == ZERO
```

True

```
defin[E₃[1], X₂[e₅]] // Expand
```

$$\left(\begin{array}{ccc}\begin{pmatrix}0 & \{0,0,0\}\\ \{0,0,0\} & 0\end{pmatrix} & \begin{pmatrix}0 & \{0,0,0\}\\ \{0,0,0\} & 0\end{pmatrix} & \begin{pmatrix}0 & \{0,0,0\}\\ \{0,0,0\} & 0\end{pmatrix}\\ \begin{pmatrix}0 & \{0,0,0\}\\ \{0,0,0\} & 0\end{pmatrix} & \begin{pmatrix}0 & \{0,0,0\}\\ \{0,0,0\} & 0\end{pmatrix} & \begin{pmatrix}0 & \{0,0,0\}\\ \{0,0,0\} & 0\end{pmatrix}\\ \begin{pmatrix}0 & \{0,0,0\}\\ \{0,0,0\} & 0\end{pmatrix} & \begin{pmatrix}0 & \{0,0,0\}\\ \{0,0,0\} & 0\end{pmatrix} & \begin{pmatrix}\gamma_8+\lambda_3 & \{0,0,0\}\\ \{0,0,0\} & \gamma_8+\lambda_3\end{pmatrix}\end{array}\right)$$

```
Δ[X₂[e₅]] = Δ[X₂[e₅]] //. {λ₃ → -γ₈}
```

$$\left(\begin{array}{ccc}\begin{pmatrix}0 & \{0,0,0\}\\ \{0,0,0\} & 0\end{pmatrix} & \begin{pmatrix}0 & \{0,0,\lambda_8\}\\ \{\lambda_9,\lambda_{10},0\} & \lambda_5\end{pmatrix} & \begin{pmatrix}\lambda_{20} & \{\lambda_{22},\lambda_{23},\lambda_{24}\}\\ \{\lambda_{25},\lambda_{26},\lambda_{27}\} & \lambda_{21}\end{pmatrix}\\ \begin{pmatrix}\lambda_5 & \{0,0,-\lambda_8\}\\ \{-\lambda_9,-\lambda_{10},0\} & 0\end{pmatrix} & \begin{pmatrix}\gamma_8 & \{0,0,0\}\\ \{0,0,0\} & \gamma_8\end{pmatrix} & \begin{pmatrix}\lambda_{12} & \{\lambda_{14},\lambda_{15},\lambda_{16}\}\\ \{\lambda_{17},\lambda_{18},0\} & \lambda_{13}\end{pmatrix}\\ \begin{pmatrix}\lambda_{21} & \{-\lambda_{22},-\lambda_{23},-\lambda_{24}\}\\ \{-\lambda_{25},-\lambda_{26},-\lambda_{27}\} & \lambda_{20}\end{pmatrix} & \begin{pmatrix}\lambda_{13} & \{-\lambda_{14},-\lambda_{15},-\lambda_{16}\}\\ \{-\lambda_{17},-\lambda_{18},0\} & \lambda_{12}\end{pmatrix} & \begin{pmatrix}-\gamma_8 & \{0,0,0\}\\ \{0,0,0\} & -\gamma_8\end{pmatrix}\end{array}\right)$$

```
U_{X₂[e₅]}[E₃[1]] == ZERO
```

True

```
defin[X₂[e₅], E₃[1]] // Expand
```

$$\left(\begin{array}{ccc}\begin{pmatrix}0 & \{0,0,0\}\\ \{0,0,0\} & 0\end{pmatrix} & \begin{pmatrix}0 & \{0,0,-\lambda_{20}\}\\ \{-\lambda_{23},\lambda_{22},0\} & -\lambda_{27}\end{pmatrix} & \begin{pmatrix}0 & \{0,0,0\}\\ \{0,0,0\} & 0\end{pmatrix}\\ \begin{pmatrix}-\lambda_{27} & \{0,0,\lambda_{20}\}\\ \{\lambda_{23},-\lambda_{22},0\} & 0\end{pmatrix} & \begin{pmatrix}0 & \{0,0,0\}\\ \{0,0,0\} & 0\end{pmatrix} & \begin{pmatrix}0 & \{0,0,0\}\\ \{0,0,0\} & 0\end{pmatrix}\\ \begin{pmatrix}0 & \{0,0,0\}\\ \{0,0,0\} & 0\end{pmatrix} & \begin{pmatrix}0 & \{0,0,0\}\\ \{0,0,0\} & 0\end{pmatrix} & \begin{pmatrix}0 & \{0,0,0\}\\ \{0,0,0\} & 0\end{pmatrix}\end{array}\right)$$

**Δ[X₂[e₅]] = Δ[X₂[e₅]] //. {λ₂₀ → 0, λ₂₂ → 0, λ₂₃ → 0, λ₂₇ → 0}**

$$\begin{pmatrix} \begin{pmatrix} 0 & \{0,0,0\} \\ \{0,0,0\} & 0 \end{pmatrix} & \begin{pmatrix} 0 & \{0,0,\lambda_8\} \\ \{\lambda_9,\lambda_{10},0\} & \lambda_5 \end{pmatrix} & \begin{pmatrix} 0 & \{0,0,\lambda_{24}\} \\ \{\lambda_{25},\lambda_{26},0\} & \lambda_{21} \end{pmatrix} \\ \begin{pmatrix} \lambda_5 & \{0,0,-\lambda_8\} \\ \{-\lambda_9,-\lambda_{10},0\} & 0 \end{pmatrix} & \begin{pmatrix} \gamma_8 & \{0,0,0\} \\ \{0,0,0\} & \gamma_8 \end{pmatrix} & \begin{pmatrix} \lambda_{12} & \{\lambda_{14},\lambda_{15},\lambda_{16}\} \\ \{\lambda_{17},\lambda_{18},0\} & \lambda_{13} \end{pmatrix} \\ \begin{pmatrix} \lambda_{21} & \{0,0,-\lambda_{24}\} \\ \{-\lambda_{25},-\lambda_{26},0\} & 0 \end{pmatrix} & \begin{pmatrix} \lambda_{13} & \{-\lambda_{14},-\lambda_{15},-\lambda_{16}\} \\ \{-\lambda_{17},-\lambda_{18},0\} & \lambda_{12} \end{pmatrix} & \begin{pmatrix} -\gamma_8 & \{0,0,0\} \\ \{0,0,0\} & -\gamma_8 \end{pmatrix} \end{pmatrix}$$

**U_{X₂[e₅]}[X₁[e₁]] == ZERO**

True

**defin[X₂[e₅], X₁[e₁]] // Expand**

$$\begin{pmatrix} \begin{pmatrix} 0 & \{0,0,0\} \\ \{0,0,0\} & 0 \end{pmatrix} & \begin{pmatrix} 0 & \{0,0,0\} \\ \{0,0,0\} & 0 \end{pmatrix} & \begin{pmatrix} 0 & \{0,0,0\} \\ \{0,0,0\} & 0 \end{pmatrix} \\ \begin{pmatrix} 0 & \{0,0,0\} \\ \{0,0,0\} & 0 \end{pmatrix} & \begin{pmatrix} 0 & \{0,0,0\} \\ \{0,0,0\} & 0 \end{pmatrix} & \begin{pmatrix} 0 & \{0,0,\lambda_5-\beta_8\} \\ \{0,0,0\} & 0 \end{pmatrix} \\ \begin{pmatrix} 0 & \{0,0,0\} \\ \{0,0,0\} & 0 \end{pmatrix} & \begin{pmatrix} 0 & \{0,0,\beta_8-\lambda_5\} \\ \{0,0,0\} & 0 \end{pmatrix} & \begin{pmatrix} 0 & \{0,0,0\} \\ \{0,0,0\} & 0 \end{pmatrix} \end{pmatrix}$$

**Δ[X₂[e₅]] = Δ[X₂[e₅]] //. {λ₅ → β₈}**

$$\begin{pmatrix} \begin{pmatrix} 0 & \{0,0,0\} \\ \{0,0,0\} & 0 \end{pmatrix} & \begin{pmatrix} 0 & \{0,0,\lambda_8\} \\ \{\lambda_9,\lambda_{10},0\} & \beta_8 \end{pmatrix} & \begin{pmatrix} 0 & \{0,0,\lambda_{24}\} \\ \{\lambda_{25},\lambda_{26},0\} & \lambda_{21} \end{pmatrix} \\ \begin{pmatrix} \beta_8 & \{0,0,-\lambda_8\} \\ \{-\lambda_9,-\lambda_{10},0\} & 0 \end{pmatrix} & \begin{pmatrix} \gamma_8 & \{0,0,0\} \\ \{0,0,0\} & \gamma_8 \end{pmatrix} & \begin{pmatrix} \lambda_{12} & \{\lambda_{14},\lambda_{15},\lambda_{16}\} \\ \{\lambda_{17},\lambda_{18},0\} & \lambda_{13} \end{pmatrix} \\ \begin{pmatrix} \lambda_{21} & \{0,0,-\lambda_{24}\} \\ \{-\lambda_{25},-\lambda_{26},0\} & 0 \end{pmatrix} & \begin{pmatrix} \lambda_{13} & \{-\lambda_{14},-\lambda_{15},-\lambda_{16}\} \\ \{-\lambda_{17},-\lambda_{18},0\} & \lambda_{12} \end{pmatrix} & \begin{pmatrix} -\gamma_8 & \{0,0,0\} \\ \{0,0,0\} & -\gamma_8 \end{pmatrix} \end{pmatrix}$$

**U_{X₂[e₅]}[X₁[e₃]] == ZERO**

True

**defin[X₂[e₅], X₁[e₃]] // Expand**

$$\begin{pmatrix} \begin{pmatrix} 0 & \{0,0,0\} \\ \{0,0,0\} & 0 \end{pmatrix} & \begin{pmatrix} 0 & \{0,0,0\} \\ \{0,0,0\} & 0 \end{pmatrix} & \begin{pmatrix} 0 & \{0,0,0\} \\ \{0,0,0\} & 0 \end{pmatrix} \\ \begin{pmatrix} 0 & \{0,0,0\} \\ \{0,0,0\} & 0 \end{pmatrix} & \begin{pmatrix} 0 & \{0,0,0\} \\ \{0,0,0\} & 0 \end{pmatrix} & \begin{pmatrix} 0 & \{0,0,\beta_4-\lambda_9\} \\ \{0,0,0\} & 0 \end{pmatrix} \\ \begin{pmatrix} 0 & \{0,0,0\} \\ \{0,0,0\} & 0 \end{pmatrix} & \begin{pmatrix} 0 & \{0,0,\lambda_9-\beta_4\} \\ \{0,0,0\} & 0 \end{pmatrix} & \begin{pmatrix} 0 & \{0,0,0\} \\ \{0,0,0\} & 0 \end{pmatrix} \end{pmatrix}$$

**Δ[X₂[e₅]] = Δ[X₂[e₅]] //. {λ₉ → β₄}**

$$\begin{pmatrix} \begin{pmatrix} 0 & \{0,0,0\} \\ \{0,0,0\} & 0 \end{pmatrix} & \begin{pmatrix} 0 & \{0,0,\lambda_8\} \\ \{\beta_4,\lambda_{10},0\} & \beta_8 \end{pmatrix} & \begin{pmatrix} 0 & \{0,0,\lambda_{24}\} \\ \{\lambda_{25},\lambda_{26},0\} & \lambda_{21} \end{pmatrix} \\ \begin{pmatrix} \beta_8 & \{0,0,-\lambda_8\} \\ \{-\beta_4,-\lambda_{10},0\} & 0 \end{pmatrix} & \begin{pmatrix} \gamma_8 & \{0,0,0\} \\ \{0,0,0\} & \gamma_8 \end{pmatrix} & \begin{pmatrix} \lambda_{12} & \{\lambda_{14},\lambda_{15},\lambda_{16}\} \\ \{\lambda_{17},\lambda_{18},0\} & \lambda_{13} \end{pmatrix} \\ \begin{pmatrix} \lambda_{21} & \{0,0,-\lambda_{24}\} \\ \{-\lambda_{25},-\lambda_{26},0\} & 0 \end{pmatrix} & \begin{pmatrix} \lambda_{13} & \{-\lambda_{14},-\lambda_{15},-\lambda_{16}\} \\ \{-\lambda_{17},-\lambda_{18},0\} & \lambda_{12} \end{pmatrix} & \begin{pmatrix} -\gamma_8 & \{0,0,0\} \\ \{0,0,0\} & -\gamma_8 \end{pmatrix} \end{pmatrix}$$

**U_{X₁[e₈]}[X₂[e₅]] == ZERO**

True

**defin[X₁[e₈], X₂[e₅]] // Expand**

$$\begin{pmatrix} \begin{pmatrix} 0 & \{0,0,0\} \\ \{0,0,0\} & 0 \end{pmatrix} & \begin{pmatrix} 0 & \{0,0,0\} \\ \{0,0,\beta_1-\lambda_8\} & 0 \end{pmatrix} & \begin{pmatrix} 0 & \{0,0,0\} \\ \{0,0,0\} & 0 \end{pmatrix} \\ \begin{pmatrix} 0 & \{0,0,0\} \\ \{0,0,\lambda_8-\beta_1\} & 0 \end{pmatrix} & \begin{pmatrix} 0 & \{0,0,0\} \\ \{0,0,0\} & 0 \end{pmatrix} & \begin{pmatrix} 0 & \{0,0,0\} \\ \{0,0,0\} & 0 \end{pmatrix} \\ \begin{pmatrix} 0 & \{0,0,0\} \\ \{0,0,0\} & 0 \end{pmatrix} & \begin{pmatrix} 0 & \{0,0,0\} \\ \{0,0,0\} & 0 \end{pmatrix} & \begin{pmatrix} 0 & \{0,0,0\} \\ \{0,0,0\} & 0 \end{pmatrix} \end{pmatrix}$$

**Δ[X₂[e₅]] = Δ[X₂[e₅]] //. {λ₈ → β₁}**

$$\begin{pmatrix} \begin{pmatrix} 0 & \{0,0,0\} \\ \{0,0,0\} & 0 \end{pmatrix} & \begin{pmatrix} 0 & \{0,0,\beta_1\} \\ \{\beta_4,\lambda_{10},0\} & \beta_8 \end{pmatrix} & \begin{pmatrix} 0 & \{0,0,\lambda_{24}\} \\ \{\lambda_{25},\lambda_{26},0\} & \lambda_{21} \end{pmatrix} \\ \begin{pmatrix} \beta_8 & \{0,0,-\beta_1\} \\ \{-\beta_4,-\lambda_{10},0\} & 0 \end{pmatrix} & \begin{pmatrix} \gamma_8 & \{0,0,0\} \\ \{0,0,0\} & \gamma_8 \end{pmatrix} & \begin{pmatrix} \lambda_{12} & \{\lambda_{14},\lambda_{15},\lambda_{16}\} \\ \{\lambda_{17},\lambda_{18},0\} & \lambda_{13} \end{pmatrix} \\ \begin{pmatrix} \lambda_{21} & \{0,0,-\lambda_{24}\} \\ \{-\lambda_{25},-\lambda_{26},0\} & 0 \end{pmatrix} & \begin{pmatrix} \lambda_{13} & \{-\lambda_{14},-\lambda_{15},-\lambda_{16}\} \\ \{-\lambda_{17},-\lambda_{18},0\} & \lambda_{12} \end{pmatrix} & \begin{pmatrix} -\gamma_8 & \{0,0,0\} \\ \{0,0,0\} & -\gamma_8 \end{pmatrix} \end{pmatrix}$$

**T[E₁[1], E₁[1], X₂[e₅]] == ZERO**

True

**Leib[E₁[1], E₁[1], X₂[e₅]] // Expand**

$$\begin{pmatrix} \begin{pmatrix} 0 & \{0,0,0\} \\ \{0,0,0\} & 0 \end{pmatrix} & \begin{pmatrix} 0 & \{0,0,0\} \\ \{0,\beta_3+\lambda_{10},0\} & 0 \end{pmatrix} & \begin{pmatrix} 0 & \{0,0,\alpha_1+\lambda_{24}\} \\ \{\alpha_4+\lambda_{25},\lambda_{26}-\alpha_3,0\} & \alpha_8+\lambda_{21} \end{pmatrix} \\ \begin{pmatrix} 0 & \{0,0,0\} \\ \{0,-\beta_3-\lambda_{10},0\} & 0 \end{pmatrix} & \begin{pmatrix} 0 & \{0,0,0\} \\ \{0,0,0\} & 0 \end{pmatrix} & \begin{pmatrix} 0 & \{0,0,0\} \\ \{0,0,0\} & 0 \end{pmatrix} \\ \begin{pmatrix} \alpha_8+\lambda_{21} & \{0,0,-\alpha_1-\lambda_{24}\} \\ \{-\alpha_4-\lambda_{25},\alpha_3-\lambda_{26},0\} & 0 \end{pmatrix} & \begin{pmatrix} 0 & \{0,0,0\} \\ \{0,0,0\} & 0 \end{pmatrix} & \begin{pmatrix} 0 & \{0,0,0\} \\ \{0,0,0\} & 0 \end{pmatrix} \end{pmatrix}$$

$\Delta[X_2[e_5]] =$
$\quad \Delta[X_2[e_5]] \,//.\, \{\lambda_{10} \to -\beta_3,\, \lambda_{21} \to -\alpha_8,\, \lambda_{24} \to -\alpha_1,\, \lambda_{25} \to -\alpha_4,\, \lambda_{26} \to \alpha_3\}$

$$\begin{pmatrix} \begin{pmatrix} 0 & \{0,0,0\} \\ \{0,0,0\} & 0 \end{pmatrix} & \begin{pmatrix} 0 & \{0,0,\beta_1\} \\ \{\beta_4,-\beta_3,0\} & \beta_8 \end{pmatrix} & \begin{pmatrix} 0 & \{0,0,-\alpha_1\} \\ \{-\alpha_4,\alpha_3,0\} & -\alpha_8 \end{pmatrix} \\ \begin{pmatrix} \beta_8 & \{0,0,-\beta_1\} \\ \{-\beta_4,\beta_3,0\} & 0 \end{pmatrix} & \begin{pmatrix} \gamma_8 & \{0,0,0\} \\ \{0,0,0\} & \gamma_8 \end{pmatrix} & \begin{pmatrix} \lambda_{12} & \{\lambda_{14},\lambda_{15},\lambda_{16}\} \\ \{\lambda_{17},\lambda_{18},0\} & \lambda_{13} \end{pmatrix} \\ \begin{pmatrix} -\alpha_8 & \{0,0,\alpha_1\} \\ \{\alpha_4,-\alpha_3,0\} & 0 \end{pmatrix} & \begin{pmatrix} \lambda_{13} & \{-\lambda_{14},-\lambda_{15},-\lambda_{16}\} \\ \{-\lambda_{17},-\lambda_{18},0\} & \lambda_{12} \end{pmatrix} & \begin{pmatrix} -\gamma_8 & \{0,0,0\} \\ \{0,0,0\} & -\gamma_8 \end{pmatrix} \end{pmatrix}$$

**T[X$_1$[e$_1$], X$_1$[e$_2$], X$_2$[e$_5$]] == ZERO**

True

**Leib[X$_1$[e$_1$], X$_1$[e$_2$], X$_2$[e$_5$]] // Expand**

$$\begin{pmatrix} \begin{pmatrix} 0 & \{0,0,0\} \\ \{0,0,0\} & 0 \end{pmatrix} & \begin{pmatrix} 0 & \{0,0,0\} \\ \{0,0,0\} & 0 \end{pmatrix} & \begin{pmatrix} 0 & \{0,0,0\} \\ \{0,0,0\} & 0 \end{pmatrix} \\ \begin{pmatrix} 0 & \{0,0,0\} \\ \{0,0,0\} & 0 \end{pmatrix} & \begin{pmatrix} 0 & \{0,0,0\} \\ \{0,0,0\} & 0 \end{pmatrix} & \begin{pmatrix} 0 & \{0,0,0\} \\ \{\lambda_{17}+\rho_2,\lambda_{18}-\rho_1,0\} & \lambda_{13}+\rho_6 \end{pmatrix} \\ \begin{pmatrix} 0 & \{0,0,0\} \\ \{0,0,0\} & 0 \end{pmatrix} & \begin{pmatrix} \lambda_{13}+\rho_6 & \{0,0,0\} \\ \{-\lambda_{17}-\rho_2,\rho_1-\lambda_{18},0\} & 0 \end{pmatrix} & \begin{pmatrix} 0 & \{0,0,0\} \\ \{0,0,0\} & 0 \end{pmatrix} \end{pmatrix}$$

$\Delta[X_2[e_5]] = \Delta[X_2[e_5]] \,//.\, \{\lambda_{13} \to -\rho_6,\, \lambda_{17} \to -\rho_2,\, \lambda_{18} \to \rho_1\}$

$$\begin{pmatrix} \begin{pmatrix} 0 & \{0,0,0\} \\ \{0,0,0\} & 0 \end{pmatrix} & \begin{pmatrix} 0 & \{0,0,\beta_1\} \\ \{\beta_4,-\beta_3,0\} & \beta_8 \end{pmatrix} & \begin{pmatrix} 0 & \{0,0,-\alpha_1\} \\ \{-\alpha_4,\alpha_3,0\} & -\alpha_8 \end{pmatrix} \\ \begin{pmatrix} \beta_8 & \{0,0,-\beta_1\} \\ \{-\beta_4,\beta_3,0\} & 0 \end{pmatrix} & \begin{pmatrix} \gamma_8 & \{0,0,0\} \\ \{0,0,0\} & \gamma_8 \end{pmatrix} & \begin{pmatrix} \lambda_{12} & \{\lambda_{14},\lambda_{15},\lambda_{16}\} \\ \{-\rho_2,\rho_1,0\} & -\rho_6 \end{pmatrix} \\ \begin{pmatrix} -\alpha_8 & \{0,0,\alpha_1\} \\ \{\alpha_4,-\alpha_3,0\} & 0 \end{pmatrix} & \begin{pmatrix} -\rho_6 & \{-\lambda_{14},-\lambda_{15},-\lambda_{16}\} \\ \{\rho_2,-\rho_1,0\} & \lambda_{12} \end{pmatrix} & \begin{pmatrix} -\gamma_8 & \{0,0,0\} \\ \{0,0,0\} & -\gamma_8 \end{pmatrix} \end{pmatrix}$$

**T[X$_1$[e$_1$], X$_1$[e$_5$], X$_2$[e$_5$]] == ZERO**

True

**Leib[X$_1$[e$_1$], X$_1$[e$_5$], X$_2$[e$_5$]] // Expand**

$$\begin{pmatrix} \begin{pmatrix} 0 & \{0,0,0\} \\ \{0,0,0\} & 0 \end{pmatrix} & \begin{pmatrix} 0 & \{0,0,0\} \\ \{0,0,0\} & 0 \end{pmatrix} & \begin{pmatrix} 0 & \{0,0,0\} \\ \{0,0,0\} & 0 \end{pmatrix} \\ \begin{pmatrix} 0 & \{0,0,0\} \\ \{0,0,0\} & 0 \end{pmatrix} & \begin{pmatrix} 0 & \{0,0,0\} \\ \{0,0,0\} & 0 \end{pmatrix} & \begin{pmatrix} 0 & \{0,0,0\} \\ \{\phi_2-\lambda_{15},\lambda_{14}-\phi_1,0\} & 0 \end{pmatrix} \\ \begin{pmatrix} 0 & \{0,0,0\} \\ \{0,0,0\} & 0 \end{pmatrix} & \begin{pmatrix} 0 & \{0,0,0\} \\ \{\lambda_{15}-\phi_2,\phi_1-\lambda_{14},0\} & 0 \end{pmatrix} & \begin{pmatrix} 0 & \{0,0,0\} \\ \{0,0,0\} & 0 \end{pmatrix} \end{pmatrix}$$

$\Delta[X_2[e_5]] = \Delta[X_2[e_5]] \,//.\, \{\lambda_{14} \to \phi_1,\, \lambda_{15} \to \phi_2\}$

$$\begin{pmatrix}
\begin{pmatrix} 0 & \{0,0,0\} \\ \{0,0,0\} & 0 \end{pmatrix} &
\begin{pmatrix} 0 & \{0,0,\beta_1\} \\ \{\beta_4,-\beta_3,0\} & \beta_8 \end{pmatrix} &
\begin{pmatrix} 0 & \{0,0,-\alpha_1\} \\ \{-\alpha_4,\alpha_3,0\} & -\alpha_8 \end{pmatrix} \\
\begin{pmatrix} \beta_8 & \{0,0,-\beta_1\} \\ \{-\beta_4,\beta_3,0\} & 0 \end{pmatrix} &
\begin{pmatrix} \gamma_8 & \{0,0,0\} \\ \{0,0,0\} & \gamma_8 \end{pmatrix} &
\begin{pmatrix} \lambda_{12} & \{\phi_1,\phi_2,\lambda_{16}\} \\ \{-\rho_2,\rho_1,0\} & -\rho_6 \end{pmatrix} \\
\begin{pmatrix} -\alpha_8 & \{0,0,\alpha_1\} \\ \{\alpha_4,-\alpha_3,0\} & 0 \end{pmatrix} &
\begin{pmatrix} -\rho_6 & \{-\phi_1,-\phi_2,-\lambda_{16}\} \\ \{\rho_2,-\rho_1,0\} & \lambda_{12} \end{pmatrix} &
\begin{pmatrix} -\gamma_8 & \{0,0,0\} \\ \{0,0,0\} & -\gamma_8 \end{pmatrix}
\end{pmatrix}$$

**T[X₁[e₁], X₁[e₆], X₂[e₅]] == ZERO**

True

**Leib[X₁[e₁], X₁[e₆], X₂[e₅]] // Expand**

$$\begin{pmatrix}
\begin{pmatrix} 0 & \{0,0,0\} \\ \{0,0,0\} & 0 \end{pmatrix} &
\begin{pmatrix} 0 & \{0,0,0\} \\ \{0,0,0\} & 0 \end{pmatrix} &
\begin{pmatrix} 0 & \{0,0,0\} \\ \{0,0,0\} & 0 \end{pmatrix} \\
\begin{pmatrix} 0 & \{0,0,0\} \\ \{0,0,0\} & 0 \end{pmatrix} &
\begin{pmatrix} 0 & \{0,0,0\} \\ \{0,0,0\} & 0 \end{pmatrix} &
\begin{pmatrix} 0 & \{0,0,0\} \\ \{\lambda_{12}+\psi_1,0,0\} & 0 \end{pmatrix} \\
\begin{pmatrix} 0 & \{0,0,0\} \\ \{0,0,0\} & 0 \end{pmatrix} &
\begin{pmatrix} 0 & \{0,0,0\} \\ \{-\lambda_{12}-\psi_1,0,0\} & 0 \end{pmatrix} &
\begin{pmatrix} 0 & \{0,0,0\} \\ \{0,0,0\} & 0 \end{pmatrix}
\end{pmatrix}$$

$\Delta[X_2[e_5]] = \Delta[X_2[e_5]] \,//.\, \{\lambda_{12} \to -\psi_1\}$

$$\begin{pmatrix}
\begin{pmatrix} 0 & \{0,0,0\} \\ \{0,0,0\} & 0 \end{pmatrix} &
\begin{pmatrix} 0 & \{0,0,\beta_1\} \\ \{\beta_4,-\beta_3,0\} & \beta_8 \end{pmatrix} &
\begin{pmatrix} 0 & \{0,0,-\alpha_1\} \\ \{-\alpha_4,\alpha_3,0\} & -\alpha_8 \end{pmatrix} \\
\begin{pmatrix} \beta_8 & \{0,0,-\beta_1\} \\ \{-\beta_4,\beta_3,0\} & 0 \end{pmatrix} &
\begin{pmatrix} \gamma_8 & \{0,0,0\} \\ \{0,0,0\} & \gamma_8 \end{pmatrix} &
\begin{pmatrix} -\psi_1 & \{\phi_1,\phi_2,\lambda_{16}\} \\ \{-\rho_2,\rho_1,0\} & -\rho_6 \end{pmatrix} \\
\begin{pmatrix} -\alpha_8 & \{0,0,\alpha_1\} \\ \{\alpha_4,-\alpha_3,0\} & 0 \end{pmatrix} &
\begin{pmatrix} -\rho_6 & \{-\phi_1,-\phi_2,-\lambda_{16}\} \\ \{\rho_2,-\rho_1,0\} & -\psi_1 \end{pmatrix} &
\begin{pmatrix} -\gamma_8 & \{0,0,0\} \\ \{0,0,0\} & -\gamma_8 \end{pmatrix}
\end{pmatrix}$$

**T[X₁[e₁], X₁[e₈], X₂[e₅]] == X₂[e₂]**

True

**Δ[X₂[e₂]] − Leib[X₁[e₁], X₁[e₈], X₂[e₅]] // Expand**

$$\begin{pmatrix}
\begin{pmatrix} 0 & \{0,0,0\} \\ \{0,0,0\} & 0 \end{pmatrix} &
\begin{pmatrix} 0 & \{0,0,0\} \\ \{0,0,0\} & 0 \end{pmatrix} &
\begin{pmatrix} 0 & \{0,0,0\} \\ \{0,0,0\} & 0 \end{pmatrix} \\
\begin{pmatrix} 0 & \{0,0,0\} \\ \{0,0,0\} & 0 \end{pmatrix} &
\begin{pmatrix} 0 & \{0,0,0\} \\ \{0,0,0\} & 0 \end{pmatrix} &
\begin{pmatrix} 0 & \{0,0,0\} \\ \{0,0,0\} & \xi-\epsilon_2-\eta_1-\lambda_{16} \end{pmatrix} \\
\begin{pmatrix} 0 & \{0,0,0\} \\ \{0,0,0\} & 0 \end{pmatrix} &
\begin{pmatrix} \xi-\epsilon_2-\eta_1-\lambda_{16} & \{0,0,0\} \\ \{0,0,0\} & 0 \end{pmatrix} &
\begin{pmatrix} 0 & \{0,0,0\} \\ \{0,0,0\} & 0 \end{pmatrix}
\end{pmatrix}$$

$\Delta[X_2[e_5]] = \Delta[X_2[e_5]] \,//.\, \{\lambda_{16} \to \xi - \epsilon_2 - \eta_1\}$

$$\begin{pmatrix} \begin{pmatrix} 0 & \{0,0,0\} \\ \{0,0,0\} & 0 \end{pmatrix} & \begin{pmatrix} 0 & \{0,0,\beta_1\} \\ \{\beta_4,-\beta_3,0\} & \beta_8 \end{pmatrix} & \begin{pmatrix} 0 & \{0,0,-\alpha_1\} \\ \{-\alpha_4,\alpha_3,0\} & -\alpha_8 \end{pmatrix} \\ \begin{pmatrix} \beta_8 & \{0,0,-\beta_1\} \\ \{-\beta_4,\beta_3,0\} & 0 \end{pmatrix} & \begin{pmatrix} \gamma_8 & \{0,0,0\} \\ \{0,0,0\} & \gamma_8 \end{pmatrix} & \begin{pmatrix} -\psi_1 & \{\phi_1,\phi_2,\xi-\epsilon_2-\eta_1\} \\ \{-\rho_2,\rho_1,0\} & -\rho_6 \end{pmatrix} \\ \begin{pmatrix} -\alpha_8 & \{0,0,\alpha_1\} \\ \{\alpha_4,-\alpha_3,0\} & 0 \end{pmatrix} & \begin{pmatrix} -\rho_6 & \{-\phi_1,-\phi_2,-\xi+\epsilon_2+\eta_1\} \\ \{\rho_2,-\rho_1,0\} & -\psi_1 \end{pmatrix} & \begin{pmatrix} -\gamma_8 & \{0,0,0\} \\ \{0,0,0\} & -\gamma_8 \end{pmatrix} \end{pmatrix}$$

**Variables[Δ[X₂[e₅]]]**

$\{\xi, \alpha_1, \alpha_3, \alpha_4, \alpha_8, \beta_1, \beta_3, \beta_4, \beta_8, \gamma_8, \epsilon_2, \eta_1, \rho_1, \rho_2, \rho_6, \phi_1, \phi_2, \psi_1\}$

- **X₂[e₆]**

$\Delta[X_2[e_6]] = \text{generic};$

$U_{E_1[1]}[X_2[e_6]] == \text{ZERO}$

True

$(\text{defin}[E_1[1], X_2[e_6]] \,//\, \text{Expand})$

$$\begin{pmatrix} \begin{pmatrix} \lambda_1 & \{0,0,0\} \\ \{0,0,0\} & \lambda_1 \end{pmatrix} & \begin{pmatrix} 0 & \{0,0,0\} \\ \{0,0,0\} & 0 \end{pmatrix} & \begin{pmatrix} 0 & \{0,0,0\} \\ \{0,0,0\} & 0 \end{pmatrix} \\ \begin{pmatrix} 0 & \{0,0,0\} \\ \{0,0,0\} & 0 \end{pmatrix} & \begin{pmatrix} 0 & \{0,0,0\} \\ \{0,0,0\} & 0 \end{pmatrix} & \begin{pmatrix} 0 & \{0,0,0\} \\ \{0,0,0\} & 0 \end{pmatrix} \\ \begin{pmatrix} 0 & \{0,0,0\} \\ \{0,0,0\} & 0 \end{pmatrix} & \begin{pmatrix} 0 & \{0,0,0\} \\ \{0,0,0\} & 0 \end{pmatrix} & \begin{pmatrix} 0 & \{0,0,0\} \\ \{0,0,0\} & 0 \end{pmatrix} \end{pmatrix}$$

$\Delta[X_2[e_6]] = \Delta[X_2[e_6]] \,//.\, \{\lambda_1 \to 0\}$

$$\begin{pmatrix} \begin{pmatrix} 0 & \{0,0,0\} \\ \{0,0,0\} & 0 \end{pmatrix} & \begin{pmatrix} \lambda_4 & \{\lambda_6,\lambda_7,\lambda_8\} \\ \{\lambda_9,\lambda_{10},\lambda_{11}\} & \lambda_5 \end{pmatrix} & \begin{pmatrix} \lambda_{20} & \{\lambda_{22},\lambda_{23},\lambda_{24}\} \\ \{\lambda_{25},\lambda_{26},\lambda_{27}\} & \lambda_{21} \end{pmatrix} \\ \begin{pmatrix} \lambda_5 & \{-\lambda_6,-\lambda_7,-\lambda_8\} \\ \{-\lambda_9,-\lambda_{10},-\lambda_{11}\} & \lambda_4 \end{pmatrix} & \begin{pmatrix} \lambda_2 & \{0,0,0\} \\ \{0,0,0\} & \lambda_2 \end{pmatrix} & \begin{pmatrix} \lambda_{12} & \{\lambda_{14},\lambda_{15},\lambda_{16}\} \\ \{\lambda_{17},\lambda_{18},\lambda_{19}\} & \lambda_{13} \end{pmatrix} \\ \begin{pmatrix} \lambda_{21} & \{-\lambda_{22},-\lambda_{23},-\lambda_{24}\} \\ \{-\lambda_{25},-\lambda_{26},-\lambda_{27}\} & \lambda_{20} \end{pmatrix} & \begin{pmatrix} \lambda_{13} & \{-\lambda_{14},-\lambda_{15},-\lambda_{16}\} \\ \{-\lambda_{17},-\lambda_{18},-\lambda_{19}\} & \lambda_{12} \end{pmatrix} & \begin{pmatrix} \lambda_3 & \{0,0,0\} \\ \{0,0,0\} & \lambda_3 \end{pmatrix} \end{pmatrix}$$

$U_{X_2[e_6]}[E_2[1]] == \text{ZERO}$

True

$\text{defin}[X_2[e_6], E_2[1]] \,//\, \text{Expand}$

$$\begin{pmatrix} \begin{pmatrix} 0 & \{0,0,0\} \\ \{0,0,0\} & 0 \end{pmatrix} & \begin{pmatrix} 0 & \{0,0,0\} \\ \{0,0,0\} & 0 \end{pmatrix} & \begin{pmatrix} \lambda_6 & \{0,-\lambda_{11},\lambda_{10}\} \\ \{\lambda_5,0,0\} & 0 \end{pmatrix} \\ \begin{pmatrix} 0 & \{0,0,0\} \\ \{0,0,0\} & 0 \end{pmatrix} & \begin{pmatrix} 0 & \{0,0,0\} \\ \{0,0,0\} & 0 \end{pmatrix} & \begin{pmatrix} 0 & \{0,0,0\} \\ \{\lambda_2-\gamma_3,0,0\} & 0 \end{pmatrix} \\ \begin{pmatrix} 0 & \{0,\lambda_{11},-\lambda_{10}\} \\ \{-\lambda_5,0,0\} & \lambda_6 \end{pmatrix} & \begin{pmatrix} 0 & \{0,0,0\} \\ \{\gamma_3-\lambda_2,0,0\} & 0 \end{pmatrix} & \begin{pmatrix} -\lambda_{14} & \{0,0,0\} \\ \{0,0,0\} & -\lambda_{14} \end{pmatrix} \end{pmatrix}$$

**Δ[X₂[e₆]] =**
**Δ[X₂[e₆]] //. {λ₂ → γ₃, λ₅ → 0, λ₆ → 0, λ₁₀ → 0, λ₁₁ → 0, λ₁₄ → 0}**

$$\left( \begin{array}{ccc} \begin{pmatrix} 0 & \{0,0,0\} \\ \{0,0,0\} & 0 \end{pmatrix} & \begin{pmatrix} \lambda_4 & \{0,\lambda_7,\lambda_8\} \\ \{\lambda_9,0,0\} & 0 \end{pmatrix} & \begin{pmatrix} \lambda_{20} & \{\lambda_{22},\lambda_{23},\lambda_{24}\} \\ \{\lambda_{25},\lambda_{26},\lambda_{27}\} & \lambda_{21} \end{pmatrix} \\ \begin{pmatrix} 0 & \{0,-\lambda_7,-\lambda_8\} \\ \{-\lambda_9,0,0\} & \lambda_4 \end{pmatrix} & \begin{pmatrix} \gamma_3 & \{0,0,0\} \\ \{0,0,0\} & \gamma_3 \end{pmatrix} & \begin{pmatrix} \lambda_{12} & \{0,\lambda_{15},\lambda_{16}\} \\ \{\lambda_{17},\lambda_{18},\lambda_{19}\} & \lambda_{13} \end{pmatrix} \\ \begin{pmatrix} \lambda_{21} & \{-\lambda_{22},-\lambda_{23},-\lambda_{24}\} \\ \{-\lambda_{25},-\lambda_{26},-\lambda_{27}\} & \lambda_{20} \end{pmatrix} & \begin{pmatrix} \lambda_{13} & \{0,-\lambda_{15},-\lambda_{16}\} \\ \{-\lambda_{17},-\lambda_{18},-\lambda_{19}\} & \lambda_{12} \end{pmatrix} & \begin{pmatrix} \lambda_3 & \{0,0,0\} \\ \{0,0,0\} & \lambda_3 \end{pmatrix} \end{array} \right)$$

**U_{E₃[1]}[X₂[e₆]] == ZERO**

True

**defin[E₃[1], X₂[e₆]] // Expand**

$$\left( \begin{array}{ccc} \begin{pmatrix} 0 & \{0,0,0\} \\ \{0,0,0\} & 0 \end{pmatrix} & \begin{pmatrix} 0 & \{0,0,0\} \\ \{0,0,0\} & 0 \end{pmatrix} & \begin{pmatrix} 0 & \{0,0,0\} \\ \{0,0,0\} & 0 \end{pmatrix} \\ \begin{pmatrix} 0 & \{0,0,0\} \\ \{0,0,0\} & 0 \end{pmatrix} & \begin{pmatrix} 0 & \{0,0,0\} \\ \{0,0,0\} & 0 \end{pmatrix} & \begin{pmatrix} 0 & \{0,0,0\} \\ \{0,0,0\} & 0 \end{pmatrix} \\ \begin{pmatrix} 0 & \{0,0,0\} \\ \{0,0,0\} & 0 \end{pmatrix} & \begin{pmatrix} 0 & \{0,0,0\} \\ \{0,0,0\} & 0 \end{pmatrix} & \begin{pmatrix} \gamma_3+\lambda_3 & \{0,0,0\} \\ \{0,0,0\} & \gamma_3+\lambda_3 \end{pmatrix} \end{array} \right)$$

**Δ[X₂[e₆]] = Δ[X₂[e₆]] //. {λ₃ → -γ₃}**

$$\left( \begin{array}{ccc} \begin{pmatrix} 0 & \{0,0,0\} \\ \{0,0,0\} & 0 \end{pmatrix} & \begin{pmatrix} \lambda_4 & \{0,\lambda_7,\lambda_8\} \\ \{\lambda_9,0,0\} & 0 \end{pmatrix} & \begin{pmatrix} \lambda_{20} & \{\lambda_{22},\lambda_{23},\lambda_{24}\} \\ \{\lambda_{25},\lambda_{26},\lambda_{27}\} & \lambda_{21} \end{pmatrix} \\ \begin{pmatrix} 0 & \{0,-\lambda_7,-\lambda_8\} \\ \{-\lambda_9,0,0\} & \lambda_4 \end{pmatrix} & \begin{pmatrix} \gamma_3 & \{0,0,0\} \\ \{0,0,0\} & \gamma_3 \end{pmatrix} & \begin{pmatrix} \lambda_{12} & \{0,\lambda_{15},\lambda_{16}\} \\ \{\lambda_{17},\lambda_{18},\lambda_{19}\} & \lambda_{13} \end{pmatrix} \\ \begin{pmatrix} \lambda_{21} & \{-\lambda_{22},-\lambda_{23},-\lambda_{24}\} \\ \{-\lambda_{25},-\lambda_{26},-\lambda_{27}\} & \lambda_{20} \end{pmatrix} & \begin{pmatrix} \lambda_{13} & \{0,-\lambda_{15},-\lambda_{16}\} \\ \{-\lambda_{17},-\lambda_{18},-\lambda_{19}\} & \lambda_{12} \end{pmatrix} & \begin{pmatrix} -\gamma_3 & \{0,0,0\} \\ \{0,0,0\} & -\gamma_3 \end{pmatrix} \end{array} \right)$$

**U_{X₂[e₆]}[E₃[1]] == ZERO**

True

**defin[X₂[e₆], E₃[1]] // Expand**

$$\left( \begin{array}{ccc} \begin{pmatrix} 0 & \{0,0,0\} \\ \{0,0,0\} & 0 \end{pmatrix} & \begin{pmatrix} -\lambda_{22} & \{0,\lambda_{27},-\lambda_{26}\} \\ \{-\lambda_{21},0,0\} & 0 \end{pmatrix} & \begin{pmatrix} 0 & \{0,0,0\} \\ \{0,0,0\} & 0 \end{pmatrix} \\ \begin{pmatrix} 0 & \{0,-\lambda_{27},\lambda_{26}\} \\ \{\lambda_{21},0,0\} & -\lambda_{22} \end{pmatrix} & \begin{pmatrix} 0 & \{0,0,0\} \\ \{0,0,0\} & 0 \end{pmatrix} & \begin{pmatrix} 0 & \{0,0,0\} \\ \{0,0,0\} & 0 \end{pmatrix} \\ \begin{pmatrix} 0 & \{0,0,0\} \\ \{0,0,0\} & 0 \end{pmatrix} & \begin{pmatrix} 0 & \{0,0,0\} \\ \{0,0,0\} & 0 \end{pmatrix} & \begin{pmatrix} 0 & \{0,0,0\} \\ \{0,0,0\} & 0 \end{pmatrix} \end{array} \right)$$

$\Delta[X_2[e_6]] = \Delta[X_2[e_6]] \; //. \; \{\lambda_{21} \to 0, \; \lambda_{22} \to 0, \; \lambda_{26} \to 0, \; \lambda_{27} \to 0\}$

$$\begin{pmatrix} \begin{pmatrix} 0 & \{0,0,0\} \\ \{0,0,0\} & 0 \end{pmatrix} & \begin{pmatrix} \lambda_4 & \{0,\lambda_7,\lambda_8\} \\ \{\lambda_9,0,0\} & 0 \end{pmatrix} & \begin{pmatrix} \lambda_{20} & \{0,\lambda_{23},\lambda_{24}\} \\ \{\lambda_{25},0,0\} & 0 \end{pmatrix} \\ \begin{pmatrix} 0 & \{0,-\lambda_7,-\lambda_8\} \\ \{-\lambda_9,0,0\} & \lambda_4 \end{pmatrix} & \begin{pmatrix} \gamma_3 & \{0,0,0\} \\ \{0,0,0\} & \gamma_3 \end{pmatrix} & \begin{pmatrix} \lambda_{12} & \{0,\lambda_{15},\lambda_{16}\} \\ \{\lambda_{17},\lambda_{18},\lambda_{19}\} & \lambda_{13} \end{pmatrix} \\ \begin{pmatrix} 0 & \{0,-\lambda_{23},-\lambda_{24}\} \\ \{-\lambda_{25},0,0\} & \lambda_{20} \end{pmatrix} & \begin{pmatrix} \lambda_{13} & \{0,-\lambda_{15},-\lambda_{16}\} \\ \{-\lambda_{17},-\lambda_{18},-\lambda_{19}\} & \lambda_{12} \end{pmatrix} & \begin{pmatrix} -\gamma_3 & \{0,0,0\} \\ \{0,0,0\} & -\gamma_3 \end{pmatrix} \end{pmatrix}$$

$U_{X_1[e_2]}[X_2[e_6]] == \text{ZERO}$

True

$\text{defin}[X_1[e_2], X_2[e_6]] \; // \; \text{Expand}$

$$\begin{pmatrix} \begin{pmatrix} 0 & \{0,0,0\} \\ \{0,0,0\} & 0 \end{pmatrix} & \begin{pmatrix} 0 & \{0,0,0\} \\ \{0,0,0\} & \lambda_4-\beta_3 \end{pmatrix} & \begin{pmatrix} 0 & \{0,0,0\} \\ \{0,0,0\} & 0 \end{pmatrix} \\ \begin{pmatrix} \lambda_4-\beta_3 & \{0,0,0\} \\ \{0,0,0\} & 0 \end{pmatrix} & \begin{pmatrix} 0 & \{0,0,0\} \\ \{0,0,0\} & 0 \end{pmatrix} & \begin{pmatrix} 0 & \{0,0,0\} \\ \{0,0,0\} & 0 \end{pmatrix} \\ \begin{pmatrix} 0 & \{0,0,0\} \\ \{0,0,0\} & 0 \end{pmatrix} & \begin{pmatrix} 0 & \{0,0,0\} \\ \{0,0,0\} & 0 \end{pmatrix} & \begin{pmatrix} 0 & \{0,0,0\} \\ \{0,0,0\} & 0 \end{pmatrix} \end{pmatrix}$$

$\Delta[X_2[e_6]] = \Delta[X_2[e_6]] \; //. \; \{\lambda_4 \to \beta_3\}$

$$\begin{pmatrix} \begin{pmatrix} 0 & \{0,0,0\} \\ \{0,0,0\} & 0 \end{pmatrix} & \begin{pmatrix} \beta_3 & \{0,\lambda_7,\lambda_8\} \\ \{\lambda_9,0,0\} & 0 \end{pmatrix} & \begin{pmatrix} \lambda_{20} & \{0,\lambda_{23},\lambda_{24}\} \\ \{\lambda_{25},0,0\} & 0 \end{pmatrix} \\ \begin{pmatrix} 0 & \{0,-\lambda_7,-\lambda_8\} \\ \{-\lambda_9,0,0\} & \beta_3 \end{pmatrix} & \begin{pmatrix} \gamma_3 & \{0,0,0\} \\ \{0,0,0\} & \gamma_3 \end{pmatrix} & \begin{pmatrix} \lambda_{12} & \{0,\lambda_{15},\lambda_{16}\} \\ \{\lambda_{17},\lambda_{18},\lambda_{19}\} & \lambda_{13} \end{pmatrix} \\ \begin{pmatrix} 0 & \{0,-\lambda_{23},-\lambda_{24}\} \\ \{-\lambda_{25},0,0\} & \lambda_{20} \end{pmatrix} & \begin{pmatrix} \lambda_{13} & \{0,-\lambda_{15},-\lambda_{16}\} \\ \{-\lambda_{17},-\lambda_{18},-\lambda_{19}\} & \lambda_{12} \end{pmatrix} & \begin{pmatrix} -\gamma_3 & \{0,0,0\} \\ \{0,0,0\} & -\gamma_3 \end{pmatrix} \end{pmatrix}$$

$U_{X_2[e_6]}[X_1[e_3]] == \text{ZERO}$

True

$\text{defin}[X_2[e_6], X_1[e_3]] \; // \; \text{Expand}$

$$\begin{pmatrix} \begin{pmatrix} 0 & \{0,0,0\} \\ \{0,0,0\} & 0 \end{pmatrix} & \begin{pmatrix} 0 & \{0,0,0\} \\ \{0,0,0\} & 0 \end{pmatrix} & \begin{pmatrix} 0 & \{0,0,0\} \\ \{0,0,0\} & 0 \end{pmatrix} \\ \begin{pmatrix} 0 & \{0,0,0\} \\ \{0,0,0\} & 0 \end{pmatrix} & \begin{pmatrix} 0 & \{0,0,0\} \\ \{0,0,0\} & 0 \end{pmatrix} & \begin{pmatrix} 0 & \{0,0,0\} \\ \{\beta_2-\lambda_9,0,0\} & 0 \end{pmatrix} \\ \begin{pmatrix} 0 & \{0,0,0\} \\ \{0,0,0\} & 0 \end{pmatrix} & \begin{pmatrix} 0 & \{0,0,0\} \\ \{\lambda_9-\beta_2,0,0\} & 0 \end{pmatrix} & \begin{pmatrix} 0 & \{0,0,0\} \\ \{0,0,0\} & 0 \end{pmatrix} \end{pmatrix}$$

$\Delta[X_2[e_6]] = \Delta[X_2[e_6]] \,//\, \{\lambda_9 \to \beta_2\}$

$$\left(\begin{array}{ccc} \begin{pmatrix} 0 & \{0,0,0\} \\ \{0,0,0\} & 0 \end{pmatrix} & \begin{pmatrix} \beta_3 & \{0,\lambda_7,\lambda_8\} \\ \{\beta_2,0,0\} & 0 \end{pmatrix} & \begin{pmatrix} \lambda_{20} & \{0,\lambda_{23},\lambda_{24}\} \\ \{\lambda_{25},0,0\} & 0 \end{pmatrix} \\ \begin{pmatrix} 0 & \{0,-\lambda_7,-\lambda_8\} \\ \{-\beta_2,0,0\} & \beta_3 \end{pmatrix} & \begin{pmatrix} \gamma_3 & \{0,0,0\} \\ \{0,0,0\} & \gamma_3 \end{pmatrix} & \begin{pmatrix} \lambda_{12} & \{0,\lambda_{15},\lambda_{16}\} \\ \{\lambda_{17},\lambda_{18},\lambda_{19}\} & \lambda_{13} \end{pmatrix} \\ \begin{pmatrix} 0 & \{0,-\lambda_{23},-\lambda_{24}\} \\ \{-\lambda_{25},0,0\} & \lambda_{20} \end{pmatrix} & \begin{pmatrix} \lambda_{13} & \{0,-\lambda_{15},-\lambda_{16}\} \\ \{-\lambda_{17},-\lambda_{18},-\lambda_{19}\} & \lambda_{12} \end{pmatrix} & \begin{pmatrix} -\gamma_3 & \{0,0,0\} \\ \{0,0,0\} & -\gamma_3 \end{pmatrix} \end{array}\right)$$

$U_{X_1[e_8]}[X_2[e_6]] == \text{ZERO}$

True

$\text{defin}[X_1[e_8], X_2[e_6]] \,//\, \text{Expand}$

$$\left(\begin{array}{ccc} \begin{pmatrix} 0 & \{0,0,0\} \\ \{0,0,0\} & 0 \end{pmatrix} & \begin{pmatrix} 0 & \{0,0,0\} \\ \{0,0,\beta_7-\lambda_8\} & 0 \end{pmatrix} & \begin{pmatrix} 0 & \{0,0,0\} \\ \{0,0,0\} & 0 \end{pmatrix} \\ \begin{pmatrix} 0 & \{0,0,0\} \\ \{0,0,\lambda_8-\beta_7\} & 0 \end{pmatrix} & \begin{pmatrix} 0 & \{0,0,0\} \\ \{0,0,0\} & 0 \end{pmatrix} & \begin{pmatrix} 0 & \{0,0,0\} \\ \{0,0,0\} & 0 \end{pmatrix} \\ \begin{pmatrix} 0 & \{0,0,0\} \\ \{0,0,0\} & 0 \end{pmatrix} & \begin{pmatrix} 0 & \{0,0,0\} \\ \{0,0,0\} & 0 \end{pmatrix} & \begin{pmatrix} 0 & \{0,0,0\} \\ \{0,0,0\} & 0 \end{pmatrix} \end{array}\right)$$

$\Delta[X_2[e_6]] = \Delta[X_2[e_6]] \,//\, \{\lambda_8 \to \beta_7\}$

$$\left(\begin{array}{ccc} \begin{pmatrix} 0 & \{0,0,0\} \\ \{0,0,0\} & 0 \end{pmatrix} & \begin{pmatrix} \beta_3 & \{0,\lambda_7,\beta_7\} \\ \{\beta_2,0,0\} & 0 \end{pmatrix} & \begin{pmatrix} \lambda_{20} & \{0,\lambda_{23},\lambda_{24}\} \\ \{\lambda_{25},0,0\} & 0 \end{pmatrix} \\ \begin{pmatrix} 0 & \{0,-\lambda_7,-\beta_7\} \\ \{-\beta_2,0,0\} & \beta_3 \end{pmatrix} & \begin{pmatrix} \gamma_3 & \{0,0,0\} \\ \{0,0,0\} & \gamma_3 \end{pmatrix} & \begin{pmatrix} \lambda_{12} & \{0,\lambda_{15},\lambda_{16}\} \\ \{\lambda_{17},\lambda_{18},\lambda_{19}\} & \lambda_{13} \end{pmatrix} \\ \begin{pmatrix} 0 & \{0,-\lambda_{23},-\lambda_{24}\} \\ \{-\lambda_{25},0,0\} & \lambda_{20} \end{pmatrix} & \begin{pmatrix} \lambda_{13} & \{0,-\lambda_{15},-\lambda_{16}\} \\ \{-\lambda_{17},-\lambda_{18},-\lambda_{19}\} & \lambda_{12} \end{pmatrix} & \begin{pmatrix} -\gamma_3 & \{0,0,0\} \\ \{0,0,0\} & -\gamma_3 \end{pmatrix} \end{array}\right)$$

$T[E_1[1], E_1[1], X_2[e_6]] == \text{ZERO}$

True

$\text{Leib}[E_1[1], E_1[1], X_2[e_6]] \,//\, \text{Expand}$

$$\left(\begin{array}{ccc} \begin{pmatrix} 0 & \{0,0,0\} \\ \{0,0,0\} & 0 \end{pmatrix} & \begin{pmatrix} 0 & \{0,\beta_8+\lambda_7,0\} \\ \{0,0,0\} & 0 \end{pmatrix} & \begin{pmatrix} \alpha_3+\lambda_{20} & \{0,\lambda_{23}-\alpha_8,\alpha_7+\lambda_{24}\} \\ \{\alpha_2+\lambda_{25},0,0\} & 0 \end{pmatrix} \\ \begin{pmatrix} 0 & \{0,-\beta_8-\lambda_7,0\} \\ \{0,0,0\} & 0 \end{pmatrix} & \begin{pmatrix} 0 & \{0,0,0\} \\ \{0,0,0\} & 0 \end{pmatrix} & \begin{pmatrix} 0 & \{0,0,0\} \\ \{0,0,0\} & 0 \end{pmatrix} \\ \begin{pmatrix} 0 & \{0,\alpha_8-\lambda_{23},-\alpha_7-\lambda_{24}\} \\ \{-\alpha_2-\lambda_{25},0,0\} & \alpha_3+\lambda_{20} \end{pmatrix} & \begin{pmatrix} 0 & \{0,0,0\} \\ \{0,0,0\} & 0 \end{pmatrix} & \begin{pmatrix} 0 & \{0,0,0\} \\ \{0,0,0\} & 0 \end{pmatrix} \end{array}\right)$$

$\Delta[X_2[e_6]] =$
$\quad \Delta[X_2[e_6]] \ //. \ \{\lambda_7 \to -\beta_8, \ \lambda_{20} \to -\alpha_3, \ \lambda_{23} \to \alpha_8, \ \lambda_{24} \to -\alpha_7, \ \lambda_{25} \to -\alpha_2\}$

$$\begin{pmatrix}
\begin{pmatrix} 0 & \{0,0,0\} \\ \{0,0,0\} & 0 \end{pmatrix} & \begin{pmatrix} \beta_3 & \{0,-\beta_8,\beta_7\} \\ \{\beta_2,0,0\} & 0 \end{pmatrix} & \begin{pmatrix} -\alpha_3 & \{0,\alpha_8,-\alpha_7\} \\ \{-\alpha_2,0,0\} & 0 \end{pmatrix} \\
\begin{pmatrix} 0 & \{0,\beta_8,-\beta_7\} \\ \{-\beta_2,0,0\} & \beta_3 \end{pmatrix} & \begin{pmatrix} \gamma_3 & \{0,0,0\} \\ \{0,0,0\} & \gamma_3 \end{pmatrix} & \begin{pmatrix} \lambda_{12} & \{0,\lambda_{15},\lambda_{16}\} \\ \{\lambda_{17},\lambda_{18},\lambda_{19}\} & \lambda_{13} \end{pmatrix} \\
\begin{pmatrix} 0 & \{0,-\alpha_8,\alpha_7\} \\ \{\alpha_2,0,0\} & -\alpha_3 \end{pmatrix} & \begin{pmatrix} \lambda_{13} & \{0,-\lambda_{15},-\lambda_{16}\} \\ \{-\lambda_{17},-\lambda_{18},-\lambda_{19}\} & \lambda_{12} \end{pmatrix} & \begin{pmatrix} -\gamma_3 & \{0,0,0\} \\ \{0,0,0\} & -\gamma_3 \end{pmatrix}
\end{pmatrix}$$

**T[X₁[e₁], X₁[e₃], X₂[e₆]] == ZERO**

True

**Leib[X₁[e₁], X₁[e₃], X₂[e₆]] // Expand**

$$\begin{pmatrix}
\begin{pmatrix} 0 & \{0,0,0\} \\ \{0,0,0\} & 0 \end{pmatrix} & \begin{pmatrix} 0 & \{0,0,0\} \\ \{0,0,0\} & 0 \end{pmatrix} & \begin{pmatrix} 0 & \{0,0,0\} \\ \{0,0,0\} & 0 \end{pmatrix} \\
\begin{pmatrix} 0 & \{0,0,0\} \\ \{0,0,0\} & 0 \end{pmatrix} & \begin{pmatrix} 0 & \{0,0,0\} \\ \{0,0,0\} & 0 \end{pmatrix} & \begin{pmatrix} 0 & \{0,0,0\} \\ \{0,-\delta_6-\lambda_{16},\lambda_{15}-\delta_7\} & 0 \end{pmatrix} \\
\begin{pmatrix} 0 & \{0,0,0\} \\ \{0,0,0\} & 0 \end{pmatrix} & \begin{pmatrix} 0 & \{0,0,0\} \\ \{0,\delta_6+\lambda_{16},\delta_7-\lambda_{15}\} & 0 \end{pmatrix} & \begin{pmatrix} 0 & \{0,0,0\} \\ \{0,0,0\} & 0 \end{pmatrix}
\end{pmatrix}$$

$\Delta[X_2[e_6]] = \Delta[X_2[e_6]] \ //. \ \{\lambda_{15} \to \delta_7, \ \lambda_{16} \to -\delta_6\}$

$$\begin{pmatrix}
\begin{pmatrix} 0 & \{0,0,0\} \\ \{0,0,0\} & 0 \end{pmatrix} & \begin{pmatrix} \beta_3 & \{0,-\beta_8,\beta_7\} \\ \{\beta_2,0,0\} & 0 \end{pmatrix} & \begin{pmatrix} -\alpha_3 & \{0,\alpha_8,-\alpha_7\} \\ \{-\alpha_2,0,0\} & 0 \end{pmatrix} \\
\begin{pmatrix} 0 & \{0,\beta_8,-\beta_7\} \\ \{-\beta_2,0,0\} & \beta_3 \end{pmatrix} & \begin{pmatrix} \gamma_3 & \{0,0,0\} \\ \{0,0,0\} & \gamma_3 \end{pmatrix} & \begin{pmatrix} \lambda_{12} & \{0,\delta_7,-\delta_6\} \\ \{\lambda_{17},\lambda_{18},\lambda_{19}\} & \lambda_{13} \end{pmatrix} \\
\begin{pmatrix} 0 & \{0,-\alpha_8,\alpha_7\} \\ \{\alpha_2,0,0\} & -\alpha_3 \end{pmatrix} & \begin{pmatrix} \lambda_{13} & \{0,-\delta_7,\delta_6\} \\ \{-\lambda_{17},-\lambda_{18},-\lambda_{19}\} & \lambda_{12} \end{pmatrix} & \begin{pmatrix} -\gamma_3 & \{0,0,0\} \\ \{0,0,0\} & -\gamma_3 \end{pmatrix}
\end{pmatrix}$$

**T[X₁[e₁], X₁[e₆], X₂[e₆]] == ZERO**

True

**Leib[X₁[e₁], X₁[e₆], X₂[e₆]] // Expand**

$$\begin{pmatrix}
\begin{pmatrix} 0 & \{0,0,0\} \\ \{0,0,0\} & 0 \end{pmatrix} & \begin{pmatrix} 0 & \{0,0,0\} \\ \{0,0,0\} & 0 \end{pmatrix} & \begin{pmatrix} 0 & \{0,0,0\} \\ \{0,0,0\} & 0 \end{pmatrix} \\
\begin{pmatrix} 0 & \{0,0,0\} \\ \{0,0,0\} & 0 \end{pmatrix} & \begin{pmatrix} 0 & \{0,0,0\} \\ \{0,0,0\} & 0 \end{pmatrix} & \begin{pmatrix} 0 & \{0,0,0\} \\ \{\delta_2+\lambda_{12},0,0\} & 0 \end{pmatrix} \\
\begin{pmatrix} 0 & \{0,0,0\} \\ \{0,0,0\} & 0 \end{pmatrix} & \begin{pmatrix} 0 & \{0,0,0\} \\ \{-\delta_2-\lambda_{12},0,0\} & 0 \end{pmatrix} & \begin{pmatrix} 0 & \{0,0,0\} \\ \{0,0,0\} & 0 \end{pmatrix}
\end{pmatrix}$$

$\Delta[X_2[e_6]] = \Delta[X_2[e_6]] \,/\!/.\, \{\lambda_{12} \to -\delta_2\}$

$$\begin{pmatrix}
\begin{pmatrix} 0 & \{0,0,0\} \\ \{0,0,0\} & 0 \end{pmatrix} & \begin{pmatrix} \beta_3 & \{0,-\beta_8,\beta_7\} \\ \{\beta_2,0,0\} & 0 \end{pmatrix} & \begin{pmatrix} -\alpha_3 & \{0,\alpha_8,-\alpha_7\} \\ \{-\alpha_2,0,0\} & 0 \end{pmatrix} \\
\begin{pmatrix} 0 & \{0,\beta_8,-\beta_7\} \\ \{-\beta_2,0,0\} & \beta_3 \end{pmatrix} & \begin{pmatrix} \gamma_3 & \{0,0,0\} \\ \{0,0,0\} & \gamma_3 \end{pmatrix} & \begin{pmatrix} -\delta_2 & \{0,\delta_7,-\delta_6\} \\ \{\lambda_{17},\lambda_{18},\lambda_{19}\} & \lambda_{13} \end{pmatrix} \\
\begin{pmatrix} 0 & \{0,-\alpha_8,\alpha_7\} \\ \{\alpha_2,0,0\} & -\alpha_3 \end{pmatrix} & \begin{pmatrix} \lambda_{13} & \{0,-\delta_7,\delta_6\} \\ \{-\lambda_{17},-\lambda_{18},-\lambda_{19}\} & -\delta_2 \end{pmatrix} & \begin{pmatrix} -\gamma_3 & \{0,0,0\} \\ \{0,0,0\} & -\gamma_3 \end{pmatrix}
\end{pmatrix}$$

**T[X$_1$[e$_2$], X$_1$[e$_8$], X$_2$[e$_6$]] == -X$_2$[e$_4$]**

True

**$\Delta$[X$_2$[e$_4$]] + Leib[X$_1$[e$_2$], X$_1$[e$_8$], X$_2$[e$_6$]] // Expand**

$$\begin{pmatrix}
\begin{pmatrix} 0 & \{0,0,0\} \\ \{0,0,0\} & 0 \end{pmatrix} & \begin{pmatrix} 0 & \{0,0,0\} \\ \{0,0,0\} & 0 \end{pmatrix} & \begin{pmatrix} 0 & \{0,0,0\} \\ \{0,0,0\} & 0 \end{pmatrix} \\
\begin{pmatrix} 0 & \{0,0,0\} \\ \{0,0,0\} & 0 \end{pmatrix} & \begin{pmatrix} 0 & \{0,0,0\} \\ \{0,0,0\} & 0 \end{pmatrix} & \begin{pmatrix} 0 & \{\epsilon_1+\lambda_{18},\xi+\delta_1-\eta_1-\lambda_{17},0\} \\ \{0,0,0\} & 0 \end{pmatrix} \\
\begin{pmatrix} 0 & \{0,0,0\} \\ \{0,0,0\} & 0 \end{pmatrix} & \begin{pmatrix} 0 & \{-\epsilon_1-\lambda_{18},-\xi-\delta_1+\eta_1+\lambda_{17},0\} \\ \{0,0,0\} & 0 \end{pmatrix} & \begin{pmatrix} 0 & \{0,0,0\} \\ \{0,0,0\} & 0 \end{pmatrix}
\end{pmatrix}$$

$\Delta[X_2[e_6]] = \Delta[X_2[e_6]] \,/\!/.\, \{\lambda_{18} \to -\epsilon_1,\, \lambda_{17} \to \xi+\delta_1-\eta_1\}$

$$\begin{pmatrix}
\begin{pmatrix} 0 & \{0,0,0\} \\ \{0,0,0\} & 0 \end{pmatrix} & \begin{pmatrix} \beta_3 & \{0,-\beta_8,\beta_7\} \\ \{\beta_2,0,0\} & 0 \end{pmatrix} & \begin{pmatrix} -\alpha_3 & \{0,\alpha_8,-\alpha_7\} \\ \{-\alpha_2,0,0\} & 0 \end{pmatrix} \\
\begin{pmatrix} 0 & \{0,\beta_8,-\beta_7\} \\ \{-\beta_2,0,0\} & \beta_3 \end{pmatrix} & \begin{pmatrix} \gamma_3 & \{0,0,0\} \\ \{0,0,0\} & \gamma_3 \end{pmatrix} & \begin{pmatrix} -\delta_2 & \{0,\delta_7,-\delta_6 \\ \{\xi+\delta_1-\eta_1,-\epsilon_1,\lambda_{19}\} & \lambda_{13} \end{pmatrix} \\
\begin{pmatrix} 0 & \{0,-\alpha_8,\alpha_7\} \\ \{\alpha_2,0,0\} & -\alpha_3 \end{pmatrix} & \begin{pmatrix} \lambda_{13} & \{0,-\delta_7,\delta_6\} \\ \{-\xi-\delta_1+\eta_1,\epsilon_1,-\lambda_{19}\} & -\delta_2 \end{pmatrix} & \begin{pmatrix} -\gamma_3 & \{0,0,0\} \\ \{0,0,0\} & -\gamma_3 \end{pmatrix}
\end{pmatrix}$$

**T[X$_1$[e$_5$], X$_1$[e$_8$], X$_2$[e$_6$]] == -X$_2$[e$_6$]**

True

**$\Delta$[X$_2$[e$_6$]] + Leib[X$_1$[e$_5$], X$_1$[e$_8$], X$_2$[e$_6$]] // Expand**

$$\begin{pmatrix}
\begin{pmatrix} 0 & \{0,0,0\} \\ \{0,0,0\} & 0 \end{pmatrix} & \begin{pmatrix} 0 & \{0,0,0\} \\ \{0,0,0\} & 0 \end{pmatrix} & \begin{pmatrix} 0 & \{0,0,0\} \\ \{0,0,0\} & 0 \end{pmatrix} \\
\begin{pmatrix} 0 & \{0,0,0\} \\ \{0,0,0\} & 0 \end{pmatrix} & \begin{pmatrix} 0 & \{0,0,0\} \\ \{0,0,0\} & 0 \end{pmatrix} & \begin{pmatrix} 0 & \{0,0,0\} \\ \{0,0,\lambda_{19}+\phi_1\} & \lambda_{13}-\epsilon_4 \end{pmatrix} \\
\begin{pmatrix} 0 & \{0,0,0\} \\ \{0,0,0\} & 0 \end{pmatrix} & \begin{pmatrix} \lambda_{13}-\epsilon_4 & \{0,0,0\} \\ \{0,0,-\lambda_{19}-\phi_1\} & 0 \end{pmatrix} & \begin{pmatrix} 0 & \{0,0,0\} \\ \{0,0,0\} & 0 \end{pmatrix}
\end{pmatrix}$$

**Δ[X₂[e₆]] = Δ[X₂[e₆]] //. {λ₁₉ → -ϕ₁, λ₁₃ → ϵ₄}**

$$\begin{pmatrix} \begin{pmatrix} 0 & \{0,0,0\} \\ \{0,0,0\} & 0 \end{pmatrix} & \begin{pmatrix} \beta_3 & \{0,-\beta_8,\beta_7\} \\ \{\beta_2,0,0\} & 0 \end{pmatrix} & \begin{pmatrix} -\alpha_3 & \{0,\alpha_8,-\alpha_7\} \\ \{-\alpha_2,0,0\} & 0 \end{pmatrix} \\ \begin{pmatrix} 0 & \{0,\beta_8,-\beta_7\} \\ \{-\beta_2,0,0\} & \beta_3 \end{pmatrix} & \begin{pmatrix} \gamma_3 & \{0,0,0\} \\ \{0,0,0\} & \gamma_3 \end{pmatrix} & \begin{pmatrix} -\delta_2 & \{0,\delta_7,-\delta_6\} \\ \{\xi+\delta_1-\eta_1,-\epsilon_1,-\phi_1\} & \epsilon_4 \end{pmatrix} \\ \begin{pmatrix} 0 & \{0,-\alpha_8,\alpha_7\} \\ \{\alpha_2,0,0\} & -\alpha_3 \end{pmatrix} & \begin{pmatrix} \epsilon_4 & \{0,-\delta_7,\delta_6\} \\ \{-\xi-\delta_1+\eta_1,\epsilon_1,\phi_1\} & -\delta_2 \end{pmatrix} & \begin{pmatrix} -\gamma_3 & \{0,0,0\} \\ \{0,0,0\} & -\gamma_3 \end{pmatrix} \end{pmatrix}$$

**Variables[Δ[X₂[e₆]]]**

$\{\xi, \alpha_2, \alpha_3, \alpha_7, \alpha_8, \beta_2, \beta_3, \beta_7, \beta_8, \gamma_3, \delta_1, \delta_2, \delta_6, \delta_7, \epsilon_1, \epsilon_4, \eta_1, \phi_1\}$

- **X₂[e₇]**

    **Δ[X₂[e₇]] = generic;**

    **U_{E₁[1]}[X₂[e₇]] == ZERO**

    True

    **(defin[E₁[1], X₂[e₇]] // Expand)**

$$\begin{pmatrix} \begin{pmatrix} \lambda_1 & \{0,0,0\} \\ \{0,0,0\} & \lambda_1 \end{pmatrix} & \begin{pmatrix} 0 & \{0,0,0\} \\ \{0,0,0\} & 0 \end{pmatrix} & \begin{pmatrix} 0 & \{0,0,0\} \\ \{0,0,0\} & 0 \end{pmatrix} \\ \begin{pmatrix} 0 & \{0,0,0\} \\ \{0,0,0\} & 0 \end{pmatrix} & \begin{pmatrix} 0 & \{0,0,0\} \\ \{0,0,0\} & 0 \end{pmatrix} & \begin{pmatrix} 0 & \{0,0,0\} \\ \{0,0,0\} & 0 \end{pmatrix} \\ \begin{pmatrix} 0 & \{0,0,0\} \\ \{0,0,0\} & 0 \end{pmatrix} & \begin{pmatrix} 0 & \{0,0,0\} \\ \{0,0,0\} & 0 \end{pmatrix} & \begin{pmatrix} 0 & \{0,0,0\} \\ \{0,0,0\} & 0 \end{pmatrix} \end{pmatrix}$$

    **Δ[X₂[e₇]] = Δ[X₂[e₇]] //. {λ₁ → 0}**

$$\begin{pmatrix} \begin{pmatrix} 0 & \{0,0,0\} \\ \{0,0,0\} & 0 \end{pmatrix} & \begin{pmatrix} \lambda_4 & \{\lambda_6,\lambda_7,\lambda_8\} \\ \{\lambda_9,\lambda_{10},\lambda_{11}\} & \lambda_5 \end{pmatrix} & \begin{pmatrix} \lambda_{20} & \{\lambda_{22},\lambda_{23},\lambda_{24}\} \\ \{\lambda_{25},\lambda_{26},\lambda_{27}\} & \lambda_{21} \end{pmatrix} \\ \begin{pmatrix} \lambda_5 & \{-\lambda_6,-\lambda_7,-\lambda_8\} \\ \{-\lambda_9,-\lambda_{10},-\lambda_{11}\} & \lambda_4 \end{pmatrix} & \begin{pmatrix} \lambda_2 & \{0,0,0\} \\ \{0,0,0\} & \lambda_2 \end{pmatrix} & \begin{pmatrix} \lambda_{12} & \{\lambda_{14},\lambda_{15},\lambda_{16}\} \\ \{\lambda_{17},\lambda_{18},\lambda_{19}\} & \lambda_{13} \end{pmatrix} \\ \begin{pmatrix} \lambda_{21} & \{-\lambda_{22},-\lambda_{23},-\lambda_{24}\} \\ \{-\lambda_{25},-\lambda_{26},-\lambda_{27}\} & \lambda_{20} \end{pmatrix} & \begin{pmatrix} \lambda_{13} & \{-\lambda_{14},-\lambda_{15},-\lambda_{16}\} \\ \{-\lambda_{17},-\lambda_{18},-\lambda_{19}\} & \lambda_{12} \end{pmatrix} & \begin{pmatrix} \lambda_3 & \{0,0,0\} \\ \{0,0,0\} & \lambda_3 \end{pmatrix} \end{pmatrix}$$

    **U_{X₂[e₇]}[E₂[1]] == ZERO**

    True

    **defin[X₂[e₇], E₂[1]] // Expand**

$$\begin{pmatrix} \begin{pmatrix} 0 & \{0,0,0\} \\ \{0,0,0\} & 0 \end{pmatrix} & \begin{pmatrix} 0 & \{0,0,0\} \\ \{0,0,0\} & 0 \end{pmatrix} & \begin{pmatrix} \lambda_7 & \{\lambda_{11},0,-\lambda_9\} \\ \{0,\lambda_5,0\} & 0 \end{pmatrix} \\ \begin{pmatrix} 0 & \{0,0,0\} \\ \{0,0,0\} & 0 \end{pmatrix} & \begin{pmatrix} 0 & \{0,0,0\} \\ \{0,0,0\} & 0 \end{pmatrix} & \begin{pmatrix} 0 & \{0,0,0\} \\ \{0,\lambda_2-\gamma_4,0\} & 0 \end{pmatrix} \\ \begin{pmatrix} 0 & \{-\lambda_{11},0,\lambda_9\} \\ \{0,-\lambda_5,0\} & \lambda_7 \end{pmatrix} & \begin{pmatrix} 0 & \{0,0,0\} \\ \{0,\gamma_4-\lambda_2,0\} & 0 \end{pmatrix} & \begin{pmatrix} -\lambda_{15} & \{0,0,0\} \\ \{0,0,0\} & -\lambda_{15} \end{pmatrix} \end{pmatrix}$$

```
Δ[X₂[e₇]] =
  Δ[X₂[e₇]] //. {λ₂ → γ₄, λ₅ → 0, λ₇ → 0, λ₉ → 0, λ₁₁ → 0, λ₁₅ → 0}
```

$$\left(\begin{array}{ccc} \begin{pmatrix} 0 & \{0,0,0\} \\ \{0,0,0\} & 0 \end{pmatrix} & \begin{pmatrix} \lambda_4 & \{\lambda_6, 0, \lambda_8\} \\ \{0, \lambda_{10}, 0\} & 0 \end{pmatrix} & \begin{pmatrix} \lambda_{20} & \{\lambda_{22}, \lambda_{23}, \lambda_{24}\} \\ \{\lambda_{25}, \lambda_{26}, \lambda_{27}\} & \lambda_{21} \end{pmatrix} \\ \begin{pmatrix} 0 & \{-\lambda_6, 0, -\lambda_8\} \\ \{0, -\lambda_{10}, 0\} & \lambda_4 \end{pmatrix} & \begin{pmatrix} \gamma_4 & \{0,0,0\} \\ \{0,0,0\} & \gamma_4 \end{pmatrix} & \begin{pmatrix} \lambda_{12} & \{\lambda_{14}, 0, \lambda_{16}\} \\ \{\lambda_{17}, \lambda_{18}, \lambda_{19}\} & \lambda_{13} \end{pmatrix} \\ \begin{pmatrix} \lambda_{21} & \{-\lambda_{22}, -\lambda_{23}, -\lambda_{24}\} \\ \{-\lambda_{25}, -\lambda_{26}, -\lambda_{27}\} & \lambda_{20} \end{pmatrix} & \begin{pmatrix} \lambda_{13} & \{-\lambda_{14}, 0, -\lambda_{16}\} \\ \{-\lambda_{17}, -\lambda_{18}, -\lambda_{19}\} & \lambda_{12} \end{pmatrix} & \begin{pmatrix} \lambda_3 & \{0,0,0\} \\ \{0,0,0\} & \lambda_3 \end{pmatrix} \end{array}\right)$$

**U_{E₃[1]}[X₂[e₇]] == ZERO**

True

**defin[E₃[1], X₂[e₇]] // Expand**

$$\left(\begin{array}{ccc} \begin{pmatrix} 0 & \{0,0,0\} \\ \{0,0,0\} & 0 \end{pmatrix} & \begin{pmatrix} 0 & \{0,0,0\} \\ \{0,0,0\} & 0 \end{pmatrix} & \begin{pmatrix} 0 & \{0,0,0\} \\ \{0,0,0\} & 0 \end{pmatrix} \\ \begin{pmatrix} 0 & \{0,0,0\} \\ \{0,0,0\} & 0 \end{pmatrix} & \begin{pmatrix} 0 & \{0,0,0\} \\ \{0,0,0\} & 0 \end{pmatrix} & \begin{pmatrix} 0 & \{0,0,0\} \\ \{0,0,0\} & 0 \end{pmatrix} \\ \begin{pmatrix} 0 & \{0,0,0\} \\ \{0,0,0\} & 0 \end{pmatrix} & \begin{pmatrix} 0 & \{0,0,0\} \\ \{0,0,0\} & 0 \end{pmatrix} & \begin{pmatrix} \gamma_4+\lambda_3 & \{0,0,0\} \\ \{0,0,0\} & \gamma_4+\lambda_3 \end{pmatrix} \end{array}\right)$$

**Δ[X₂[e₇]] = Δ[X₂[e₇]] //. {λ₃ → −γ₄}**

$$\left(\begin{array}{ccc} \begin{pmatrix} 0 & \{0,0,0\} \\ \{0,0,0\} & 0 \end{pmatrix} & \begin{pmatrix} \lambda_4 & \{\lambda_6, 0, \lambda_8\} \\ \{0, \lambda_{10}, 0\} & 0 \end{pmatrix} & \begin{pmatrix} \lambda_{20} & \{\lambda_{22}, \lambda_{23}, \lambda_{24}\} \\ \{\lambda_{25}, \lambda_{26}, \lambda_{27}\} & \lambda_{21} \end{pmatrix} \\ \begin{pmatrix} 0 & \{-\lambda_6, 0, -\lambda_8\} \\ \{0, -\lambda_{10}, 0\} & \lambda_4 \end{pmatrix} & \begin{pmatrix} \gamma_4 & \{0,0,0\} \\ \{0,0,0\} & \gamma_4 \end{pmatrix} & \begin{pmatrix} \lambda_{12} & \{\lambda_{14}, 0, \lambda_{16}\} \\ \{\lambda_{17}, \lambda_{18}, \lambda_{19}\} & \lambda_{13} \end{pmatrix} \\ \begin{pmatrix} \lambda_{21} & \{-\lambda_{22}, -\lambda_{23}, -\lambda_{24}\} \\ \{-\lambda_{25}, -\lambda_{26}, -\lambda_{27}\} & \lambda_{20} \end{pmatrix} & \begin{pmatrix} \lambda_{13} & \{-\lambda_{14}, 0, -\lambda_{16}\} \\ \{-\lambda_{17}, -\lambda_{18}, -\lambda_{19}\} & \lambda_{12} \end{pmatrix} & \begin{pmatrix} -\gamma_4 & \{0,0,0\} \\ \{0,0,0\} & -\gamma_4 \end{pmatrix} \end{array}\right)$$

**U_{X₂[e₇]}[E₃[1]] == ZERO**

True

**defin[X₂[e₇], E₃[1]] // Expand**

$$\left(\begin{array}{ccc} \begin{pmatrix} 0 & \{0,0,0\} \\ \{0,0,0\} & 0 \end{pmatrix} & \begin{pmatrix} -\lambda_{23} & \{-\lambda_{27}, 0, \lambda_{25}\} \\ \{0, -\lambda_{21}, 0\} & 0 \end{pmatrix} & \begin{pmatrix} 0 & \{0,0,0\} \\ \{0,0,0\} & 0 \end{pmatrix} \\ \begin{pmatrix} 0 & \{\lambda_{27}, 0, -\lambda_{25}\} \\ \{0, \lambda_{21}, 0\} & -\lambda_{23} \end{pmatrix} & \begin{pmatrix} 0 & \{0,0,0\} \\ \{0,0,0\} & 0 \end{pmatrix} & \begin{pmatrix} 0 & \{0,0,0\} \\ \{0,0,0\} & 0 \end{pmatrix} \\ \begin{pmatrix} 0 & \{0,0,0\} \\ \{0,0,0\} & 0 \end{pmatrix} & \begin{pmatrix} 0 & \{0,0,0\} \\ \{0,0,0\} & 0 \end{pmatrix} & \begin{pmatrix} 0 & \{0,0,0\} \\ \{0,0,0\} & 0 \end{pmatrix} \end{array}\right)$$

$\Delta[X_2[e_7]] = \Delta[X_2[e_7]] //. \{\lambda_{21} \to 0, \lambda_{23} \to 0, \lambda_{25} \to 0, \lambda_{27} \to 0\}$

$$\begin{pmatrix} \begin{pmatrix} 0 & \{0,0,0\} \\ \{0,0,0\} & 0 \end{pmatrix} & \begin{pmatrix} \lambda_4 & \{\lambda_6,0,\lambda_8\} \\ \{0,\lambda_{10},0\} & 0 \end{pmatrix} & \begin{pmatrix} \lambda_{20} & \{\lambda_{22},0,\lambda_{24}\} \\ \{0,\lambda_{26},0\} & 0 \end{pmatrix} \\ \begin{pmatrix} 0 & \{-\lambda_6,0,-\lambda_8\} \\ \{0,-\lambda_{10},0\} & \lambda_4 \end{pmatrix} & \begin{pmatrix} \gamma_4 & \{0,0,0\} \\ \{0,0,0\} & \gamma_4 \end{pmatrix} & \begin{pmatrix} \lambda_{12} & \{\lambda_{14},0,\lambda_{16}\} \\ \{\lambda_{17},\lambda_{18},\lambda_{19}\} & \lambda_{13} \end{pmatrix} \\ \begin{pmatrix} 0 & \{-\lambda_{22},0,-\lambda_{24}\} \\ \{0,-\lambda_{26},0\} & \lambda_{20} \end{pmatrix} & \begin{pmatrix} \lambda_{13} & \{-\lambda_{14},0,-\lambda_{16}\} \\ \{-\lambda_{17},-\lambda_{18},-\lambda_{19}\} & \lambda_{12} \end{pmatrix} & \begin{pmatrix} -\gamma_4 & \{0,0,0\} \\ \{0,0,0\} & -\gamma_4 \end{pmatrix} \end{pmatrix}$$

$U_{X_1[e_2]}[X_2[e_7]] == \text{ZERO}$

True

$\text{defin}[X_1[e_2], X_2[e_7]] // \text{Expand}$

$$\begin{pmatrix} \begin{pmatrix} 0 & \{0,0,0\} \\ \{0,0,0\} & 0 \end{pmatrix} & \begin{pmatrix} 0 & \{0,0,0\} \\ \{0,0,0\} & \lambda_4-\beta_4 \end{pmatrix} & \begin{pmatrix} 0 & \{0,0,0\} \\ \{0,0,0\} & 0 \end{pmatrix} \\ \begin{pmatrix} \lambda_4-\beta_4 & \{0,0,0\} \\ \{0,0,0\} & 0 \end{pmatrix} & \begin{pmatrix} 0 & \{0,0,0\} \\ \{0,0,0\} & 0 \end{pmatrix} & \begin{pmatrix} 0 & \{0,0,0\} \\ \{0,0,0\} & 0 \end{pmatrix} \\ \begin{pmatrix} 0 & \{0,0,0\} \\ \{0,0,0\} & 0 \end{pmatrix} & \begin{pmatrix} 0 & \{0,0,0\} \\ \{0,0,0\} & 0 \end{pmatrix} & \begin{pmatrix} 0 & \{0,0,0\} \\ \{0,0,0\} & 0 \end{pmatrix} \end{pmatrix}$$

$\Delta[X_2[e_7]] = \Delta[X_2[e_7]] //. \{\lambda_4 \to \beta_4\}$

$$\begin{pmatrix} \begin{pmatrix} 0 & \{0,0,0\} \\ \{0,0,0\} & 0 \end{pmatrix} & \begin{pmatrix} \beta_4 & \{\lambda_6,0,\lambda_8\} \\ \{0,\lambda_{10},0\} & 0 \end{pmatrix} & \begin{pmatrix} \lambda_{20} & \{\lambda_{22},0,\lambda_{24}\} \\ \{0,\lambda_{26},0\} & 0 \end{pmatrix} \\ \begin{pmatrix} 0 & \{-\lambda_6,0,-\lambda_8\} \\ \{0,-\lambda_{10},0\} & \beta_4 \end{pmatrix} & \begin{pmatrix} \gamma_4 & \{0,0,0\} \\ \{0,0,0\} & \gamma_4 \end{pmatrix} & \begin{pmatrix} \lambda_{12} & \{\lambda_{14},0,\lambda_{16}\} \\ \{\lambda_{17},\lambda_{18},\lambda_{19}\} & \lambda_{13} \end{pmatrix} \\ \begin{pmatrix} 0 & \{-\lambda_{22},0,-\lambda_{24}\} \\ \{0,-\lambda_{26},0\} & \lambda_{20} \end{pmatrix} & \begin{pmatrix} \lambda_{13} & \{-\lambda_{14},0,-\lambda_{16}\} \\ \{-\lambda_{17},-\lambda_{18},-\lambda_{19}\} & \lambda_{12} \end{pmatrix} & \begin{pmatrix} -\gamma_4 & \{0,0,0\} \\ \{0,0,0\} & -\gamma_4 \end{pmatrix} \end{pmatrix}$$

$U_{X_1[e_4]}[X_2[e_7]] == \text{ZERO}$

True

$\text{defin}[X_1[e_4], X_2[e_7]] // \text{Expand}$

$$\begin{pmatrix} \begin{pmatrix} 0 & \{0,0,0\} \\ \{0,0,0\} & 0 \end{pmatrix} & \begin{pmatrix} 0 & \{0,\beta_2-\lambda_{10},0\} \\ \{0,0,0\} & 0 \end{pmatrix} & \begin{pmatrix} 0 & \{0,0,0\} \\ \{0,0,0\} & 0 \end{pmatrix} \\ \begin{pmatrix} 0 & \{0,\lambda_{10}-\beta_2,0\} \\ \{0,0,0\} & 0 \end{pmatrix} & \begin{pmatrix} 0 & \{0,0,0\} \\ \{0,0,0\} & 0 \end{pmatrix} & \begin{pmatrix} 0 & \{0,0,0\} \\ \{0,0,0\} & 0 \end{pmatrix} \\ \begin{pmatrix} 0 & \{0,0,0\} \\ \{0,0,0\} & 0 \end{pmatrix} & \begin{pmatrix} 0 & \{0,0,0\} \\ \{0,0,0\} & 0 \end{pmatrix} & \begin{pmatrix} 0 & \{0,0,0\} \\ \{0,0,0\} & 0 \end{pmatrix} \end{pmatrix}$$

$\Delta[X_2[e_7]] = \Delta[X_2[e_7]] //. \{\lambda_{10} \to \beta_2\}$

$$\begin{pmatrix} \begin{pmatrix} 0 & \{0,0,0\} \\ \{0,0,0\} & 0 \end{pmatrix} & \begin{pmatrix} \beta_4 & \{\lambda_6,0,\lambda_8\} \\ \{0,\beta_2,0\} & 0 \end{pmatrix} & \begin{pmatrix} \lambda_{20} & \{\lambda_{22},0,\lambda_{24}\} \\ \{0,\lambda_{26},0\} & 0 \end{pmatrix} \\ \begin{pmatrix} 0 & \{-\lambda_6,0,-\lambda_8\} \\ \{0,-\beta_2,0\} & \beta_4 \end{pmatrix} & \begin{pmatrix} \gamma_4 & \{0,0,0\} \\ \{0,0,0\} & \gamma_4 \end{pmatrix} & \begin{pmatrix} \lambda_{12} & \{\lambda_{14},0,\lambda_{16}\} \\ \{\lambda_{17},\lambda_{18},\lambda_{19}\} & \lambda_{13} \end{pmatrix} \\ \begin{pmatrix} 0 & \{-\lambda_{22},0,-\lambda_{24}\} \\ \{0,-\lambda_{26},0\} & \lambda_{20} \end{pmatrix} & \begin{pmatrix} \lambda_{13} & \{-\lambda_{14},0,-\lambda_{16}\} \\ \{-\lambda_{17},-\lambda_{18},-\lambda_{19}\} & \lambda_{12} \end{pmatrix} & \begin{pmatrix} -\gamma_4 & \{0,0,0\} \\ \{0,0,0\} & -\gamma_4 \end{pmatrix} \end{pmatrix}$$

```
U_{X_1[e_8]}[X_2[e_7]] == ZERO
```

True

```
defin[X_1[e_8], X_2[e_7]] // Expand
```

$$\begin{pmatrix} \begin{pmatrix} 0 & \{0,0,0\} \\ \{0,0,0\} & 0 \end{pmatrix} & \begin{pmatrix} 0 & \{0,0,0\} \\ \{0,0,-\beta_6-\lambda_8\} & 0 \end{pmatrix} & \begin{pmatrix} 0 & \{0,0,0\} \\ \{0,0,0\} & 0 \end{pmatrix} \\ \begin{pmatrix} 0 & \{0,0,0\} \\ \{0,0,\beta_6+\lambda_8\} & 0 \end{pmatrix} & \begin{pmatrix} 0 & \{0,0,0\} \\ \{0,0,0\} & 0 \end{pmatrix} & \begin{pmatrix} 0 & \{0,0,0\} \\ \{0,0,0\} & 0 \end{pmatrix} \\ \begin{pmatrix} 0 & \{0,0,0\} \\ \{0,0,0\} & 0 \end{pmatrix} & \begin{pmatrix} 0 & \{0,0,0\} \\ \{0,0,0\} & 0 \end{pmatrix} & \begin{pmatrix} 0 & \{0,0,0\} \\ \{0,0,0\} & 0 \end{pmatrix} \end{pmatrix}$$

```
Δ[X_2[e_7]] = Δ[X_2[e_7]] //. {λ_8 → -β_6}
```

$$\begin{pmatrix} \begin{pmatrix} 0 & \{0,0,0\} \\ \{0,0,0\} & 0 \end{pmatrix} & \begin{pmatrix} \beta_4 & \{\lambda_6,0,-\beta_6\} \\ \{0,\beta_2,0\} & 0 \end{pmatrix} & \begin{pmatrix} \lambda_{20} & \{\lambda_{22},0,\lambda_{24}\} \\ \{0,\lambda_{26},0\} & 0 \end{pmatrix} \\ \begin{pmatrix} 0 & \{-\lambda_6,0,\beta_6\} \\ \{0,-\beta_2,0\} & \beta_4 \end{pmatrix} & \begin{pmatrix} \gamma_4 & \{0,0,0\} \\ \{0,0,0\} & \gamma_4 \end{pmatrix} & \begin{pmatrix} \lambda_{12} & \{\lambda_{14},0,\lambda_{16}\} \\ \{\lambda_{17},\lambda_{18},\lambda_{19}\} & \lambda_{13} \end{pmatrix} \\ \begin{pmatrix} 0 & \{-\lambda_{22},0,-\lambda_{24}\} \\ \{0,-\lambda_{26},0\} & \lambda_{20} \end{pmatrix} & \begin{pmatrix} \lambda_{13} & \{-\lambda_{14},0,-\lambda_{16}\} \\ \{-\lambda_{17},-\lambda_{18},-\lambda_{19}\} & \lambda_{12} \end{pmatrix} & \begin{pmatrix} -\gamma_4 & \{0,0,0\} \\ \{0,0,0\} & -\gamma_4 \end{pmatrix} \end{pmatrix}$$

```
T[E_1[1], E_1[1], X_2[e_7]] == ZERO
```

True

```
Leib[E_1[1], E_1[1], X_2[e_7]] // Expand
```

$$\begin{pmatrix} \begin{pmatrix} 0 & \{0,0,0\} \\ \{0,0,0\} & 0 \end{pmatrix} & \begin{pmatrix} 0 & \{\lambda_6-\beta_8,0,0\} \\ \{0,0,0\} & 0 \end{pmatrix} & \begin{pmatrix} \alpha_4+\lambda_{20} & \{\alpha_8+\lambda_{22},0,\lambda_{24}-\alpha_6\} \\ \{0,\alpha_2+\lambda_{26},0\} & 0 \end{pmatrix} \\ \begin{pmatrix} 0 & \{\beta_8-\lambda_6,0,0\} \\ \{0,0,0\} & 0 \end{pmatrix} & \begin{pmatrix} 0 & \{0,0,0\} \\ \{0,0,0\} & 0 \end{pmatrix} & \begin{pmatrix} 0 & \{0,0,0\} \\ \{0,0,0\} & 0 \end{pmatrix} \\ \begin{pmatrix} 0 & \{-\alpha_8-\lambda_{22},0,\alpha_6-\lambda_{24}\} \\ \{0,-\alpha_2-\lambda_{26},0\} & \alpha_4+\lambda_{20} \end{pmatrix} & \begin{pmatrix} 0 & \{0,0,0\} \\ \{0,0,0\} & 0 \end{pmatrix} & \begin{pmatrix} 0 & \{0,0,0\} \\ \{0,0,0\} & 0 \end{pmatrix} \end{pmatrix}$$

```
Δ[X_2[e_7]] =
  Δ[X_2[e_7]] //. {λ_6 → β_8, λ_20 → -α_4, λ_22 → -α_8, λ_24 → α_6, λ_26 → -α_2}
```

$$\begin{pmatrix} \begin{pmatrix} 0 & \{0,0,0\} \\ \{0,0,0\} & 0 \end{pmatrix} & \begin{pmatrix} \beta_4 & \{\beta_8,0,-\beta_6\} \\ \{0,\beta_2,0\} & 0 \end{pmatrix} & \begin{pmatrix} -\alpha_4 & \{-\alpha_8,0,\alpha_6\} \\ \{0,-\alpha_2,0\} & 0 \end{pmatrix} \\ \begin{pmatrix} 0 & \{-\beta_8,0,\beta_6\} \\ \{0,-\beta_2,0\} & \beta_4 \end{pmatrix} & \begin{pmatrix} \gamma_4 & \{0,0,0\} \\ \{0,0,0\} & \gamma_4 \end{pmatrix} & \begin{pmatrix} \lambda_{12} & \{\lambda_{14},0,\lambda_{16}\} \\ \{\lambda_{17},\lambda_{18},\lambda_{19}\} & \lambda_{13} \end{pmatrix} \\ \begin{pmatrix} 0 & \{\alpha_8,0,-\alpha_6\} \\ \{0,\alpha_2,0\} & -\alpha_4 \end{pmatrix} & \begin{pmatrix} \lambda_{13} & \{-\lambda_{14},0,-\lambda_{16}\} \\ \{-\lambda_{17},-\lambda_{18},-\lambda_{19}\} & \lambda_{12} \end{pmatrix} & \begin{pmatrix} -\gamma_4 & \{0,0,0\} \\ \{0,0,0\} & -\gamma_4 \end{pmatrix} \end{pmatrix}$$

```
T[X_1[e_1], X_1[e_3], X_2[e_7]] == ZERO
```

True

**Leib[X$_1$[e$_1$], X$_1$[e$_3$], X$_2$[e$_7$]] // Expand**

$$\left(\begin{array}{ccc} \begin{pmatrix} 0 & \{0,0,0\} \\ \{0,0,0\} & 0 \end{pmatrix} & \begin{pmatrix} 0 & \{0,0,0\} \\ \{0,0,0\} & 0 \end{pmatrix} & \begin{pmatrix} 0 & \{0,0,0\} \\ \{0,0,0\} & 0 \end{pmatrix} \\ \begin{pmatrix} 0 & \{0,0,0\} \\ \{0,0,0\} & 0 \end{pmatrix} & \begin{pmatrix} 0 & \{0,0,0\} \\ \{0,0,0\} & 0 \end{pmatrix} & \begin{pmatrix} 0 & \{0,0,0\} \\ \{0, \delta_5 - \lambda_{16}, 0\} & 0 \end{pmatrix} \\ \begin{pmatrix} 0 & \{0,0,0\} \\ \{0,0,0\} & 0 \end{pmatrix} & \begin{pmatrix} 0 & \{0,0,0\} \\ \{0, \lambda_{16} - \delta_5, 0\} & 0 \end{pmatrix} & \begin{pmatrix} 0 & \{0,0,0\} \\ \{0,0,0\} & 0 \end{pmatrix} \end{array}\right)$$

**Δ[X$_2$[e$_7$]] = Δ[X$_2$[e$_7$]] //. {λ$_{16}$ → δ$_5$}**

$$\left(\begin{array}{ccc} \begin{pmatrix} 0 & \{0,0,0\} \\ \{0,0,0\} & 0 \end{pmatrix} & \begin{pmatrix} \beta_4 & \{\beta_8, 0, -\beta_6\} \\ \{0, \beta_2, 0\} & 0 \end{pmatrix} & \begin{pmatrix} -\alpha_4 & \{-\alpha_8, 0, \alpha_6\} \\ \{0, -\alpha_2, 0\} & 0 \end{pmatrix} \\ \begin{pmatrix} 0 & \{-\beta_8, 0, \beta_6\} \\ \{0, -\beta_2, 0\} & \beta_4 \end{pmatrix} & \begin{pmatrix} \gamma_4 & \{0,0,0\} \\ \{0,0,0\} & \gamma_4 \end{pmatrix} & \begin{pmatrix} \lambda_{12} & \{\lambda_{14}, 0, \delta_5\} \\ \{\lambda_{17}, \lambda_{18}, \lambda_{19}\} & \lambda_{13} \end{pmatrix} \\ \begin{pmatrix} 0 & \{\alpha_8, 0, -\alpha_6\} \\ \{0, \alpha_2, 0\} & -\alpha_4 \end{pmatrix} & \begin{pmatrix} \lambda_{13} & \{-\lambda_{14}, 0, -\delta_5\} \\ \{-\lambda_{17}, -\lambda_{18}, -\lambda_{19}\} & \lambda_{12} \end{pmatrix} & \begin{pmatrix} -\gamma_4 & \{0,0,0\} \\ \{0,0,0\} & -\gamma_4 \end{pmatrix} \end{array}\right)$$

**T[X$_1$[e$_1$], X$_1$[e$_6$], X$_2$[e$_7$]] == ZERO**

True

**Leib[X$_1$[e$_1$], X$_1$[e$_6$], X$_2$[e$_7$]] // Expand**

$$\left(\begin{array}{ccc} \begin{pmatrix} 0 & \{0,0,0\} \\ \{0,0,0\} & 0 \end{pmatrix} & \begin{pmatrix} 0 & \{0,0,0\} \\ \{0,0,0\} & 0 \end{pmatrix} & \begin{pmatrix} 0 & \{0,0,0\} \\ \{0,0,0\} & 0 \end{pmatrix} \\ \begin{pmatrix} 0 & \{0,0,0\} \\ \{0,0,0\} & 0 \end{pmatrix} & \begin{pmatrix} 0 & \{0,0,0\} \\ \{0,0,0\} & 0 \end{pmatrix} & \begin{pmatrix} 0 & \{0,0,0\} \\ \{\delta_3 + \lambda_{12}, 0, 0\} & \delta_7 + \lambda_{14} \end{pmatrix} \\ \begin{pmatrix} 0 & \{0,0,0\} \\ \{0,0,0\} & 0 \end{pmatrix} & \begin{pmatrix} \delta_7 + \lambda_{14} & \{0,0,0\} \\ \{-\delta_3 - \lambda_{12}, 0, 0\} & 0 \end{pmatrix} & \begin{pmatrix} 0 & \{0,0,0\} \\ \{0,0,0\} & 0 \end{pmatrix} \end{array}\right)$$

**Δ[X$_2$[e$_7$]] = Δ[X$_2$[e$_7$]] //. {λ$_{12}$ → -δ$_3$}**

$$\left(\begin{array}{ccc} \begin{pmatrix} 0 & \{0,0,0\} \\ \{0,0,0\} & 0 \end{pmatrix} & \begin{pmatrix} \beta_4 & \{\beta_8, 0, -\beta_6\} \\ \{0, \beta_2, 0\} & 0 \end{pmatrix} & \begin{pmatrix} -\alpha_4 & \{-\alpha_8, 0, \alpha_6\} \\ \{0, -\alpha_2, 0\} & 0 \end{pmatrix} \\ \begin{pmatrix} 0 & \{-\beta_8, 0, \beta_6\} \\ \{0, -\beta_2, 0\} & \beta_4 \end{pmatrix} & \begin{pmatrix} \gamma_4 & \{0,0,0\} \\ \{0,0,0\} & \gamma_4 \end{pmatrix} & \begin{pmatrix} -\delta_3 & \{\lambda_{14}, 0, \delta_5\} \\ \{\lambda_{17}, \lambda_{18}, \lambda_{19}\} & \lambda_{13} \end{pmatrix} \\ \begin{pmatrix} 0 & \{\alpha_8, 0, -\alpha_6\} \\ \{0, \alpha_2, 0\} & -\alpha_4 \end{pmatrix} & \begin{pmatrix} \lambda_{13} & \{-\lambda_{14}, 0, -\delta_5\} \\ \{-\lambda_{17}, -\lambda_{18}, -\lambda_{19}\} & -\delta_3 \end{pmatrix} & \begin{pmatrix} -\gamma_4 & \{0,0,0\} \\ \{0,0,0\} & -\gamma_4 \end{pmatrix} \end{array}\right)$$

**T[X$_1$[e$_2$], X$_1$[e$_8$], X$_2$[e$_7$]] == X$_2$[e$_3$]**

True

**Δ[X₂[e₃]] - Leib[X₁[e₂], X₁[e₈], X₂[e₇]] // Expand**

$$\left(\begin{pmatrix} 0 & \{0,0,0\} \\ \{0,0,0\} & 0 \end{pmatrix} \quad \begin{pmatrix} 0 & \{0,0,0\} \\ \{0,0,0\} & 0 \end{pmatrix} \quad \begin{pmatrix} 0 & \{0,0,0\} \\ \{0,0,0\} & 0 \end{pmatrix}\right.$$
$$\begin{pmatrix} 0 & \{0,0,0\} \\ \{0,0,0\} & 0 \end{pmatrix} \quad \begin{pmatrix} 0 & \{0,0,0\} \\ \{0,0,0\} & 0 \end{pmatrix} \quad \begin{pmatrix} 0 & \{\xi+\delta_1-\epsilon_2-\lambda_{18},\eta_2+\lambda_{17}, \\ \{0,0,0\} & 0 \end{pmatrix}$$
$$\left.\begin{pmatrix} 0 & \{0,0,0\} \\ \{0,0,0\} & 0 \end{pmatrix} \quad \begin{pmatrix} 0 & \{-\xi-\delta_1+\epsilon_2+\lambda_{18},-\eta_2-\lambda_{17},0\} \\ \{0,0,0\} & 0 \end{pmatrix} \quad \begin{pmatrix} 0 & \{0,0,0\} \\ \{0,0,0\} & 0 \end{pmatrix}\right)$$

**Δ[X₂[e₇]] = Δ[X₂[e₇]] //. {λ₁₈ → ξ + δ₁ − ϵ₂, λ₁₇ → −η₂}**

$$\left(\begin{pmatrix} 0 & \{0,0,0\} \\ \{0,0,0\} & 0 \end{pmatrix} \quad \begin{pmatrix} \beta_4 & \{\beta_8,0,-\beta_6\} \\ \{0,\beta_2,0\} & 0 \end{pmatrix} \quad \begin{pmatrix} -\alpha_4 & \{-\alpha_8,0,\alpha_6\} \\ \{0,-\alpha_2,0\} & 0 \end{pmatrix}\right.$$
$$\begin{pmatrix} 0 & \{-\beta_8,0,\beta_6\} \\ \{0,-\beta_2,0\} & \beta_4 \end{pmatrix} \quad \begin{pmatrix} \gamma_4 & \{0,0,0\} \\ \{0,0,0\} & \gamma_4 \end{pmatrix} \quad \begin{pmatrix} -\delta_3 & \{\lambda_{14},0, \\ \{-\eta_2,\xi+\delta_1-\epsilon_2,\lambda_{19}\} & \lambda_{13} \end{pmatrix}$$
$$\left.\begin{pmatrix} 0 & \{\alpha_8,0,-\alpha_6\} \\ \{0,\alpha_2,0\} & -\alpha_4 \end{pmatrix} \quad \begin{pmatrix} \lambda_{13} & \{-\lambda_{14},0,-\delta_5\} \\ \{\eta_2,-\xi-\delta_1+\epsilon_2,-\lambda_{19}\} & -\delta_3 \end{pmatrix} \quad \begin{pmatrix} -\gamma_4 & \{0,0,0\} \\ \{0,0,0\} & -\gamma_4 \end{pmatrix}\right)$$

**T[X₁[e₅], X₁[e₈], X₂[e₇]] == -X₂[e₇]**

True

**Δ[X₂[e₇]] + Leib[X₁[e₅], X₁[e₈], X₂[e₇]] // Expand**

$$\left(\begin{pmatrix} 0 & \{0,0,0\} \\ \{0,0,0\} & 0 \end{pmatrix} \quad \begin{pmatrix} 0 & \{0,0,0\} \\ \{0,0,0\} & 0 \end{pmatrix} \quad \begin{pmatrix} 0 & \{0,0,0\} \\ \{0,0,0\} & 0 \end{pmatrix}\right.$$
$$\begin{pmatrix} 0 & \{0,0,0\} \\ \{0,0,0\} & 0 \end{pmatrix} \quad \begin{pmatrix} 0 & \{0,0,0\} \\ \{0,0,0\} & 0 \end{pmatrix} \quad \begin{pmatrix} 0 & \{\delta_7+\lambda_{14},0,0\} \\ \{0,0,\lambda_{19}+\phi_2\} & \eta_5+\lambda_{13} \end{pmatrix}$$
$$\left.\begin{pmatrix} 0 & \{0,0,0\} \\ \{0,0,0\} & 0 \end{pmatrix} \quad \begin{pmatrix} \eta_5+\lambda_{13} & \{-\delta_7-\lambda_{14},0,0\} \\ \{0,0,-\lambda_{19}-\phi_2\} & 0 \end{pmatrix} \quad \begin{pmatrix} 0 & \{0,0,0\} \\ \{0,0,0\} & 0 \end{pmatrix}\right)$$

**Δ[X₂[e₇]] = Δ[X₂[e₇]] //. {λ₁₉ → −ϕ₂, λ₁₃ → −η₅, λ₁₄ → −δ₇}**

$$\left(\begin{pmatrix} 0 & \{0,0,0\} \\ \{0,0,0\} & 0 \end{pmatrix} \quad \begin{pmatrix} \beta_4 & \{\beta_8,0,-\beta_6\} \\ \{0,\beta_2,0\} & 0 \end{pmatrix} \quad \begin{pmatrix} -\alpha_4 & \{-\alpha_8,0,\alpha_6\} \\ \{0,-\alpha_2,0\} & 0 \end{pmatrix}\right.$$
$$\begin{pmatrix} 0 & \{-\beta_8,0,\beta_6\} \\ \{0,-\beta_2,0\} & \beta_4 \end{pmatrix} \quad \begin{pmatrix} \gamma_4 & \{0,0,0\} \\ \{0,0,0\} & \gamma_4 \end{pmatrix} \quad \begin{pmatrix} -\delta_3 & \{-\delta_7,0,\delta_5\} \\ \{-\eta_2,\xi+\delta_1-\epsilon_2,-\phi_2\} & -\eta_5 \end{pmatrix}$$
$$\left.\begin{pmatrix} 0 & \{\alpha_8,0,-\alpha_6\} \\ \{0,\alpha_2,0\} & -\alpha_4 \end{pmatrix} \quad \begin{pmatrix} -\eta_5 & \{\delta_7,0,-\delta_5\} \\ \{\eta_2,-\xi-\delta_1+\epsilon_2,\phi_2\} & -\delta_3 \end{pmatrix} \quad \begin{pmatrix} -\gamma_4 & \{0,0,0\} \\ \{0,0,0\} & -\gamma_4 \end{pmatrix}\right)$$

**Variables[Δ[X₂[e₇]]]**

$\{\xi, \alpha_2, \alpha_4, \alpha_6, \alpha_8, \beta_2, \beta_4, \beta_6, \beta_8, \gamma_4, \delta_1, \delta_3, \delta_5, \delta_7, \epsilon_2, \eta_2, \eta_5, \phi_2\}$

- **X₂[e₈]**

  **Δ[X₂[e₈]] = generic;**

  **U_{E₁[1]}[X₂[e₈]] == ZERO**

  True

**(defin[E₁[1], X₂[e₈]] // Expand)**

$$\begin{pmatrix} \begin{pmatrix} \lambda_1 & \{0,0,0\} \\ \{0,0,0\} & \lambda_1 \end{pmatrix} & \begin{pmatrix} 0 & \{0,0,0\} \\ \{0,0,0\} & 0 \end{pmatrix} & \begin{pmatrix} 0 & \{0,0,0\} \\ \{0,0,0\} & 0 \end{pmatrix} \\ \begin{pmatrix} 0 & \{0,0,0\} \\ \{0,0,0\} & 0 \end{pmatrix} & \begin{pmatrix} 0 & \{0,0,0\} \\ \{0,0,0\} & 0 \end{pmatrix} & \begin{pmatrix} 0 & \{0,0,0\} \\ \{0,0,0\} & 0 \end{pmatrix} \\ \begin{pmatrix} 0 & \{0,0,0\} \\ \{0,0,0\} & 0 \end{pmatrix} & \begin{pmatrix} 0 & \{0,0,0\} \\ \{0,0,0\} & 0 \end{pmatrix} & \begin{pmatrix} 0 & \{0,0,0\} \\ \{0,0,0\} & 0 \end{pmatrix} \end{pmatrix}$$

**Δ[X₂[e₈]] = Δ[X₂[e₈]] //. {λ₁ → 0}**

$$\begin{pmatrix} \begin{pmatrix} 0 & \{0,0,0\} \\ \{0,0,0\} & 0 \end{pmatrix} & \begin{pmatrix} \lambda_4 & \{\lambda_6,\lambda_7,\lambda_8\} \\ \{\lambda_9,\lambda_{10},\lambda_{11}\} & \lambda_5 \end{pmatrix} & \begin{pmatrix} \lambda_{20} & \{\lambda_{22},\lambda_{23},\lambda_{24}\} \\ \{\lambda_{25},\lambda_{26},\lambda_{27}\} & \lambda_{21} \end{pmatrix} \\ \begin{pmatrix} \lambda_5 & \{-\lambda_6,-\lambda_7,-\lambda_8\} \\ \{-\lambda_9,-\lambda_{10},-\lambda_{11}\} & \lambda_4 \end{pmatrix} & \begin{pmatrix} \lambda_2 & \{0,0,0\} \\ \{0,0,0\} & \lambda_2 \end{pmatrix} & \begin{pmatrix} \lambda_{12} & \{\lambda_{14},\lambda_{15},\lambda_{16}\} \\ \{\lambda_{17},\lambda_{18},\lambda_{19}\} & \lambda_{13} \end{pmatrix} \\ \begin{pmatrix} \lambda_{21} & \{-\lambda_{22},-\lambda_{23},-\lambda_{24}\} \\ \{-\lambda_{25},-\lambda_{26},-\lambda_{27}\} & \lambda_{20} \end{pmatrix} & \begin{pmatrix} \lambda_{13} & \{-\lambda_{14},-\lambda_{15},-\lambda_{16}\} \\ \{-\lambda_{17},-\lambda_{18},-\lambda_{19}\} & \lambda_{12} \end{pmatrix} & \begin{pmatrix} \lambda_3 & \{0,0,0\} \\ \{0,0,0\} & \lambda_3 \end{pmatrix} \end{pmatrix}$$

**U_{X₂[e₈]}[E₂[1]] == ZERO**

True

**defin[X₂[e₈], E₂[1]] // Expand**

$$\begin{pmatrix} \begin{pmatrix} 0 & \{0,0,0\} \\ \{0,0,0\} & 0 \end{pmatrix} & \begin{pmatrix} 0 & \{0,0,0\} \\ \{0,0,0\} & 0 \end{pmatrix} & \begin{pmatrix} \lambda_8 & \{-\lambda_{10},\lambda_9,0\} \\ \{0,0,\lambda_5\} & 0 \end{pmatrix} \\ \begin{pmatrix} 0 & \{0,0,0\} \\ \{0,0,0\} & 0 \end{pmatrix} & \begin{pmatrix} 0 & \{0,0,0\} \\ \{0,0,0\} & 0 \end{pmatrix} & \begin{pmatrix} 0 & \{0,0,0\} \\ \{0,0,\lambda_2-\gamma_5\} & 0 \end{pmatrix} \\ \begin{pmatrix} 0 & \{\lambda_{10},-\lambda_9,0\} \\ \{0,0,-\lambda_5\} & \lambda_8 \end{pmatrix} & \begin{pmatrix} 0 & \{0,0,0\} \\ \{0,0,\gamma_5-\lambda_2\} & 0 \end{pmatrix} & \begin{pmatrix} -\lambda_{16} & \{0,0,0\} \\ \{0,0,0\} & -\lambda_{16} \end{pmatrix} \end{pmatrix}$$

**Δ[X₂[e₈]] =**
 **Δ[X₂[e₈]] //. {λ₂ → γ₅, λ₅ → 0, λ₈ → 0, λ₉ → 0, λ₁₀ → 0, λ₁₆ → 0}**

$$\begin{pmatrix} \begin{pmatrix} 0 & \{0,0,0\} \\ \{0,0,0\} & 0 \end{pmatrix} & \begin{pmatrix} \lambda_4 & \{\lambda_6,\lambda_7,0\} \\ \{0,0,\lambda_{11}\} & 0 \end{pmatrix} & \begin{pmatrix} \lambda_{20} & \{\lambda_{22},\lambda_{23},\lambda_{24}\} \\ \{\lambda_{25},\lambda_{26},\lambda_{27}\} & \lambda_{21} \end{pmatrix} \\ \begin{pmatrix} 0 & \{-\lambda_6,-\lambda_7,0\} \\ \{0,0,-\lambda_{11}\} & \lambda_4 \end{pmatrix} & \begin{pmatrix} \gamma_5 & \{0,0,0\} \\ \{0,0,0\} & \gamma_5 \end{pmatrix} & \begin{pmatrix} \lambda_{12} & \{\lambda_{14},\lambda_{15},0\} \\ \{\lambda_{17},\lambda_{18},\lambda_{19}\} & \lambda_{13} \end{pmatrix} \\ \begin{pmatrix} \lambda_{21} & \{-\lambda_{22},-\lambda_{23},-\lambda_{24}\} \\ \{-\lambda_{25},-\lambda_{26},-\lambda_{27}\} & \lambda_{20} \end{pmatrix} & \begin{pmatrix} \lambda_{13} & \{-\lambda_{14},-\lambda_{15},0\} \\ \{-\lambda_{17},-\lambda_{18},-\lambda_{19}\} & \lambda_{12} \end{pmatrix} & \begin{pmatrix} \lambda_3 & \{0,0,0\} \\ \{0,0,0\} & \lambda_3 \end{pmatrix} \end{pmatrix}$$

**U_{E₃[1]}[X₂[e₈]] == ZERO**

True

`defin[E₃[1], X₂[e₈]] // Expand`

$$\begin{pmatrix} \begin{pmatrix} 0 & \{0,0,0\} \\ \{0,0,0\} & 0 \end{pmatrix} & \begin{pmatrix} 0 & \{0,0,0\} \\ \{0,0,0\} & 0 \end{pmatrix} & \begin{pmatrix} 0 & \{0,0,0\} \\ \{0,0,0\} & 0 \end{pmatrix} \\ \begin{pmatrix} 0 & \{0,0,0\} \\ \{0,0,0\} & 0 \end{pmatrix} & \begin{pmatrix} 0 & \{0,0,0\} \\ \{0,0,0\} & 0 \end{pmatrix} & \begin{pmatrix} 0 & \{0,0,0\} \\ \{0,0,0\} & 0 \end{pmatrix} \\ \begin{pmatrix} 0 & \{0,0,0\} \\ \{0,0,0\} & 0 \end{pmatrix} & \begin{pmatrix} 0 & \{0,0,0\} \\ \{0,0,0\} & 0 \end{pmatrix} & \begin{pmatrix} \gamma_5+\lambda_3 & \{0,0,0\} \\ \{0,0,0\} & \gamma_5+\lambda_3 \end{pmatrix} \end{pmatrix}$$

`Δ[X₂[e₈]] = Δ[X₂[e₈]] //. {λ₃ → -γ₅}`

$$\begin{pmatrix} \begin{pmatrix} 0 & \{0,0,0\} \\ \{0,0,0\} & 0 \end{pmatrix} & \begin{pmatrix} \lambda_4 & \{\lambda_6,\lambda_7,0\} \\ \{0,0,\lambda_{11}\} & 0 \end{pmatrix} & \begin{pmatrix} \lambda_{20} & \{\lambda_{22},\lambda_{23},\lambda_{24}\} \\ \{\lambda_{25},\lambda_{26},\lambda_{27}\} & \lambda_{21} \end{pmatrix} \\ \begin{pmatrix} 0 & \{-\lambda_6,-\lambda_7,0\} \\ \{0,0,-\lambda_{11}\} & \lambda_4 \end{pmatrix} & \begin{pmatrix} \gamma_5 & \{0,0,0\} \\ \{0,0,0\} & \gamma_5 \end{pmatrix} & \begin{pmatrix} \lambda_{12} & \{\lambda_{14},\lambda_{15},0\} \\ \{\lambda_{17},\lambda_{18},\lambda_{19}\} & \lambda_{13} \end{pmatrix} \\ \begin{pmatrix} \lambda_{21} & \{-\lambda_{22},-\lambda_{23},-\lambda_{24}\} \\ \{-\lambda_{25},-\lambda_{26},-\lambda_{27}\} & \lambda_{20} \end{pmatrix} & \begin{pmatrix} \lambda_{13} & \{-\lambda_{14},-\lambda_{15},0\} \\ \{-\lambda_{17},-\lambda_{18},-\lambda_{19}\} & \lambda_{12} \end{pmatrix} & \begin{pmatrix} -\gamma_5 & \{0,0,0\} \\ \{0,0,0\} & -\gamma_5 \end{pmatrix} \end{pmatrix}$$

`U_{X₂[e₈]}[E₃[1]] == ZERO`

True

`defin[X₂[e₈], E₃[1]] // Expand`

$$\begin{pmatrix} \begin{pmatrix} 0 & \{0,0,0\} \\ \{0,0,0\} & 0 \end{pmatrix} & \begin{pmatrix} -\lambda_{24} & \{\lambda_{26},-\lambda_{25},0\} \\ \{0,0,-\lambda_{21}\} & 0 \end{pmatrix} & \begin{pmatrix} 0 & \{0,0,0\} \\ \{0,0,0\} & 0 \end{pmatrix} \\ \begin{pmatrix} 0 & \{-\lambda_{26},\lambda_{25},0\} \\ \{0,0,\lambda_{21}\} & -\lambda_{24} \end{pmatrix} & \begin{pmatrix} 0 & \{0,0,0\} \\ \{0,0,0\} & 0 \end{pmatrix} & \begin{pmatrix} 0 & \{0,0,0\} \\ \{0,0,0\} & 0 \end{pmatrix} \\ \begin{pmatrix} 0 & \{0,0,0\} \\ \{0,0,0\} & 0 \end{pmatrix} & \begin{pmatrix} 0 & \{0,0,0\} \\ \{0,0,0\} & 0 \end{pmatrix} & \begin{pmatrix} 0 & \{0,0,0\} \\ \{0,0,0\} & 0 \end{pmatrix} \end{pmatrix}$$

`Δ[X₂[e₈]] = Δ[X₂[e₈]] //. {λ₂₁ → 0, λ₂₄ → 0, λ₂₅ → 0, λ₂₆ → 0}`

$$\begin{pmatrix} \begin{pmatrix} 0 & \{0,0,0\} \\ \{0,0,0\} & 0 \end{pmatrix} & \begin{pmatrix} \lambda_4 & \{\lambda_6,\lambda_7,0\} \\ \{0,0,\lambda_{11}\} & 0 \end{pmatrix} & \begin{pmatrix} \lambda_{20} & \{\lambda_{22},\lambda_{23},0\} \\ \{0,0,\lambda_{27}\} & 0 \end{pmatrix} \\ \begin{pmatrix} 0 & \{-\lambda_6,-\lambda_7,0\} \\ \{0,0,-\lambda_{11}\} & \lambda_4 \end{pmatrix} & \begin{pmatrix} \gamma_5 & \{0,0,0\} \\ \{0,0,0\} & \gamma_5 \end{pmatrix} & \begin{pmatrix} \lambda_{12} & \{\lambda_{14},\lambda_{15},0\} \\ \{\lambda_{17},\lambda_{18},\lambda_{19}\} & \lambda_{13} \end{pmatrix} \\ \begin{pmatrix} 0 & \{-\lambda_{22},-\lambda_{23},0\} \\ \{0,0,-\lambda_{27}\} & \lambda_{20} \end{pmatrix} & \begin{pmatrix} \lambda_{13} & \{-\lambda_{14},-\lambda_{15},0\} \\ \{-\lambda_{17},-\lambda_{18},-\lambda_{19}\} & \lambda_{12} \end{pmatrix} & \begin{pmatrix} -\gamma_5 & \{0,0,0\} \\ \{0,0,0\} & -\gamma_5 \end{pmatrix} \end{pmatrix}$$

`U_{X₁[e₂]}[X₂[e₈]] == ZERO`

True

**defin[X$_1$[e$_2$], X$_2$[e$_8$]] // Expand**

$$\left(\begin{pmatrix} \begin{pmatrix} 0 & \{0,0,0\} \\ \{0,0,0\} & 0 \end{pmatrix} & \begin{pmatrix} 0 & \{0,0,0\} \\ \{0,0,0\} & \lambda_4 - \beta_5 \end{pmatrix} & \begin{pmatrix} 0 & \{0,0,0\} \\ \{0,0,0\} & 0 \end{pmatrix} \\ \begin{pmatrix} \lambda_4 - \beta_5 & \{0,0,0\} \\ \{0,0,0\} & 0 \end{pmatrix} & \begin{pmatrix} 0 & \{0,0,0\} \\ \{0,0,0\} & 0 \end{pmatrix} & \begin{pmatrix} 0 & \{0,0,0\} \\ \{0,0,0\} & 0 \end{pmatrix} \\ \begin{pmatrix} 0 & \{0,0,0\} \\ \{0,0,0\} & 0 \end{pmatrix} & \begin{pmatrix} 0 & \{0,0,0\} \\ \{0,0,0\} & 0 \end{pmatrix} & \begin{pmatrix} 0 & \{0,0,0\} \\ \{0,0,0\} & 0 \end{pmatrix} \end{pmatrix}\right)$$

**Δ[X$_2$[e$_8$]] = Δ[X$_2$[e$_8$]] //. {λ$_4$ → β$_5$}**

$$\begin{pmatrix} \begin{pmatrix} 0 & \{0,0,0\} \\ \{0,0,0\} & 0 \end{pmatrix} & \begin{pmatrix} \beta_5 & \{\lambda_6, \lambda_7, 0\} \\ \{0,0,\lambda_{11}\} & 0 \end{pmatrix} & \begin{pmatrix} \lambda_{20} & \{\lambda_{22}, \lambda_{23}, 0\} \\ \{0,0,\lambda_{27}\} & 0 \end{pmatrix} \\ \begin{pmatrix} 0 & \{-\lambda_6, -\lambda_7, 0\} \\ \{0,0,-\lambda_{11}\} & \beta_5 \end{pmatrix} & \begin{pmatrix} \gamma_5 & \{0,0,0\} \\ \{0,0,0\} & \gamma_5 \end{pmatrix} & \begin{pmatrix} \lambda_{12} & \{\lambda_{14}, \lambda_{15}, 0\} \\ \{\lambda_{17}, \lambda_{18}, \lambda_{19}\} & \lambda_{13} \end{pmatrix} \\ \begin{pmatrix} 0 & \{-\lambda_{22}, -\lambda_{23}, 0\} \\ \{0,0,-\lambda_{27}\} & \lambda_{20} \end{pmatrix} & \begin{pmatrix} \lambda_{13} & \{-\lambda_{14}, -\lambda_{15}, 0\} \\ \{-\lambda_{17}, -\lambda_{18}, -\lambda_{19}\} & \lambda_{12} \end{pmatrix} & \begin{pmatrix} -\gamma_5 & \{0,0,0\} \\ \{0,0,0\} & -\gamma_5 \end{pmatrix} \end{pmatrix}$$

**U$_{X_2[e_8]}$[X$_1$[e$_5$]] == ZERO**

True

**defin[X$_2$[e$_8$], X$_1$[e$_5$]] // Expand**

$$\left(\begin{pmatrix} \begin{pmatrix} 0 & \{0,0,0\} \\ \{0,0,0\} & 0 \end{pmatrix} & \begin{pmatrix} 0 & \{0,0,0\} \\ \{0,0,0\} & 0 \end{pmatrix} & \begin{pmatrix} 0 & \{0,0,0\} \\ \{0,0,0\} & 0 \end{pmatrix} \\ \begin{pmatrix} 0 & \{0,0,0\} \\ \{0,0,0\} & 0 \end{pmatrix} & \begin{pmatrix} 0 & \{0,0,0\} \\ \{0,0,0\} & 0 \end{pmatrix} & \begin{pmatrix} 0 & \{0,0,0\} \\ \{0,0,\beta_2 - \lambda_{11}\} & 0 \end{pmatrix} \\ \begin{pmatrix} 0 & \{0,0,0\} \\ \{0,0,0\} & 0 \end{pmatrix} & \begin{pmatrix} 0 & \{0,0,0\} \\ \{0,0,\lambda_{11} - \beta_2\} & 0 \end{pmatrix} & \begin{pmatrix} 0 & \{0,0,0\} \\ \{0,0,0\} & 0 \end{pmatrix} \end{pmatrix}\right)$$

**Δ[X$_2$[e$_8$]] = Δ[X$_2$[e$_8$]] //. {λ$_{11}$ → β$_2$}**

$$\begin{pmatrix} \begin{pmatrix} 0 & \{0,0,0\} \\ \{0,0,0\} & 0 \end{pmatrix} & \begin{pmatrix} \beta_5 & \{\lambda_6, \lambda_7, 0\} \\ \{0,0,\beta_2\} & 0 \end{pmatrix} & \begin{pmatrix} \lambda_{20} & \{\lambda_{22}, \lambda_{23}, 0\} \\ \{0,0,\lambda_{27}\} & 0 \end{pmatrix} \\ \begin{pmatrix} 0 & \{-\lambda_6, -\lambda_7, 0\} \\ \{0,0,-\beta_2\} & \beta_5 \end{pmatrix} & \begin{pmatrix} \gamma_5 & \{0,0,0\} \\ \{0,0,0\} & \gamma_5 \end{pmatrix} & \begin{pmatrix} \lambda_{12} & \{\lambda_{14}, \lambda_{15}, 0\} \\ \{\lambda_{17}, \lambda_{18}, \lambda_{19}\} & \lambda_{13} \end{pmatrix} \\ \begin{pmatrix} 0 & \{-\lambda_{22}, -\lambda_{23}, 0\} \\ \{0,0,-\lambda_{27}\} & \lambda_{20} \end{pmatrix} & \begin{pmatrix} \lambda_{13} & \{-\lambda_{14}, -\lambda_{15}, 0\} \\ \{-\lambda_{17}, -\lambda_{18}, -\lambda_{19}\} & \lambda_{12} \end{pmatrix} & \begin{pmatrix} -\gamma_5 & \{0,0,0\} \\ \{0,0,0\} & -\gamma_5 \end{pmatrix} \end{pmatrix}$$

**U$_{X_2[e_8]}$[X$_1$[e$_6$]] == ZERO**

True

**defin[X$_2$[e$_8$], X$_1$[e$_6$]] // Expand**

$$\left(\begin{array}{ccc} \begin{pmatrix} 0 & \{0,0,0\} \\ \{0,0,0\} & 0 \end{pmatrix} & \begin{pmatrix} 0 & \{0,0,0\} \\ \{0,0,0\} & 0 \end{pmatrix} & \begin{pmatrix} 0 & \{0,0,0\} \\ \{0,0,0\} & 0 \end{pmatrix} \\ \begin{pmatrix} 0 & \{0,0,0\} \\ \{0,0,0\} & 0 \end{pmatrix} & \begin{pmatrix} 0 & \{0,0,0\} \\ \{0,0,0\} & 0 \end{pmatrix} & \begin{pmatrix} 0 & \{0,0,0\} \\ \{0,0,-\beta_7-\lambda_6\} & 0 \end{pmatrix} \\ \begin{pmatrix} 0 & \{0,0,0\} \\ \{0,0,0\} & 0 \end{pmatrix} & \begin{pmatrix} 0 & \{0,0,0\} \\ \{0,0,\beta_7+\lambda_6\} & 0 \end{pmatrix} & \begin{pmatrix} 0 & \{0,0,0\} \\ \{0,0,0\} & 0 \end{pmatrix} \end{array}\right)$$

**Δ[X$_2$[e$_8$]] = Δ[X$_2$[e$_8$]] //. {λ$_6$ → −β$_7$}**

$$\left(\begin{array}{ccc} \begin{pmatrix} 0 & \{0,0,0\} \\ \{0,0,0\} & 0 \end{pmatrix} & \begin{pmatrix} \beta_5 & \{-\beta_7,\lambda_7,0\} \\ \{0,0,\beta_2\} & 0 \end{pmatrix} & \begin{pmatrix} \lambda_{20} & \{\lambda_{22},\lambda_{23},0\} \\ \{0,0,\lambda_{27}\} & 0 \end{pmatrix} \\ \begin{pmatrix} 0 & \{\beta_7,-\lambda_7,0\} \\ \{0,0,-\beta_2\} & \beta_5 \end{pmatrix} & \begin{pmatrix} \gamma_5 & \{0,0,0\} \\ \{0,0,0\} & \gamma_5 \end{pmatrix} & \begin{pmatrix} \lambda_{12} & \{\lambda_{14},\lambda_{15},0\} \\ \{\lambda_{17},\lambda_{18},\lambda_{19}\} & \lambda_{13} \end{pmatrix} \\ \begin{pmatrix} 0 & \{-\lambda_{22},-\lambda_{23},0\} \\ \{0,0,-\lambda_{27}\} & \lambda_{20} \end{pmatrix} & \begin{pmatrix} \lambda_{13} & \{-\lambda_{14},-\lambda_{15},0\} \\ \{-\lambda_{17},-\lambda_{18},-\lambda_{19}\} & \lambda_{12} \end{pmatrix} & \begin{pmatrix} -\gamma_5 & \{0,0,0\} \\ \{0,0,0\} & -\gamma_5 \end{pmatrix} \end{array}\right)$$

**T[E$_1$[1], E$_1$[1], X$_2$[e$_8$]] == ZERO**

True

**Leib[E$_1$[1], E$_1$[1], X$_2$[e$_8$]] // Expand**

$$\left(\begin{array}{ccc} \begin{pmatrix} 0 & \{0,0,0\} \\ \{0,0,0\} & 0 \end{pmatrix} & \begin{pmatrix} 0 & \{0,\lambda_7-\beta_6,0\} \\ \{0,0,0\} & 0 \end{pmatrix} & \begin{pmatrix} \alpha_5+\lambda_{20} & \{\lambda_{22}-\alpha_7,\alpha_6+\lambda_{23},0\} \\ \{0,0,\alpha_2+\lambda_{27}\} & 0 \end{pmatrix} \\ \begin{pmatrix} 0 & \{0,\beta_6-\lambda_7,0\} \\ \{0,0,0\} & 0 \end{pmatrix} & \begin{pmatrix} 0 & \{0,0,0\} \\ \{0,0,0\} & 0 \end{pmatrix} & \begin{pmatrix} 0 & \{0,0,0\} \\ \{0,0,0\} & 0 \end{pmatrix} \\ \begin{pmatrix} 0 & \{\alpha_7-\lambda_{22},-\alpha_6-\lambda_{23},0\} \\ \{0,0,-\alpha_2-\lambda_{27}\} & \alpha_5+\lambda_{20} \end{pmatrix} & \begin{pmatrix} 0 & \{0,0,0\} \\ \{0,0,0\} & 0 \end{pmatrix} & \begin{pmatrix} 0 & \{0,0,0\} \\ \{0,0,0\} & 0 \end{pmatrix} \end{array}\right)$$

**Δ[X$_2$[e$_8$]] =**
  **Δ[X$_2$[e$_8$]] //. {λ$_7$ → β$_6$, λ$_{20}$ → −α$_5$, λ$_{22}$ → α$_7$, λ$_{23}$ → −α$_6$, λ$_{27}$ → −α$_2$}**

$$\left(\begin{array}{ccc} \begin{pmatrix} 0 & \{0,0,0\} \\ \{0,0,0\} & 0 \end{pmatrix} & \begin{pmatrix} \beta_5 & \{-\beta_7,\beta_6,0\} \\ \{0,0,\beta_2\} & 0 \end{pmatrix} & \begin{pmatrix} -\alpha_5 & \{\alpha_7,-\alpha_6,0\} \\ \{0,0,-\alpha_2\} & 0 \end{pmatrix} \\ \begin{pmatrix} 0 & \{\beta_7,-\beta_6,0\} \\ \{0,0,-\beta_2\} & \beta_5 \end{pmatrix} & \begin{pmatrix} \gamma_5 & \{0,0,0\} \\ \{0,0,0\} & \gamma_5 \end{pmatrix} & \begin{pmatrix} \lambda_{12} & \{\lambda_{14},\lambda_{15},0\} \\ \{\lambda_{17},\lambda_{18},\lambda_{19}\} & \lambda_{13} \end{pmatrix} \\ \begin{pmatrix} 0 & \{-\alpha_7,\alpha_6,0\} \\ \{0,0,\alpha_2\} & -\alpha_5 \end{pmatrix} & \begin{pmatrix} \lambda_{13} & \{-\lambda_{14},-\lambda_{15},0\} \\ \{-\lambda_{17},-\lambda_{18},-\lambda_{19}\} & \lambda_{12} \end{pmatrix} & \begin{pmatrix} -\gamma_5 & \{0,0,0\} \\ \{0,0,0\} & -\gamma_5 \end{pmatrix} \end{array}\right)$$

**T[X$_1$[e$_1$], X$_1$[e$_3$], X$_2$[e$_8$]] == ZERO**

True

**Leib[X₁[e₁], X₁[e₃], X₂[e₈]] // Expand**

$$\left( \begin{pmatrix} 0 & \{0,0,0\} \\ \{0,0,0\} & 0 \end{pmatrix} \begin{pmatrix} 0 & \{0,0,0\} \\ \{0,0,0\} & 0 \end{pmatrix} \begin{pmatrix} 0 & \{0,0,0\} \\ \{0,0,0\} & 0 \end{pmatrix} \\ \begin{pmatrix} 0 & \{0,0,0\} \\ \{0,0,0\} & 0 \end{pmatrix} \begin{pmatrix} 0 & \{0,0,0\} \\ \{0,0,0\} & 0 \end{pmatrix} \begin{pmatrix} 0 & \{0,0,0\} \\ \{0,0,\delta_5+\lambda_{15}\} & 0 \end{pmatrix} \\ \begin{pmatrix} 0 & \{0,0,0\} \\ \{0,0,0\} & 0 \end{pmatrix} \begin{pmatrix} 0 & \{0,0,0\} \\ \{0,0,-\delta_5-\lambda_{15}\} & 0 \end{pmatrix} \begin{pmatrix} 0 & \{0,0,0\} \\ \{0,0,0\} & 0 \end{pmatrix} \right)$$

**Δ[X₂[e₈]] = Δ[X₂[e₈]] //. {λ₁₅ → −δ₅}**

$$\left( \begin{pmatrix} 0 & \{0,0,0\} \\ \{0,0,0\} & 0 \end{pmatrix} \begin{pmatrix} \beta_5 & \{-\beta_7,\beta_6,0\} \\ \{0,0,\beta_2\} & 0 \end{pmatrix} \begin{pmatrix} -\alpha_5 & \{\alpha_7,-\alpha_6,0\} \\ \{0,0,-\alpha_2\} & 0 \end{pmatrix} \\ \begin{pmatrix} 0 & \{\beta_7,-\beta_6,0\} \\ \{0,0,-\beta_2\} & \beta_5 \end{pmatrix} \begin{pmatrix} \gamma_5 & \{0,0,0\} \\ \{0,0,0\} & \gamma_5 \end{pmatrix} \begin{pmatrix} \lambda_{12} & \{\lambda_{14},-\delta_5,0\} \\ \{\lambda_{17},\lambda_{18},\lambda_{19}\} & \lambda_{13} \end{pmatrix} \\ \begin{pmatrix} 0 & \{-\alpha_7,\alpha_6,0\} \\ \{0,0,\alpha_2\} & -\alpha_5 \end{pmatrix} \begin{pmatrix} \lambda_{13} & \{-\lambda_{14},\delta_5,0\} \\ \{-\lambda_{17},-\lambda_{18},-\lambda_{19}\} & \lambda_{12} \end{pmatrix} \begin{pmatrix} -\gamma_5 & \{0,0,0\} \\ \{0,0,0\} & -\gamma_5 \end{pmatrix} \right)$$

**T[X₁[e₁], X₁[e₆], X₂[e₈]] == ZERO**

True

**Leib[X₁[e₁], X₁[e₆], X₂[e₈]] // Expand**

$$\left( \begin{pmatrix} 0 & \{0,0,0\} \\ \{0,0,0\} & 0 \end{pmatrix} \begin{pmatrix} 0 & \{0,0,0\} \\ \{0,0,0\} & 0 \end{pmatrix} \begin{pmatrix} 0 & \{0,0,0\} \\ \{0,0,0\} & 0 \end{pmatrix} \\ \begin{pmatrix} 0 & \{0,0,0\} \\ \{0,0,0\} & 0 \end{pmatrix} \begin{pmatrix} 0 & \{0,0,0\} \\ \{0,0,0\} & 0 \end{pmatrix} \begin{pmatrix} 0 & \{0,0,0\} \\ \{\delta_4+\lambda_{12},0,0\} & \lambda_{14}-\delta_6 \end{pmatrix} \\ \begin{pmatrix} 0 & \{0,0,0\} \\ \{0,0,0\} & 0 \end{pmatrix} \begin{pmatrix} \lambda_{14}-\delta_6 & \{0,0,0\} \\ \{-\delta_4-\lambda_{12},0,0\} & 0 \end{pmatrix} \begin{pmatrix} 0 & \{0,0,0\} \\ \{0,0,0\} & 0 \end{pmatrix} \right)$$

**Δ[X₂[e₈]] = Δ[X₂[e₈]] //. {λ₁₂ → −δ₄, λ₁₄ → δ₆}**

$$\left( \begin{pmatrix} 0 & \{0,0,0\} \\ \{0,0,0\} & 0 \end{pmatrix} \begin{pmatrix} \beta_5 & \{-\beta_7,\beta_6,0\} \\ \{0,0,\beta_2\} & 0 \end{pmatrix} \begin{pmatrix} -\alpha_5 & \{\alpha_7,-\alpha_6,0\} \\ \{0,0,-\alpha_2\} & 0 \end{pmatrix} \\ \begin{pmatrix} 0 & \{\beta_7,-\beta_6,0\} \\ \{0,0,-\beta_2\} & \beta_5 \end{pmatrix} \begin{pmatrix} \gamma_5 & \{0,0,0\} \\ \{0,0,0\} & \gamma_5 \end{pmatrix} \begin{pmatrix} -\delta_4 & \{\delta_6,-\delta_5,0\} \\ \{\lambda_{17},\lambda_{18},\lambda_{19}\} & \lambda_{13} \end{pmatrix} \\ \begin{pmatrix} 0 & \{-\alpha_7,\alpha_6,0\} \\ \{0,0,\alpha_2\} & -\alpha_5 \end{pmatrix} \begin{pmatrix} \lambda_{13} & \{-\delta_6,\delta_5,0\} \\ \{-\lambda_{17},-\lambda_{18},-\lambda_{19}\} & -\delta_4 \end{pmatrix} \begin{pmatrix} -\gamma_5 & \{0,0,0\} \\ \{0,0,0\} & -\gamma_5 \end{pmatrix} \right)$$

**T[X₁[e₂], X₁[e₈], X₂[e₈]] == ZERO**

True

**Leib[X₁[e₂], X₁[e₈], X₂[e₈]] // Expand**

$$\begin{pmatrix}
\begin{pmatrix} 0 & \{0,0,0\} \\ \{0,0,0\} & 0 \end{pmatrix} & \begin{pmatrix} 0 & \{0,0,0\} \\ \{0,0,0\} & 0 \end{pmatrix} & \begin{pmatrix} 0 & \{0,0,0\} \\ \{0,0,0\} & 0 \end{pmatrix} \\
\begin{pmatrix} 0 & \{0,0,0\} \\ \{0,0,0\} & 0 \end{pmatrix} & \begin{pmatrix} 0 & \{0,0,0\} \\ \{0,0,0\} & 0 \end{pmatrix} & \begin{pmatrix} 0 & \{\epsilon_3+\lambda_{18}, -\eta_3-\lambda_{17}, 0\} \\ \{0,0,0\} & 0 \end{pmatrix} \\
\begin{pmatrix} 0 & \{0,0,0\} \\ \{0,0,0\} & 0 \end{pmatrix} & \begin{pmatrix} 0 & \{-\epsilon_3-\lambda_{18}, \eta_3+\lambda_{17}, 0\} \\ \{0,0,0\} & 0 \end{pmatrix} & \begin{pmatrix} 0 & \{0,0,0\} \\ \{0,0,0\} & 0 \end{pmatrix}
\end{pmatrix}$$

**Δ[X₂[e₈]] = Δ[X₂[e₈]] //. {λ₁₈ → −ε₃, λ₁₇ → −η₃}**

$$\begin{pmatrix}
\begin{pmatrix} 0 & \{0,0,0\} \\ \{0,0,0\} & 0 \end{pmatrix} & \begin{pmatrix} \beta_5 & \{-\beta_7, \beta_6, 0\} \\ \{0,0,\beta_2\} & 0 \end{pmatrix} & \begin{pmatrix} -\alpha_5 & \{\alpha_7, -\alpha_6, 0\} \\ \{0,0,-\alpha_2\} & 0 \end{pmatrix} \\
\begin{pmatrix} 0 & \{\beta_7, -\beta_6, 0\} \\ \{0,0,-\beta_2\} & \beta_5 \end{pmatrix} & \begin{pmatrix} \gamma_5 & \{0,0,0\} \\ \{0,0,0\} & \gamma_5 \end{pmatrix} & \begin{pmatrix} -\delta_4 & \{\delta_6, -\delta_5, 0\} \\ \{-\eta_3, -\epsilon_3, \lambda_{19}\} & \lambda_{13} \end{pmatrix} \\
\begin{pmatrix} 0 & \{-\alpha_7, \alpha_6, 0\} \\ \{0,0,\alpha_2\} & -\alpha_5 \end{pmatrix} & \begin{pmatrix} \lambda_{13} & \{-\delta_6, \delta_5, 0\} \\ \{\eta_3, \epsilon_3, -\lambda_{19}\} & -\delta_4 \end{pmatrix} & \begin{pmatrix} -\gamma_5 & \{0,0,0\} \\ \{0,0,0\} & -\gamma_5 \end{pmatrix}
\end{pmatrix}$$

**T[X₁[e₇], X₂[e₆], X₂[e₈]] == −X₁[e₁]**

True

**Δ[X₁[e₁]] + Leib[X₁[e₇], X₂[e₆], X₂[e₈]] // Expand**

$$\begin{pmatrix}
\begin{pmatrix} 0 & \{0,0,0\} \\ \{0,0,0\} & 0 \end{pmatrix} & \begin{pmatrix} -\xi+\epsilon_2+\eta_1-\lambda_{19} & \{0,0,0\} \\ \{0,0,0\} & 0 \end{pmatrix} & \begin{pmatrix} 0 & \{0,0,0\} \\ \{0,0,0\} & 0 \end{pmatrix} \\
\begin{pmatrix} 0 & \{0,0,0\} \\ \{0,0,0\} & -\xi+\epsilon_2+\eta_1-\lambda_{19} \end{pmatrix} & \begin{pmatrix} 0 & \{0,0,0\} \\ \{0,0,0\} & 0 \end{pmatrix} & \begin{pmatrix} 0 & \{0,0,0\} \\ \{0,0,0\} & 0 \end{pmatrix} \\
\begin{pmatrix} 0 & \{0,0,0\} \\ \{0,0,0\} & 0 \end{pmatrix} & \begin{pmatrix} 0 & \{0,0,0\} \\ \{0,0,0\} & 0 \end{pmatrix} & \begin{pmatrix} 0 & \{0,0,0\} \\ \{0,0,0\} & 0 \end{pmatrix}
\end{pmatrix}$$

**Δ[X₂[e₈]] = Δ[X₂[e₈]] //. {λ₁₉ → −ξ + ε₂ + η₁}**

$$\begin{pmatrix}
\begin{pmatrix} 0 & \{0,0,0\} \\ \{0,0,0\} & 0 \end{pmatrix} & \begin{pmatrix} \beta_5 & \{-\beta_7, \beta_6, 0\} \\ \{0,0,\beta_2\} & 0 \end{pmatrix} & \begin{pmatrix} -\alpha_5 & \{\alpha_7, -\alpha_6, 0\} \\ \{0,0,-\alpha_2\} & 0 \end{pmatrix} \\
\begin{pmatrix} 0 & \{\beta_7, -\beta_6, 0\} \\ \{0,0,-\beta_2\} & \beta_5 \end{pmatrix} & \begin{pmatrix} \gamma_5 & \{0,0,0\} \\ \{0,0,0\} & \gamma_5 \end{pmatrix} & \begin{pmatrix} -\delta_4 & \{\delta_6, -\delta_5, 0\} \\ \{-\eta_3, -\epsilon_3, -\xi+\epsilon_2+\eta_1\} & \lambda_{13} \end{pmatrix} \\
\begin{pmatrix} 0 & \{-\alpha_7, \alpha_6, 0\} \\ \{0,0,\alpha_2\} & -\alpha_5 \end{pmatrix} & \begin{pmatrix} \lambda_{13} & \{-\delta_6, \delta_5, 0\} \\ \{\eta_3, \epsilon_3, \xi-\epsilon_2-\eta_1\} & -\delta_4 \end{pmatrix} & \begin{pmatrix} -\gamma_5 & \{0,0,0\} \\ \{0,0,0\} & -\gamma_5 \end{pmatrix}
\end{pmatrix}$$

**T[X₁[e₂], X₂[e₇], X₂[e₈]] == X₁[e₃]**

True

```
Δ[X₁[e₃]] - Leib[X₁[e₂], X₂[e₇], X₂[e₈]] // Expand
```

$$\left(\begin{array}{ccc} \begin{pmatrix} 0 & \{0,0,0\} \\ \{0,0,0\} & 0 \end{pmatrix} & \begin{pmatrix} 0 & \{0,0,0\} \\ \{0, \eta_4 - \lambda_{13}, 0\} & 0 \end{pmatrix} & \begin{pmatrix} 0 & \{0,0,0\} \\ \{0,0,0\} & 0 \end{pmatrix} \\ \begin{pmatrix} 0 & \{0,0,0\} \\ \{0, \lambda_{13} - \eta_4, 0\} & 0 \end{pmatrix} & \begin{pmatrix} 0 & \{0,0,0\} \\ \{0,0,0\} & 0 \end{pmatrix} & \begin{pmatrix} 0 & \{0,0,0\} \\ \{0,0,0\} & 0 \end{pmatrix} \\ \begin{pmatrix} 0 & \{0,0,0\} \\ \{0,0,0\} & 0 \end{pmatrix} & \begin{pmatrix} 0 & \{0,0,0\} \\ \{0,0,0\} & 0 \end{pmatrix} & \begin{pmatrix} 0 & \{0,0,0\} \\ \{0,0,0\} & 0 \end{pmatrix} \end{array}\right)$$

```
Δ[X₂[e₈]] = Δ[X₂[e₈]] //. {λ₁₃ → η₄}
```

$$\left(\begin{array}{ccc} \begin{pmatrix} 0 & \{0,0,0\} \\ \{0,0,0\} & 0 \end{pmatrix} & \begin{pmatrix} \beta_5 & \{-\beta_7, \beta_6, 0\} \\ \{0,0,\beta_2\} & 0 \end{pmatrix} & \begin{pmatrix} -\alpha_5 & \{\alpha_7, -\alpha_6, 0\} \\ \{0,0,-\alpha_2\} & 0 \end{pmatrix} \\ \begin{pmatrix} 0 & \{\beta_7, -\beta_6, 0\} \\ \{0,0,-\beta_2\} & \beta_5 \end{pmatrix} & \begin{pmatrix} \gamma_5 & \{0,0,0\} \\ \{0,0,0\} & \gamma_5 \end{pmatrix} & \begin{pmatrix} -\delta_4 & \{\delta_6, -\delta_5, 0\} \\ \{-\eta_3, -\epsilon_3, -\xi+\epsilon_2+\eta_1\} & \eta_4 \end{pmatrix} \\ \begin{pmatrix} 0 & \{-\alpha_7, \alpha_6, 0\} \\ \{0,0,\alpha_2\} & -\alpha_5 \end{pmatrix} & \begin{pmatrix} \eta_4 & \{-\delta_6, \delta_5, 0\} \\ \{\eta_3, \epsilon_3, \xi-\epsilon_2-\eta_1\} & -\delta_4 \end{pmatrix} & \begin{pmatrix} -\gamma_5 & \{0,0,0\} \\ \{0,0,0\} & -\gamma_5 \end{pmatrix} \end{array}\right)$$

```
Variables[Δ[X₂[e₈]]]
```

$\{\xi, \alpha_2, \alpha_5, \alpha_6, \alpha_7, \beta_2, \beta_5, \beta_6, \beta_7, \gamma_5, \delta_4, \delta_5, \delta_6, \epsilon_2, \epsilon_3, \eta_1, \eta_3, \eta_4\}$

## ■ Looking for identities for $X_3[e_i]$ and its image

```
T[E₁[1], X₁[e₁], X₂[e₁]] == X₃[e₁]
```

True

```
Δ[X₃[e₁]] = Leib[E₁[1], X₁[e₁], X₂[e₁]]
```

$$\left(\begin{array}{ccc} \begin{pmatrix} -\beta_2 & \{0,0,0\} \\ \{0,0,0\} & -\beta_2 \end{pmatrix} & \begin{pmatrix} -\gamma_2 & \{\gamma_3, \gamma_4, \gamma_5\} \\ \{0,0,0\} & 0 \end{pmatrix} & \begin{pmatrix} \xi+\delta_1 & \{\epsilon_4, -\eta_5, \eta_4\} \\ \{\delta_5, \delta_6, \delta_7\} & 0 \end{pmatrix} \\ \begin{pmatrix} 0 & \{-\gamma_3, -\gamma_4, -\gamma_5\} \\ \{0,0,0\} & -\gamma_2 \end{pmatrix} & \begin{pmatrix} 0 & \{0,0,0\} \\ \{0,0,0\} & 0 \end{pmatrix} & \begin{pmatrix} \alpha_2 & \{0,0,0\} \\ \{-\alpha_6, -\alpha_7, -\alpha_8\} & 0 \end{pmatrix} \\ \begin{pmatrix} 0 & \{-\epsilon_4, \eta_5, -\eta_4\} \\ \{-\delta_5, -\delta_6, -\delta_7\} & \xi+\delta_1 \end{pmatrix} & \begin{pmatrix} 0 & \{0,0,0\} \\ \{\alpha_6, \alpha_7, \alpha_8\} & \alpha_2 \end{pmatrix} & \begin{pmatrix} \beta_2 & \{0,0,0\} \\ \{0,0,0\} & \beta_2 \end{pmatrix} \end{array}\right)$$

```
T[E₁[1], X₁[e₂], X₂[e₂]] == X₃[e₂]
```

True

```
Δ[X₃[e₂]] = Leib[E₁[1], X₁[e₂], X₂[e₂]]
```

$$\left(\begin{array}{ccc} \begin{pmatrix} -\beta_1 & \{0,0,0\} \\ \{0,0,0\} & -\beta_1 \end{pmatrix} & \begin{pmatrix} 0 & \{0,0,0\} \\ \{\gamma_6, \gamma_7, \gamma_8\} & -\gamma_1 \end{pmatrix} & \begin{pmatrix} 0 & \{\rho_1, \rho_2, \rho_3\} \\ \{-\chi_1, \psi_2, -\psi_1\} & -\xi-\delta_1 \end{pmatrix} \\ \begin{pmatrix} -\gamma_1 & \{0,0,0\} \\ \{-\gamma_6, -\gamma_7, -\gamma_8\} & 0 \end{pmatrix} & \begin{pmatrix} 0 & \{0,0,0\} \\ \{0,0,0\} & 0 \end{pmatrix} & \begin{pmatrix} 0 & \{-\alpha_3, -\alpha_4, -\alpha_5\} \\ \{0,0,0\} & \alpha_1 \end{pmatrix} \\ \begin{pmatrix} -\xi-\delta_1 & \{-\rho_1, -\rho_2, -\rho_3\} \\ \{\chi_1, -\psi_2, \psi_1\} & 0 \end{pmatrix} & \begin{pmatrix} \alpha_1 & \{\alpha_3, \alpha_4, \alpha_5\} \\ \{0,0,0\} & 0 \end{pmatrix} & \begin{pmatrix} \beta_1 & \{0,0,0\} \\ \{0,0,0\} & \beta_1 \end{pmatrix} \end{array}\right)$$

**T[E₁[1], X₁[e₁], X₂[e₃]] == X₃[e₃]**

True

**Δ[X₃[e₃]] = Leib[E₁[1], X₁[e₁], X₂[e₃]]**

$$\left(\begin{pmatrix} \beta_6 & \{0,0,0\} \\ \{0,0,0\} & \beta_6 \end{pmatrix} \begin{pmatrix} \gamma_6 & \{-\gamma_1,0,0\} \\ \{0,-\gamma_5,\gamma_4\} & 0 \end{pmatrix} \begin{pmatrix} -\chi_1 & \{\eta_1-\xi,\eta_2,\eta_3\} \\ \{0,\delta_4,-\delta_3\} & \delta_5 \end{pmatrix}\right.$$
$$\begin{pmatrix} 0 & \{\gamma_1,0,0\} \\ \{0,\gamma_5,-\gamma_4\} & \gamma_6 \end{pmatrix} \begin{pmatrix} 0 & \{0,0,0\} \\ \{0,0,0\} & 0 \end{pmatrix} \begin{pmatrix} 0 & \{\alpha_2,0,0\} \\ \{0,-\alpha_5,\alpha_4\} & -\alpha_6 \end{pmatrix}$$
$$\left.\begin{pmatrix} \delta_5 & \{\xi-\eta_1,-\eta_2,-\eta_3\} \\ \{0,-\delta_4,\delta_3\} & -\chi_1 \end{pmatrix} \begin{pmatrix} -\alpha_6 & \{-\alpha_2,0,0\} \\ \{0,\alpha_5,-\alpha_4\} & 0 \end{pmatrix} \begin{pmatrix} -\beta_6 & \{0,0,0\} \\ \{0,0,0\} & -\beta_6 \end{pmatrix}\right)$$

**T[E₁[1], X₁[e₁], X₂[e₄]] == X₃[e₄]**

True

**Δ[X₃[e₄]] = Leib[E₁[1], X₁[e₁], X₂[e₄]]**

$$\left(\begin{pmatrix} \beta_7 & \{0,0,0\} \\ \{0,0,0\} & \beta_7 \end{pmatrix} \begin{pmatrix} \gamma_7 & \{0,-\gamma_1,0\} \\ \{\gamma_5,0,-\gamma_3\} & 0 \end{pmatrix} \begin{pmatrix} \psi_2 & \{\epsilon_1,\epsilon_2-\xi,\epsilon_3\} \\ \{-\delta_4,0,\delta_2\} & \delta_6 \end{pmatrix}\right.$$
$$\begin{pmatrix} 0 & \{0,\gamma_1,0\} \\ \{-\gamma_5,0,\gamma_3\} & \gamma_7 \end{pmatrix} \begin{pmatrix} 0 & \{0,0,0\} \\ \{0,0,0\} & 0 \end{pmatrix} \begin{pmatrix} 0 & \{0,\alpha_2,0\} \\ \{\alpha_5,0,-\alpha_3\} & -\alpha_7 \end{pmatrix}$$
$$\left.\begin{pmatrix} \delta_6 & \{-\epsilon_1,\xi-\epsilon_2,-\epsilon_3\} \\ \{\delta_4,0,-\delta_2\} & \psi_2 \end{pmatrix} \begin{pmatrix} -\alpha_7 & \{0,-\alpha_2,0\} \\ \{-\alpha_5,0,\alpha_3\} & 0 \end{pmatrix} \begin{pmatrix} -\beta_7 & \{0,0,0\} \\ \{0,0,0\} & -\beta_7 \end{pmatrix}\right)$$

**T[E₁[1], X₁[e₁], X₂[e₅]] == X₃[e₅]**

True

**Δ[X₃[e₅]] = Leib[E₁[1], X₁[e₁], X₂[e₅]]**

$$\left(\begin{pmatrix} \beta_8 & \{0,0,0\} \\ \{0,0,0\} & \beta_8 \end{pmatrix} \begin{pmatrix} \gamma_8 & \{0,0,-\gamma_1\} \\ \{-\gamma_4,\gamma_3,0\} & 0 \end{pmatrix} \begin{pmatrix} -\psi_1 & \{\phi_1,\phi_2,\xi+\delta_1-\epsilon_2-\eta_1\} \\ \{\delta_3,-\delta_2,0\} & \delta_7 \end{pmatrix}\right.$$
$$\begin{pmatrix} 0 & \{0,0,\gamma_1\} \\ \{\gamma_4,-\gamma_3,0\} & \gamma_8 \end{pmatrix} \begin{pmatrix} 0 & \{0,0,0\} \\ \{0,0,0\} & 0 \end{pmatrix} \begin{pmatrix} 0 & \{0,0,\alpha_2\} \\ \{-\alpha_4,\alpha_3,0\} & -\alpha_8 \end{pmatrix}$$
$$\left.\begin{pmatrix} \delta_7 & \{-\phi_1,-\phi_2,-\xi-\delta_1+\epsilon_2+\eta_1\} \\ \{-\delta_3,\delta_2,0\} & -\psi_1 \end{pmatrix} \begin{pmatrix} -\alpha_8 & \{0,0,-\alpha_2\} \\ \{\alpha_4,-\alpha_3,0\} & 0 \end{pmatrix} \begin{pmatrix} -\beta_8 & \{0,0,0\} \\ \{0,0,0\} & -\beta_8 \end{pmatrix}\right)$$

**T[E₁[1], X₁[e₆], X₂[e₁]] == X₃[e₆]**

True

**Δ[X₃[e₆]] = Leib[E₁[1], X₁[e₆], X₂[e₁]]**

$$\left(\begin{pmatrix} \beta_3 & \{0,0,0\} \\ \{0,0,0\} & \beta_3 \end{pmatrix} \begin{pmatrix} 0 & \{0,\gamma_8,-\gamma_7\} \\ \{-\gamma_2,0,0\} & \gamma_3 \end{pmatrix} \begin{pmatrix} \rho_1 & \{0,-\rho_6,\rho_5\} \\ \{\xi-\eta_1,-\epsilon_1,-\phi_1\} & \epsilon_4 \end{pmatrix}\right.$$
$$\begin{pmatrix} \gamma_3 & \{0,-\gamma_8,\gamma_7\} \\ \{\gamma_2,0,0\} & 0 \end{pmatrix} \begin{pmatrix} 0 & \{0,0,0\} \\ \{0,0,0\} & 0 \end{pmatrix} \begin{pmatrix} -\alpha_3 & \{0,\alpha_8,-\alpha_7\} \\ \{\alpha_1,0,0\} & 0 \end{pmatrix}$$
$$\left.\begin{pmatrix} \epsilon_4 & \{0,\rho_6,-\rho_5\} \\ \{\eta_1-\xi,\epsilon_1,\phi_1\} & \rho_1 \end{pmatrix} \begin{pmatrix} 0 & \{0,-\alpha_8,\alpha_7\} \\ \{-\alpha_1,0,0\} & -\alpha_3 \end{pmatrix} \begin{pmatrix} -\beta_3 & \{0,0,0\} \\ \{0,0,0\} & -\beta_3 \end{pmatrix}\right)$$

**T[E₁[1], X₁[e₇], X₂[e₁]] == X₃[e₇]**

True

**Δ[X₃[e₇]] = Leib[E₁[1], X₁[e₇], X₂[e₁]]**

$$\left(\begin{array}{ccc} \begin{pmatrix} \beta_4 & \{0,0,0\} \\ \{0,0,0\} & \beta_4 \end{pmatrix} & \begin{pmatrix} 0 & \{-\gamma_8,0,\gamma_6\} \\ \{0,-\gamma_2,0\} & \gamma_4 \end{pmatrix} & \begin{pmatrix} \rho_2 & \{\rho_6,0,-\rho_4\} \\ \{-\eta_2,\xi-\epsilon_2,-\phi_2\} & -\eta_5 \end{pmatrix} \\ \begin{pmatrix} \gamma_4 & \{\gamma_8,0,-\gamma_6\} \\ \{0,\gamma_2,0\} & 0 \end{pmatrix} & \begin{pmatrix} 0 & \{0,0,0\} \\ \{0,0,0\} & 0 \end{pmatrix} & \begin{pmatrix} -\alpha_4 & \{-\alpha_8,0,\alpha_6\} \\ \{0,\alpha_1,0\} & 0 \end{pmatrix} \\ \begin{pmatrix} -\eta_5 & \{-\rho_6,0,\rho_4\} \\ \{\eta_2,\epsilon_2-\xi,\phi_2\} & \rho_2 \end{pmatrix} & \begin{pmatrix} 0 & \{\alpha_8,0,-\alpha_6\} \\ \{0,-\alpha_1,0\} & -\alpha_4 \end{pmatrix} & \begin{pmatrix} -\beta_4 & \{0,0,0\} \\ \{0,0,0\} & -\beta_4 \end{pmatrix} \end{array}\right)$$

**T[E₁[1], X₁[e₈], X₂[e₁]] == X₃[e₈]**

True

**Δ[X₃[e₈]] = Leib[E₁[1], X₁[e₈], X₂[e₁]]**

$$\left(\begin{array}{ccc} \begin{pmatrix} \beta_5 & \{0,0,0\} \\ \{0,0,0\} & \beta_5 \end{pmatrix} & \begin{pmatrix} 0 & \{\gamma_7,-\gamma_6,0\} \\ \{0,0,-\gamma_2\} & \gamma_5 \end{pmatrix} & \begin{pmatrix} \rho_3 & \{-\rho_5,\rho_4,0\} \\ \{-\eta_3,-\epsilon_3,-\xi-\delta_1+\epsilon_2+\eta_1\} & \eta_4 \end{pmatrix} \\ \begin{pmatrix} \gamma_5 & \{-\gamma_7,\gamma_6,0\} \\ \{0,0,\gamma_2\} & 0 \end{pmatrix} & \begin{pmatrix} 0 & \{0,0,0\} \\ \{0,0,0\} & 0 \end{pmatrix} & \begin{pmatrix} -\alpha_5 & \{\alpha_7,-\alpha_6,0\} \\ \{0,0,\alpha_1\} & 0 \end{pmatrix} \\ \begin{pmatrix} \eta_4 & \{\rho_5,-\rho_4,0\} \\ \{\eta_3,\epsilon_3,\xi+\delta_1-\epsilon_2-\eta_1\} & \rho_3 \end{pmatrix} & \begin{pmatrix} 0 & \{-\alpha_7,\alpha_6,0\} \\ \{0,0,-\alpha_1\} & -\alpha_5 \end{pmatrix} & \begin{pmatrix} -\beta_5 & \{0,0,0\} \\ \{0,0,0\} & -\beta_5 \end{pmatrix} \end{array}\right)$$

## ■ Generic matrix and a basis

```
paramet = Union[Variables[Δ[E₁[1]]],
    Variables[Δ[E₂[1]]], Variables[Δ[E₃[1]]]];
Do[paramet = Union[paramet, Variables[Δ[Xᵢ[eⱼ]]]], {i, 3}, {j, 8}]
```

**paramet**

$\{\xi, \alpha_1, \alpha_2, \alpha_3, \alpha_4, \alpha_5, \alpha_6, \alpha_7, \alpha_8, \beta_1, \beta_2, \beta_3, \beta_4, \beta_5, \beta_6, \beta_7,$
$\beta_8, \gamma_1, \gamma_2, \gamma_3, \gamma_4, \gamma_5, \gamma_6, \gamma_7, \gamma_8, \delta_1, \delta_2, \delta_3, \delta_4, \delta_5, \delta_6, \delta_7, \epsilon_1, \epsilon_2,$
$\epsilon_3, \epsilon_4, \eta_1, \eta_2, \eta_3, \eta_4, \eta_5, \rho_1, \rho_2, \rho_3, \rho_4, \rho_5, \rho_6, \phi_1, \phi_2, \chi_1, \psi_1, \psi_2\}$

**Length[paramet]**

```
li = {Δ[E₁[1]], Δ[E₂[1]], Δ[E₃[1]], Δ[X₁[e₁]],
    Δ[X₁[e₂]], Δ[X₁[e₃]], Δ[X₁[e₄]], Δ[X₁[e₅]],
    Δ[X₁[e₆]], Δ[X₁[e₇]], Δ[X₁[e₈]], Δ[X₂[e₁]], Δ[X₂[e₂]],
    Δ[X₂[e₃]], Δ[X₂[e₄]], Δ[X₂[e₅]], Δ[X₂[e₆]], Δ[X₂[e₇]],
    Δ[X₂[e₈]], Δ[X₃[e₁]], Δ[X₃[e₂]], Δ[X₃[e₃]], Δ[X₃[e₄]],
    Δ[X₃[e₅]], Δ[X₃[e₆]], Δ[X₃[e₇]], Δ[X₃[e₈]]};
```

```
supermatrix = {};
Do[se = li[[i]]; f = {se[[1, 1]][[1, 1]],
   se[[2, 2]][[1, 1]], se[[3, 3]][[1, 1]], se[[1, 2]][[1, 1]],
      se[[1, 2]][[2, 2]], se[[1, 2]][[1, 2]][[1]],
   se[[1, 2]][[1, 2]][[2]], se[[1, 2]][[1, 2]][[3]],
   se[[1, 2]][[2, 1]][[1]], se[[1, 2]][[2, 1]][[2]],
   se[[1, 2]][[2, 1]][[3]], se[[2, 3]][[1, 1]],
      se[[2, 3]][[2, 2]], se[[2, 3]][[1, 2]][[1]],
   se[[2, 3]][[1, 2]][[2]], se[[2, 3]][[1, 2]][[3]],
   se[[2, 3]][[2, 1]][[1]], se[[2, 3]][[2, 1]][[2]],
   se[[2, 3]][[2, 1]][[3]],
      se[[1, 3]][[1, 1]], se[[1, 3]][[2, 2]],
   se[[1, 3]][[1, 2]][[1]], se[[1, 3]][[1, 2]][[2]],
   se[[1, 3]][[1, 2]][[3]], se[[1, 3]][[2, 1]][[1]],
   se[[1, 3]][[2, 1]][[2]], se[[1, 3]][[2, 1]][[3]]};
   AppendTo[supermatrix, f], {i, 27}]

supermatrix
```

[27×27 matrix displayed with entries in terms of $\alpha_i, \beta_i, \gamma_i, \delta_i, \epsilon_i, \eta_i, \rho_i, \phi_i, \psi_i, \xi, \chi_i$]

```
δ_{i_,j_} := If[i == j, 1, 0];
var = Variables[supermatrix]; num = Length[var];
reglas = Table[var[[j]] -> δ_{i,j}, {i, num}, {j, num}];
Do[x_i = supermatrix //. reglas[[i]], {i, num}]

Dimensions[supermatrix]
```

{27, 27}

## ■ supermatrix $\subset f_4(O_s)$

We will prove that any matrix as the previous one is a derivation:

```
Coor[(α_  a_  c_ )
     (_   β_  b_)
     (_   _   χ_)] :=
 {{α[[1, 1]], β[[1, 1]], χ[[1, 1]], a[[1, 1]], a[[2, 2]],
   a[[1, 2]][[1]], a[[1, 2]][[2]], a[[1, 2]][[3]],
   a[[2, 1]][[1]], a[[2, 1]][[2]], a[[2, 1]][[3]], b[[1, 1]],
   b[[2, 2]], b[[1, 2]][[1]], b[[1, 2]][[2]], b[[1, 2]][[3]],
   b[[2, 1]][[1]], b[[2, 1]][[2]], b[[2, 1]][[3]], c[[1, 1]],
   c[[2, 2]], c[[1, 2]][[1]], c[[1, 2]][[2]], c[[1, 2]][[3]],
   c[[2, 1]][[1]], c[[2, 1]][[2]], c[[2, 1]][[3]]}}

InvCoor[{{a1_, a2_, a3_, a4_, a5_, a6_, a7_, a8_,
   a9_, a10_, a11_, a12_, a13_, a14_, a15_, a16_, a17_, a18_,
   a19_, a20_, a21_, a22_, a23_, a24_, a25_, a26_, a27_}}] :=
```

$$\begin{pmatrix} \begin{pmatrix} a1 & \{0,0,0\} \\ \{0,0,0\} & a1 \end{pmatrix} & \begin{pmatrix} a4 & \{a6,a7,a8\} \\ \{a9,a10,a11\} & a5 \end{pmatrix} & \begin{pmatrix} a20 & \cdots \\ \{a25,a26 & \end{pmatrix} \\ \sigma\left[\begin{pmatrix} a4 & \{a6,a7,a8\} \\ \{a9,a10,a11\} & a5 \end{pmatrix}\right] & \begin{pmatrix} a2 & \{0,0,0\} \\ \{0,0,0\} & a2 \end{pmatrix} & \begin{pmatrix} a12 & \cdots \\ \{a17,a18 & \end{pmatrix} \\ \sigma\left[\begin{pmatrix} a20 & \{a22,a23,a24\} \\ \{a25,a26,a27\} & a21 \end{pmatrix}\right] & \sigma\left[\begin{pmatrix} a12 & \{a14,a15,a16\} \\ \{a17,a18,a19\} & a13 \end{pmatrix}\right] & \begin{pmatrix} \cdots \\ \{0 \end{pmatrix} \end{pmatrix}$$

D(1)=0

```
InvCoor[Coor[ONE].supermatrix] == ZERO
```

True

$D(U_x(y)) = T(D(x), y, x) + U_x(D(y))$

```
InvCoor[Coor[U_generic[E_1[1]]].supermatrix] -
  T[InvCoor[Coor[generic].supermatrix], E_1[1], generic] -
  U_generic[InvCoor[Coor[E_1[1]].supermatrix]] // Expand
```

$$\begin{pmatrix} \begin{pmatrix} 0 & \{0,0,0\} \\ \{0,0,0\} & 0 \end{pmatrix} & \begin{pmatrix} 0 & \{0,0,0\} \\ \{0,0,0\} & 0 \end{pmatrix} & \begin{pmatrix} 0 & \{0,0,0\} \\ \{0,0,0\} & 0 \end{pmatrix} \\ \begin{pmatrix} 0 & \{0,0,0\} \\ \{0,0,0\} & 0 \end{pmatrix} & \begin{pmatrix} 0 & \{0,0,0\} \\ \{0,0,0\} & 0 \end{pmatrix} & \begin{pmatrix} 0 & \{0,0,0\} \\ \{0,0,0\} & 0 \end{pmatrix} \\ \begin{pmatrix} 0 & \{0,0,0\} \\ \{0,0,0\} & 0 \end{pmatrix} & \begin{pmatrix} 0 & \{0,0,0\} \\ \{0,0,0\} & 0 \end{pmatrix} & \begin{pmatrix} 0 & \{0,0,0\} \\ \{0,0,0\} & 0 \end{pmatrix} \end{pmatrix}$$

```
InvCoor[Coor[U_generic[E_2[1]]].supermatrix] -
  T[InvCoor[Coor[generic].supermatrix], E_2[1], generic] -
  U_generic[InvCoor[Coor[E_2[1]].supermatrix]] // Expand
```

$$\begin{pmatrix} \begin{pmatrix} 0 & \{0,0,0\} \\ \{0,0,0\} & 0 \end{pmatrix} & \begin{pmatrix} 0 & \{0,0,0\} \\ \{0,0,0\} & 0 \end{pmatrix} & \begin{pmatrix} 0 & \{0,0,0\} \\ \{0,0,0\} & 0 \end{pmatrix} \\ \begin{pmatrix} 0 & \{0,0,0\} \\ \{0,0,0\} & 0 \end{pmatrix} & \begin{pmatrix} 0 & \{0,0,0\} \\ \{0,0,0\} & 0 \end{pmatrix} & \begin{pmatrix} 0 & \{0,0,0\} \\ \{0,0,0\} & 0 \end{pmatrix} \\ \begin{pmatrix} 0 & \{0,0,0\} \\ \{0,0,0\} & 0 \end{pmatrix} & \begin{pmatrix} 0 & \{0,0,0\} \\ \{0,0,0\} & 0 \end{pmatrix} & \begin{pmatrix} 0 & \{0,0,0\} \\ \{0,0,0\} & 0 \end{pmatrix} \end{pmatrix}$$

```
InvCoor[Coor[U_generic[E_3[1]]].supermatrix] -
  T[InvCoor[Coor[generic].supermatrix], E_3[1], generic] -
  U_generic[InvCoor[Coor[E_3[1]].supermatrix]] // Expand
```

$$\left(\begin{array}{cc}\begin{pmatrix} 0 & \{0,0,0\} \\ \{0,0,0\} & 0 \end{pmatrix} & \begin{pmatrix} 0 & \{0,0,0\} \\ \{0,0,0\} & 0 \end{pmatrix} & \begin{pmatrix} 0 & \{0,0,0\} \\ \{0,0,0\} & 0 \end{pmatrix} \\ \begin{pmatrix} 0 & \{0,0,0\} \\ \{0,0,0\} & 0 \end{pmatrix} & \begin{pmatrix} 0 & \{0,0,0\} \\ \{0,0,0\} & 0 \end{pmatrix} & \begin{pmatrix} 0 & \{0,0,0\} \\ \{0,0,0\} & 0 \end{pmatrix} \\ \begin{pmatrix} 0 & \{0,0,0\} \\ \{0,0,0\} & 0 \end{pmatrix} & \begin{pmatrix} 0 & \{0,0,0\} \\ \{0,0,0\} & 0 \end{pmatrix} & \begin{pmatrix} 0 & \{0,0,0\} \\ \{0,0,0\} & 0 \end{pmatrix}\end{array}\right)$$

```
Do[Print[Expand[
    InvCoor[Coor[U_generic[X_{1,2}[e_i]]].supermatrix] - T[InvCoor[
        Coor[generic].supermatrix], X_{1,2}[e_i], generic] - U_generic[
      InvCoor[Coor[X_{1,2}[e_i]].supermatrix]]] == ZERO], {i, 8}]
```

True

True

True

True

True

True

True

True

```
Do[Print[Expand[
    InvCoor[Coor[U_generic[X_{2,3}[e_i]]].supermatrix] - T[InvCoor[
        Coor[generic].supermatrix], X_{2,3}[e_i], generic] - U_generic[
      InvCoor[Coor[X_{2,3}[e_i]].supermatrix]]] == ZERO], {i, 8}]
```

True

True

True

True

True

True

True

True

```
Do[Print[Expand[
    InvCoor[Coor[U_generic[X_{1,3}[e_i]]].supermatrix] - T[InvCoor[
      Coor[generic].supermatrix], X_{1,3}[e_i], generic] - U_generic[
      InvCoor[Coor[X_{1,3}[e_i]].supermatrix]]] == ZERO], {i, 8}]
```

True

True

True

True

True

True

True

True